# Model-Based Derivative-Free Optimization Methods and Software

Tom M. Ragonneau

*Doctor of Philosophy*

*The Hong Kong Polytechnic University*

2022

The Hong Kong Polytechnic University
Department of Applied Mathematics

# Model-Based Derivative-Free Optimization Methods and Software

Tom M. Ragonneau

*A thesis submitted in partial fulfilment of the requirements for the degree of Doctor of Philosophy*

August, 2022

*This page intentionally left blank*

# CERTIFICATE OF ORIGINALITY

I hereby declare that this thesis is my own work and that, to the best of my knowledge and belief, it reproduces no material previously published or written, nor material that has been accepted for the award of any other degree or diploma, except where due acknowledgement has been made in the text.

____________________ (Signed)

*Tom M. Ragonneau* (Name of student)

*This page intentionally left blank*

To

Marie-Thérèse Bretzner,
Isabelle and Franck Ragonneau, and
Julie Venturini,

with thanks and love

*This page intentionally left blank*

# Abstract

This thesis studies *derivative-free optimization* (DFO), particularly model-based methods and software. These methods are motivated by optimization problems for which it is impossible or prohibitively expensive to access the first-order information of the objective function and possibly the constraint functions. Such problems frequently arise in engineering and industrial applications, and keep emerging due to recent advances in data science and machine learning.

We first provide an overview of DFO and interpolation models for DFO methods. In addition, we show that an interpolation set employed by Powell for underdetermined quadratic interpolation is optimal in terms of well-poisedness.

We then present PDFO, a package that we develop to provide both MATLAB and Python interfaces to Powell's model-based DFO solvers, namely COBYLA, UOBYQA, NEWUOA, BOBYQA, and LINCOA. They were implemented by Powell in Fortran 77, and hence, are becoming inaccessible to many users nowadays. PDFO provides user-friendly interfaces to these solvers, so that users do not need to deal with the Fortran code. In addition, it patches bugs in the original Fortran implementation. We also share some observations about the behavior of Powell's solvers.

A major part of this thesis is devoted to the development of a new DFO method based on the *sequential quadratic programming* (SQP) method. Therefore, we first present an overview of the SQP method and provide some perspectives on its theory and practice. In particular, we show that the objective function of the SQP subproblem is a natural quadratic approximation of the original objective function in the tangent space of a surface. Moreover, we propose an extension of the Byrd-Omojokun approach for solving trust-region SQP subproblems with inequality constraints. This extension works quite well in our experiments.

Finally, we elaborate on the development of our new DFO method, named COBYQA after *Constrained Optimization BY Quadratic Approximations*. This derivative-free trust-region SQP method is designed to tackle nonlinearly constrained optimization problems that admit equality and inequality constraints. An important feature of COBYQA is that it always respects bound constraints, if any, which is motivated by applications where the objective function is undefined when bounds are violated. COBYQA builds quadratic trust-region models based on the derivative-free symmetric Broyden update proposed by Powell. We provide a detailed description of COBYQA, including its subproblem solvers, and introduce its Python implementation. We expose extensive numerical experiments of COBYQA, showing evident advantages of COBYQA compared with Powell's DFO solvers. These experiments demonstrate that COBYQA is an excellent successor to COBYLA as a general-purpose DFO solver.

# Acknowledgments

First, I would like to express my deepest gratitude to Prof. Xiaojun Chen and Dr. Zaikun Zhang from the Department of Applied Mathematics of The Hong Kong Polytechnic University for supervising this thesis. Their patience and motivation helped me to overcome the challenges I faced during my study.

The support and encouragement provided by Prof. Serge Gratton from the Department of Applied Mathematics of Toulouse INP/ENSEEIHT and Prof. Ya-xiang Yuan from the Institute of Computational Mathematics and Scientific/Engineering Computing of the Chinese Academy of Sciences are greatly appreciated. They shared beneficial comments on my research projects and helped me through my research studies.

Words cannot express my gratitude to my parents and my sister for their wholesome counsel, wise guidance, and sympathetic ears during the research study, and my friends, for always being by my side, providing invigorating discussions along with resourceful and necessary distractions.

Finally, I would like to extend my sincere thanks to the Research Grants Council of Hong Kong, under the aegis of the University Grants Committee of Hong Kong, for providing financial support under the Hong Kong Ph.D. Fellowship Scheme. I also thank The Hong Kong Polytechnic University for supporting the research study and, more specifically, the Department of Applied Mathematics for providing the necessary equipment and research facilities.

*This page intentionally left blank*

# Contents

# List of Figures

*This page intentionally left blank*

# List of Tables

*This page intentionally left blank*

# List of Algorithms

*This page intentionally left blank*

# Listings

# Acronyms

| | |
|---|---|
| **ACQ** | Abadie constraint qualification. 8 |
| **AUC** | Area under the ROC curve. 5 |
| **BBO** | Black-box optimization. 1 |
| **BFO** | Brute force optimizer. 5 |
| **CDT** | Celis-Dennis-Tapia. 57 |
| **CPLD** | Constant positive linear dependence constraint qualification. 8 |
| **CPU** | Central processing unit. 4 |
| **CRCQ** | Constant rank constraint qualification. 8 |
| **DFO** | Derivative-free optimization. vii, 1 |
| **DFP** | Davidon-Fletcher-Powell. 45 |
| **GCQ** | Guignard constraint qualification. 8 |
| **GPS** | Generalized pattern search. 12 |
| **KKT** | Karush-Kuhn-Tucker. 9 |
| **LICQ** | Linear independence constraint qualification. 8 |
| **MADS** | Mesh adaptive direct search. 5 |
| **MDO** | Multi-disciplinary design optimization. 6 |
| **MFCQ** | Mangasarian-Fromovitz constraint qualification. 8 |
| **PSB** | Powell symmetric Broyden. 33 |
| **QNCQ** | Quasi-normality constraint qualification. 8 |
| **RBF** | Radial basis function. 12 |
| **RS** | Random search. 54 |
| **S$\ell_p$QP** | Sequential $\ell_p$-quadratic programming. 81 |
| **SC** | Slater's condition. 8 |
| **SLQP** | Sequential linear-quadratic programming. 81 |
| **SQP** | Sequential quadratic programming. vii, 22 |
| **SVC** | Support vector classification. 5 |
| **SVM** | Support vector machine. 5 |
| **TCG** | Truncated conjugate gradient. 48 |

*This page intentionally left blank*

# List of Symbols

In what follows, the descriptions of the mathematical symbols are suffixed with the page numbers of their introductory use for reference.

## Functions

| | |
|---|---|
| $f$ | Real-valued objective function defined on $\mathbb{R}^n$. 1 |
| $c_i$ | Real-valued constraint function defined on $\mathbb{R}^n$, with $i \in \mathcal{I} \cup \mathcal{E}$. 2 |
| $\mathcal{L}$ | Lagrangian function. 9 |
| $\mathcal{L}_\mathsf{A}$ | Augmented Lagrangian function. 66 |
| $\varphi$ | Merit function ($\ell_2$-penalty function, unless otherwise stated). 69 |
| $\hat{f}$ | Quadratic model of $f$. 24 |
| $\widehat{\mathcal{L}}$ | Quadratic model of $\mathcal{L}$. 87 |
| $\widehat{\varphi}$ | Quadratic model of $\varphi$. 71 |

## Notations

| | |
|---|---|
| $\in$ | Set membership notation. 2 |
| $\subseteq$ | Set inclusion notation. 1 |
| $\square$ | Halmos symbol. 38 |
| $A^\mathsf{T}, v^\mathsf{T}$ | Transpose of a matrix or a vector. 3 |
| $I_n$ | Identity matrix on $\mathbb{R}^{n \times n}$. 3 |
| $e_i$ | Standard coordinate vector of $\mathbb{R}^n$ ($i$th column of $I_n$), with $1 \le i \le n$. 3 |
| $\mathcal{O}(\cdot)$ | Big-O notation. 13 |
| $\mathrm{o}(\cdot)$ | Little-O notation. 67 |
| $N(\mu, \sigma^2)$ | Gaussian distribution with mean $\mu$ and variance $\sigma^2$. 18 |

## Operators

| | |
|---|---|
| $[\cdot]_+$ | Positive-part operator. 47 |
| $[\cdot]_-$ | Negative-part operator. 98 |
| $\lvert\cdot\rvert$ | Modulus operator. 3 |
| $\lfloor\cdot\rfloor$ | Floor operator. 87 |
| $\lVert\cdot\rVert$ | Norm of a vector or a matrix (may be subscripted for clarity). 3 |
| $\nabla$ | Gradient operator (elements $\partial / \partial x_i$, with $i \in \{1, 2, \dots, n\}$). 8 |
| $\nabla^2$ | Hessian operator (elements $\partial^2 / \partial x_i \partial x_j$, with $i, j \in \{1, 2, \dots, n\}$). 11 |
| $\odot$ | Hadamard product. 31 |
| card | Cardinal function. 16 |
| rank | Rank operator. 51 |
| sgn | Sign operator. 5 |
| arg min | Global minimizer operator. 74 |
| arg max | Global maximizer operator. 102 |

## Sets

| | |
|---|---|
| $\emptyset$ | Empty set. 4 |
| $[a, b]$ | Closed set $\{x \in \mathbb{R} : a \leq x \leq b\}$ with $a \leq b$. 54 |
| $(a, b)$ | Open set $\{x \in \mathbb{R} : a < x < b\}$ with $a < b$. 16 |
| $(a, b]$ | Semi-open set $\{x \in \mathbb{R} : a < x \leq b\}$ with $a < b$. 75 |
| $\mathbb{R}$ | Set of real numbers. 5 |
| $\mathbb{R}^n$ | Real coordinate space of dimension $n$. 1 |
| $\mathbb{R}^{m \times n}$ | Real matrix space of dimension $m \times n$. 31 |
| $\mathscr{L}(\mathbb{R}^n)$ | Space of linear polynomials in $\mathbb{R}^n$. 25 |
| $\mathbb{Q}(\mathbb{R}^n)$ | Space of quadratic polynomials in $\mathbb{R}^n$. 25 |
| $\Omega$ | Feasible set, included in $\mathbb{R}^n$. 1 |
| $\mathcal{E}$ | Set of indices of the equality constraints. 2 |
| $\mathcal{I}$ | Set of indices of the inequality constraints. 2 |

# 1 Introduction

This thesis will develop model-based *derivative-free optimization* (DFO) methods and software. In this chapter, we provide an overview of DFO, introduce some applications, and present the concepts and tools that we will use throughout this thesis. In particular, Section 1.1 includes basic concepts of DFO and its motivation. Section 1.2 provides several examples of DFO problems from academic, engineering, and industrial applications. Section 1.3 introduces the optimality conditions for smooth nonlinear optimization, which will be used in the development of our algorithms. We then summarize existing methodologies for the design of DFO algorithms in Section 1.4. Finally, Section 1.5 presents two benchmarking tools that will be used to compare and assess the performance of DFO solvers.

## 1.1 Overview of derivative-free optimization

Optimization is the study of extremal points and values of mathematical functions. It aims at minimizing (or maximizing) a real-valued function $f$, referred to as the *objective function*, within a given set of points $\Omega \subseteq \mathbb{R}^n$, referred to as the *feasible set*. It is well known that essential information for optimization is embraced in the (possibly generalized) derivatives of the functions involved. However, in practice, evaluations of such derivatives may be unreliable or prohibitively expensive, if not impossible. It motivates the study of DFO [47, 7, 55, 123], where problems are solved using only function values. This thesis focuses on methods and software for DFO.

DFO problems arise naturally when the objective function or the feasible set results from complex experiments or simulations. Regarding these functions as black boxes, people often refer to those problems as *black-box optimization* (BBO) problems [7], which constitute a significant type of DFO problem in practice. Note that DFO differs from nonsmooth optimization [41, 53], which studies problems involving nonsmooth functions. In DFO, the major difficulty is not the possible nonsmoothness of the functions involved but the lack of knowledge about the structures of the problems. In theoretical analysis of DFO methods, it is not uncommon to assume that the underlying functions enjoy some smoothness, although algorithms cannot retrieve their (classical or generalized) derivatives. We emphasize that if any derivative information can be evaluated at an affordable cost or approximated well enough, DFO methods are not recommended, as they are very unlikely to outperform methods that use derivatives. Consider, for example, minimizing an objective function defined by a sophisticated simulation whose source code is available. One may then attempt to evaluate derivatives using automatic differentiation tools [97, 98] and apply derivative-based methods.

For DFO methods, the leading complexity measure we consider is the number of function evaluations. In practice, each function evaluation may require several minutes or even several days to complete [7, § 1.4]. For instance, a recent application of DFO is hyperparameter tuning in machine learning [82], for which every objective function evaluation necessitates training a machine learning model (see Section 1.2.3). Hence, in DFO methods, the expense of numerical linear algebra is less of a concern, although we will maintain it acceptable.

In this introduction, we consider the nonlinearly constrained problem

$$\min_{x \in \mathbb{R}^n} \quad f(x) \tag{1.1.1a}$$

$$\text{s.t.} \quad c_i(x) \le 0, \; i \in \mathcal{I}, \tag{1.1.1b}$$

$$c_i(x) = 0, \; i \in \mathcal{E}, \tag{1.1.1c}$$

where the *objective* and *constraint functions* $f$ and $c_i$, with $i \in \mathcal{I} \cup \mathcal{E}$, are real-valued functions on $\mathbb{R}^n$, and where the sets of indices $\mathcal{I}$ and $\mathcal{E}$ are finite (perhaps empty) and disjoint. The feasible set of this problem is

$$\Omega \stackrel{\text{def}}{=} \{x \in \mathbb{R}^n : c_i(x) \le 0, \; i \in \mathcal{I} \text{ and } c_i(x) = 0, \; i \in \mathcal{E}\}.$$

If $f$ is convex, while $c_i$ is convex for all $i \in \mathcal{I}$ and affine for all $i \in \mathcal{E}$, then problem (1.1.1) is *convex*. However, this thesis does *not* assume convexity.

We emphasize that (1.1.1b) and (1.1.1c) may include bound constraints. We do not extract them explicitly in this chapter, but they may need to be handled differently from other constraints because they often represent inalienable physical or theoretical restrictions. This will be considered by our new DFO method (see Chapter 5).

## 1.2 Examples of applications

### 1.2.1 Automatic error analysis

A typical example of DFO applications is automatic error analysis [105, 106], which formulates numerical computation's accuracies and stabilities using optimization problems. Consider, for instance, the Gaussian elimination with partial pivoting of a matrix $A \in \mathbb{R}^{n \times n}$, given in Algorithm 1.1, where the superscripts denote iteration numbers.

Wilkinson's backward error analysis (see, e.g., equation (25.14) of chapter 3 in [209], where $t$ is introduced at the beginning of Paragraph 10 and $g$ at the end of p. 97) demonstrates that the growth factor of the Gaussian elimination, defined as

$$\rho_n(A) \stackrel{\text{def}}{=} \frac{1}{\|A\|_{\max}} \max_{0 \le k \le n-1} \|A^{(k)}\|_{\max}, \tag{1.2.1}$$

**Algorithm 1.1:** Gaussian elimination with partial pivoting

**Data:** Matrix $A \in \mathbb{R}^{n \times n}$.

**Result:** Factorized matrix $A^{(n-1)} \in \mathbb{R}^{n \times n}$.

Initialize $A^{(0)} \leftarrow A$

**for** $k = 1, 2, \ldots, n-1$ **do**

  Determine a pivot index $j$ that solves max $\left\{ \left| A_{i,k}^{(k-1)} \right| : k \leq i \leq n \right\}$

  **if** $A_{j,k}^{(k-1)} \neq 0$ **then**

    Exchange the $k$th and the $j$th rows of $A^{(k-1)}$

    Evaluate the multiplier $\tau^k \in \mathbb{R}^n$ with components

$$\tau_i^k = \begin{cases} A_{i,k}^{(k-1)} / A_{k,k}^{(k-1)}, & \text{if } i > k, \\ 0, & \text{otherwise} \end{cases}$$

    Update $A^{(k)} \leftarrow (I_n - \tau^k e_k^\top) A^{(k-1)}$

  **else**

    Set $A^{(k)} \leftarrow A^{(k-1)}$

  **end if**

**end for**

determinates the numerical stability of Algorithm 1.1, where $\|\cdot\|_{\max}$ denotes the max norm of a matrix, i.e., the largest absolute value of the matrix's entries. More specifically, the $\ell_\infty$-norm of the backward error of the computed solution is bounded from above by a term proportional to $\rho_n(A)$. To study the worst-case scenario, we wish to determine how large $\rho_n$ can be and hence, to solve

$$\max_{A \in \mathbb{R}^{n \times n}} \rho_n(A). \tag{1.2.2}$$

Note that $\mathbb{R}^{n \times n}$ is isomorphic to $\mathbb{R}^{n^2}$; hence, problem (1.2.2) can be formulated as problem (1.1.1). Besides, although the growth factor is defined everywhere, it may not be continuous at the points yielding a tie in selecting the pivot element. Moreover, it is not differentiable at the points yielding a tie in any maximum operator in equation (1.2.1). Hence, optimization methods based on derivative information cannot be used for this problem. In such a case, DFO methods can help solve problem (1.2.2). Note that the optimal value and all local solutions to problem (1.2.2) are known [107], but DFO methods can be used to help the theoretical development [105].

### 1.2.2 Tuning nonlinear optimization methods

Another well-known example of DFO applications is the parameter tuning of nonlinear optimization methods [9]. For example, consider Algorithm 1.2, a basic trust-region method for solving the problem (1.1.1) when $\mathcal{I} = \mathcal{E} = \emptyset$, where $\|\cdot\|$ can be any norm.

**Algorithm 1.2:** Basic trust-region method for unconstrained optimization

```
    Data: Objective function f, initial guess x^0 ∈ ℝ^n, initial trust-region
          radius Δ^0 > 0, and parameters 0 < η_1 ≤ η_2 < 1 and 0 < θ_1 < 1 < θ_2.
1   for k = 0, 1, ... do
2       Define a simple function m_k such that m_k(d) ≈ f(x^k + d) for ||d|| ≤ Δ^k
3       Set the trial step d^k to an approximate solution to
            min_{d ∈ ℝ^n}  m_k(d)
            s.t.           ||d|| ≤ Δ^k
4       Evaluate the trust-region ratio
            ρ^k ← (f(x^k) − f(x^k + d^k)) / (m_k(0) − m_k(d^k))
5       if ρ^k ≥ η_1 then
6           Update the trial point x^{k+1} ← x^k + d^k
7       else
8           Retain the trial point x^{k+1} ← x^k
9       end if
10      Update the trust-region radius
            Δ^{k+1} ← { θ_1 Δ^k,  if ρ^k ≤ η_1,
                       { Δ^k,      if η_1 < ρ^k ≤ η_2,
                       { θ_2 Δ^k,  otherwise
11  end for
```

A critical simplification in Algorithm 1.2 lies in Line 5. A complete framework includes a parameter $\eta_0 \geq 0$ satisfying $\eta_0 \leq \eta_1$, and the condition in Line 5 is replaced by $\rho^k \geq \eta_0$. In practice, it is not uncommon to set $\eta_0 = 0$.

To choose the parameters $\eta_1$, $\eta_2$, $\theta_1$, and $\theta_2$, we minimize some measure of the expense of the method (e.g., the CPU time to solve a given set of problems), which will be denoted by $C$. In other words, we wish to solve

$$\min \quad C(\eta_1, \eta_2, \theta_1, \theta_2) \tag{1.2.3a}$$

$$\text{s.t.} \quad 0 \leq \eta_1 \leq \eta_2 < 1, \tag{1.2.3b}$$

$$0 < \theta_1 < 1 < \theta_2. \tag{1.2.3c}$$

Derivatives of $C$ cannot be evaluated, if they even exist. Moreover, observe that the objective function $C$ is likely to be undefined if the bound constraints are violated. This problem can be solved using DFO methods, such as the *mesh adaptive direct search* (MADS) method [9, 8].

Interestingly, DFO methods can be self-tuned following this methodology. The *brute force optimizer* (BFO) [151], a method for bound-constrained problems mixing continuous and discrete variables, is an example of self-tuned DFO methods.

### 1.2.3 Hyperparameter tuning in machine learning

A more recent example of DFO applications is hyperparameter tuning in machine learning [82]. Google solves for instance hyperparameter tuning problems using Google Vizier [88], the Google-internal service for performing black-box optimization. To illustrate this example, we consider the following hyperparameter tuning problem of a *support vector machine* (SVM) for binary classification. Given a binary-labeled data-set $\{(x_i, y_i)\}_{i=1,2,\dots,m} \subseteq \mathbb{R}^n \times \{\pm 1\}$, we build an SVM to classify the data with their respective labels. Binary classification is obtained using a *$C$-support vector classification* (SVC) [37] by solving the optimization problem

$$\begin{aligned}
\min_{(\omega,\beta,\xi)\in\mathbb{R}^\ell\times\mathbb{R}\times\mathbb{R}^m} \quad & \frac{1}{2}\|\omega\|_2^2 + C\|\xi\|_1 && \text{(1.2.4a)}\\
\text{s.t.} \quad & y_i\bigl(\beta + \omega^\mathsf{T}\psi_\gamma(x_i)\bigr) \ge 1 - \xi_i,\ i \in \{1, 2, \dots, m\}, && \text{(1.2.4b)}\\
& \xi \ge 0, && \text{(1.2.4c)}
\end{aligned}$$

where $\psi_\gamma$ is a function mapping the data to a higher-dimensional space $\mathbb{R}^\ell$ and $\gamma > 0$ and $C > 0$ are parameters. Given a solution $(\omega^*, \beta^*, \xi^*) \in \mathbb{R}^\ell \times \mathbb{R} \times \mathbb{R}^m$ to the problem (1.2.4), the $C$-SVC classifies any data $x \in \mathbb{R}^n$ according to

$$\delta(x) \overset{\text{def}}{=} \operatorname{sgn}\bigl(\beta^* + (\omega^*)^\mathsf{T}\psi_\gamma(x)\bigr), \tag{1.2.5}$$

which maps an observation $x \in \mathbb{R}^n$ to a label in $\{\pm 1\}$. It is clear that $\delta$ depends on the two parameters $C$ and $\gamma$, which can be chosen by solving an optimization problem. The objective function $P$ of this problem can be a 5-fold cross-validation based on some performance measure of model (1.2.5), presented in Algorithm 1.3.

Typical example of model's performances are the model's accuracy, i.e., the ratio of data correctly classified, and the *area under the ROC curve* (AUC) [102], particularly effective for imbalanced datasets [27]. The hyperparameter tuning problem is

$$\begin{aligned}
\min \quad & P(C, \gamma)\\
\text{s.t.} \quad & C, \gamma > 0,
\end{aligned}$$

It is clear that derivatives of the objective function of such a problem cannot be easily

**Algorithm 1.3:** $k$-fold cross-validation of an SVC with parameters $C$ and $\gamma$

**Data:** Labeled dataset $\{(x_i, y_i)\}_{i=1,2,\dots,m} \subseteq \mathbb{R}^n \times \{\pm 1\}$, fold number $k > 0$, and parameters $C$ and $\gamma$.

**Result:** Performance measure $P(C, \gamma)$.

```
Split the dataset into k balanced groups
for i = 1, 2, ..., k do
    Calculate (ω*, β*, ξ*) with the all the data except those in the ith group
    Evaluate the performance p_i of the model (1.2.5) on the data in the ith group
end for
Define P(C, γ) by summarizing the performances in {p_1, p_2, ..., p_k}
```

evaluated and may even not exist. This problem can be solved using DFO methods.

As in [176], DFO can also be applied to reinforcement learning. Instead of training a model on a fixed labeled dataset, reinforcement learning bases the training process on rewarding expected behaviors and punishing undesired ones. Hence, it often consists in finding optimal parameters that maximize a reward. However, these reward derivatives often cannot be evaluated, and DFO methods can be an approach to solving such problems. This concept is often referred to as derivative-free reinforcement learning.

### 1.2.4 Other industrial and engineering applications

DFO methods are also widely used in industry and engineering, especially for solving problems that involve heavy simulations. Such problems arise from helicopter rotor blade manufacturing [22, 23, 191], aeroacoustic shape design [131, 132], computational fluid dynamics [62], worst-case analysis of analog circuits [124], rapid-cycling synchrotron accelerator modeling [63], nuclear energy engineering [117, 119, 118], reservoir engineering and engine calibration [122], and groundwater supply and bioremediation engineering [76, 140, 215], to name but a few. In general, problems that involve sophisticated models, simulations, or experiments, induce DFO problems.

A particular application of DFO comes from *multi-disciplinary design optimization* (MDO) in the industry. It is a field that uses optimization methods to solve design problems defined by multiple disciplines. The objective and constraint functions of an MDO problem can be provided by different departments of the same company or even by different companies. It is the case in aircraft engineering [81], where the design problem of one component is solved while taking into account constraints imposed by other components handled by different departments. MDO problems often involve simulations or experiments, and therefore, DFO methods are often needed. In Chapter 3 of this thesis, we will present a piece of software we implemented for solving DFO problems based on methods by Powell [161, 163, 166, 168, 171]. It has been included in GEMSEO [78], an engine for MDO initiated by a team from IRT Saint Exupéry[1] in France.

## 1.3 Optimality conditions for smooth optimization

We discuss in this section optimality conditions for problem (1.1.1). We do not assume any structure on the objective and constraint functions, except some smoothness. More specialized results can be obtained by assuming that problem (1.1.1) is convex, for example, but it is out of the scope of this work.

### 1.3.1 Local and global solutions

Before solving problem (1.1.1), we must define what a solution is. Definition 1.3.1 presents the most natural understanding of a solution.

**Definition 1.3.1** (Global solution)**.** To problem (1.1.1), a point $x^* \in \mathbb{R}^n$ is a *global solution* if $x^* \in \Omega$ and $f(x) \geq f(x^*)$ for all $x \in \Omega$.

The following relaxed concepts of solutions are of interest in theory and practice.

**Definition 1.3.2** (Local solution)**.** To problem (1.1.1), a point $x^* \in \mathbb{R}^n$ is

- a *local solution* if $x^* \in \Omega$ and there exists a neighborhood $\mathcal{N} \subseteq \mathbb{R}^n$ of $x^*$ such that $f(x) \geq f(x^*)$ for all $x \in \mathcal{N} \cap \Omega$,
- a *strict local solution* if $x^* \in \Omega$ and there exists a neighborhood $\mathcal{N} \subseteq \mathbb{R}^n$ of $x^*$ such that $f(x) > f(x^*)$ for all $x \in \mathcal{N} \cap \Omega \setminus \{x^*\}$, and
- an *isolated local solution* if there exists a neighborhood $\mathcal{N} \subseteq \mathbb{R}^n$ of $x^*$ such that it is the only local solution in $\mathcal{N} \cap \Omega$.

It is known that finding a local solution to a general nonconvex problem is NP-hard. In many cases, it is already NP-hard to check whether a given point is a local solution. There also exist convex optimization problems that are NP-hard, such as copositive programming, for which checking the feasibility of a given point is indeed NP-hard. See [141] for more discussions.

The methods we consider in this thesis are local. They attempt to find approximate local solutions to problem (1.1.1). However, in general, theoretical analyses of these methods can only guarantee approximations of stationary points, which will be introduced hereafter.

[1] See `https://www.irt-saintexupery.com`.

### 1.3.2 Constraint qualifications

Before introducing any necessary and sufficient conditions for local optimality, we discuss some regularity conditions on the constraints (1.1.1b) and (1.1.1c), referred to as *constraint qualifications*. They will be required for the necessary conditions to hold. We first introduce the notion of *active sets*.

**Definition 1.3.3** (Active set)**.** The *active set* $\mathcal{A}(x) \subseteq \mathcal{I} \cup \mathcal{E}$ for problem (1.1.1) at a point $x \in \mathbb{R}^n$ is defined by

$$\mathcal{A}(x) \overset{\text{def}}{=} \mathcal{E} \cup \{i \in \mathcal{I} : c_i(x) \geq 0\}.$$

If a constraint belongs to the active set[2] at a given point, it is said to be *active* at this point and *inactive* otherwise. Note that a violated constraint is always considered active.

We introduce hereafter two classical constraint qualifications.

**Definition 1.3.4** (Constraint qualifications)**.** Given $x \in \Omega$, denote $\mathcal{A}(x)$ the active set for problem (1.1.1) at $x$, and assume that the constraint functions $c_i$ are differentiable at $x$ for all $i \in \mathcal{A}(x)$. We say that

- the *linear independence constraint qualification* (LICQ) holds at $x$ if the gradients $\nabla c_i(x)$ are linearly independent for all $i \in \mathcal{A}(x)$, and
- the *Mangasarian-Fromovitz constraint qualification* (MFCQ) holds at $x$ if the gradients $\nabla c_i(x)$ are linearly independent for all $i \in \mathcal{E}$ and there exists a vector $z \in \mathbb{R}^n$ such that

$$\begin{cases} \nabla c_i(x)^\mathsf{T} z < 0, \ i \in \mathcal{A}(x) \cap \mathcal{I}, & \text{(1.3.1a)} \\ \nabla c_i(x)^\mathsf{T} z = 0, \ i \in \mathcal{E}. & \text{(1.3.1b)} \end{cases}$$

The LICQ is stronger than the MFCQ. If the LICQ holds at $x \in \Omega$, then system (1.3.1) is consistent because of the linear independence of all $\nabla c_i(x)$ for $i \in \mathcal{A}(x)$.

Many other constraint qualifications exist. Examples include the *Abadie constraint qualification* (ACQ) and the *Guignard constraint qualification* (GCQ), which are formulated using tangent and linearized cones of the feasible set. There also exist several traditional constraint qualifications weaker than the MFCQ, such as the *constant rank constraint qualification* (CRCQ), the *constant positive linear dependence constraint qualification* (CPLD), or the *quasi-normality constraint qualification* (QNCQ). We also note that dedicated constraint qualifications may exist when the problem has a particular structure, such as the *Slater's condition* (SC) for convex problems.

[2]For simplicity, we do not distinguish a constraint from its index.

### 1.3.3 First-order optimality conditions

**Statement of the optimality conditions**

Let $\mathcal{L}$ be the *Lagrangian* of problem (1.1.1), defined by

$$\mathcal{L}(x,\lambda) \overset{\text{def}}{=} f(x) + \sum_{i\in\mathcal{I}\cup\mathcal{E}} \lambda_i c_i(x), \quad \text{for } x \in \mathbb{R}^n \text{ and } \lambda_i \in \mathbb{R} \text{ for } i \in \mathcal{I}\cup\mathcal{E},$$

where $\lambda = [\lambda_i]^{\mathsf{T}}_{i\in\mathcal{I}\cup\mathcal{E}}$.

**Theorem 1.3.1** (First-order necessary conditions [144, Thm. 12.1]). *Let $x^* \in \Omega$ be a local solution to (1.1.1), assume that the functions $f$ and $c_i$ are continuously differentiable in a neighborhood of $x^*$ for all $i \in \mathcal{I}\cup\mathcal{E}$, and that the LICQ holds at $x^*$. Then there exists a Lagrange multiplier $\lambda^* = [\lambda_i^*]^{\mathsf{T}}_{i\in\mathcal{I}\cup\mathcal{E}}$ with $\lambda_i^* \in \mathbb{R}$ for all $i \in \mathcal{I}\cup\mathcal{E}$ such that*

$$\begin{cases} \nabla_x \mathcal{L}(x^*,\lambda^*) = 0, & \text{(1.3.2a)} \\ c_i(x^*) \le 0,\ i \in \mathcal{I} & \text{(1.3.2b)} \\ c_i(x^*) = 0,\ i \in \mathcal{E} & \text{(1.3.2c)} \\ \lambda_i^* c_i(x^*) = 0,\ i \in \mathcal{I} & \text{(1.3.2d)} \\ \lambda_i^* \ge 0,\ i \in \mathcal{I} & \text{(1.3.2e)} \end{cases}$$

We present Theorem 1.3.1 with the LICQ as an example, but a similar conclusion can be established with other constraint qualifications, such as the MFCQ presented above (see, e.g., [144, p. 339]). The conditions (1.3.2) are commonly referred to as the *Karush-Kuhn-Tucker* (KKT) conditions [113, 120], and the pair $(x^*, \lambda^*)$ in Theorem 1.3.1 is referred to as a *KKT pair*. More specifically, condition (1.3.2a) is the *stationarity* condition, conditions (1.3.2b) and (1.3.2c) as the *primal feasibility* conditions, condition (1.3.2d) as the *complementary slackness* condition, and condition (1.3.2e) as the *dual feasibility* condition. A *first-order stationary point* is a point $x \in \mathbb{R}^n$ satisfying the KKT conditions (1.3.2). Such a point may not be a local solution.

**An illustration of the first-order optimality conditions**

We do not provide a proof of Theorem 1.3.1, but we illustrate graphically the main idea on the simple 2-dimensional example

$$\begin{aligned} \min_{x\in\mathbb{R}^2} \quad & f(x) = x_1 + x_2 & \text{(1.3.3a)} \\ \text{s.t.} \quad & c_1(x) = x_1^2 + x_2^2 - 2 \le 0, & \text{(1.3.3b)} \\ & c_2(x) = -x_2 \le 0, & \text{(1.3.3c)} \end{aligned}$$

whose solution is $x^* = [-\sqrt{2}, 0]^\mathsf{T}$. A graphical representation of problem (1.3.3) is given in Figure 1.1, where the white area represents the feasible set.

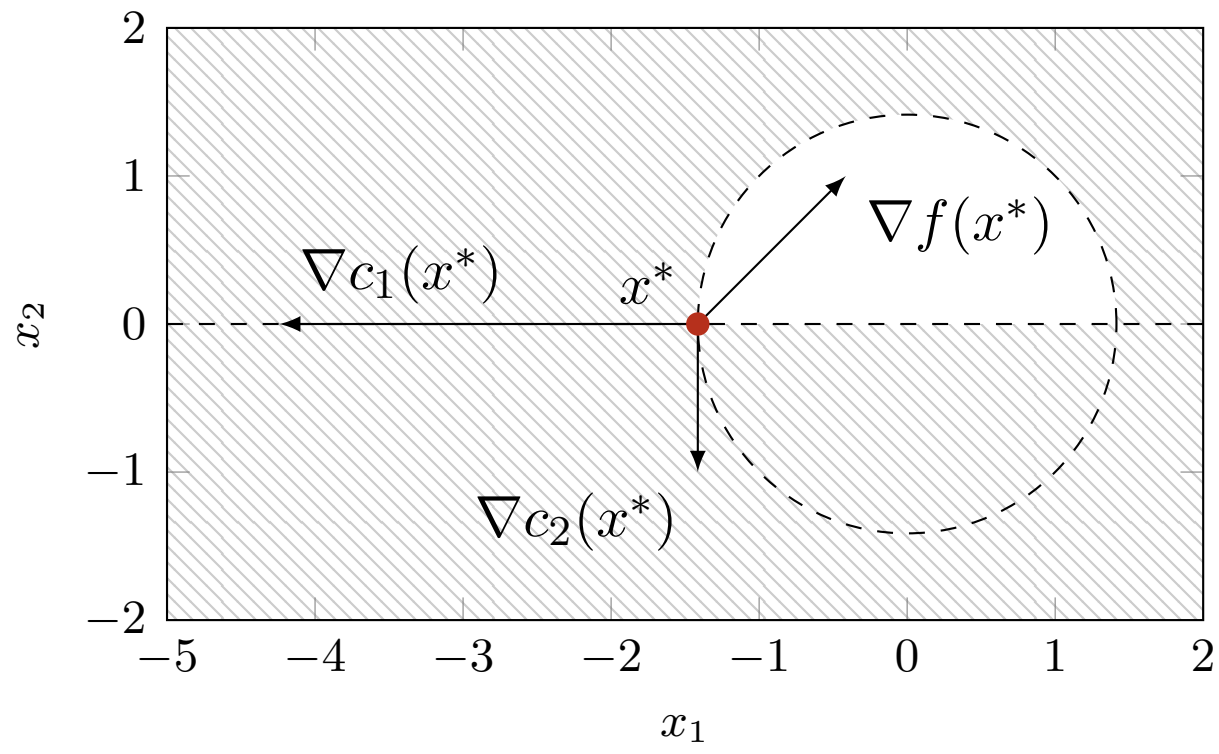

Figure 1.1: Graphical representation of the problem (1.3.3)

As manifest in Figure 1.1, there does not exist any direction $d \in \mathbb{R}^2$ satisfying

$$\begin{cases} \nabla f(x^*)^\mathsf{T} d < 0, & (1.3.4a) \\ \nabla c_1(x^*)^\mathsf{T} d \leq 0, & (1.3.4b) \\ \nabla c_2(x^*)^\mathsf{T} d \leq 0. & (1.3.4c) \end{cases}$$

The Farkas' lemma [64] ensures, therefore, that there exists a nonnegative Lagrange multiplier $\lambda^* = [\lambda_1^*, \lambda_2^*]^\mathsf{T}$ such that

$$\nabla f(x^*) + \lambda_2^* \nabla c_1(x^*) + \lambda_2^* \nabla c_2(x^*) = 0.$$

This validates the condition (1.3.2a), while (1.3.2b) to (1.3.2e) are conspicuous. A solution to (1.3.4) is both a descent direction for $f$ and a linearized feasible direction for the constraints. In the context of Theorem 1.3.1, the nonexistence of such a direction is guaranteed by the LICQ. See [144, § 12.4] for a complete proof of Theorem 1.3.1.

### 1.3.4 Second-order optimality conditions

It is known that at a local solution of a smooth unconstrained optimization problem, the gradient of the objective function is zero and its Hessian matrix is positive semidefinite. Theorem 1.3.2 generalizes this fact to problem (1.1.1).

**Theorem 1.3.2** (Second-order necessary conditions [144, Thm. 12.5])**.** *Let $x^* \in \Omega$ be a local solution to (1.1.1). Assume that the functions $f$ and $c_i$, with $i \in \mathcal{I} \cup \mathcal{E}$, are twice continuously differentiable in a neighborhood or $x^*$, and that the LICQ holds at $x^*$. Denote the active set for problem (1.1.1) at $x^*$ by $\mathcal{A}(x^*)$. Let $\lambda^* = [\lambda_i^*]_{i \in \mathcal{I} \cup \mathcal{E}}^\mathsf{T}$ be a Lagrange multiplier with $\lambda_i^* \in \mathbb{R}$ for all $i \in \mathcal{I} \cup \mathcal{E}$ satisfying the KKT conditions (1.3.2), and* *let $z \in \mathbb{R}^n$ be any vector such that*

$$\begin{cases} \nabla c_i(x^*)^\mathsf{T} z = 0, \ i \in \mathcal{E}, & (1.3.5a) \\ \nabla c_i(x^*)^\mathsf{T} z = 0, \ i \in \mathcal{A}(x^*) \cap \mathcal{I}, \ \lambda_i^* > 0, & (1.3.5b) \\ \nabla c_i(x^*)^\mathsf{T} z \geq 0, \ i \in \mathcal{A}(x^*) \cap \mathcal{I}, \ \lambda_i^* = 0. & (1.3.5c) \end{cases}$$

*Then $z^\mathsf{T} \nabla^2_{x,x} \mathcal{L}(x^*, \lambda^*) z \geq 0$.*

A first-order stationary point $x \in \Omega$ is called a *second-order stationary point* if it satisfies the conclusion of Theorem 1.3.2. We emphasize that a second-order stationary point may not be a local solution. However, it is known in smooth unconstrained optimization that a point at which the gradient of the objective function is zero and its Hessian matrix is positive definite is a strict local solution. Theorem 1.3.3 generalizes this fact to the optimization problem (1.1.1).

**Theorem 1.3.3** (Second-order sufficient conditions [144, Thm. 12.6])**.** *Let $x^* \in \Omega$ be a given point. Assume that the functions $f$ and $c_i$, with $i \in \mathcal{I} \cup \mathcal{E}$, are twice continuously differentiable in a neighborhood of $x^*$, and that at this point there exists a Lagrange multiplier $\lambda^* = [\lambda_i^*]^\mathsf{T}_{i \in \mathcal{I} \cup \mathcal{E}}$ with $\lambda_i^* \in \mathbb{R}$ for all $i \in \mathcal{I} \cup \mathcal{E}$ satisfying the KKT condition (1.3.2). If for all vector $z \in \mathbb{R}^n \setminus \{0\}$ satisfying the conditions (1.3.5) we have $z^\mathsf{T} \nabla^2_{x,x} \mathcal{L}(x^*, \lambda^*) z > 0$, then $x^*$ is a strict local solution to problem (1.1.1).*

## 1.4 Methodology of derivative-free optimization

As summarized in [47], two main strategies have been developed for solving DFO problems. One strategy consists of sampling the objective function around the current iterate and choosing the next iterate among the sampled points based on simple comparisons. The *direct-search* methods [116] are based on this framework. The other strategy builds iteratively models that approximate the problems (e.g., using polynomials) around the current iterate and choose the next iterate according to the approximated problems. These methods are referred to as *model-based* methods. Hybrid methods also exist, such as implicit filtering [114] and SID-PSM [54], which combine direct search and models.

We remark that there exist methods that can solve DFO problems, but that are not covered by the above categories. Examples include random search methods [221], simulated annealing methods [115], the genetic algorithm [57, 108], and Bayesian optimization methods [135, 192]. For more information on these methods, we refer to the review [123] and the reference therein.

### 1.4.1 Direct-search methods

An early example of DFO is a method from Fermi and Metropolis [66], who developed in 1952 a nonlinear least-squares solver on MANIAC, a computer based on the von Neumann architecture. From a modern viewpoint, this method is a coordinate search method, a particular example of direct-search methods, where the search directions are defined to be the coordinate axes.

Several direct-search methods appeared after that. Rosenbrock designed a direct-search method for unconstrained problems. It constructs an orthogonal basis using previous steps and searching along the directions in this basis [184]. Another famous early direct-search method is the Hooke-Jeeves method [109]. This method combines exploratory moves along coordinates axes with pattern moves, to exploit the pattern revealed by previous successful directions. In 1965, Nelder and Mead [142] introduced the simplex method[3], which evaluates the objective function at the vertices of a simplex, and updates this simplex according to these function values, one vertex at each iteration. It is arguably the most widely used DFO method, available as the `fminsearch` function in MATLAB, and many variations exist [213]. Many other works on direct search appeared in the same period, and summaries can be found in [67, 26]. Nowadays, the direct-search paradigm offers an abundance of algorithms [116], such as the *generalized pattern search* (GPS) methods [24], later extended to the MADS methods [5, 1, 2, 6, 126]. A recent example of direct-search methods is the BFO [151, 152], which handles integer and categorical variables and can be self-tuned.

Gratton et al. recently proposed incorporating stochastic strategies in direct-search methods [93, 95]. This improves both the performance and the worst-case complexity bounds compared with traditional direct-search methods.

### 1.4.2 Model-based methods

Unlike direct-search methods, model-based methods approximate locally the functions involved in the optimization problems by simple functions called *models* or *surrogates*. To make model-based methods globally convergent, the models are exploited by a globalization strategy, such as a *line-search* framework [144, Ch. 3] or a *trust-region* framework [44, 218]. Most of the existing model-based methods use linear or quadratic approximations [161, 45, 46], although other models have also been successfully used, such as *radial basis functions* (RBFs) [170, 145].

Model-based methods are highly appealing in practice, as they provide excellent performances in real applications. Compared with direct-search methods, the information contained in the function values is better exploited due to the use of models. The first trust-region DFO method was developed by Winfield [211, 212] in 1969. It is also regarded as the first trust-region method with or without derivatives [44, § 1.2]. A similar method was later developed by Powell, namely UOBYQA [163]. Both methods use fully-determined quadratic interpolation models (see Section 2.2 for detailed discussions). By updating the Lagrange functions of the interpolation problem, UOBYQA needs only $\mathcal{O}(n^4)$ computer operations to obtain each quadratic interpolant, whereas the complexity is $\mathcal{O}(n^6)$ in Winfield's method.

[3] Note that it is different from the simplex method in linear programming.

The methods of both Winfield and Powell require $\mathcal{O}(n^2)$ function evaluations to establish each quadratic model. Such an amount of function values are needed to initialize the methods, although most of them will be reused at subsequent iterations. However, this amount of function evaluations is not scalable to moderately large problems. It motivates methods that use underdetermined quadratic interpolation models. Such models typically require $\mathcal{O}(n)$ function values to be built, and the remaining freedom bequeathed by the interpolation conditions is taken up by minimizing a functional that reflects the regularity of the model. For example, the method of Conn and Toint [51] uses such models with the least $\ell_2$-norm of the coefficients. Powell's methods, namely NEWUOA [166], BOBYQA [168], and LINCOA [171], use models that minimize the Frobenius norm of the change to their Hessian matrices. Another example of such algorithms is MNH [207], whose models are underdetermined quadratic interpolants with the least Frobenius norm of their Hessian matrices.

Particular care must be given to the geometry of the interpolation set in order to maintain reasonable accuracies of the models (see Section 2.3). Model-based methods usually make explicit geometry-improving steps when necessary [45, 46]. The Wedge method [130] of Marazzi and Nocedal, however, does not include such steps. Instead, it adds a so-called wedge constraint to the trust-region subproblems to prevent the trial steps from lying in a region that is likely to worsen the geometry of the interpolation set. The method of Fasano, Morales, and Nocedal [65] does not include geometry-improving steps either but performs surprisingly well. The authors conjecture that a self-correcting mechanism may be at play and prevents the geometry from deteriorating. Scheinberg and Toint [185] point out that the geometry improvement cannot be dispensed in general, but a slight modification of the algorithm in [65] does enjoy a self-correcting mechanism that guarantees the convergence without taking explicit geometry-improving steps.

There are many other trust-region DFO methods, such as CSV2 [17], which determines its models by regression and not interpolation. Examples of solvers that use nonpolynomial models include ORBIT [208], CONORBIT [177], and BOOSTERS [146], which use cubic RBF. There also exist methods for more specific problems, such as DFLS [219] and DFO-LS [33], aiming at solving nonlinear least-squares DFO problems.

Randomization is also exploited to improve the performance of trust-region DFO methods. This idea was first proposed by Bandeira, Scheinberg, and Vicente [10], who established the global convergence of a trust-region method based on random mod-els. They only require approximating the objective function well enough with a certain probability, and such models can be obtained by interpolation on randomly selected points [11]. The global convergence rate of this method is established by Gratton et al. [94], using an idea elaborated earlier in [93, § 6]. Similar results are established independently by Cartis and Scheinberg [34] but for a more general class of methods. The work of Bandeira, Scheinberg, and Vicente has motivated many investigations on methods that use randomized models, such as [38].

### 1.4.3 Comments on methods based on finite differences

Perhaps the most straightforward approach to solving DFO problems is to make finite-difference approximations of the derivatives and then employ derivative-based methods. Here, we briefly discuss the advantages and disadvantages of this approach.

The $i$th coordinate of the gradient of a smooth function $f : \mathbb{R}^n \to \mathbb{R}$ at a point $x \in \mathbb{R}^n$ can be approximated by the forward difference

$$\frac{\partial f}{\partial x_i}(x) = \frac{f(x + he_i) - f(x)}{h} + \mathcal{O}(h), \tag{1.4.1}$$

or the central difference

$$\frac{\partial f}{\partial x_i}(x) = \frac{f(x + he_i) - f(x - he_i)}{2h} + \mathcal{O}(h^2), \tag{1.4.2}$$

where $h > 0$ is the difference parameter, $e_i \in \mathbb{R}^n$ is the $i$th standard coordinate vector of $\mathbb{R}^n$, and the order of precision require standard assumptions (see, e.g., [144, § 8.1]).

Methods based on finite differences have several advantages. Firstly, they are relatively easy to implement, whereas the implementation of many commonly used DFO methods is highly nontrivial and challenging[4]. Moreover, there is a profusion of well-established derivative-based methods that can be explored with finite-difference approximations. Finally, we highlight that the function evaluations needed by the approximated derivatives in (1.4.1) and (1.4.2) can be straightforwardly parallelized, making the resulting algorithm scalable in many situations.

Meanwhile, such methods come with disadvantages as well. First of all, if the problem is noisy, it is nontrivial to choose the difference parameter $h$. From a theoretical standpoint, the optimal choice for the difference parameter depends on the noise level and the Lipschitz constants of the derivatives of the function $f$ (see, e.g., [144, § 8.1] and [194, Eqs. (2.13) and (2.14)]). There exist procedures to estimate these quantities (see, e.g., [139, § 3] and [194, Proc. I]). However, these procedures are often costly in terms of function evaluations in practice. In addition, it is difficult to reuse the function values in methods based on finite differences. A finite difference at $x \in \mathbb{R}^n$ requires the function to be evaluated on a mesh of $\mathcal{O}(n)$ points around $x$ oriented along the coordinate directions. When $x$ changes, most of the mesh points change completely, and hence, the function needs to be evaluated at another batch of $\mathcal{O}(n)$ points.

[4] For example, Powell [166] wrote "The development of NEWUOA has taken nearly three years. The work was very frustrating [...]"

There is no sharp division between methods based on finite difference and those based on interpolation models. Indeed, the forward difference (1.4.1) produces the gradient of the linear function that interpolates $f$ at

$$\{x, x + he_1, x + he_2, \dots, x + he_n\}.$$

Similarly, the central difference (1.4.2) generates the gradient of the unique quadratic function with diagonal Hessian matrix that interpolates $f$ at

$$\{x, x + he_1, x - he_1, x + he_2, x - he_2, \dots, x + he_n, x - he_n\}.$$

Interpolation can be regarded as a generalization of finite difference because there is more freedom in choosing the interpolation set. This freedom allows us to reuse most of the interpolation points, which is essential for the efficiency of model-based methods.

Section 1.5 will present some simple numerical experiments after introducing the benchmarking tools we use in this thesis for comparing DFO solvers. These results will demonstrate that the methods based on finite differences are highly sensitive to noise. Unsurprisingly, the methods suffer if the difference parameter is not adapted to the noise level.

## 1.5 Benchmarking tools for derivative-free optimization methods

We introduce in this section the benchmarking tools that we will use throughout this thesis to compare DFO solvers. We use the performance and data profiles [60, 138], developed by Dolan, Moré, and Wild. Most of the information in this section can be found in the aforementioned articles.

### 1.5.1 Expense measure and convergence test

We denote by $\mathcal{S}$ a set of solvers to benchmark and $\mathcal{P}$ a set of test problems, assumed to represent the problems for which the solvers have been designed. Let $t_{p,s}$ be the expense for the solver $s \in \mathcal{S}$ to achieve a given convergence test on the problem $p \in \mathcal{P}$. As mentioned in Section 1.1, the major cost of DFO in practice is the function evaluations. Therefore, in a DFO context, $t_{p,s}$ measures the number of function evaluations required by $s \in \mathcal{S}$ to solve $p \in \mathcal{P}$ up to the convergence test. When $s$ fails to satisfy the convergence test for $p$ within a given budget (e.g., a maximal number of function evaluations), we define $t_{p,s} = \infty$ by convention.

In this thesis, the numerical experiments select $\mathcal{P}$ from the CUTEst set [90]. We consider the following convergence test for DFO solvers, following [138, § 2]. We define for each problem $p \in \mathcal{P}$ a merit function $\varphi_p$, i.e., a function that measures the quality of a point, taking into account the values of both the constraint and the objective functions. The smaller value of $\varphi_p$, the better. Let $x_p^0$ be the initial guess for a given problem $p \in \mathcal{P}$, and $\varphi_p^*$ be the least value of $\varphi_p$ obtained by the solvers in $\mathcal{S}$. Given a tolerance $\tau \in (0, 1)$, a point $x$ satisfies the convergence test if

$$\varphi_p(x) \le \varphi_p^* + \tau[\varphi_p(x_p^0) - \varphi_p^*], \tag{1.5.1}$$

in which case we say that $x$ solves problem $p$ up to the tolerance $\tau$. This convergence test can be interpreted as follows. A point $x$ satisfies the test if the reduction $\varphi_p(x_p^0) - \varphi_p(x)$ is at least $1 - \tau$ times the maximal reduction $\varphi_p(x_p^0) - \varphi_p^*$ achieved by all solvers in $\mathcal{S}$. We say that a solver solves the problem $p$ up to the tolerance $\tau$ whenever it produces an iterate that achieves the convergence test (1.5.1).

In the convergence test (1.5.1), it is tempting to define $\varphi_p^*$ as the merit function value at a minimizer of $p$, which is the case when testing derivative-based solvers [60]. However, according to [138], this may not be appropriate in DFO because it may happen that no solver in $\mathcal{S}$ achieves the test within the given computational budget if the function evaluations are expensive.

### 1.5.2 Performance profile

Fix a tolerance $\tau \in (0, 1)$ in the convergence test (1.5.1). For a solver $s \in \mathcal{S}$ and a problem $p \in \mathcal{P}$, define the *performance ratio* by

$$r_{p,s} \overset{\text{def}}{=} \frac{t_{p,s}}{\min\{t_{p,u} : u \in \mathcal{S}\}}, \tag{1.5.2}$$

which is the *relative* expense for $s$ to solve $p$ compared with the most efficient solver in $\mathcal{S}$ for this problem. The *performance profile* of $s$ is defined as

$$\rho_s(\alpha) \overset{\text{def}}{=} \frac{1}{\operatorname{card}(\mathcal{P})} \operatorname{card}(\{p \in \mathcal{P} : r_{p,s} \le \alpha\}), \quad \text{for } \alpha \ge 1,$$

where $\operatorname{card}(\cdot)$ denotes the cardinal number of a set. Clearly, $\rho_s(\alpha)$ is the proportion of problems in $\mathcal{P}$ that are solved by $s$ with a performance ratio at most $\alpha$. It can also be interpreted as the probability for the solver $s$ to solve a random problem from $\mathcal{P}$ under the restriction on the performance ratio. In particular, $\rho_s(1)$ is the proportion of problems that $s$ solves faster than any other solver in $\mathcal{S}$. Meanwhile,

$$\lim_{\alpha \to \infty} \rho_s(\alpha)$$

is the proportion of problems that are solved by $s$ (within the budget restriction). Given two solvers $s_1$ and $s_2$, $\rho_{s_1}(\alpha) > \rho_{s_2}(\alpha)$ means that $s_1$ solves more problems than $s_2$ under the constraint that $r_{p,s_1} \leq \alpha$ and $r_{p,s_2} \leq \alpha$ for $p \in \mathcal{P}$. Therefore, a larger value of $\rho_s$ indicates a better performance of $s$.

### 1.5.3 Data profile

We now introduce the data profile, another benchmarking tool proposed by [138]. The *data profile* of a solver $s \in \mathcal{S}$ is defined as

$$d_s(\alpha) \overset{\text{def}}{=} \frac{1}{\operatorname{card}(\mathcal{P})} \operatorname{card}\left(\left\{p \in \mathcal{P} : \frac{t_{p,s}}{n_p + 1} \leq \alpha\right\}\right), \quad \text{for } \alpha \geq 0, \tag{1.5.3}$$

where $n_p$ is the dimension of the problem $p$. Therefore, $d_s(\alpha)$ is the proportion of problems in $\mathcal{P}$ that are solved by $s$ with at most $\alpha(n_p+1)$ function evaluations. It can also be interpreted as the probability for the solver $s$ to solve a random problem from $\mathcal{P}$, under the budget described above. In particular, $d_s(0) = 0$ and, the same as $\rho_s$,

$$\lim_{\alpha \to \infty} d_s(\alpha)$$

is the proportion of problems that are solved by $s$. In equation (1.5.3), the denominator $n_p + 1$ is a unit cost that serves to normalize $t_{p,s}$, so that the computational expenses for different problems are comparable even if their dimensions are quite different. Note that $n_p + 1$ is the number of function evaluations needed for evaluating a simplex gradient [25]. Therefore, $d_s(\alpha)$ is the proportion of problems solved by $s$ within a budget equivalent to $\alpha$ simplex gradient estimates. Given two solvers $s_1$ and $s_2$, $d_{s_1}(\alpha) > d_{s_2}(\alpha)$ means that $s_1$ solves more problems than $s_2$ using at most $\alpha(n_p+1)$ function evaluations for each $p \in \mathcal{P}$. Therefore, a larger value of $d_s$ indicates a better performance of $s$.

### 1.5.4 An illustrative example

As an example, we now compare three solvers using performance and data profiles. These solvers, constituting the set $\mathcal{S}$, are

1. NEWUOA (see Section 3.2.4), a model-based DFO method by Powell [166]; and
2. BFGS and CG, two gradient-based solvers provided by SciPy 1.0 [204]. When no derivatives are provided, they use forward finite-difference to approximate gradients, with the difference parameter $h = \sqrt{u}$, where $u$ is the unit roundoff.

The set $\mathcal{P}$ contains $154$ smooth unconstrained CUTEst problems of dimension at most $50$, and the convergence tolerance in (1.5.1) is $\tau = 10^{-3}$. Since these problems are unconstrained, the merit function $\varphi_p$ of each problem $p \in \mathcal{P}$ is set to be the objective function $f_p$. Therefore, the convergence test (1.5.1) becomes

$$f_p(x) \leq f_p^* + \tau[f_p(x_p^0) - f_p^*]. \tag{1.5.4}$$

We conduct two experiments as follows.

1. The first experiment is made without modifying the problems in $\mathcal{P}$. In (1.5.4), $f_p^*$ is set to the least value of $f_p$ obtained by all solvers in $\mathcal{S}$. The value of $t_{p,s}$ is defined as the number of function evaluations that the solver $s \in \mathcal{S}$ needs to solve the problem $p \in \mathcal{P}$ up to the tolerance $\tau$.

2. The second experiment is a noisy variation of the previous one. Let $\sigma > 0$ be the noise level. For each problem $p \in \mathcal{P}$, the objective function is evaluated by

$$\tilde{f}_p(x) \stackrel{\text{def}}{=} [1 + \epsilon(x)] f_p(x), \tag{1.5.5}$$

where $\epsilon(x) \sim N(0, \sigma^2)$. Each problem is solved ten times with each solver. In the convergence test (1.5.4), $f_p^*$ is either the least value of $f_p$ obtained by all solvers during these ten runs, or the value obtained in the previous noise-free experiment, whichever is smaller. The value of $t_{p,s}$ is set to the average number of function evaluations that $s \in \mathcal{S}$ needs to solve $p \in \mathcal{P}$ up to the tolerance $\tau$. Note that the convergence test (1.5.4) uses the values of $f_p$ and not those of $\tilde{f}_p$. This means that we evaluate the solvers according to the true objective function, even though the solvers can access only the noisy values produced by (1.5.5).

In these experiments, we invoke NEWUOA via the PDFO package (see Chapter 3) under the default settings. In particular, the initial trust-region radius is $1$, the final trust-region radius is $10^{-6}$, and the number of interpolation points is $2n + 1$, with $n$ being the dimension of the problem being solved. BFGS and CG are called with the default configurations in SciPy 1.0. For each testing problem, the starting point is set to the one provided by CUTEst, and the maximal number of function evaluations is set to $500$ times the number of variables.

### Performance profiles

Figure 1.2 plots the performance profiles obtained by these experiments. In general, the best solver is indicated by the highest curve. Note that $\alpha$ is displayed in $\log_2$-scale. Figure 1.2a corresponds to the noise-free experiment and Figures 1.2b to 1.2d correspond to the noisy experiment, with $\sigma$ being $10^{-10}$, $10^{-8}$, and $10^{-6}$, respectively.

We can make the following observations on the noise-free experiment in Figure 1.2a.

1. BFGS solves more than $90\,\%$ of the problems, while the two other solvers solve slightly fewer problems. According to [138, § 2], BFGS is said more reliable than the other two solvers on the test problems in $\mathcal{P}$, but the margin is narrow.

2. NEWUOA uses the least number of function evaluations on about $50\,\%$ of the problems, while the percentages for BFGS and CG are about $30\,\%$ and $25\,\%$ respectively.

We observe the following for the noisy experiments according to Figures 1.2b to 1.2d.

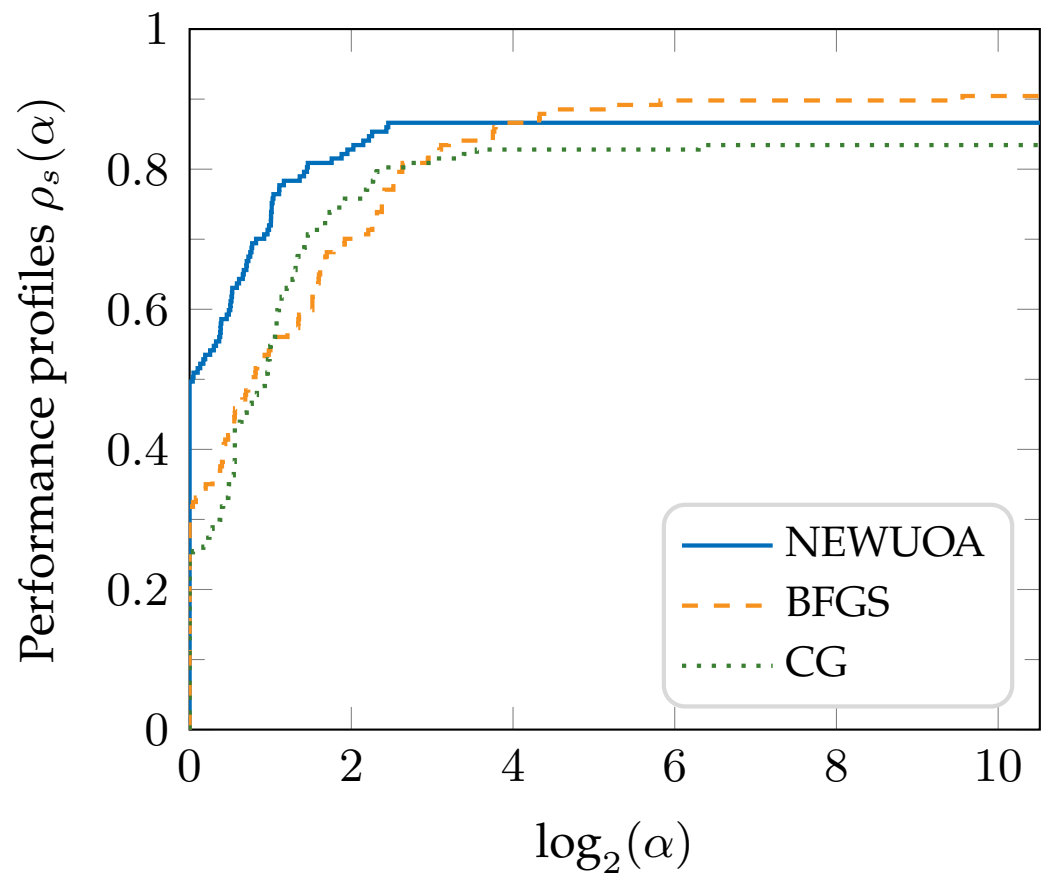

(a) Noise-free experiment

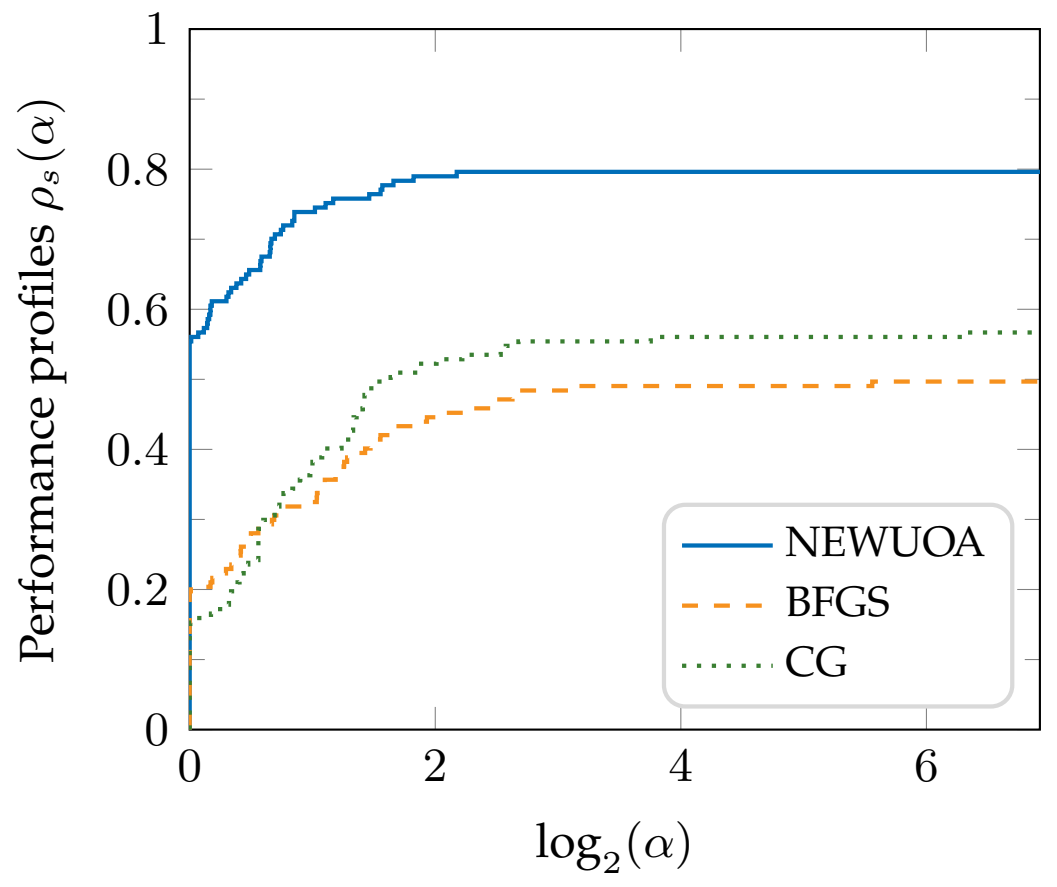

(b) Noisy experiment with $\sigma = 10^{-10}$

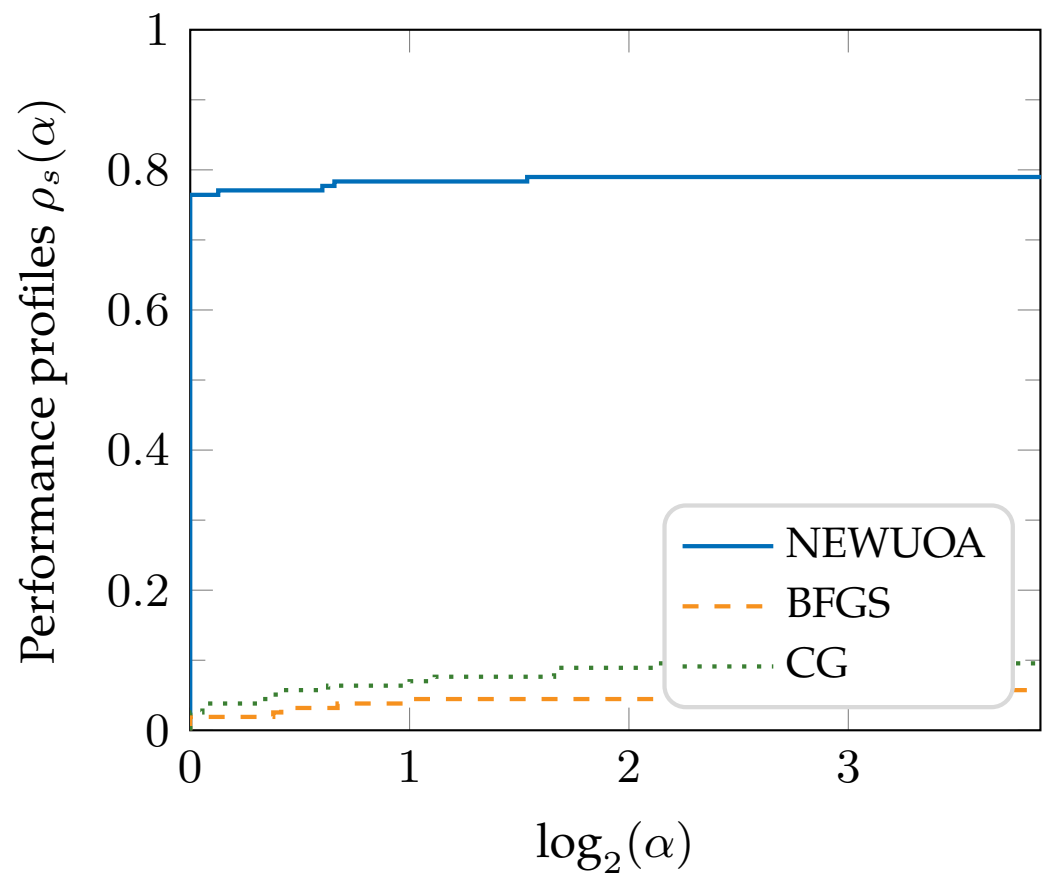

(c) Noisy experiment with $\sigma = 10^{-8}$

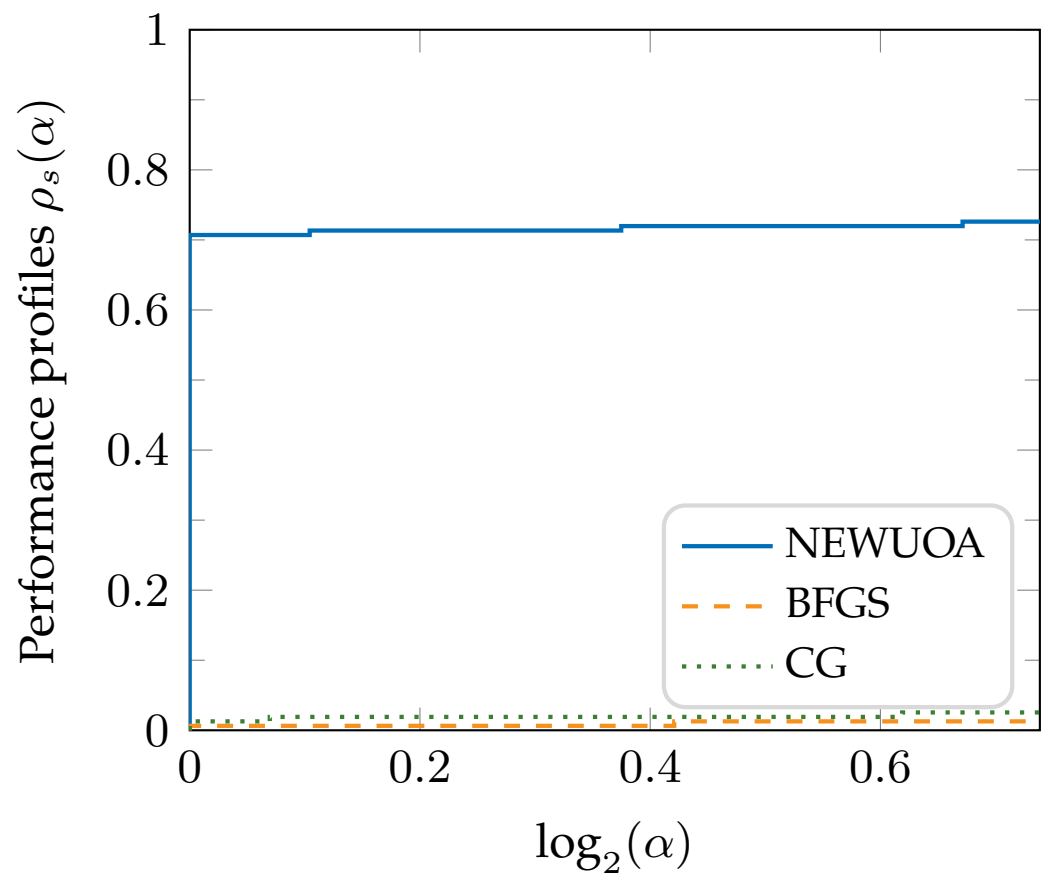

(d) Noisy experiment with $\sigma = 10^{-6}$

Figure 1.2: Performance profiles of NEWUOA, BFGS, and CG with $\tau = 10^{-3}$

1. NEWUOA always solves more than 70 % of the problems while the two other solvers solve significantly fewer problems.
2. When $\sigma = 10^{-10}$, NEWUOA uses the least number of function evaluations on more than 55 % of the problems, while the percentages for BFGS and CG are about 20 % and 15 % respectively. For $\sigma \in \{10^{-8}, 10^{-6}\}$, NEWUOA is the most efficient on almost all the problems that it can solve, while the other two solvers fail on most of the problems in $\mathcal{P}$.

Overall, on the noise-free problems, BFGS is slightly more reliable than NEWUOA and CG, although NEWUOA is more efficient on about half of the problems. In contrast, NEWUOA outperforms the two other solvers on noisy problems. The performances of BFGS and CG deteriorate significantly when Gaussian noise is imposed as in (1.5.5), even though the noise level is not high and the convergence tolerance is not demanding.

Note that our results do not contradict the observations of Shi et al. [194], where they choose the difference parameter $h$ more carefully, according to the noise levels and the smoothness of the problems. We rather keep the default value provided for the solvers BFGS and CG in SciPy 1.0 [204]. As commented in [194], the performance of methods based on finite differences is encouraging when there is no noise, yet much more care is needed when the problems are noisy. NEWUOA performs robustly with the noise according to this experiment, which agrees with the observations in [194].

We also mention that no conclusive comparison can be made between BFGS and CG according to Figures 1.2b to 1.2d. In fact, as pointed out in [91], when there are more than two solvers to compare, performance profiles have limitations when ranking the solvers, except for the one that outperforms the others (if any), such as NEWUOA in Figures 1.2b to 1.2d.

### Data profiles

Figure 1.3 displays the data profiles obtained in these experiments. A higher curve indicates a better performance. Figure 1.3a corresponds to the noise-free experiment, and Figures 1.3b to 1.3d correspond to the noisy experiment, with $\sigma$ being $10^{-10}$, $10^{-8}$, and $10^{-6}$, respectively.

We can make the following observations on the noise-free experiment on Figure 1.3a.

1. BFGS solves more than $90\,\%$ of the problems while the two other solvers solve slightly fewer problems. It agrees with what is indicated by Figure 1.2a, as it should be the case according to theory.

2. BFGS uses less than $300(n_p + 1)$ function evaluations to solve most of the problems $p \in \mathcal{P}$ that it can solve, while NEWUOA and CG use about $400(n_p + 1)$.

We observe the following for the noisy experiments according to Figures 1.3b to 1.3d.

1. NEWUOA always solves more than $70\,\%$ of the problems while the two other solvers solve significantly fewer problems. It agrees with what is indicated by Figures 1.2b to 1.2d, as it should be the case.

2. When $\sigma = 10^{-10}$, NEWUOA uses less than $200(n_p + 1)$ function evaluations to solve most of the problems $p \in \mathcal{P}$ that it can solve, while BFGS and CG use less than $100(n_p+1)$, although they solve much fewer problems. For $\sigma \in \{10^{-8}, 10^{-6}\}$, NEWUOA uses less than $100(n_p + 1)$ function evaluations[5] to solve most of the problems that it can solve, while the other two solvers fail on most of the problems.

[5] Even though the corresponding value is $200(n_p + 1)$ when $\sigma = 10^{-10}$, we should not interpret that NEWUOA performs better with higher noise. When $\sigma$ is larger, NEWUOA solves fewer problems and $f_p^*$ in (1.5.4) is often larger, leading to a weaker convergence test.

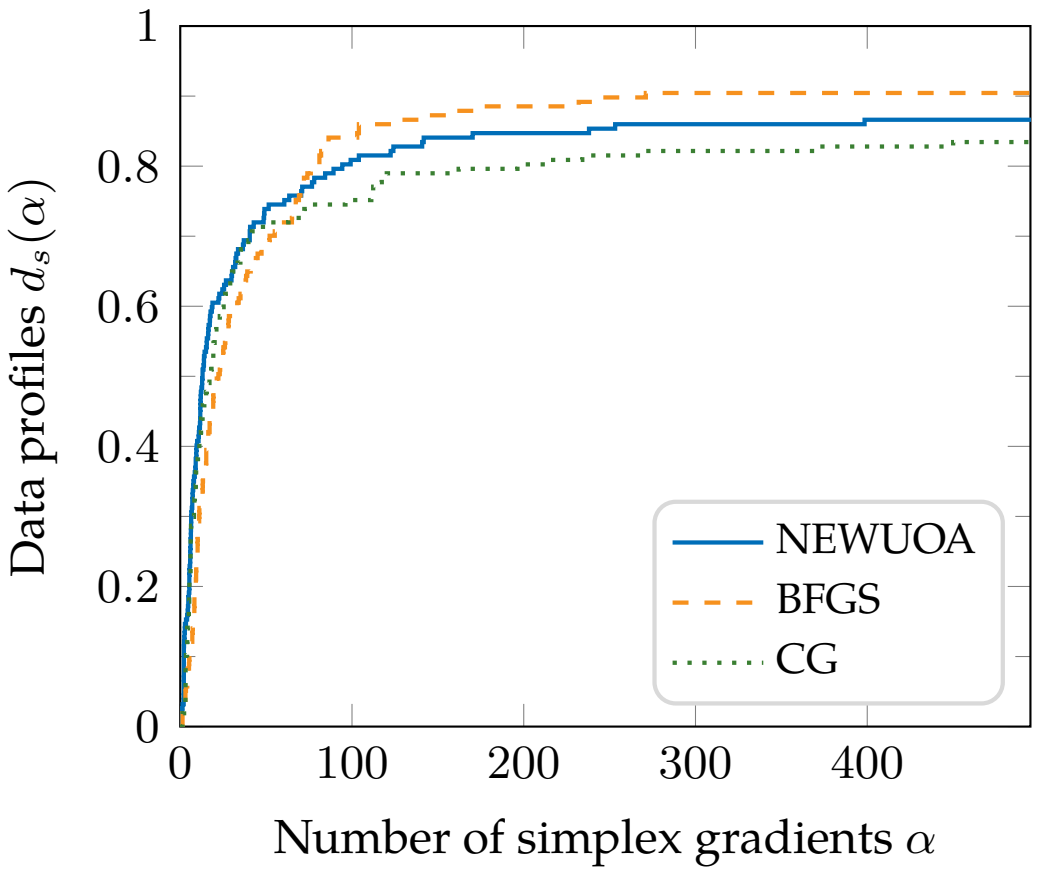

(a) Noise-free experiment

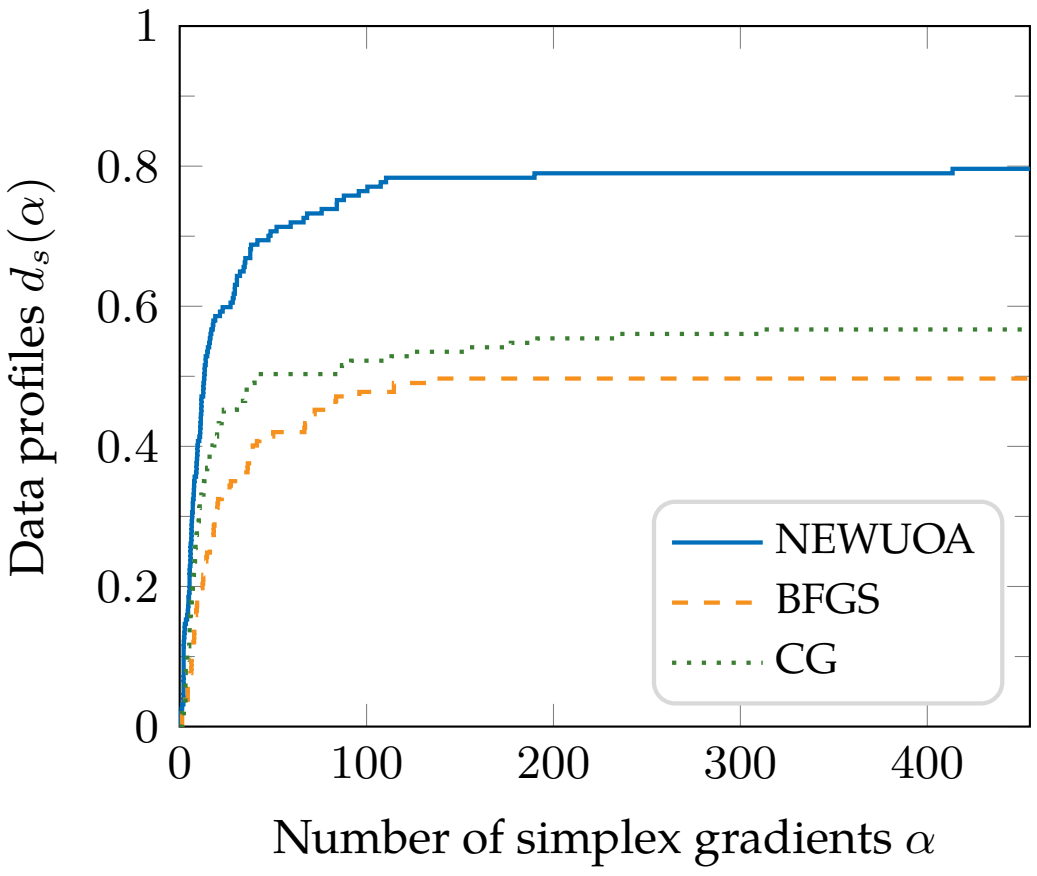

(b) Noisy experiment with $\sigma = 10^{-10}$

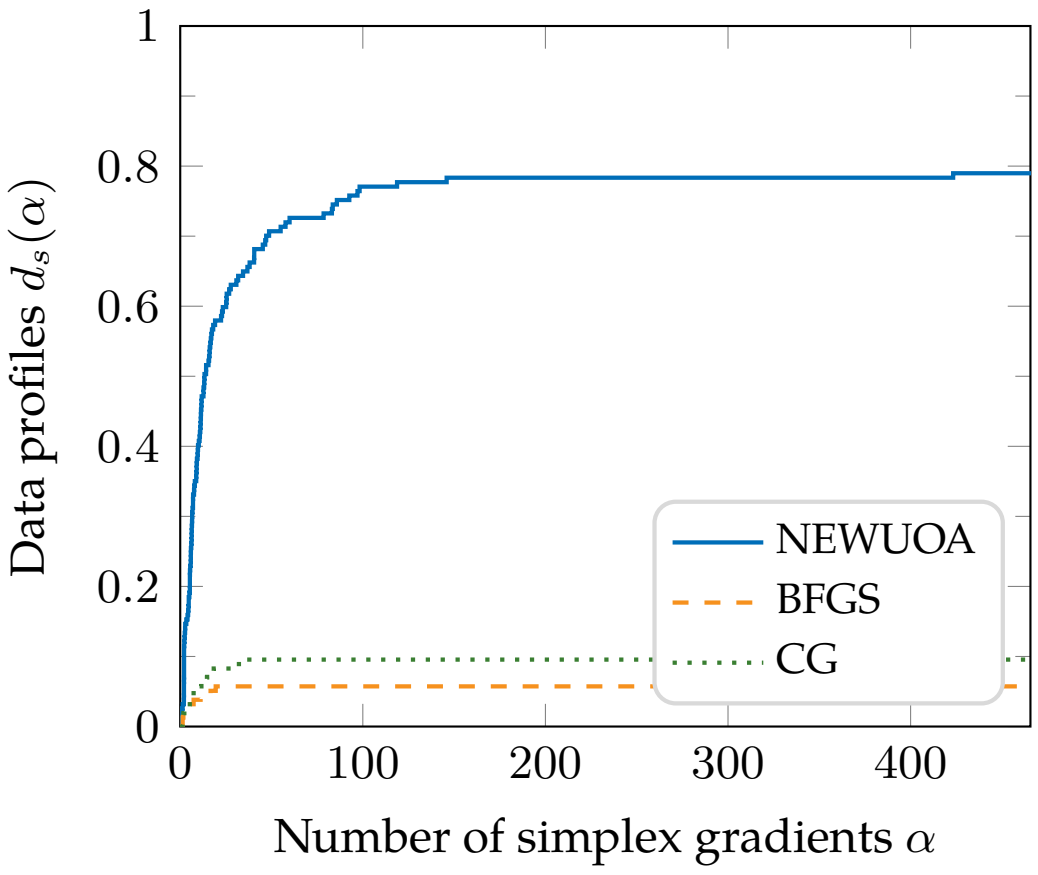

(c) Noisy experiment with $\sigma = 10^{-8}$

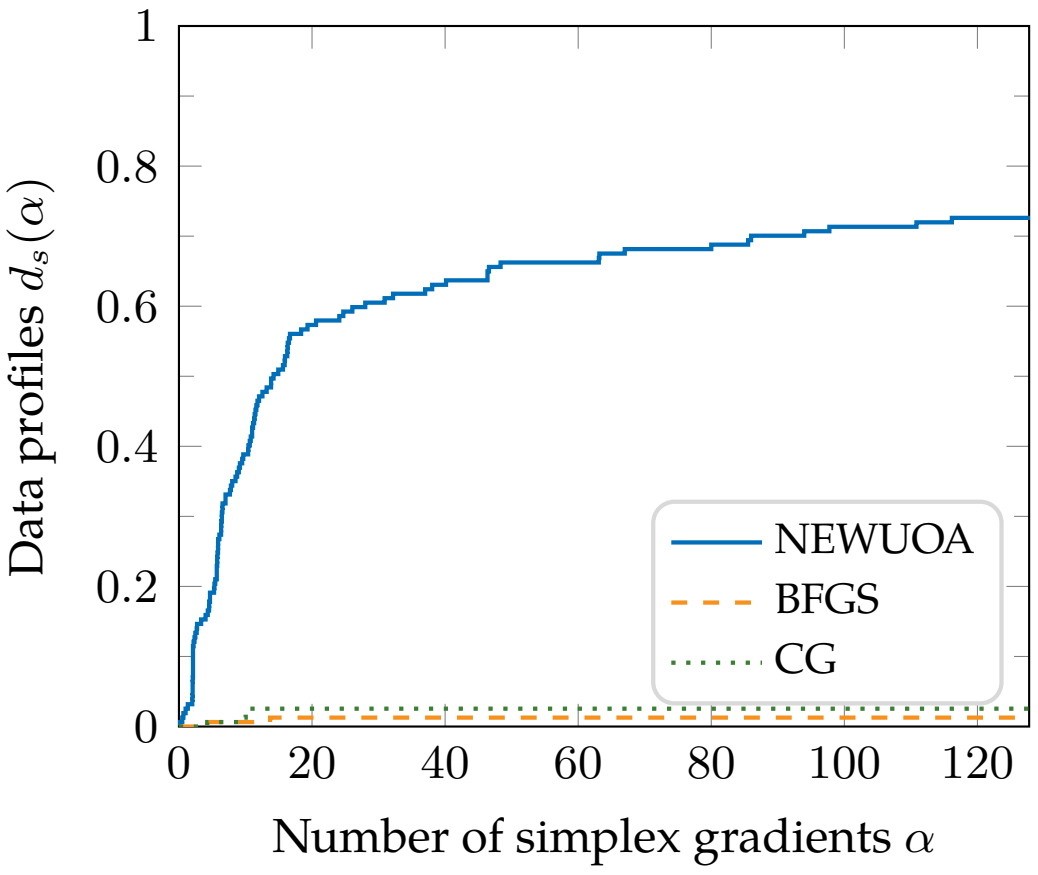

(d) Noisy experiment with $\sigma = 10^{-6}$

Figure 1.3: Data profiles of NEWUOA, BFGS, and CG with $\tau = 10^{-3}$

### 1.5.5 Absoluteness and relativeness of the profiles

Before concluding this section, we mention a significant difference between performance and data profiles. Data profiles are absolute because (1.5.3) compares the expense $t_{p,s}$ with a quantity independent of the other solvers, namely, $n_p + 1$ for each $p \in \mathcal{P}$. The data profile of a solver does depend weakly on the other solvers, the only dependence being the value of $\varphi_p^*$ in the convergence test (1.5.1). On the other hand, performance profiles are relative because the value of $r_{p,s}$ in (1.5.2) compares the expense of $s$ with that of the most efficient solver on $p$, which depends on all the other solvers. For example, the data profiles in Figure 1.3b suggest that CG outperforms BFGS during the corresponding experiment, and the profiles are unlikely to change dramatically if NEWUOA is removed from the comparison. In contrast, the comparison between BFGS and CG is not as conclusive based on the performance profiles in Figure 1.2b because the profiles may change significantly if NEWUOA is removed from the comparison. This phenomenon is stud-ied in [91], which presents examples on which comparisons using performance profiles should be made with care.

## 1.6 Contributions of the thesis

The main contributions of this thesis are as follows.

1. *An analysis on an interpolation set for underdetermined quadratic interpolation.* We show in Section 2.5 that an interpolation set introduced by Powell [166] for underdetermined quadratic interpolation is optimal in terms of the $\Lambda$-poisedness in the smallest ball that contains the set. This analysis justifies the choice of the default interpolation set employed by Powell in NEWUOA [166]. It also provides theoretical guidance for us to choose the initial interpolation set of our new DFO solver presented in Chapter 5.
2. *A cross-platform software package for Powell's DFO solvers.* We introduce in Chapter 3 our MATLAB and Python package named PDFO, which is open-source and distributed on all the major platforms (Windows, macOS, and Linux), available at `https://www.pdfo.net/`. It provides user-friendly interfaces for Powell's five model-based DFO solvers that are widely used in industrial and engineering applications. Powell originally implemented these solvers in Fortran 77, which is becoming an obstacle for many users. PDFO simplifies the usage of the solvers, and users do not need to deal with Fortran when using PDFO. As of August 2022, the package has been downloaded more than 30,000 times. It has been included as a major optimization engine in GEMSEO (`https://gemseo.readthedocs.io/`), an industrial software package for MDO.
3. *New perspectives on theory and practice of the SQP method.*
   (a) Theorem 4.1.1 provides a new interpretation of the objective function of the *sequential quadratic programming* (SQP) subproblem, showing that it is a natural quadratic approximation of the original objective function in the tangent space of a surface.
   (b) We reveal in Proposition 4.3.1 a relation between the Byrd-Omojokun and the Vardi approaches for solving the trust-region SQP subproblem in the equality-constrained case, interpreting the former as an approximation of the latter.
   (c) We extend the Byrd-Omojokun and the Vardi approaches to the inequality-constrained case in Section 4.3.4. They differ from those given in [44, § 15.4.4]. In particular, our Byrd-Omojokun approach is employed in the new solver described in Chapter 5, and performs much better than the one in [44, § 15.4.4].
4. *A new model-based DFO algorithm and its implementation.*

(a) We present in Chapter 5 a new derivative-free trust-region SQP algorithm for nonlinearly constrained optimization, named COBYQA. To the best of our knowledge, this is the first trust-region SQP algorithm that can solve DFO problems with nonlinear equality and inequality constraints without introducing slack variables.

(b) We detail in Chapter 6 the methods for solving the subproblems of COBYQA in our implementation. In particular, we propose in Chapter 6 a new method for approximately solving convex piecewise quadratic programming problems within a trust region. We use this method to deal with the normal subproblem of COBYQA.

(c) Chapter 7 presents the Python implementation of COBYQA, which is open-source and publicly available at `https://www.cobyqa.com/`. Section 7.3 details numerical experiments that compare COBYQA with Powell's solvers. COBYQA is competitive with NEWUOA on unconstrained problems and outperforms BOBYQA on bound-constrained ones, while being able to handle more general problems. Compared with LINCOA and COBYLA, a strength of COBYQA is that it always respects the bound constraints, while LINCOA and COBYLA do not. Such a strength is crucial in many applications, as is requested by engineers from IRT Saint Exupéry. On linearly constrained problems, COBYQA performs much better than LINCOA if the problems also contain bound constraints that cannot be violated. Most importantly, COBYQA outperforms COBYLA on all types of problems, no matter whether bound constraints (if any) can be violated or not. As a general-purpose solver, COBYQA is a good successor to COBYLA.

## 1.7 Organization of the thesis

The remaining of this thesis is organized as follows.

Chapter 2 presents interpolation models that are useful to DFO. We will focus on linear and quadratic models, as these are the models employed by the DFO methods presented later in this thesis.

Chapter 3 introduces PDFO, a package providing MATLAB and Python interfaces for using late Prof. M. J. D. Powell's DFO solvers. It also provides an overview of Powell's DFO solvers.

Chapter 4 presents an overview of the SQP framework and provides some perspectives. In particular, we elaborate an interpretation of the SQP subproblem that is new to the best of our knowledge. We also present a trust-region SQP method, discuss various approaches to solving its subproblem, and show an interesting connection between two of them.

Chapter 5 presents a new model-based DFO method named COBYQA. It is a trust-region SQP method, designed to tackle nonlinearly constrained problems that include equality and inequality constraints. It builds the trust-region quadratic models using derivative-free symmetric Broyden update presented in Chapter 2.

Chapter 6 introduces in details the methods employed by COBYQA to solve its subproblems, including the trust-region subproblem, the multiplier-estimating subproblem, and the geometry-improving subproblem.

Chapter 7 provides details on the Python implementation of COBYQA, in particular its preprocessing procedure, its stopping criteria, and some algorithmic details. It illustrates the usage of the Python package by several examples. It finally presents some numerical experiments that compare COBYQA with NEWUOA, BOBYQA, LINCOA, and COBYLA.

Chapter 8 summarizes the thesis, and discusses several future research directions.

## 1.8 Notations

Throughout this thesis, the functions $f$ and $c_i$, with $i \in \mathcal{I} \cup \mathcal{E}$, always denote the objective and constraint functions of an optimization problem. The Lagrangian function of this problem is denoted by $(x, \lambda) \mapsto \mathcal{L}(x, \lambda)$. Its first argument $x$ is the decision variable, while its second argument $\lambda = [\lambda_i]_{i \in \mathcal{I} \cup \mathcal{E}}^{\mathsf{T}}$ is the dual variable of the considered problem. A quadratic approximation of a function $u$ is always denoted by $\hat{u}$, and may be superscripted if this model is used in an iterative process. For example, if an optimization method maintains a quadratic approximation of the objective function, the $k$th model of the objective function is denoted by $\hat{f}^k$.

In general, an iteration number is written as a *superscript*. When an exponent is used in a formula that may lead to confusion, we will explicitly state the meaning of the notation in the context. Moreover, a *subscript* usually signifies a component index of a vector, a matrix, etc. An exception to this rule is the notation $e_i \in \mathbb{R}^n$, which represents the $i$th standard coordinate vector, i.e., the $i$th column of the identity matrix in $\mathbb{R}^{n \times n}$.

We will mainly use the $\ell_2$-norm in this thesis. Hence, $\|\cdot\|$ denotes the $\ell_2$-norm, unless otherwise stated. If any other norm is needed, its choice is clarified using subscripts. For example, $\|\cdot\|_{\mathsf{F}}$ denotes the Frobenius norm of a matrix, and $\|\cdot\|_1$ denotes the $\ell_1$-norm of a vector or the corresponding induced norm of a matrix.

# 2 Interpolation models for derivative-free optimization

As mentioned in Chapter 1, model-based DFO methods necessitate approximating locally the functions involved in optimization problems by simple functions, referred to as *models* or *surrogates*. These models are used to construct subproblems that are, in turn, approximately solved. Examples of such functions include polynomials and RBFs. In this thesis, we focus on linear and quadratic polynomial models, i.e., polynomials of degree at most one and two, respectively.

This chapter introduces the basic concepts of polynomial interpolation that we will need to present and establish DFO methods. Section 2.1 motivates the use of linear and quadratic polynomials in DFO methods. We then introduce the poisedness of interpolation sets and the quadratic Lagrange polynomials in Section 2.2. Section 2.3 presents a notion of well-poisedness referred to as the $\Lambda$-poisedness. Further, we summarize in Section 2.4 the methodology of underdetermined quadratic interpolation that will be employed later in this thesis. Finally, in Section 2.5, we study a particular interpolation set proposed by Powell [166], showing that this set is, in some sense, optimal. This set is essential to us because we will use a variation of it as the initial interpolation set of a new DFO method we will present in Chapter 5.

## 2.1 Introduction and motivation

Let $\mathscr{L}(\mathbb{R}^n)$ and $\mathbb{Q}(\mathbb{R}^n)$ denote the spaces of linear and quadratic polynomials on $\mathbb{R}^n$, respectively. In a DFO context, models from $\mathscr{L}(\mathbb{R}^n)$ or $\mathbb{Q}(\mathbb{R}^n)$ are built for a real-valued function $f$ without using derivatives. This can be done by interpolation schemes based on function values. Given a finite set of points $\mathcal{Y} = \{y^1, y^2, \dots, y^m\} \subseteq \mathbb{R}^n$, we construct a model $\hat{f}$ that interpolates the function $f$ on $\mathcal{Y}$, i.e.,

$$\hat{f}(y^i) = f(y^i), \quad \text{for } i \in \{1, 2, \dots, m\}. \tag{2.1.1}$$

The conditions (2.1.1) may be inconsistent. In such a case, models can be built using regression schemes. For example, a least-square regression model $\hat{f}$ minimizes

$$\sum_{i=1}^{m} [f(y^i) - \hat{f}(y^i)]^2.$$

Although there are successful methods that use regression models (see, e.g., [17, 46]), the DFO methods we present in this thesis use interpolation models and ensure that the interpolation conditions are consistent and well-conditioned (see Chapters 3 and 5).

It is also possible to use polynomials of degree higher than two. However, we do not consider such models in this thesis due to the following observations.

1. As shown in [206, Thm. 2.5], the space of polynomials on $\mathbb{R}^n$ of degree at most $k$ has a dimension of
$$\binom{n+k}{n} = \frac{1}{k!}\prod_{i=1}^{k}(n+i) \ge \frac{n^k}{k!}.$$
Therefore, to determine a model from this space merely by the interpolation conditions (2.1.1), we need in general $\mathcal{O}(n^k)$ function values. This amount is unacceptable in a DFO context unless $k$ is small. It is possible to reduce this number with underdetermined interpolation, which is used by several optimization methods for $k \le 2$ (see Section 2.4), including the method COBYQA developed in this thesis (see Chapter 5). Using underdetermined interpolation models with $k \ge 3$ is out of the scope of this thesis, although it is an interesting research direction.

2. The DFO methods need to solve approximately subproblems (e.g., trust-region subproblems) built upon these models. Sophisticated models usually lead to complicated subproblems to solve. On the other hand, even with models that are not quadratic, practical algorithms normally solve the subproblems based on first- or second-order approximations of the models, e.g., calculate the approximate Cauchy point discussed in [44, § 6.3.3]. Therefore, building polynomial models of degrees higher than two may not be necessary in practice.

Although we do not study RBF models, we mention that there also exist many DFO methods based on these models. Examples of such methods include ORBIT [208], CONORBIT [177], and BOOSTERS [146].

## 2.2 Elementary concepts of multivariate interpolation

In this section, the space in which $\hat{f}$ lies is either $\mathscr{L}(\mathbb{R}^n)$ or $\mathbb{Q}(\mathbb{R}^n)$. Moreover, we assume that the number of interpolation points $m$ equals the dimension of the chosen space, i.e., either $n+1$ or $(n+1)(n+2)/2$, respectively.

### 2.2.1 Poisedness of interpolation sets

Before studying properties of multivariate interpolation, we must introduce the following notion of poisedness.

**Definition 2.2.1** (Poisedness)**.** The set $\mathcal{Y}$ is *poised* for interpolation on $\mathscr{L}(\mathbb{R}^n)$ if the interpolation system (2.1.1) has a unique solution in $\mathscr{L}(\mathbb{R}^n)$ for any real-valued function $f$. Similarly, it is *poised* for interpolation on $\mathbb{Q}(\mathbb{R}^n)$ if the system has a unique solution in $\mathbb{Q}(\mathbb{R}^n)$ for any real-valued function $f$.

Let us first consider the problem of finding a linear model $\hat{f} \in \mathscr{L}(\mathbb{R}^n)$ satisfying the interpolation system (2.1.1) whenever $m = \dim \mathscr{L}(\mathbb{R}^n) = n + 1$. In the natural basis of $\mathscr{L}(\mathbb{R}^n)$, system (2.1.1) can be reformulated as

$$\alpha + g^{\mathsf{T}} y^i = f(y^i), \quad \text{for } i \in \{1, 2, \dots, m\}, \tag{2.2.1}$$

where $\alpha \in \mathbb{R}$ and $g \in \mathbb{R}^n$ are the coefficients to determine. If $\mathcal{Y}$ is poised for linear interpolation, given $(\alpha^*, g^*)$ the solution to (2.2.1), the linear model $\hat{f}$ is defined by

$$\hat{f}(x) = \alpha^* + (g^*)^{\mathsf{T}} x, \quad \text{for } x \in \mathbb{R}^n,$$

and the vector $\nabla \hat{f} \equiv g^*$ is referred to as the *simplex gradient* of $f$ for the interpolation set $\mathcal{Y}$. Such models are used, for instance, by COBYLA [161], a DFO method for nonlinearly constrained optimization detailed in Section 3.2.2.

The problem of finding a quadratic model $\hat{f} \in \mathbb{Q}(\mathbb{R}^n)$ satisfying the interpolation system (2.1.1) when $m = \dim \mathbb{Q}(\mathbb{R}^n) = (n+1)(n+2)/2$ is very similar. This kind of models are used, for instance, by UOBYQA [163], a DFO method for unconstrained optimization detailed in Section 3.2.3. The advantage of quadratic models is that they capture curvature information of the function $f$, and are more precise than linear models. In the remaining of this chapter, we focus our discussions on quadratic models.

### 2.2.2 Lagrange polynomials

Assume that $\mathcal{Y} = \{y^1, y^2, \dots, y^m\} \subseteq \mathbb{R}^n$ is a poised interpolation set of $m = (n + 1)(n + 2)/2$ points. We introduce the Lagrange polynomials for the interpolation problem

$$\hat{f}(y^i) = f(y^i), \quad \text{for } i \in \{1, 2, \dots, m\} \text{ with } \hat{f} \in \mathbb{Q}(\mathbb{R}^n). \tag{2.2.2}$$

**Definition 2.2.2** (Lagrange polynomials for (2.2.2))**.** For each $i \in \{1, 2, \dots, m\}$, the $i$th Lagrange polynomial $L_i$ for the interpolation problem (2.2.2) is the unique quadratic polynomial that satisfies

$$L_i(y^j) = \begin{cases} 1, & j = i, \\ 0, & j \in \{1, 2, \dots, m\} \setminus \{i\}. \end{cases}$$

We note that in Definition 2.2.2, the poisedness of the interpolation set ensures the existence and the uniqueness of the Lagrange polynomials. Further, it is well known that the interpolant $\hat{f} \in \mathbb{Q}(\mathbb{R}^n)$ of $f$ on $\mathcal{Y}$ can be formulated as a linear combination of the Lagrange polynomials, as detailed in Theorem 2.2.1

**Theorem 2.2.1.** *The Lagrange polynomials* $\{L_1, L_2, \ldots, L_m\}$ *form a basis of* $\mathbb{Q}(\mathbb{R}^n)$. *Moreover, the quadratic interpolant* $\hat{f}$ *of* $f$ *on* $\mathcal{Y}$ *is given by*

$$\hat{f}(x) = \sum_{i=1}^{m} f(y^i)L_i(x), \quad \textit{for } x \in \mathbb{R}^n.$$

One of the most important application of the Lagrange polynomials is to measure the well-poisedness of the interpolation set $\mathcal{Y}$, as we discuss in the next section.

## 2.3 Well-poisedness of interpolation sets

As we mentioned above, the interpolation set $\mathcal{Y}$ is poised if it admits a unique interpolant in $\mathbb{Q}(\mathbb{R}^n)$. However, given a poised set, how good is this interpolation set? This section briefly discusses this topic and introduces the concept of $\Lambda$-poisedness.

Suppose that $f$ is thrice differentiable on $\mathbb{R}^n$. It is well known that for any $x \in \mathbb{R}^n$, we have

$$|f(x) - \hat{f}(x)| \le \frac{\theta}{6} \sum_{i=1}^{m} |L_i(x)| \|x - y^i\|^3, \tag{2.3.1}$$

where $\theta$ is an upper bound on the absolute value of the third-order directional derivatives of $f$ [162, Thm. 2]. A similar bound is established in [40, Thm. 2]. For any compact set $\mathcal{C} \subseteq \mathbb{R}^n$, we then clearly have

$$\max_{x \in \mathcal{C}} |f(x) - \hat{f}(x)| \le \frac{m\theta\Lambda}{6} \max_{1 \le i \le m} \max_{x \in \mathcal{C}} \|x - y^i\|^3, \tag{2.3.2}$$

where $\Lambda$ is defined by

$$\Lambda \overset{\text{def}}{=} \max_{1 \le i \le m} \max_{x \in \mathcal{C}} |L_i(x)|. \tag{2.3.3}$$

In this sense, $\Lambda$ measures the well-poisedness of the interpolation set $\mathcal{Y}$ with respect to the set $\mathcal{C}$. This motivates the following concept of $\Lambda$-poisedness.

**Definition 2.3.1** ($\Lambda$-poisedness [47, Def. 3.6])**.** A poised interpolation set $\mathcal{Y} \subseteq \mathbb{R}^n$ is said to be $\Lambda$-poised in a compact set $\mathcal{C} \subseteq \mathbb{R}^n$, for some $\Lambda > 0$, if

$$\Lambda \ge \max_{1 \le i \le m} \max_{x \in \mathcal{C}} |L_i(x)|.$$

Note that in this definition, the compact set $\mathcal{C} \subseteq \mathbb{R}^n$ may or may not contain the interpolation points in $\mathcal{Y}$. The compactness of $\mathcal{C}$ can be relaxed if we take the supremum instead of the maximum over $\mathcal{C}$. If $\mathcal{Y}$ is $\Lambda_0$-poised in $\mathcal{C}$, it is obviously $\Lambda$-poised in $\mathcal{C}$ for any $\Lambda \geq \Lambda_0$. Alternative definitions of the $\Lambda$-poisedness of sets are given in [47, § 3.3].

As discussed in [47, § 3.3], the $\Lambda$-poisedness measures how good an interpolation set is. A lower $\Lambda$ indicates a better interpolation set. Indeed, if $\mathcal{Y}$ is $\Lambda$-poised in a compact set $\mathcal{C} \subseteq \mathbb{R}^n$ and $\Lambda$ is reasonably low, then (2.3.2) shows that the quadratic interpolant $\hat{f}$ represents $f$ reasonably well in $\mathcal{C}$.

The notion of $\Lambda$-poisedness is intrinsically related to the well-conditioning of the interpolation system. Note that the interpolation system (2.2.2) is linear with respect to the coefficients of the interpolant $\hat{f}$. One can show that the set $\mathcal{Y}$ is $\Lambda$-poised if and only if the condition number of the coefficient matrix of this linear system is bounded by some terms proportional to $\Lambda$ [47, Thm. 3.14]. Therefore, the lower $\Lambda$ is, the better the conditioning of the interpolation system is, and hence, the better the interpolation set is from a numerical standpoint.

In practice, however, we do not use the $\Lambda$-poisedness of interpolation sets directly. This is because determining whether a set is $\Lambda$-poised is a difficult problem, even if the compact set $\mathcal{C} \subseteq \mathbb{R}^n$ is simple (e.g., a ball). It is a theoretical tool that we mostly use to justify design choices of algorithms.

## 2.4 Underdetermined quadratic interpolation

As before, we focus on quadratic interpolation in this section. Building a quadratic model by the interpolation system (2.2.2) requires $\mathcal{O}(n^2)$ function evaluations. Therefore, if a DFO method employs such models, its initialization already costs $\mathcal{O}(n^2)$ function evaluations, which are needed to build the first model. We use underdetermined interpolation to reduce this amount of function evaluations, as presented hereinafter.

Underdetermined interpolation works as follows. Assume that we are given an interpolation set $\mathcal{Y} \subseteq \mathbb{R}^n$ such that the interpolation conditions (2.1.1) are consistent. However, the number $m$ of interpolation points may be lower than $\dim \mathbb{Q}(\mathbb{R}^n)$. In such a case, the interpolation conditions may *not* define a unique interpolant. To select one, we choose a functional $\mathcal{F} : \mathbb{Q}(\mathbb{R}^n) \to \mathbb{R}$ that reflects a desired property or regularity of the interpolants. An interpolant $\hat{f}$ of $f$ is then defined as a solution of

$$\min_{Q \in \mathbb{Q}(\mathbb{R}^n)} \quad \mathcal{F}(Q) \tag{2.4.1a}$$

$$\text{s.t.} \quad Q(y^i) = f(y^i), \ i \in \{1, 2, \dots, m\}. \tag{2.4.1b}$$

Two examples of the functional $\mathcal{F}$ are given in Sections 2.4.1 and 2.4.2. If $m = \dim \mathbb{Q}(\mathbb{R}^n)$ and if $\mathcal{Y}$ is poised in the sense of Definition 2.4.1, then the underdetermined interpolant reduces to the unique model that satisfies the interpolation conditions.

In what follows, we adapt the definitions and some properties presented in the previous section to special cases of underdetermined interpolation. We briefly present two underdetermined schemes: the least Frobenius norms scheme and the derivative-free symmetric Broyden update. The latter is of particular interest to us, because it underlies the three most recent DFO solvers of Powell and our new model-based DFO method, presented in Chapter 5.

### 2.4.1 Least Frobenius norm quadratic models

One of the simplest functional $\mathcal{F}$ mentioned above is the Frobenius norm of the Hessian matrix of the interpolants. In this context, a quadratic model $\hat{f}$ solves

$$\min_{Q \in \mathbb{Q}(\mathbb{R}^n)} \quad \|\nabla^2 Q\|_{\mathsf{F}} \tag{2.4.2a}$$

$$\text{s.t.} \quad Q(y^i) = f(y^i),\ i \in \{1, 2, \dots, m\}, \tag{2.4.2b}$$

where $\|\cdot\|_{\mathsf{F}}$ denotes the Frobenius norm. This norm is chosen because it can be easily evaluated, unlike the spectral norm, for example. The Frobenius norm also brings computational advantages because (2.4.2) is essentially a quadratic programming problem, which can be tackled by solving a linear system. Examples of DFO methods that use least Frobenius norm quadratic models include MNH [207] and DFO [48–50]. We also assume hereafter that $m \geq n + 2$, as otherwise, any solution to (2.4.2) would be linear.

We first present the notion of poisedness for problem (2.4.2) as follows.

> **Definition 2.4.1** (Poisedness)**.** The set $\mathcal{Y}$ is *poised* in the minimum Frobenius norm sense if the solution to problem (2.4.2) exists and is unique for any real-valued function $f$.

To get an in-depth view of the poisedness of $\mathcal{Y}$, let us investigate how to solve the problem (2.4.2). Given any fix $\bar{x} \in \mathbb{R}^n$, we expand a quadratic polynomial $Q \in \mathbb{Q}(\mathbb{R}^n)$ for any $x \in \mathbb{R}^n$ as

$$Q(x) = \alpha + g^{\mathsf{T}}(x - \bar{x}) + \frac{1}{2}(x - \bar{x})^{\mathsf{T}} H (x - \bar{x}),$$

where $\alpha \in \mathbb{R}$, $g \in \mathbb{R}^n$, and $H \in \mathbb{R}^{n \times n}$ are the coefficients of $Q$, with $H$ being symmetric. Note that the constraint (2.4.2b) is linear with respect to $(\alpha, g, H)$. Consequently, the problem (2.4.2) is a convex quadratic problem with equality constraints. We then consider its KKT system. Let $\lambda = [\lambda_1, \lambda_2, \dots, \lambda_m]^{\mathsf{T}}$ be the Lagrange multiplier of (2.4.2). The stationarity condition (1.3.2a) can be formulated as

$$\sum_{i=1}^m \lambda_i (y^i - \bar{x})(y^i - \bar{x})^{\mathsf{T}} = H, \quad \sum_{i=1}^m \lambda_i = 0, \quad \text{and} \quad \sum_{i=1}^m \lambda_i (y^i - \bar{x}) = 0. \tag{2.4.3}$$

Combining (2.4.3) with the primal feasibility condition (1.3.2c), we obtain

$$\begin{bmatrix} Y_H & e & Y_L^{\mathsf{T}} \\ e^{\mathsf{T}} & 0 & 0 \\ Y_L & 0 & 0 \end{bmatrix} \begin{bmatrix} \lambda \\ \alpha \\ g \end{bmatrix} = \begin{bmatrix} f(\mathcal{Y}) \\ 0 \\ 0 \end{bmatrix}, \tag{2.4.4}$$

where $f(\mathcal{Y}) \overset{\text{def}}{=} [f(y^1), f(y^2), \dots, f(y^m)]^{\mathsf{T}}$, $e \in \mathbb{R}^m$ is the all-one vector, while $Y_L \in \mathbb{R}^{n \times m}$ and $Y_H \in \mathbb{R}^{m \times m}$ are defined by

$$Y_L \overset{\text{def}}{=} \begin{bmatrix} y^1 - \bar{x} & y^2 - \bar{x} & \cdots & y^m - \bar{x} \end{bmatrix} \quad \text{and} \quad Y_H \overset{\text{def}}{=} \frac{1}{2}(Y_L^{\mathsf{T}} Y_L) \odot (Y_L^{\mathsf{T}} Y_L). \tag{2.4.5}$$

In (2.4.5), $\odot$ denotes the Hadamard product, and hence, the $(i, j)$th entry of $Y_H$ is given by $[(y^i - \bar{x})^{\mathsf{T}}(y^j - \bar{x})]^2/2$. In fact, the problem (2.4.2) is equivalent to the linear system (2.4.4), where the matrix $H$ can be retrieved from the value of $\lambda$ by (2.4.3) [164, § 2], and the interpolation set $\mathcal{Y}$ is poised if and only if the coefficient matrix in (2.4.4) is nonsingular [47, § 5.3]. Since the poisedness of $\mathcal{Y}$ does not depend on the choice of $\bar{x}$, neither does the nonsingularity of the coefficient matrix in (2.4.4).

Any interpolation set that is poised in the sense of Definition 2.4.1 must be affinely independent, i.e., it does not lie in a low-dimensional affine subset of $\mathbb{R}^n$, because the matrix $[e, Y_L^{\mathsf{T}}]$ has full-column rank according to the nonsingularity of the coefficient matrix in (2.4.4). Even if $\mathcal{Y}$ is not poised, the full-rankness of $[e, Y_L^{\mathsf{T}}]$ holds as long as the solution to (2.4.2) is unique when $f \equiv 0$. This is because $\alpha + g^{\mathsf{T}}(x - \bar{x})$ is a nonzero linear polynomial that interpolates the zero function on $\mathcal{Y}$ if $[\alpha, g^{\mathsf{T}}]^{\mathsf{T}} \in \mathbb{R}^{n+1}$ is a nonzero vector in $\ker([e, Y_L^{\mathsf{T}}])$. Such a polynomial should not exist because of the aforementioned uniqueness, since the zero polynomial should be the only solution.

### Minimum Frobenius norm Lagrange polynomials

Let $\mathcal{Y} \subseteq \mathbb{R}^n$ be a poised interpolation set of $m$ points in the minimum Frobenius norm sense, with $n + 2 \le m \le \dim \mathbb{Q}(\mathbb{R}^n)$. We now present the minimum Frobenius norm Lagrange polynomials for the interpolation system (2.2.2).

**Definition 2.4.2** (Minimum Frobenius norm Lagrange polynomials for (2.2.2))**.** For each $i \in \{1, 2, \dots, m\}$, the $i$th minimum Frobenius norm Lagrange polynomial $L_i$ for the interpolation problem (2.2.2) is the unique quadratic polynomial that solves

$$\begin{aligned} \min_{Q \in \mathbb{Q}(\mathbb{R}^n)} \quad & \|\nabla^2 Q\|_{\mathsf{F}} \\ \text{s.t.} \quad & Q(y^i) = 1, \\ & Q(y^j) = 0, \ j \in \{1, 2, \dots, m\} \setminus \{i\}. \end{aligned}$$

As in the case where $m = \dim \mathbb{Q}(\mathbb{R}^n)$, any minimum Frobenius norm interpolant can be expressed as a linear combination of minimum Frobenius norm Lagrange polynomials, as described by Theorem 2.4.1 (see [164, § 3]), even though the set of all minimum Frobenius norm Lagrange polynomials is not a basis of $\mathbb{Q}(\mathbb{R}^n)$ if $m < \dim \mathbb{Q}(\mathbb{R}^n)$.

**Theorem 2.4.1.** *The minimum Frobenius norm quadratic interpolant $\hat{f}$ of $f$ on $\mathcal{Y}$ is given for any $x \in \mathbb{R}^n$ by*

$$\hat{f}(x) = \sum_{i=1}^{m} f(y^i) L_i(x).$$

#### Poisedness of interpolation sets in the minimum Frobenius norm sense

We now present the notion of $\Lambda$-poisedness in the minimum Frobenius norm sense, which is a straightforward adaptation of Definition 2.3.1.

**Definition 2.4.3** ($\Lambda$-poisedness in the minimum Frobenius norm sense [47, Def. 5.6])**.** A poised interpolation set $\mathcal{Y} \subseteq \mathbb{R}^n$ in the minimum Frobenius norm sense is said to be $\Lambda$-poised in the minimum Frobenius norm sense in a compact set $\mathcal{C} \subseteq \mathbb{R}^n$, for some $\Lambda > 0$, if

$$\Lambda \geq \max_{1 \leq i \leq m} \max_{x \in \mathcal{C}} |L_i(x)|.$$

Similarly to what we mentioned before, if $\mathcal{Y}$ is $\Lambda_0$-poised in the minimum Frobenius norm sense in $\mathcal{C}$, it is obviously $\Lambda$-poised in the minimum Frobenius norm sense in $\mathcal{C}$ for all $\Lambda \geq \Lambda_0$. Moreover, one can show that the set $\mathcal{Y}$ is $\Lambda$-poised if and only if the condition number of the coefficient matrix in (2.4.4) is bounded by some terms proportional to $\Lambda$ [47, Thm. 5.8].

### 2.4.2 Quadratic models based on the symmetric Broyden update

Now, we present another functional $\mathcal{F}$ that can be used in the variational problem (2.4.1). It is introduced by Powell [164]. It takes up the freedom bequeathed by the interpolation conditions (2.2.2) by minimizing the difference between the Hessian matrix of the interpolant and a prior estimation $H_{\text{old}} \in \mathbb{R}^{n \times n}$ of this matrix instead of minimizing the Frobenius norm of the Hessian matrix. More specifically, we build a quadratic model $\hat{f}$

by solving

$$\min_{Q \in \mathbb{Q}(\mathbb{R}^n)} \quad \|\nabla^2 Q - H_{\text{old}}\|_{\mathsf{F}} \tag{2.4.6a}$$

$$\text{s.t.} \quad Q(y^i) = f(y^i),\ i \in \{1, 2, \dots, m\}. \tag{2.4.6b}$$

In a DFO method, the model $\hat{f}^k \in \mathbb{Q}(\mathbb{R}^n)$ of $f$ at the $k$th iteration can be constructed by solving (2.4.6) with $H_{\text{old}} = \nabla^2 \hat{f}^{k-1}$, and we define $\hat{f}^{-1} \equiv 0$. This is referred to as the *least Frobenius norm updating of quadratic models* [164]. Examples of DFO methods that use this update include NEWUOA [166], BOBYQA [168], and LINCOA, presented in Sections 3.2.4 to 3.2.6. The new method we introduce in Chapter 5 is also based on this update.

The above discussions we had on the minimum Frobenius norm models can be easily adapted to the least Frobenius norm updating. If we are given a function $\hat{f}^{\text{old}}$ such that $H_{\text{old}} = \nabla^2 \hat{f}^{\text{old}}$, then $\hat{f} - \hat{f}^{\text{old}}$ solves

$$\min_{Q \in \mathbb{Q}(\mathbb{R}^n)} \quad \|\nabla^2 Q\|_{\mathsf{F}} \tag{2.4.7a}$$

$$\text{s.t.} \quad Q(y^i) = f(y^i) - \hat{f}^{\text{old}}(y^i),\ i \in \{1, 2, \dots, m\}. \tag{2.4.7b}$$

In a DFO method, we have $\hat{f}^{\text{old}} = \hat{f}^{k-1}$. Therefore, the coefficients of the model can be evaluated by solving a linear system similar to (2.4.4).

As mentioned in [220, § 3.6], the minimum Frobenius norm update can be regarded as a derivative-free variation of the *Powell symmetric Broyden* (PSB) update in the quasi-Newton method. The PSB method [155] approximates the function $f$ at a point $y^1 \in \mathbb{R}^n$ by a quadratic function

$$\hat{f}(x) = \hat{f}(y^1) + \nabla \hat{f}(y^1)^{\mathsf{T}}(x - y^1) + \frac{1}{2}(x - y^1)^{\mathsf{T}} \hat{H}(x - y^1),$$

and, as shown in [58, Thm. 4.2], $\hat{H} \in \mathbb{R}^{n \times n}$ is the unique solution to

$$\begin{aligned} \min_{H \in \mathbb{R}^{n \times n}} \quad & \|H - H_{\text{old}}\|_{\mathsf{F}} \\ \text{s.t.} \quad & H(y^1 - y^2) = \nabla \hat{f}(y^1) - \nabla \hat{f}(y^2), \\ & H^{\mathsf{T}} = H, \end{aligned}$$

where $H_{\text{old}} \in \mathbb{R}^{n \times n}$ is a prior estimation of $\nabla^2 f(y^1)$, and $y^2 \in \mathbb{R}^n$ is a point that differs from $y^1$. In an optimization context, $y^1$ and $y^2$ are respectively the current and the previous iterates of the optimization method. Clearly, the function $\hat{f}$ solves

$$\min_{Q \in \mathbb{Q}(\mathbb{R}^n)} \quad \|\nabla^2 Q - H_{\text{old}}\|_{\mathsf{F}} \tag{2.4.8a}$$

$$\text{s.t.} \quad Q(y^1) = f(y^1), \tag{2.4.8b}$$

$$\nabla Q(y^i) = \nabla f(y^i), \; i \in \{1, 2\}. \tag{2.4.8c}$$

The variational problem (2.4.6) is a derivative-free analog of (2.4.8). Therefore, following Powell [169] and Zhang [220], we will also refer to the minimum Frobenius norm update as the *derivative-free symmetric Broyden update*.

### 2.4.3 Implementation in derivative-free optimization methods

We consider a model-based DFO algorithm that builds models by the derivative-free symmetric Broyden update, and we summarize the implementation of this update, as proposed by Powell [164, 165]. The implementation can straightforwardly be adapted to construct minimum Frobenius norm models.

For each $k \geq 0$, we denote by $x^k \in \mathbb{R}^n$ and the $k$th iterate, and the $k$th interpolation set is $\mathcal{Y}^k = \{y^{k,1}, y^{k,2}, \dots, y^{k,m}\} \subseteq \mathbb{R}^n$, with $x^k \in \mathcal{Y}^k$. The DFO algorithm builds a quadratic model $\hat{f}^k \in \mathbb{Q}(\mathbb{R}^n)$ of the objective function $f$ by the derivative-free symmetric Broyden update (2.4.6), with $H_{\text{old}} = \nabla^2 \hat{f}^{k-1}$ and $\mathcal{Y} = \mathcal{Y}^k$. The model $\hat{f}^k$ is

$$\hat{f}^k(x) = \alpha^k + (g^k)^{\mathsf{T}}(x - x^k) + \frac{1}{2}(x - x^k)^{\mathsf{T}} H^k (x - x^k), \quad \text{for } x \in \mathbb{R}^n,$$

where $\alpha^k \in \mathbb{R}$, $g^k \in \mathbb{R}^n$, and $H^k \in \mathbb{R}^{n \times n}$ are the coefficients, $H^k$ being symmetric. As discussed in Section 2.4.1, thanks to (2.4.7), $\alpha^k$, $g^k$, and $H^k$ can be calculated by solving

$$\begin{bmatrix} Y_H & e & Y_L^{\mathsf{T}} \\ e^{\mathsf{T}} & 0 & 0 \\ Y_L & 0 & 0 \end{bmatrix} \begin{bmatrix} \lambda \\ \alpha \\ g \end{bmatrix} = \begin{bmatrix} v \\ 0 \\ 0 \end{bmatrix}, \tag{2.4.9}$$

where $v \in \mathbb{R}^m$ is defined by

$$v \overset{\text{def}}{=} \begin{bmatrix} f(y^{k,1}) - \hat{f}^{k-1}(y^{k,1}) \\ f(y^{k,2}) - \hat{f}^{k-1}(y^{k,2}) \\ \vdots \\ f(y^{k,m}) - \hat{f}^{k-1}(y^{k,m}) \end{bmatrix}.$$

If $\lambda$, $\alpha$, and $g$ solve the system (2.4.9), then the model's coefficients can be evaluated by

$$\begin{cases} \alpha^k = \hat{f}^{k-1}(x^k) + \alpha, & (2.4.10a) \\ g^k = \nabla \hat{f}^{k-1}(x^k) + g, & (2.4.10b) \\ H^k = \nabla^2 \hat{f}^{k-1} + \sum_{i=1}^{m} \lambda_i (y^{k,i} - x^k)(y^{k,i} - x^k)^\mathsf{T}. & (2.4.10c) \end{cases}$$

The methods we consider in this thesis build $\mathcal{Y}^{k+1}$ from $\mathcal{Y}^k$ by modifying only one point. Taking advantage of this fact, we can reduce the computational cost of these models and improve the stability of the computations, as is presented in the sequel.

### Storage and update of the quadratic model

The intercept $\alpha^k$ and the gradient vector $g^k$ of the model $\hat{f}^k$ can be stored directly. However, forming its Hessian matrix $H^k$ using (2.4.10c) requires making $\mathcal{O}(mn^2)$ basic operations. To reduce this cost, as suggested by Powell [164, § 3], we decompose the Hessian matrix as

$$H^k = \Gamma^k + \sum_{i=1}^{m} \gamma_i^k (y^{k,i} - x^k)(y^{k,i} - x^k)^\mathsf{T}, \qquad (2.4.11)$$

where $\Gamma^k \in \mathbb{R}^{n \times n}$ is referred to as an explicit part of $H^k$, and $[\gamma_1^k, \gamma_2^k, \dots, \gamma_m^k]^\mathsf{T} \in \mathbb{R}^m$ is referred to as an implicit part of $H^k$, which are both stored. Recall that $\mathcal{Y}^{k+1}$ differs from $\mathcal{Y}^k$ by only one point. Suppose that $y^{k+1,i} = y^{k,i}$ for all $i$ except for $i = \ell$. The explicit and implicit parts of $H^k$ can then be updated as

$$\begin{cases} \Gamma^{k+1} = \Gamma^k + \gamma_\ell^k (y^{k,\ell} - x^k)(y^{k,\ell} - x^k)^\mathsf{T}, & (2.4.12a) \\ \gamma_i^{k+1} = \gamma_i^k + \lambda_i, \ i \in \{1, 2, \dots, m\} \setminus \{\ell\}, & (2.4.12b) \\ \gamma_\ell^{k+1} = \lambda_\ell. & (2.4.12c) \end{cases}$$

In doing so, updating a quadratic model requires only $\mathcal{O}(n^2)$ basic operations, after solving the linear system (2.4.9).

Powell applied (2.4.11) and (2.4.12) in his three recent solvers NEWUOA [166], BOBYQA [168], and LINCOA, presented in Sections 3.2.4 to 3.2.6. The new method that we introduce in Chapter 5 also uses this technique.

### Storage and update of the inverse of the coefficient matrix

Since the methods we consider build $\mathcal{Y}^{k+1}$ from $\mathcal{Y}^k$ by changing only one point, only one row and one column of the coefficient matrix in (2.4.9) is modified from one iteration to another, which is a rank-2 modification of this matrix. Therefore, we can update its inverse matrix using the Sherman-Morrison-Woodbury formula from one iteration to another. Such an update requires only $\mathcal{O}(n^2)$ operations, much less than what would be required to compute the inverse from scratch (see [165, § 2] for details).

Suppose that the above-mentioned inverse is

$$\begin{bmatrix} \Omega & \Xi^\mathsf{T} \\ \Xi & \Upsilon \end{bmatrix},$$

where $\Omega \in \mathbb{R}^{m\times m}$, $\Xi \in \mathbb{R}^{(n+1)\times m}$, and $\Upsilon \in \mathbb{R}^{(n+1)\times(n+1)}$. The submatrices $\Xi$ and $\Upsilon$ can be stored directly, although the first row of $\Xi$ can be discarded, as well as the first row and the first column of $\Upsilon$ (see [166, § 4]). However, $\Omega$ should not be stored as is, due to the following reason. As mentioned by Powell [165], the submatrix $\Omega$ is positive semidefinite with rank $m-n-1$. If $\Omega$ is stored as an $m \times m$ matrix, the rank property may be lost in practice due to computer rounding errors, as observed by Powell during the development of NEWUOA (see [165, § 1] and the last paragraph of [166, § 8]). To maintain the structure of $\Omega$, we store the factorized form $\Omega = ZDZ^\mathsf{T}$, with $Z \in \mathbb{R}^{m\times(m-n-1)}$ and $D \in \mathbb{R}^{(m-n-1)\times(m-n-1)}$, where $D$ is a diagonal matrix with elements $\pm 1$ on the diagonal. Theoretically, the matrix $D$ should always be the identity. However, in practice, $D$ may have negative diagonal entries. This is because the positive semidefiniteness of $\Omega$ can be lost due to computer rounding errors.

## 2.5 An optimal interpolation set

We study in this section an interpolation set that we will use in Chapter 5 of this thesis, where we introduce a new model-based DFO method. This interpolation set is adapted from [166] as follows. Let $\delta > 0$ be fixed and for $j \in \{1, 2, \dots, 2n+1\}$, let $z^j \in \mathbb{R}^n$ be

$$z^j \stackrel{\text{def}}{=} \begin{cases} 0, & \text{if } j = 1, \quad (2.5.1a) \\ \delta e_{j-1}, & \text{if } 2 \le j \le n+1, \quad (2.5.1b) \\ -\delta e_{j-n-1}, & \text{otherwise}, \quad (2.5.1c) \end{cases}$$

where $e_j \in \mathbb{R}^n$ denotes the $j$th standard coordinate vector. We then define the interpolation set $\mathcal{Z}_m \subseteq \mathbb{R}^n$ for each $m \in \{n+2, n+3, \dots, 2n+1\}$ by

$$\mathcal{Z}_m \stackrel{\text{def}}{=} \{z^1, z^2, \dots, z^m\}. \tag{2.5.2}$$

Further, we denote by $\mathcal{B}_p^m$ the smallest $\ell_p$-norm ball containing $\mathcal{Z}_m$, for $p \in [1, \infty]$. Note that we allow $p = \infty$. By the construction of $\mathcal{Z}_m$, we observe that

$$\mathcal{B}_p^m \equiv \mathcal{B}_p(\delta) \stackrel{\text{def}}{=} \{x \in \mathbb{R}^n : \|x\|_p \le \delta\}. \tag{2.5.3}$$

In the sequel, we study the $\Lambda$-poisedness of the set $\mathcal{Z}_m$ in $\mathcal{B}_p^m$ for $p \ge 1$.

We will show that this set is 1-poised in $\mathcal{B}_p^m$ for all $p \in [1, 2]$ when $m = 2n+1$. In this sense, $\mathcal{Z}_{2n+1}$ is an optimal interpolation set.

### 2.5.1 Formulation of the Lagrange polynomials

According to Definition 2.4.3 and the equation (2.5.3), the set $\mathcal{Z}_m$ is $\Lambda_p$-poised in $\mathcal{B}_p^m$ in the minimum Frobenius norm sense with

$$\Lambda_p \overset{\text{def}}{=} \max_{1 \le i \le m} \max_{x \in \mathcal{B}_p(\delta)} |L_i(x)|,$$

where $L_i$ is the $i$th minimum Frobenius norm Lagrange polynomial associated with $\mathcal{Z}_m$ for $i \in \{1, 2, \dots, m\}$.

To study $\Lambda_p$, we first present explicit formulae for $L_i$ for all $i \in \{1, 2, \dots, m\}$. These formulae are given in [166, § 3], without a proof.

**Lemma 2.5.1.** *For each $m \in \{n+2, n+3, \dots, 2n+1\}$ and all $x \in \mathbb{R}^n$, we have*

$$L_i(x) = \begin{cases} 1 - \delta^{-2} \displaystyle\sum_{j=1}^{m-n-1} x_j^2 - \delta^{-1} \sum_{j=m-n}^{n} x_j, & \text{if } i = 1, \\ (\sqrt{2}\delta)^{-2} x_{i-1}^2 + (2\delta)^{-1} x_{i-1}, & \text{if } 2 \le i \le m-n, \\ \delta^{-1} x_{i-1}, & \text{if } m-n+1 \le i \le n+1, \\ (\sqrt{2}\delta)^{-2} x_{i-n-1}^2 - (2\delta)^{-1} x_{i-n-1}. & \text{otherwise.} \end{cases}$$

*In the formulation of $L_1$, if $m = 2n+1$, we define $\sum_{j=m-n}^{n} x_j = 0$.*

*Proof.* Let $i \in \{1, 2, \dots, m\}$ be fixed and let $L$ be a quadratic polynomial that satisfies

$$L(z^j) = \begin{cases} 1, & \text{if } j = i, & \text{(2.5.4a)} \\ 0, & \text{otherwise.} & \text{(2.5.4b)} \end{cases}$$

First, it is straightforward to verify that $L_i$ satisfies the interpolation conditions (2.5.4). Hence, it suffices to show that the Frobenius norm of its Hessian matrix is least. According to the equation (2.5.1), for any $j \in \{1, 2, \dots, m-n-1\}$, we have $z^1 = 0$, $z^{j+1} = \delta e_j$, and $z^{n+j+1} = -\delta e_j$. Therefore,

$$\begin{cases} L(z^{j+1}) = L(z^1) + \delta \nabla L(z^1)^\mathsf{T} e_j + \dfrac{\delta^2}{2} e_j^\mathsf{T} (\nabla^2 L) e_j, \\ L(z^{n+j+1}) = L(z^1) - \delta \nabla L(z^1)^\mathsf{T} e_j + \dfrac{\delta^2}{2} e_j^\mathsf{T} (\nabla^2 L) e_j, \end{cases}$$

and hence,

$$e_j^\mathsf{T} (\nabla^2 L) e_j = \frac{L(z^{j+1}) + L(z^{n+j+1}) - 2L(z^1)}{\delta^2}.$$

This fixes the first $m-n-1$ diagonal entries of $\nabla^2 L$, which are exactly those of $\nabla^2 L_i$. Since all the other entries of $\nabla^2 L_i$ are zero, we have

$$\|\nabla^2 L_i\|_{\mathsf{F}}^2 \le \|\nabla^2 L\|_{\mathsf{F}}^2,$$

which completes the proof. □

### 2.5.2 Bounds for the $\Lambda$-poisedness

The next lemma simplifies the value of $\Lambda_p$ for further computations.

**Lemma 2.5.2.** *For any $m \in \{n+2, n+3, \dots, 2n+1\}$ and any $p \in [1, \infty]$, we have*

$$\Lambda_p = \max_{x \in \mathcal{B}_p(\delta)} |L_1(x)|. \tag{2.5.5}$$

*Proof.* For each $i \in \{2, 3, \dots, n+1\}$, according to Lemma 2.5.1, $L_i(x)$ depends only on $x_{i-1}$ for all $x \in \mathbb{R}^n$, and hence

$$\max_{x \in \mathcal{B}_p(\delta)} |L_i(x)| = \max_{t \in [-\delta, \delta]} |L_i(te_{i-1})| = 1.$$

Similarly, for each $i \in \{n+2, n+3, \dots, m\}$, since $L_i(x)$ depends only on $x_{i-n-1}$ for all $x \in \mathbb{R}^n$, we have

$$\max_{x \in \mathcal{B}_p(\delta)} |L_i(x)| = \max_{t \in [-\delta, \delta]} |L_i(te_{i-n-1})| = 1.$$

Noting that $L_1(z^1) = 1$ and $z^1 \in \mathcal{B}_p(\delta)$, we thus have

$$\Lambda_p = \max_{x \in \mathcal{B}_p(\delta)} |L_1(x)|.$$

□

We are now equipped to develop bounds on $\Lambda_p$ in the general case.

**Theorem 2.5.1.** *For any $m \in \{n+2, n+3, \dots, 2n+1\}$ and any $p \in [1, \infty]$, we have*

$$1 + (2n+1-m)^{\frac{p-1}{p}} \le \Lambda_p \le n,$$

*where we assume that $0^0 = 0$ and that the lower bound is $2n+2-m$ for $p = \infty$.*

*Proof.* We will establish the bounds using the formulation of $\Lambda_p$ in Lemma 2.5.2. For the lower bound, by considering only the points in $\mathbb{R}^n$ whose leading $m-n-1$ components are zeros and whose remaining $2n+1-m$ components are equal, we have

$$\Lambda_p \geq \max_{t\in\mathbb{R}} \{1-\delta^{-1}(2n+1-m)t : (2n+1-m)^{\frac{1}{p}}|t| \leq \delta\} = 1+(2n+1-m)^{\frac{p-1}{p}}.$$

We now prove the upper bound. Note that for any $p \geq 1$, we have $\mathcal{B}_p(\delta) \subseteq \mathcal{B}_\infty(\delta)$, so that $\Lambda_p \leq \Lambda_\infty$. Therefore, we only need to show that $\Lambda_\infty \leq n$. Considering both $L_1$ and $-L_1$, we obtain

$$\Lambda_\infty = \max_{x\in\mathcal{B}_\infty(\delta)} |L_1(x)| = \max\{2n+2-m, n-1\} \leq n.$$

□

### 2.5.3 Some special cases

We now calculate the exact value of $\Lambda_p$ in some special cases.

**Proposition 2.5.1.** *For any $m \in \{n+2, n+3, \ldots, 2n+1\}$, we have*

$$\Lambda_1 = \begin{cases} 2, & \text{if } n+2 \leq m \leq 2n, & \text{(2.5.6a)} \\ 1, & \text{otherwise.} & \text{(2.5.6b)} \end{cases}$$

*Proof.* According to Theorem 2.5.1, $\Lambda_1$ is lower bounded by the right-hand side of (2.5.6). Therefore, we only need to prove that this right-hand side is also a lower bound for $\Lambda_1$, using the formulation in Lemma 2.5.2.

For any $x \in \mathcal{B}_1(\delta)$, we have

$$L_1(x) \leq 1 - \frac{1}{\delta}\sum_{j=m-n}^{n} x_j \leq 1 + \frac{1}{\delta}\sum_{j=m-n}^{n} |x_j|.$$

Therefore,

$$L_1(x) \leq \begin{cases} 2, & \text{if } n+2 \leq m \leq 2n, & \text{(2.5.7a)} \\ 1, & \text{otherwise.} & \text{(2.5.7b)} \end{cases}$$

On the other hand,

$$L_1(x) = 1 - \sum_{j=1}^{m-n-1} \frac{x_j^2}{\delta^2} - \sum_{j=m-n}^{n} \frac{x_j}{\delta} \tag{2.5.8a}$$

$$\geq 1 - \sum_{j=1}^{m-n-1} \frac{|x_j|}{\delta} - \sum_{j=m-n}^{n} \frac{|x_j|}{\delta} \tag{2.5.8b}$$

$$\geq 1 - \frac{\|x\|_1}{\delta} \geq 0. \tag{2.5.8c}$$

We conclude the proof by combining (2.5.7) and (2.5.8) with Lemma 2.5.2. □

**Proposition 2.5.2.** *For any* $m \in \{n+2, n+3, \ldots, 2n+1\}$, *we have*

$$\Lambda_2 = 1 + \sqrt{2n+1-m}. \tag{2.5.9}$$

*Proof.* If $m = 2n+1$, Lemma 2.5.1 tells us that $L_1(x) = 1 - \delta^{-2}\|x\|_2^2$ for $x \in \mathcal{B}_2(\delta)$. Therefore, Lemma 2.5.2 directly provides the desired result $\Lambda_2 = 1$. We now focus on the case with $n+2 \le m < 2n+1$.

According to Theorem 2.5.1, $\Lambda_2$ is lower bounded by the right-hand side of (2.5.9). Therefore, we only need to prove that this right-hand side is also an upper bound for $\Lambda_2$, using the formulation in Lemma 2.5.2.

For any $x \in \mathcal{B}_2(\delta)$, we have

$$L_1(x) = 1 - \frac{1}{\delta^2}\sum_{j=1}^{m-n-1} x_j^2 - \frac{1}{\delta}\sum_{j=m-n}^{n} x_j \tag{2.5.10a}$$

$$\ge 1 - \frac{1}{\delta^2}\left(\delta^2 - \sum_{j=m-n}^{n} x_j^2\right) - \frac{1}{\delta}\sum_{j=m-n}^{n} x_j = \sum_{j=m-n}^{n} \frac{x_j}{\delta}\left(\frac{x_j}{\delta} - 1\right) \tag{2.5.10b}$$

$$\ge \min_{y \in \mathcal{B}_2(1)} \sum_{j=m-n}^{n} y_j(y_j - 1). \tag{2.5.10c}$$

Let $y^* \in \mathcal{B}_2(1)$ be a minimizer in (2.5.10c). Theorem 1.3.1 ensures that there exists a Lagrange multiplier $\lambda^* \ge 0$ such that $2y_j^* - 1 + 2\lambda^* y_j^* = 0$ for all $j \in \{m-n, \ldots, n\}$. Therefore, the last $(2n+1-m)$ components of $y^*$ are equal, and hence,

$$\begin{aligned} L_1(x) \ge \min_{y \in \mathcal{B}_2(1)} \sum_{j=m-n}^{n} y_j(y_j - 1) &= \min_{t \in \mathbb{R}} \{(2n+1-m)t(t-1) : (2n+1-m)t^2 \le 1\} \\ &\ge -\sqrt{2n+1-m}. \end{aligned}$$

Furthermore,

$$L_1(x) \le 1 + \sum_{j=m-n}^{n} \frac{|x_j|}{\delta} \le 1 + \sqrt{2n+1-m} \sum_{j=m-n}^{n} \frac{x_j^2}{\delta^2} \le 1 + \sqrt{2n+1-m}.$$

Therefore, $|L_1(x)| \le 1 + \sqrt{2n+1-m}$, and hence, according to Lemma 2.5.2 and Theorem 2.5.1, we have

$$\Lambda_2 = \max_{x \in \mathcal{B}_2(\delta)} |L_1(x)| = 1 + \sqrt{2n+1-m},$$

which concludes the proof. □

Before developing the last special value of $\Lambda_p$ we introduce the following lemma.

**Lemma 2.5.3.** *For any $p \geq 1$ and $q \geq 1$, we have*

$$\max_{x \in \mathcal{B}_q(1)} \|x\|_p = \begin{cases} 1, & \text{if } p \geq q, \qquad (2.5.11a) \\ n^{\frac{q-p}{pq}}, & \text{otherwise.} \qquad (2.5.11b) \end{cases}$$

*Proof.* Let us first consider the case where $p \geq q$. For $x \in \mathcal{B}_q(1)$, we have $\|x\|_p \leq 1$, and this bound is attained at the first coordinate vector $e_1 \in \mathcal{B}_q(1)$, so that (2.5.11a) holds.

We now consider the case where $p < q$. Let $e \in \mathbb{R}^n$ be the all-one vector, $r = q/p$, and $s = r/(r-1) = q/(q-p)$. For $x \in \mathcal{B}_q(1)$, define $y = (|x_1|^p, |x_2|^p, \ldots, |x_n|^p)$. According to the Hölder inequality, we have

$$\|x\|_p = (e^{\mathsf{T}}y)^{\frac{1}{p}} \leq (\|e\|_s \|y\|_r)^{\frac{1}{p}} = n^{\frac{q-p}{pq}} \|x\|_q \leq n^{\frac{q-p}{pq}}.$$

Moreover, this bound is attained by $x = n^{-\frac{1}{q}} e$, which proves (2.5.11b). □

**Proposition 2.5.3.** *For any $p \geq 1$, if $m = 2n+1$, then*

$$\Lambda_p = \max\left\{1, n^{\frac{p-2}{p}} - 1\right\}.$$

*Proof.* It is clear that

$$\max_{x \in \mathcal{B}_p(\delta)} L_1(x) = \max_{x \in \mathcal{B}_p(\delta)} \left(1 - \frac{\|x\|_2^2}{\delta^2}\right) = 1. \qquad (2.5.12)$$

Moreover, according to Lemma 2.5.3, we have

$$\max_{x \in \mathcal{B}_p(\delta)} -L_1(x) = \max_{x \in \mathcal{B}_p(\delta)} \frac{\|x\|_2^2}{\delta^2} - 1 = \begin{cases} 0, & \text{if } p \leq 2, \qquad (2.5.13a) \\ n^{\frac{p-2}{p}} - 1, & \text{otherwise.} \qquad (2.5.13b) \end{cases}$$

The desired result is obtained by combining (2.5.12) and (2.5.13) with Lemma 2.5.2. □

### 2.5.4 Concluding remarks

Note that $\Lambda_p \geq 1$ for any $p \geq 1$. Proposition 2.5.3 shows that $\Lambda_p = 1$ for all $p \in [1, 2]$. In this sense, $\mathcal{Z}_{2n+1}$ in an optimal interpolation set in $\mathcal{B}_p^{2n+1}$ for $p \in [1, 2]$, because $\Lambda_p$ attains its lower bound. The fact that $\Lambda_p = 1$ for all $p \in \{1, 2\}$ is in fact shown also by Propositions 2.5.1 and 2.5.2.

In addition, $m = 2n+1$ is optimal for $\mathcal{Z}_m$ in the sense that $\Lambda_p$ attains the lower bound 1 in such a case, when $p \in [1, 2]$. We conjecture that $\Lambda_p$ is minimized at the value $m = 2n+1$ for any $p \geq 1$, but we do not have any proof of such a statement yet.

Note that the definition of $\mathcal{Z}_m$ in (2.5.1) and (2.5.2) assumes that $m \leq 2n+1$. Pow-ell [166] proposed an extension of $\mathcal{Z}_m$ in the case $m > 2n+1$, and hence, $\mathcal{Z}_m$ is defined for any $m \in \{1, 2, \dots, (n+1)(n+2)/2\}$. With such an extension, it is interesting to ask whether $m = 2n+1$ still minimizes $\Lambda_p$, which seems to be true according to some numerical experiments. We leave this problem open and expect the analysis to be more challenging than what we have done. One of the challenges is that Lemma 2.5.2 does not hold in the case $m > 2n+1$, and hence, the calculations will become more involving.

## 2.6 Summary and remarks

In this chapter, we presented concepts of polynomial interpolation that we will need later in this thesis when discussing model-based DFO methods.

We first introduced fully-determined linear and quadratic interpolation. We presented three concepts, namely the poisedness of interpolation sets, the quadratic Lagrange polynomials on poised interpolation sets, and the well-poisedness of interpolation sets.

We then introduced underdetermined quadratic interpolation. More specifically, we focused on least Frobenius quadratic models and on least Frobenius updating of quadratic models. Modern model-based DFO methods usually employ these models because they require only a few interpolation points. We presented the extended concepts of poisedness, quadratic Lagrange polynomials, and well-poisedness to underdetermined interpolation. We also presented how Powell proposed implementing such underdetermined quadratic models in DFO methods. In particular, when a method updates one point at a time in the interpolation set, it is possible and advantageous to maintain the inverse of the coefficient matrix of the interpolation problem, because modifying a point in the interpolation sets leads to an at-most rank-2 modification of this matrix. We will use this mechanism in our new DFO presented in Chapter 5, named COBYQA.

Finally, we studied in this chapter the well-poisedness of a particular interpolation set proposed by Powell for quadratic underdetermined interpolation, whose number of points varies between $n+2$ and $2n+1$. We will use a variation of this set as the initial interpolation set in COBYQA. We showed in particular that under some conditions, the optimal number of interpolation points in this set is $2n+1$. Powell uses this interpolation set in NEWUOA and encourages users to employ $2n+1$ points. This default value is then supported by our theory.

# 3 Development of the PDFO package

In this chapter, we present PDFO (*Powell's Derivative-Free Optimization solvers*), which is a MATLAB and Python package we develop to interface Powell's DFO solvers. These solvers were implemented by Powell in Fortran 77. The motivation is to provide user-friendly interfaces to them, so that users do not need to deal with the Fortran code. This package has been downloaded more than 30,000 times[1] as of August 2022.

Section 3.1 presents the motivations for the development of PDFO. We then provide an overview of Powell's five model-based methods in Section 3.2. Section 3.3 introduces the core features of PDFO and its implementation. Finally, we present in Section 3.4 some numerical experiments on Powell's solvers using PDFO. One of these experiments will reveal an interesting behavior of the solvers, namely that COBYLA performs quite well compared to other solvers on noisy problems. This behavior is intriguing, because COBYLA is the oldest of these solvers and uses the simplest models.

## 3.1 Introduction and motivation

Powell designed five model-based DFO solvers to tackle unconstrained and constrained problems, namely COBYLA [161], UOBYQA [163], NEWUOA [166], BOBYQA [168], and LINCOA. These solvers were implemented by Powell, with particular attention paid to their numerical stability and algebraic complexity. Renowned for their robustness and efficiency, these solvers are very appealing to practitioners and widely used in applications, including aeronautical engineering [78], astronomy [18, 128], computer vision [111], robotics [136], and statistics [12]. However, Powell implemented them in Fortran 77, an old-fashion language that damps the enthusiasm of many users to exploit these solvers in their projects.

There has been considerable demand from both researchers and practitioners for the availability of the five Powell's solvers in more user-friendly languages such as MATLAB, Python, Julia, and R. Our aim is to wrap Powell's Fortran code into a package named PDFO, which enables users of such languages to call Powell's solvers without any need to deal with the Fortran code. For each supported language, PDFO provides a simple subroutine that can invoke one of Powell's solvers according to the user's request (if any) or according to the type of the problem to solve. The current release (version 1.2) of PDFO supports MATLAB and Python, with more languages to be covered in the future. The signature of the MATLAB subroutine is consistent with the `fmincon` function of the MATLAB Optimization Toolbox; the signature of the Python subroutine is consis-

[1] See `https://www.pdfo.net/`.

tent with the `minimize` function of the SciPy optimization library [204]. The package is cross-platform, available on Linux, macOS, and Windows at once. It is open-source and distributed under the LGPLv3+ license. The source code of PDFO is available at

`https://github.com/pdfo/pdfo`.

A complete documentation for PDFO is available at

`https://www.pdfo.net/`.

PDFO is not the first attempt to facilitate the usage of Powell's solvers in languages other than Fortran. Various efforts have been made in this direction in response to the continual demands from both researchers and practitioners: Py-BOBYQA [33] provides a Python implementation of BOBYQA; NLopt [112] includes multi-language interfaces for COBYLA, NEWUOA, and BOBYQA; minqa [13] wraps UOBYQA, NEWUOA, and BOBYQA in R; SciPy [204] makes COBYLA available in Python under its optimization library. Nevertheless, PDFO has several features that distinguish it from others.

1. *Comprehensiveness*. To the best of our knowledge, PDFO is the only package that provides all of COBYLA, UOBYQA, NEWUOA, BOBYQA, and LINCOA with a uniform interface. In addition to homogenizing the usage, such an interface eases the comparison between these solvers in case multiple of them can tackle a given problem. By doing so, we may gain insights that cannot be obtained otherwise into the behavior of the solvers, as will be illustrated in Section 3.4.
2. *Solver selection*. When using PDFO, the user can specifically invoke one of Powell's solvers; nevertheless, if the user does not specify any solver, PDFO will automatically select a solver according to the given problem. The selection takes into consideration the performance of the solvers on the CUTEst problem set [90]. This will be elaborated on in Section 3.3.3.
3. *Code patching*. During the development of PDFO, we spotted in the original Fortran code some bugs, which led, for example, to infinite cycling or segmentation faults on some ill-conditioned problems. The bugs have been patched in PDFO. Nevertheless, we provide an option that can enforce the package to use the original code of Powell without the patches, which is not recommended except for research. In addition, PDFO provides COBYLA in double precision, whereas Powell used single precision when he implemented it in the 1990s. See Section 3.3.4 for details.
4. *Fault tolerance*. PDFO takes care of failures in the evaluations of the objective or constraint functions when NaN or infinite values are returned. In case of such failures, PDFO will not exit but will try to progress. Moreover, PDFO ensures that the returned solution is not a point where the evaluation fails, while the original code of Powell may return a point whose objective function value is numerically NaN.

5. *Problem preprocessing*. PDFO preprocesses the inputs to simplify the problem and reformulate it to meet the requirements of Powell's solvers. For instance, if the problem has linear constraints $A_{\mathcal{E}}x = b_{\mathcal{E}}$, PDFO can rewrite it into a problem on the null space of $A_{\mathcal{E}}$, eliminating such constraints and reducing the dimension. Another example is that the starting point of a linearly constrained problem is projected onto the feasible region because LINCOA needs a feasible starting point to work properly.

6. *Additional options*. PDFO includes options for the user to control the solvers in some manners that are useful in practice. For example, the user can request PDFO to scale the problem according to the bounds of the variables before solving it.

## 3.2 Overview of Powell's derivative-free optimization methods

To present Powell's DFO methods, we consider the optimization problem

$$\min_{x\in\mathbb{R}^n} \quad f(x) \tag{3.2.1a}$$

$$\text{s.t.} \quad c_i(x) \le 0, \; i \in \mathcal{I}, \tag{3.2.1b}$$

where $f$ and $c_i$, with $i \in \mathcal{I}$, are real-valued functions on $\mathbb{R}^n$, and where the set of indices $\mathcal{I}$ is finite, perhaps empty. Note that this problem is exactly the problem (1.1.1) studied in Chapter 1 with $\mathcal{E} = \emptyset$. In general, any problem of the form (1.1.1) can be reformulated into (3.2.1) by considering the equalities as two inequalities. We assume that derivatives of $f$ and $c_i$, for $i \in \mathcal{I}$, are unavailable.

Powell published his first DFO method based on conjugate directions[2] in 1964 [153]. His code for this algorithm is contained in the HSL Mathematical Software Library [110] as the subroutine `VA24`. It is not included in PDFO because the code is not in the public domain, although open-source implementations are available (see [49, Fn. 4]).

From the 1990s to the final days of his career, Powell developed five model-based DFO algorithms to solve the problem (3.2.1) in different cases, namely COBYLA [161] for nonlinearly constrained problems, UOBYQA [163] and NEWUOA [166] for unconstrained problems, BOBYQA [168] for bound constrained problems, and LINCOA for linearly constrained problems. In addition, Powell implemented these algorithms into Fortran solvers and made the code publicly available. The solvers constitute the cornerstones of PDFO.

[2]According to Google Scholar, this is Powell's second published paper and also the second most cited work. The earliest and meanwhile most cited one is his paper on the *Davidon-Fletcher-Powell* (DFP) method [74], co-authored with Fletcher and published in 1963.

### 3.2.1 A sketch of the algorithms

Powell's model-based DFO algorithms are trust-region methods. At the $k$th trust-region iteration, the algorithms construct a linear or quadratic model $\hat{f}^k$ for the objective function $f$ to meet the interpolation condition

$$\hat{f}^k(y) = f(y), \quad \text{for } y \in \mathcal{Y}^k, \tag{3.2.2}$$

where $\mathcal{Y}^k \subseteq \mathbb{R}^n$ is a finite interpolation set updated along the iterations. COBYLA models the constraints by interpolants on $\mathcal{Y}^k$ as well. With such models, the algorithms form a trust-region subproblem and solve it to find a trial point. Then the interpolation set is updated by replacing one of the existing interpolation points with the trial point, unless such a replacement turns out inappropriate according to certain criteria. In addition, the algorithms may take a step to improve the geometry of the interpolation set before the next trust-region iteration. Instead of repeating Powell's description of these algorithms, we provide a sketch of them in the sequel.

In all five algorithms, the $k$th iteration places the trust-region center $x^k$ at the "best" point where the objective function and constraints have been evaluated so far. Such a point is selected according to the objective function or a merit function that takes the constraints into account. After choosing the trust-region center $x^k$, with the trust-region model $\hat{f}^k$ constructed, a trial point is then obtained by solving approximately the trust-region subproblem

$$\min_{x \in \mathbb{R}^n} \quad \hat{f}^k(x) \tag{3.2.3a}$$

$$\text{s.t.} \quad c_i(x) \le 0, \; i \in \mathcal{I}, \tag{3.2.3b}$$

$$\|x - x^k\| \le \Delta^k. \tag{3.2.3c}$$

where $\Delta^k$ is the trust-region radius. An exception should be made for COBYLA, which replaces (3.2.3b) by

$$\hat{c}_i^k(x) \le 0, \; i \in \mathcal{I}.$$

The other solvers do not need models for the constraints, because either $\mathcal{I} = \emptyset$, or the constraints are simple constraints that do not necessitate models, namely bounds and linear constraints.

### 3.2.2 The COBYLA method

Published in 1994, COBYLA was the first model-based DFO solver of Powell. The solver is named after *Constrained Optimization BY Linear Approximations*. It aims to solve the optimization problem (3.2.1) whenever the constraint functions $c_i$, for $i \in \mathcal{I}$, are nonlinear functions whose derivatives are unknown. In other words, only function values of the constraint functions are accessible.

At the $k$th iteration, COBYLA models the objective and the constraint functions with linear interpolants on the interpolation set $\mathcal{Y}^k \subseteq \mathbb{R}^n$, which consists of $n+1$ points that are updated along the iterations.

Once the linear models $\hat{c}_i^k$ of the constraint functions $c_i$, for $i \in \mathcal{I}$, are built, the trust-region subproblem

$$\min_{x \in \mathbb{R}^n} \quad \hat{f}^k(x) \tag{3.2.4a}$$

$$\text{s.t.} \quad \hat{c}_i^k(x) \le 0, \; i \in \mathcal{I}, \tag{3.2.4b}$$

$$\|x - x^k\| \le \Delta^k, \tag{3.2.4c}$$

needs to be handled. This problem is solved by imagining that $\Delta^k$ is replaced with a constant continuously increasing from zero to $\Delta^k$, which would generate a piecewise linear path from $x^k$ to the solution of this problem. To locate this solution, the trust-region subproblem solver of COBYLA follows this path by updating the active sets of the linear constraints (3.2.4b). However, the linear constraints (3.2.4b) and the trust-region constraint (3.2.4c) may contradict each others, in which case the trial point is chosen to solve approximately

$$\begin{aligned} \min_{x \in \mathbb{R}^n} \quad & \max_{i \in \mathcal{I}} [\hat{c}_i^k(x)]_+ \\ \text{s.t.} \quad & \|x - x^k\| \le \Delta^k, \end{aligned}$$

where $[\cdot]_+$ takes the positive-part of a given number. In doing so, the method attempts to reduce the $\ell_\infty$-constraint violation of the linearized constraints within the trust region.

As we already mentioned, it is essential to maintain a good geometry of $\mathcal{Y}^k$ to ensure the accuracy of the models. When the geometry of $\mathcal{Y}^k$ turns out inadequate for producing accurate models, COBYLA removes a point from $\mathcal{Y}^k$, and adds another one chosen on the direction perpendicular to the face of $\mathcal{Y}^k$ (regarded as a simplex) that is to the opposite of the removed point. This replacement tends to increase the volume of the simplex engendered by the interpolation set, and hence, improves the conditioning of the interpolation system (3.2.2).

### 3.2.3 The UOBYQA method

In 2002, Powell published UOBYQA [163], named after *Unconstrained Optimization BY Quadratic Approximations*. It aims at solving the nonlinear optimization problem (3.2.1) in the unconstrained case, i.e., when $\mathcal{I} = \emptyset$.

At the $k$th iteration, UOBYQA models the objective function with a quadratic obtained by fully-determined interpolation on the set $\mathcal{Y}^k \subseteq \mathbb{R}^n$ containing $(n+1)(n+2)/2$ points. Then a trial point is obtained by solving the trust-region subproblem (3.2.3) with the Moré-Sorensen algorithm [137]. Further, UOBYQA replaces an interpolation point in $\mathcal{Y}^k$ by the trial point, obtaining $\mathcal{Y}^{k+1}$, unless this replacement seems to be unreasonable for maintaining a good geometry of $\mathcal{Y}^{k+1}$.

UOBYQA undertakes a geometry-improving step when the model seem not to be accurate. It first removes a point $y \in \mathcal{Y}^k$, and builds $\mathcal{Y}^{k+1}$ by including a solution to

$$\max_{x \in \mathbb{R}^n} \quad |L_y(x)| \tag{3.2.5a}$$

$$\text{s.t.} \quad \|x - x^k\| \leq \bar{\Delta}^k, \tag{3.2.5b}$$

where $L_y \colon \mathbb{R}^n \to \mathbb{R}$ denotes the Lagrange function associated with $y$ (see Section 2.2.2 for details), and $\bar{\Delta}^k > 0$ is defined according to the value of $\Delta^k$ and the diameter of $\mathcal{Y}^k$. This subproblem is motivated by the analysis in [162, § 2]. Powell developed an algorithm to solve (3.2.5) approximately, based on an estimation of $|L_y(\cdot)|$. Once such a point is calculated, the solver constructs $\mathcal{Y}^{k+1}$ and continues with a trust-region step.

### 3.2.4 The NEWUOA method

NEWUOA [166] was designed as a successor to UOBYQA. Using the derivative-free symmetric Broyden update presented in Section 2.4.2, it can build quadratic models using much fewer interpolation points than UOBYQA. Hence, it is able to solve much larger problems.

Powell did not mention the meaning of the acronym NEWUOA, but we speculate that it means *NEW Unconstrained Optimization Algorithm*. In addition to the symmetric Broyden update, NEWUOA also differs from UOBYQA by the way that it solves its trust-region subproblem (3.2.3). Instead of the Moré-Sorensen algorithm, this subproblem is handled using the Steihaug-Toint *truncated conjugate gradient* (TCG) method [196, 198] (see Algorithm 6.1) followed by a two-dimensional search.

The geometry-improving mechanism of NEWUOA first chooses a point $y \in \mathcal{Y}^k$ to remove from the interpolation set, and then selects a new interpolation point by two approaches. The first one, similar to that of UOBYQA, finds this point by solving approximately (3.2.5). To understand the second approach, recall that the implementation of the derivative-free symmetric Broyden update stores the inverse of the coefficient matrix of the interpolation problem (see Section 2.4.3). Whenever an interpolation point is replaced with a new one, this inverse will be updated by the formula detailed in [165, Eq. (2.12)], which has a scalar denominator. The second mechanism maximizes this denominator within the trust region (3.2.5b), and is employed only if the first one does not render a sufficiently big absolute value of this denominator.

### 3.2.5 The BOBYQA method

BOBYQA [168] is named after *Bound Optimization QY Quadratic Approximations*. It also builds quadratic model using the derivative-free symmetric Broyden update. An im-portant feature of BOBYQA is that the bound constraints are always respected by each iterate and each interpolation point. Therefore, both its trust-region and geometry-improving subproblems incorporate the bound constraints. Its trust-region subproblem solver is an active-set variation of the TCG method, which will be detailed in Algorithm 6.2. The geometry improving procedure of BOBYQA consists in estimating two approximate solutions to (3.2.5) subject to the bounds, and choosing the better one. Details on this procedure will be given in Section 6.5.

### 3.2.6 The LINCOA method

LINCOA is named after *LINearly Constrained Optimization Algorithm*. It also builds quadratic model using the derivative-free symmetric Broyden update. Note that Powell never published a paper introducing LINCOA, and [171] discusses only its trust-region subproblem solver. This solver is again an active-set variation of the TCG method, which will be detailed in Algorithm 6.3.

The geometry-improving procedure of LINCOA also approximately solves (3.2.5). Note that it may violate the linear constraints. LINCOA calculates three approximate solutions to (3.2.5), namely

1. the point that maximizes $|L_y(\cdot)|$ within the trust region on the lines through $x^k$ and other points in $\mathcal{Y}^k$,
2. a point obtained by a gradient step that maximizes $|L_y(\cdot)|$ within the trust region, and
3. a point obtained by a projected gradient step that maximizes $|L_y(\cdot)|$ within the trust region, the projection being made onto the null space of the constraints that are considered active at $x^k$.

The procedure first selects the point among the first two alternatives that provide the larger value of $|L_y(\cdot)|$. Further, this point is replaced with the third alternative if the latter provides a value of $|L_y(\cdot)|$ that is not too small compared with the above one, while being either feasible or with a constraint violation that is at least $0.2\bar{\Delta}^k$.

## 3.3 Core features of the PDFO package

In this section, we detail the main features of PDFO, including the signature of the main function, the problem preprocessing, the solver selection, and the bug fixes. Before starting, we emphasize that PDFO does not reimplement Powell's solvers but rather interfaces MATLAB and Python with the Fortran source code, using MEX and F2PY [150], respectively.

### 3.3.1 Signature of the main function

The philosophy of PDFO is simple: providing to users a single function to solve a DFO problem. It takes for input an optimization problem of the form

$$\min_{x\in\mathbb{R}^n} \quad f(x) \tag{3.3.1a}$$
$$\text{s.t.} \quad l \le x \le u, \tag{3.3.1b}$$
$$A_{\mathcal{I}} x \le b_{\mathcal{I}}, \tag{3.3.1c}$$
$$A_{\mathcal{E}} x = b_{\mathcal{E}}, \tag{3.3.1d}$$
$$c_i(x) \le 0, \; i \in \mathcal{I}, \tag{3.3.1e}$$
$$c_i(x) = 0, \; i \in \mathcal{E}, \tag{3.3.1f}$$

where $l, u \in (\mathbb{R} \cup \{\pm\infty\})^n$, $A_{\mathcal{I}}$ and $A_{\mathcal{E}}$ are real matrices, $b_{\mathcal{I}}$ and $b_{\mathcal{E}}$ are real vectors, and $c_i$, for $i \in \mathcal{I} \cup \mathcal{E}$, is a real-valued function. A simple example of usage is shown for MATLAB in Listing 3.1 and Python in Listing 3.2, where variable names have clear correspondences with the problem (3.3.1). For both MATLAB and Python, PDFO returns the best point calculated and (optionally) the corresponding optimal value. Additional information that describes the backend calculations can also be returned.

Listing 3.1: An elementary example of using PDFO in MATLAB

```
x = pdfo(@fun, x0, A, b, Aeq, beq, lb, ub, @nonlcon);

function fx = fun(x)
...
return
end

function [c, ceq] = nonlcon(x)
...
return
end
```

Listing 3.2: An elementary example of using PDFO in Python

```python
import numpy as np
from pdfo import *

def fun(x):
    return ...

def cub(x):
    return ...

def ceq(x):
    return ...

bounds = Bounds(lb, ub)
constraints = [
    LinearConstraint(A, -np.inf, b),
    LinearConstraint(Aeq, beq, beq),
    NonlinearConstraint(cub, -np.inf, 0.0),
    NonlinearConstraint(ceq, 0.0, 0.0),
]

res = pdfo(fun, x0, bounds=bounds, constraints=constraints)
```

### 3.3.2 Problem preprocessing

The package PDFO preprocesses the arguments provided by the user, detects the type of the problem, and then invokes the Powell's solver that best matches the given problem.

A crucial point of BOBYQA and LINCOA is that they require the initial guess to be feasible (LINCOA would otherwise increase the right-hand side of the linear constraints to make the initial guess feasible). Therefore, PDFO attempts to project the provided initial guess onto the feasible set.

Another noticeable preprocessing of the constraints made by PDFO is the treatment of the linear equality constraints (3.3.1d). As long as these constraints are consistent, we reformulate (3.3.1) into an $(n - \operatorname{rank} A_{\mathcal{E}})$-dimensional problem by eliminating these linear equality constraints. This is done using a QR factorization of $A_{\mathcal{E}}$.

### 3.3.3 Automatic selection of the solver

Another main feature of PDFO is its solver selection. When a problem is received, unless the user specifies the solver to use, PDFO selects a solver as follows.

1. If the problem is unconstrained, then UOBYQA is selected when $2 \le n \le 8$, and NEWUOA is selected when $n = 1$ or $n > 8$.
2. If the problem is bound-constrained, then BOBYQA is selected.
3. If the problem is linearly constrained, then LINCOA is selected.
4. Otherwise, COBYLA is selected.

For the unconstrained case, we select UOBYQA for small problems because it is more efficient, and the number $8$ is set according to our experiments on the CUTEst problems. We note that Powell's implementation of UOBYQA cannot handle problems with univariate objective functions, for which NEWUOA is invoked.

We also select a solver in this way if the user specifies a solver that is incapable of solving the problem received.

### 3.3.4 Bug fixes in the Fortran source code

The current version of PDFO also fixes several bugs in the Fortran source code, particularly the following ones.

1. The solvers may encounter infinite loops. This is because the exit conditions of some loops cannot be met because variables involved in these conditions become NaN values due to floating point exceptions.
2. The Fortran code may encounter memory errors due to uninitialized variables that are used as indices. This is because some variables are initialized according to conditions that can never be met due to NaN values, similar to the previous case.
3. COBYLA may not return the best point that is evaluated; sometimes, it returns a point with a large constraint violation, even though the initial guess is feasible. This is because COBYLA may discard points that are not considered good according to the current merit function, depending on a penalty parameter. However, when this penalty parameter is updated, a discarded point may turn out to be the best point according the updated merit function.

In PDFO, these bugs are fixed.

## 3.4 Numerical experiments

### 3.4.1 Comparison on the CUTEst library

With PDFO, we can easily compare the behaviors of all Powell's DFO solvers when multiple of them can solve the same problems. As an example, we compare the performance of these solvers on unconstrained problems, with and without noise.

The starting points of the problems are set to the default ones provided in CUTEst. The initial trust-region radius $\Delta^0$ is set to one, and the final value $\delta^\infty$ for the lower bound on the trust-region radius is set to $10^{-6}$. The maximal number of function evaluations is $500n$, where $n$ denotes the dimension of the problem being solved. For the methods that employ underdetermined quadratic interpolation models, the number of interpolation points is set to $2n + 1$.

Performance profiles on unconstrained problems of dimensions at most $10$ and $50$ are provided respectively in Figures 3.1a and 3.1b. According to Figure 3.1a, UOBYQA performs better than all the other solvers on small problems. This can be explained by the fact that it uses quadratic models obtained by fully-determined interpolation. However, we excluded UOBYQA from the second experiment, because the execution time was excessively long on problems with moderately big dimensions, due to the fully-determined interpolation that it does. Moreover, COBYLA is always outperformed by all other solvers. This is because it uses only linear models to approximate the objective and constraint functions of the problems, which are not as precise as the quadratic models employed by other solvers.

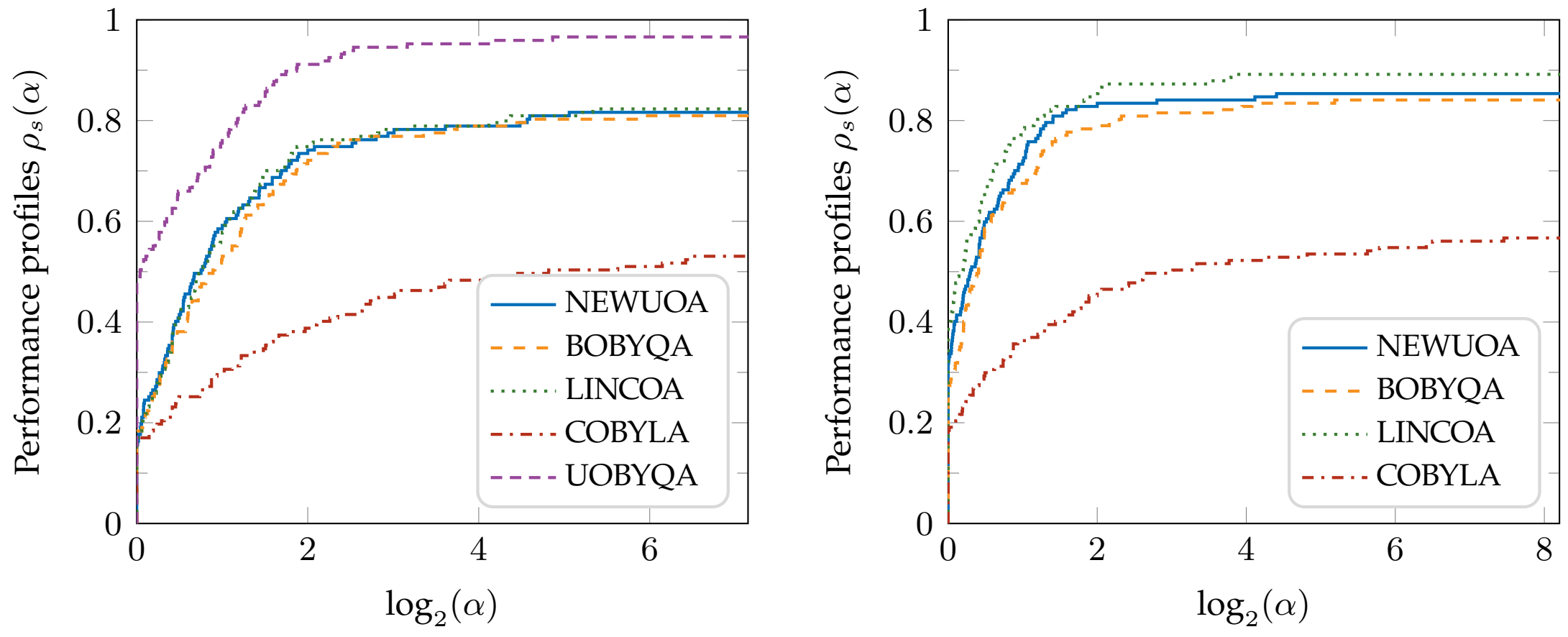

(a) Dimension at most $10$ (b) Dimension at most $50$

Figure 3.1: Performance profiles of Powell's DFO solvers on unconstrained problems with $\tau = 10^{-4}$

We now consider the experiment on noisy unconstrained problems described in Section 1.5.4. We take $\epsilon(x) \sim N(0, \sigma^2)$ with $\sigma = 10^{-2}$ in (1.5.5). Figure 3.2 presents the performance profiles on the unconstrained problems of dimension at most 50 with an error term $\sigma = 10^{-2}$.

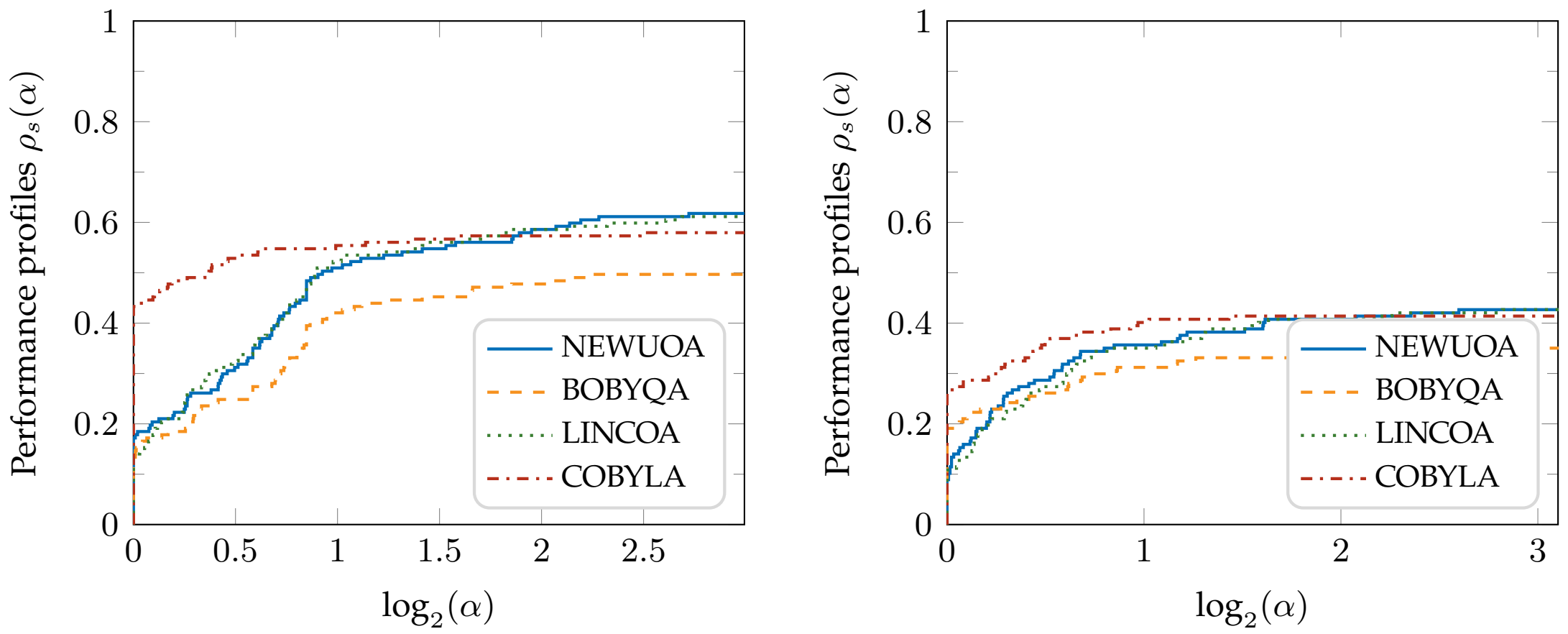

(a) Tolerance $\tau = 10^{-1}$ (b) Tolerance $\tau = 10^{-2}$

Figure 3.2: Performance profiles of Powell's DFO solvers on noisy problems

It is interesting to observe the performance of COBYLA in this experiment. Even though it is not particularly designed for such kind of problem and uses the simplest models, we observe that COBYLA defeats all other solvers on unconstrained problems, for $\tau = 10^{-1}$. It seems that the linear models of COBYLA are more robust to noise, but we have not yet derived a theory for this behavior.

### 3.4.2 An example of hyperparameter tuning problem

We now consider the hyperparameter tuning problem of an SVM described in detail in Section 1.2.3. We will compare the performance of PDFO with a Bayesian optimization method named TPE and a *random search* (RS) method. Since the hyperparameter tuning problem is a bound-constrained problem, PDFO selects BOBYQA as the solver. For the TPE and the RS methods, we use the Python package hyperopt [15].

The initial guess for the hyperparameters are the default values provided by the `SVC` class of scikit-learn [148]. For PDFO, the initial trust-region radius is 1, and the maximal number of function evaluations is 100. The maximal number of function evaluations we attempt for TPE and RS are 100 and 300.

Our experiments are based on binary classification problems from the LIBSVM datasets[3]. Table 3.1 provides a description of the datasets. The objective function of the hyperparameter tuning is given in Algorithm 1.3, where the considered performance measure is the AUC validation score on the testing dataset.

The results of this experiment are provided in Tables 3.2 to 3.5. For both the AUC score and the accuracy, the higher, the better. In terms of AUC score and accuracy, we observe that PDFO achieved a clearly better result than TPE and RS on the "splice" dataset, and they all attain comparable results on all the other datasets. However, PDFO always uses much less function evaluations, and hence, much less computation time. The difference in the computation time is particularly visible on the dataset "ijcnn1" in Table 3.5, as the size of this dataset is large, so that each function evaluation takes much time. In summary, we can conclude that PDFO performs better than TPE and RS on these problems.

[3] Available at `https://www.csie.ntu.edu.tw/~cjlin/libsvmtools/datasets/`.

Table 3.1: Description of LIBSVM datasets

| Dataset $\mathcal{D}$ | Attribute characteristic | Dimension $d$ | Dataset size |
|---|---|---|---|
| splice | $[-1, 1]$, scaled | 60 | 1000 |
| svmguide1 | $[-1, 1]$, scaled | 4 | 3088 |
| svmguide3 | $[-1, 1]$, scaled | 21 | 1242 |
| ijcnn1 | $[-1, 1]$ | 22 | 49,990 |

Table 3.2: Hyperparameter tuning problem on the dataset "splice"

| Solver | No. eval. | AUC score $(10^{-1})$ | Accuracy $(10^{-1})$ | Execution time (s) |
|---|---|---|---|---|
| PDFO | 65 | 9.568 | 9.933 | 3.697 |
| RS | 100 | 6.409 | 5.300 | 4.635 |
| RS | 300 | 7.880 | 5.300 | 13.763 |
| TPE | 100 | 5.000 | 5.033 | 4.889 |
| TPE | 300 | 7.736 | 5.300 | 15.726 |

Table 3.3: Hyperparameter tuning problem on the dataset "svmguide1"

| Solver | No. eval. | AUC score $(10^{-1})$ | Accuracy $(10^{-1})$ | Execution time (s) |
|---|---|---|---|---|
| PDFO | 68 | 9.966 | 9.730 | 4.906 |
| RS | 100 | 9.966 | 9.676 | 16.178 |
| RS | 300 | 9.966 | 9.676 | 48.404 |
| TPE | 100 | 9.966 | 9.720 | 13.057 |
| TPE | 300 | 9.966 | 9.720 | 33.392 |

Table 3.4: Hyperparameter tuning problem on the dataset "svmguide3"

| Solver | No. eval. | AUC score $(10^{-1})$ | Accuracy $(10^{-1})$ | Execution time (s) |
|---|---|---|---|---|
| PDFO | 68 | 8.241 | 8.016 | 2.793 |
| RS | 100 | 8.025 | 7.882 | 4.233 |
| RS | 300 | 8.141 | 7.775 | 12.414 |
| TPE | 100 | 7.774 | 7.453 | 4.197 |
| TPE | 300 | 8.106 | 7.989 | 12.912 |

Table 3.5: Hyperparameter tuning problem on the dataset "ijcnn1"

| Solver | No. eval. | AUC score $(10^{-1})$ | Accuracy $(10^{-1})$ | Execution time $(10^3\text{s})$ |
|---|---|---|---|---|
| PDFO | 59 | 9.940 | 9.819 | 1.892 |
| RS | 100 | 9.886 | 9.773 | 4.435 |
| RS | 300 | 9.886 | 9.773 | 13.251 |
| TPE | 100 | 9.891 | 9.791 | 4.426 |
| TPE | 300 | 9.896 | 9.786 | 12.552 |

## 3.5 Summary and remarks

We presented in this chapter PDFO, our package for interfacing Powell's five model-based DFO solvers with MATLAB and Python. To the best of our knowledge, it is the only package that provides all the five solvers with a uniform interface.

We first presented an overview of the Powell's algorithms for these solvers. In a nutshell, they are all trust-region DFO methods that use either linear or quadratic models. We then presented the core features of PDFO, including the general structure of the main function, the preprocessing of the problems, the automatic selection of the solver, and several bug fixes in the Fortran code.

Using PDFO, we provided some numerical experiments on the Powell's solvers. One of these experiments demonstrated an intriguing behavior of COBYLA. It performs quite well compared to the other four solvers on noisy unconstrained problems, even though it uses the simplest models. We plan to study this phenomenon in the future from a theoretical point of view. We also tested PDFO on a hyperparameter tuning problem of an SVM, for which PDFO outperformed two standard solvers.

Powell implemented his DFO solvers in Fortran 77. Due to the features of the language and the unique coding style, this sophisticated implementation is highly non-trivial to maintain and to extend. We are currently reimplementing these solvers using modern languages, in a structured and modularized way, so that the code is readable, maintainable, and extendable. The modernized code will be incorporated into future releases of PDFO.

# The SQP method — an overview and some perspectives

# 4

This chapter introduces the well-known SQP method, assuming that the derivatives of the objective and constraint functions are available. We first present the general SQP framework and some insights into its subproblem from four different angles in Section 4.1. Section 4.2 then introduces the concept of merit function used in globalization strategies for the SQP method. Further, Section 4.3 focuses on the trust-region SQP framework and some composite-step methods for solving its subproblem, referred to as the Byrd-Omojokun, Vardi, and *Celis-Dennis-Tapia* (CDT) approaches. Section 4.4 briefly mentions some well-known software based on the SQP method.

Most of the materials in this chapter are not new, and can be found for example in [20, 92], [44, Ch. 15], [197, Ch. 12], [144, Ch. 18], and [190, 86]. However, we highlight Theorems 4.1.1 and 4.1.4 and Proposition 4.3.1, which are new to the best of our knowledge, although the last two are easy to establish. Theorem 4.1.1 interprets the objective function of the SQP subproblem as a quadratic approximation of the original objective function in the tangent space of a surface. This explains in a new way why the SQP subproblem should not minimize second-order Taylor expansion of the original objective function unless the constraints are all linear. Theorem 4.1.4 shows that the augmented Lagrangian of the SQP subproblem is a particular quadratic approximation of the augmented Lagrangian used by Yuan et al. [143, 205]. Finally, Proposition 4.3.1 interprets the Vardi approach as an approximation to the Byrd-Omojokun approach.

In addition, we also emphasize our Byrd-Omojokun approach for the trust-region SQP subproblem in the presence of inequality constraints. This approach is different from the one given in [44, § 15.4.4], and provides evidently better performance in numerical experiments.

## 4.1 The method

Throughout this chapter, we consider a problem of the form

$$\min_{x \in \mathbb{R}^n} \quad f(x) \tag{4.1.1a}$$

$$\text{s.t.} \quad c_i(x) \le 0, \; i \in \mathcal{I}, \tag{4.1.1b}$$

$$c_i(x) = 0, \; i \in \mathcal{E}, \tag{4.1.1c}$$

Note that the problem (4.1.1) is precisely the problem (1.1.1) discussed in Chapter 1; hence, all the theory mentioned is applicable. For our later discussion, recall in particular that the Lagrangian of this problem is defined by

$$\mathcal{L}(x,\lambda) \overset{\text{def}}{=} f(x) + \sum_{i\in\mathcal{I}\cup\mathcal{E}} \lambda_i c_i(x), \quad \text{for } x\in\mathbb{R}^n \text{ and } \lambda_i\in\mathbb{R}, \text{ with } i\in\mathcal{I}\cup\mathcal{E},$$

where $\lambda = [\lambda_i]^\top_{i\in\mathcal{I}\cup\mathcal{E}}$ is the dual variable of the considered problem.

### 4.1.1 Overview of the method

The SQP method is known to be one of the most powerful methods for solving the problem (4.1.1) when derivatives of $f$ and $c_i$, with $i\in\mathcal{I}\cup\mathcal{E}$, are available. The classical SQP framework is presented in Algorithm 4.1.

**Algorithm 4.1:** Classical SQP method

**Data:** Objective function $f$, constraint functions $\{c_i\}_{i\in\mathcal{I}\cup\mathcal{E}}$, initial guess $x^0\in\mathbb{R}^n$, and estimated Lagrange multiplier $\lambda^0 = [\lambda_i^0]^\top_{i\in\mathcal{I}\cup\mathcal{E}}$.

**for** $k = 0, 1, \ldots$ **do**

  Define $H^k \approx \nabla^2_{x,x}\mathcal{L}(x^k,\lambda^k)$

  Generate a step $d^k\in\mathbb{R}^n$ by solving approximately

$$\min_{x\in\mathbb{R}^n} \quad \nabla f(x^k)^\top d + \frac{1}{2} d^\top H^k d \tag{4.1.2a}$$

$$\text{s.t.} \quad c_i(x^k) + \nabla c_i(x^k)^\top d \le 0,\ i\in\mathcal{I}, \tag{4.1.2b}$$

$$c_i(x^k) + \nabla c_i(x^k)^\top d = 0,\ i\in\mathcal{E} \tag{4.1.2c}$$

  Update the iterate $x^{k+1} \leftarrow x^k + d^k$

  Estimate the Lagrange multiplier $\lambda^{k+1} = [\lambda_i^{k+1}]^\top_{i\in\mathcal{I}\cup\mathcal{E}}$

**end for**

The earliest reference to such a method appeared in the Ph.D. thesis of Wilson [210], with $H^k = \nabla^2_{x,x}\mathcal{L}(x^k,\lambda^k)$. Robinson [181] showed the local R-quadratic convergence rate of this method. Later, Garcìa-Palomares and Mangasarian [79, 80] modified it using a quasi-Newton update for calculating $H^k$ and established a local R-superlinear convergence rate for such an algorithm. A similar method was introduced by Han [99, 100], but he only approximated $\nabla^2_{x,x}\mathcal{L}(x^k,\lambda^k)$, while Garcìa-Palomares and Mangasarian applied quasi-Newton approximations to the whole matrix $\nabla^2\mathcal{L}(x^k,\lambda^k)$. In addition, Han introduced a line-search strategy to guarantee the global convergence and local Q-superlinear convergence rate, requiring that $\nabla^2_{x,x}\mathcal{L}(x^*,\lambda^*)$ is positive definite at the solution $(x^*,\lambda^*)$. Powell [158, 157, 159] studied the method in the same direction. In particular, he proposed to apply the damped BFGS quasi-Newton formula [158, Eqs. (5.8), (5.9), and (5.10)] to update $H^k$. This formula guarantees the positive definiteness of such a matrix, which is beneficial in practice and theory (see the comments towards the end of [157, § 2]). Moreover, he introduced a practical line-search technique based on a merit function suggested by Han [99]. Furthermore, Powell established the global convergence and the local R-superlinear convergence rate for his method without requiring the positive definiteness of $\nabla^2_{x,x}\mathcal{L}(x^*, \lambda^*)$ as Han did. Recognizing the contributions of Wilson, Han, and Powell, the SQP method is also referred to as the Wilson-Han-Powell method [186, 187, 29]. See [20, 92, 190, 86] for more detailed reviews of the history, theory, and practice of the SQP method.

### 4.1.2 A simple example

In Algorithm 4.1, it is crucial that $H^k$ approximates $\nabla^2_{x,x}\mathcal{L}(x^k, \lambda^k)$. It may be tempting to set $H^k \approx \nabla f(x^k)$, because the objective function of the SQP subproblem would then be a local quadratic approximation of $f$ at $x^k$. However, such a naive idea does not work, as illustrated by the following example inspired by Boggs and Tolle [20, § 2.2].

We consider

$$\begin{aligned} \min_{x\in\mathbb{R}^2} \quad & -x_1 - \frac{(x_2)^2}{4} \\ \text{s.t.} \quad & \|x\|^2 - 1 = 0, \end{aligned}$$

whose solution is $x^* = [1, 0]^\mathsf{T}$ with the associated Lagrange multiplier $\lambda^* = 1/2$. Suppose that we have an iterate $x^k = [t, 0]^\mathsf{T}$ with $t \approx 1$, so it is already close to the solution. If $H^k = \nabla^2 f(x^k)$, then the SQP subproblem becomes

$$\min_{d\in\mathbb{R}^2} \quad -d_1 - \frac{(d_2)^2}{4} \tag{4.1.3a}$$

$$\text{s.t.} \quad d_1 = \frac{1-t^2}{2t}. \tag{4.1.3b}$$

This subproblem is unbounded from below, regardless of the value of $t$. In addition, the more $d^k$ reduces (4.1.3a), the larger $\|x^k + d^k - x^*\|$ is. If $t = 1$, we have $x^k = x^*$, but any feasible point $d^k$ for (4.1.3b) will push $x^k + d^k$ away from $x^*$, unless $d^k$ is the global maximizer of (4.1.3a) subject to (4.1.3b).

Let us now consider the SQP subproblem (4.1.2) with $H^k = \nabla^2_{x,x}\mathcal{L}(x^k, \lambda^k)$ for a dual variable $\lambda^k \approx \lambda^* = 1/2$. It is

$$\begin{aligned} \min_{d\in\mathbb{R}^2} \quad & -d_1 + \lambda^k (d_1)^2 + \left(\lambda^k - \frac{1}{4}\right)(d_2)^2 \\ \text{s.t.} \quad & d_1 = \frac{1-t^2}{2t}. \end{aligned}$$

When $\lambda^k > 1/4$, the solution to this subproblem is

$$d^k = \begin{bmatrix} \dfrac{1-t^2}{2t} & 0 \end{bmatrix}^\mathsf{T}.$$

We thus have

$$x^k + d^k = \begin{bmatrix} \frac{t^2+1}{2t} & 0 \end{bmatrix}^\mathsf{T}.$$

If we set $x^{k+1} = x^k + d^k$ and continue to iterate in this way, we will obtain a sequence of iterates that converges quadratically to $x^*$, because

$$\|x^k + d^k - x^*\| = \frac{(1-t)^2}{2|t|} = \mathcal{O}(\|x^k - x^*\|^2).$$

This is not surprising, since Robinson [181] showed the local R-quadratic convergence rate of the SQP method when $H^k$ is the exact Hessian matrix of the Lagrangian with respect to $x$.

To summarize, as indicated by this example, choosing $H^k \approx \nabla^2_{x,x}\mathcal{L}(x^k, \lambda^k)$ instead of $H^k \approx \nabla f(x^k)$ in (4.1.2) is crucial.

### 4.1.3 Interpretations of the SQP subproblem

To get some insight into the origin of the SQP method, we interpret the SQP subproblem (4.1.2). In what follows, we focus on only one iteration of Algorithm 4.1 and hence, $k$ is fixed. We will explain why it is reasonable to update $x^k$ by a solution to (4.1.2).

#### Bilinear approximation of the Karush-Kuhn-Tucker conditions

This is the most classical interpretation of the SQP subproblem. According to Theorem 1.3.1, if $x^* \in \mathbb{R}^n$ is a local solution to the problem (4.1.1), under some mild assumptions, there exists a Lagrange multiplier $\lambda^* = [\lambda_i^*]_{i\in\mathcal{I}\cup\mathcal{E}}^\mathsf{T}$ with $\lambda_i^* \in \mathbb{R}$ for all $i \in \mathcal{I}\cup\mathcal{E}$ such that

$$\begin{cases} \nabla_x \mathcal{L}(x^*, \lambda^*) = 0, & \text{(4.1.4a)} \\ c_i(x^*) \le 0,\ i \in \mathcal{I}, & \text{(4.1.4b)} \\ c_i(x^*) = 0,\ i \in \mathcal{E}, & \text{(4.1.4c)} \\ \lambda_i^* c_i(x^*) = 0,\ i \in \mathcal{I}, & \text{(4.1.4d)} \\ \lambda_i^* \ge 0,\ i \in \mathcal{I}. & \text{(4.1.4e)} \end{cases}$$

Regard (4.1.4) as a nonlinear system of inequalities and equalities, and $(x^k, \lambda^k)$ as an approximation of $(x^*, \lambda^*)$. If we want to solve this system by the Newton-Raphson method[1] starting from $(x^k, \lambda^k)$, we would seek a step $(d, \mu)$ that satisfies the system

[1] Discussions are needed on how to apply the Newton-Raphson method to systems of nonlinear inequalities and equalities. We will not go further in this direction but refer to [174, 175, 180, 56] for fundamental works on this topic.

$$
\begin{cases}
\nabla_x \mathcal{L}(x^k, \lambda^k + \mu) + \nabla^2_{x,x}\mathcal{L}(x^k, \lambda^k) d = 0, & \text{(4.1.5a)} \\
c_i(x^k) + \nabla c_i(x^k)^\mathsf{T} d \le 0, \ i \in \mathcal{I}, & \text{(4.1.5b)} \\
c_i(x^k) + \nabla c_i(x^k)^\mathsf{T} d = 0, \ i \in \mathcal{E}, & \text{(4.1.5c)} \\
\lambda_i^k [c_i(x^k) + \nabla c_i(x^k)^\mathsf{T} d] + \mu_i c_i(x^k) = 0, \ i \in \mathcal{I}, & \text{(4.1.5d)} \\
\lambda_i^k + \mu_i \ge 0, \ i \in \mathcal{I}, & \text{(4.1.5e)}
\end{cases}
$$

which is a linear approximation of (4.1.4) at $(x^k, \lambda^k)$. However, as pointed out by Robinson [179, Rem. 3], an objection to such a method is that it would not solve even a linear program in one iteration. To cope with this defect, we let $(d, \mu)$ solve instead the following bilinear approximation of (4.1.4),

$$
\begin{cases}
\nabla_x \mathcal{L}(x^k, \lambda^k + \mu) + \nabla^2_{x,x}\mathcal{L}(x^k, \lambda^k) d = 0, & \text{(4.1.6a)} \\
c_i(x^k) + \nabla c_i(x^k)^\mathsf{T} d \le 0, \ i \in \mathcal{I}, & \text{(4.1.6b)} \\
c_i(x^k) + \nabla c_i(x^k)^\mathsf{T} d = 0, \ i \in \mathcal{E}, & \text{(4.1.6c)} \\
(\lambda_i^k + \mu_i)[c_i(x^k) + \nabla c_i(x^k)^\mathsf{T} d] = 0, \ i \in \mathcal{I}, & \text{(4.1.6d)} \\
\lambda_i^k + \mu_i \ge 0, \ i \in \mathcal{I}. & \text{(4.1.6e)}
\end{cases}
$$

Its only difference with the system (4.1.5) lies in the condition (4.1.6d), which includes the bilinear term $\mu_i \nabla c_i(x^k)^\mathsf{T} d$. If the problem (4.1.1) is a linear program, then (4.1.6) is precisely its KKT system, while (4.1.5) is only an approximation. Observe that the bilinear system (4.1.6) is nothing but the KKT conditions of the SQP subproblem (4.1.2), with $\lambda^k + \mu$ being the Lagrange multiplier. Therefore, a KKT pair for the SQP subproblem (4.1.2) is similar to a Newton-Raphson step for the KKT system of the problem (4.1.1), and it is even better in the sense that the resulting method solves a linear program in one iteration.

Note that discrepancy between the systems (4.1.5) and (4.1.6) disappears if $\mathcal{I} = \emptyset$ in the problem (4.1.1) and hence, a KKT pair for the SQP subproblem (4.1.2) is exactly a Newton-Raphson step for the KKT system of the problem (4.1.1) in such a situation.

### Approximation of a modified Lagrangian

This interpretation is due to Robinson [179, Rem. 4]. Let $\widetilde{\mathcal{L}}$ be the function

$$
\widetilde{\mathcal{L}}(x, \lambda) \overset{\text{def}}{=} f(x) + \sum_{i \in \mathcal{I} \cup \mathcal{E}} \lambda_i \delta_i(x), \quad \text{for } x \in \mathbb{R}^n \text{ and } \lambda_i \in \mathbb{R}, \text{ with } i \in \mathcal{I} \cup \mathcal{E},
$$

where $\delta_i$, for $i \in \mathcal{I} \cup \mathcal{E}$, is defined by

$$
\delta_i(x) \overset{\text{def}}{=} c_i(x) - c_i(x^k) - \nabla c_i(x^k)^\mathsf{T}(x - x^k), \quad \text{for } x \in \mathbb{R}^n.
$$

The function $\delta_i$ is referred to as the departure from linearity[2] for $c_i$ at the point $x^k$ (see [86, § 2] and [84, § 2.3]). The SQP subproblem (4.1.2) with $H^k = \nabla^2_{x,x}\mathcal{L}(x^k, \lambda^k)$ can then be seen as the minimization of the second-order Taylor approximation of $\widetilde{\mathcal{L}}$ subject to the linearizations of the constraints (4.1.1b) and (4.1.1c) at $(x^k, \lambda^k)$, i.e.,

$$\begin{aligned} \min_{d \in \mathbb{R}^n} \quad & \nabla_x \widetilde{\mathcal{L}}(x^k, \lambda^k)^\mathsf{T} d + \frac{1}{2} d^\mathsf{T} \nabla^2_{x,x} \widetilde{\mathcal{L}}(x^k, \lambda^k) d \\ \text{s.t.} \quad & c_i(x^k) + \nabla c_i(x^k)^\mathsf{T} d \le 0, \ i \in \mathcal{I}, \\ & c_i(x^k) + \nabla c_i(x^k)^\mathsf{T} d = 0, \ i \in \mathcal{E}. \end{aligned}$$

By expressing its subproblem in this form, we see that the SQP method is a special case of Robinson's method [179], known to have a local R-quadratic convergence rate.

### Approximation of the objective function in the tangent space of a surface

Inspired by an observation in [86, § 2], we can also interpret the objective function of the SQP subproblem as an approximation of the original objective function in the tangent space of a particular surface described below. As shown in Theorem 4.1.1, when approximating $f$ in this space, we will naturally get the Hessian matrix of the Lagrangian in the second-order term.

For this interpretation, we focus on the equality-constrained problem

$$\min_{x \in \mathbb{R}^n} \quad f(x) \tag{4.1.7a}$$

$$\text{s.t.} \quad h(x) = 0, \tag{4.1.7b}$$

with $h : \mathbb{R}^n \to \mathbb{R}^m$. Recall that the Lagrangian function of (4.1.7) is

$$\mathcal{L}(x, \lambda) \overset{\text{def}}{=} f(x) + \lambda^\mathsf{T} h(x), \quad \text{for } x \in \mathbb{R}^n \text{ and } \lambda \in \mathbb{R}^m.$$

Let $\bar{x} \in \mathbb{R}^n$, $\bar{\lambda} \in \mathbb{R}^m$ be given, and define

$$Q(d) \overset{\text{def}}{=} f(\bar{x}) + \nabla f(\bar{x})^\mathsf{T} d + \frac{1}{2} d^\mathsf{T} \nabla^2_{x,x} \mathcal{L}(\bar{x}, \bar{\lambda}) d.$$

If $\bar{x}$ and $\bar{\lambda}$ represent the current iterate and approximate Lagrange multiplier, then $Q$ is the objective function of the SQP subproblem with an exact second-order term. At first glance, $Q$ does not seem to be a natural approximation of $f$ due to its second-order term. However, this approximation turns out to be indeed natural if we focus on the contour surface of $h$ at $\bar{x}$ and its tangent space. This is detailed in Theorem 4.1.1.

[2] When $c_i$ is strictly convex, $\delta_i$ defines the Bregman distance [28] associated with $c_i$.

**Theorem 4.1.1.** *Assume that $f$ and $h$ are twice differentiable with $\nabla^2 f$ being locally Lipschitz continuous. Consider a curve parametrized by $x : \mathbb{R} \to \mathbb{R}^n$, satisfying*

$$h\big(x(t)\big) = h(\bar{x}) \quad \textit{for all } t \in \mathbb{R}, \quad \textit{and} \quad x(0) = \bar{x}. \tag{4.1.8}$$

*If $x$ is twice differentiable and $x''$ is locally Lipschitz continuous, then there exist constants $\nu \geq 0$ and $\epsilon > 0$ such that*

$$|f(x(t)) - Q(x'(0)t)| \leq \left(\nu t + \frac{1}{2}|x''(0)^{\mathsf{T}}[\nabla f(\bar{x}) + \nabla h(\bar{x})^{\mathsf{T}}\bar{\lambda}]|\right)t^2 \quad \textit{for all } t \in (-\epsilon, \epsilon).$$

*Proof.* Define

$$\phi(t) \stackrel{\text{def}}{=} f\big(x(t)\big), \qquad \text{for } t \in \mathbb{R}.$$

By assumption, $\phi$ is twice differentiable and there exists $\epsilon > 0$ such that $\phi''$ is Lipschitz continuous in $(-\epsilon, \epsilon)$. Let $\widehat{\phi}$ be the second-order Taylor expansion of $\phi$ at $0$. Then

$$\big|f\big(x(t)\big) - Q(x'(0)t)\big| \leq \big|\phi(t) - \widehat{\phi}(t)\big| + \big|\widehat{\phi}(t) - Q(x'(0)t)\big|, \tag{4.1.9}$$

and the Lipschitz continuity of $\phi''$ ensures the existence of a constant $\nu \geq 0$ such that

$$\big|\phi(t) - \widehat{\phi}(t)\big| \leq \nu t^3 \quad \text{for } t \in (-\epsilon, \epsilon). \tag{4.1.10}$$

We now bound $\big|\widehat{\phi}(t) - Q(x'(0)t)\big|$. Noting that $x(0) = \bar{x}$, we have

$$\begin{cases} \phi'(0) = x'(0)^{\mathsf{T}}\nabla f(\bar{x}), & (4.1.11\text{a}) \\ \phi''(0) = x'(0)^{\mathsf{T}}\nabla^2 f(\bar{x})x'(0) + x''(0)^{\mathsf{T}}\nabla f(\bar{x}). & (4.1.11\text{b}) \end{cases}$$

According to (4.1.11), only the second-order terms of $\widehat{\phi}(t)$ and $Q(x'(0)t)$ differ. Hence,

$$\begin{aligned} \big|\widehat{\phi}(t) - Q(x'(0)t)\big| &= \frac{t^2}{2}\big|\phi''(0) - x'(0)^{\mathsf{T}}\nabla^2_{x,x}\mathcal{L}(\bar{x},\bar{\lambda})x'(0)\big| & (4.1.12\text{a}) \\ &= \frac{t^2}{2}\left|x''(0)^{\mathsf{T}}\nabla f(\bar{x}) - \sum_{i=1}^{m}\bar{\lambda}_i x'(0)\nabla^2 h_i(\bar{x})x'(0)\right|, & (4.1.12\text{b}) \end{aligned}$$

where the second equality uses the formula of $\phi''(0)$ in (4.1.11b) and the definition of $\mathcal{L}$. Moreover, since $h\big(x(t)\big) = h(\bar{x})$ for all $t \in \mathbb{R}$, we have for each $i \in \{1, 2, \ldots, m\}$ that

$$0 = \frac{\mathrm{d}^2}{\mathrm{d}t^2}h_i\big(x(t)\big)\bigg|_{t=0} = x'(0)^{\mathsf{T}}\nabla^2 h_i(\bar{x})x'(0) + x''(0)^{\mathsf{T}}\nabla h_i(\bar{x}). \tag{4.1.13}$$

Therefore, $x'(0)^{\mathsf{T}}\nabla^2 h_i(\bar{x})x'(0) = -x''(0)^{\mathsf{T}}\nabla h_i(\bar{x})$ for each $i$ and hence, (4.1.12) leads to

$$\big|\widehat{\phi}(t) - Q(x'(0)t)\big| \leq \frac{t^2}{2}\big|x''(0)^{\mathsf{T}}\big[\nabla f(\bar{x}) + \nabla h(\bar{x})^{\mathsf{T}}\bar{\lambda}\big]\big|. \tag{4.1.14}$$

Plugging (4.1.10) and (4.1.14) into (4.1.9), we obtain the desired result. □

Note that $x'(0)$ in Theorem 4.1.1 is a vector in

$$\mathcal{T}_{\bar{x}} = \{d \in \mathbb{R}^n : \nabla h(\bar{x})d = 0\}.$$

This is the tangent space of

$$\mathcal{S}_{\bar{x}} = \{x \in \mathbb{R}^n : h(x) = h(\bar{x})\},$$

at $\bar{x}$, and $\mathcal{S}_{\bar{x}}$ is the contour surface of the constraints at $\bar{x}$. Therefore, Theorem 4.1.1 reveals that the objective function of the SQP subproblem is a natural approximation of $f$ in this tangent space if we regard $f$ as a function on $\mathcal{S}_{\bar{x}}$.

However, the feasible region of the SQP subproblem is in general not the tangent space $\mathcal{T}_{\bar{x}}$, but the linearized feasible direction set

$$\mathcal{F}^{\mathsf{L}}_{\bar{x}} = \{d \in \mathbb{R}^n : h(\bar{x}) + \nabla h(\bar{x})d = 0\},$$

which is assumed nonempty for the current discussion. This set equals $\mathcal{T}_{\bar{x}}$ only if $\bar{x}$ is feasible so that $h(\bar{x}) = 0$. When $\bar{x}$ is infeasible, $\mathcal{F}^{\mathsf{L}}_{\bar{x}}$ is a shifted copy of $\mathcal{T}_{\bar{x}}$ as illustrated in Figure 4.1. The distance between $\mathcal{F}^{\mathsf{L}}_{\bar{x}}$ and $\mathcal{T}_{\bar{x}}$ is $\delta_{\bar{x}} = \|d^*\| = \|\nabla h(\bar{x})^\dagger h(\bar{x})\|$, where $d^*$ is the projection of $0 \in \mathcal{T}_{\bar{x}}$ onto $\mathcal{F}^{\mathsf{L}}_{\bar{x}}$, and $\nabla h(\bar{x})^\dagger$ denotes the Moore-Penrose pseudoinverse of $\nabla h(\bar{x})$. The purpose of this shift is to improve the feasibility. Indeed, if $h$ is twice continuously differentiable, then its Taylor expansion provides $\|h(\bar{x} + d)\| = \mathcal{O}(\|d\|^2)$ for all $d$ in a bounded subset of $\mathcal{F}^{\mathsf{L}}_{\bar{x}}$, which further implies $\|h(\bar{x} + d)\| = \mathcal{O}(\|h(\bar{x})\|^2)$ unless $\|d\|$ is much larger than $\|d^*\|$; in contrast, $\|h(\bar{x}+d)\| = \|h(\bar{x})\| + \mathcal{O}(\|d\|^2)$ if $d$ stays in a bounded subset of $\mathcal{T}_{\bar{x}}$.

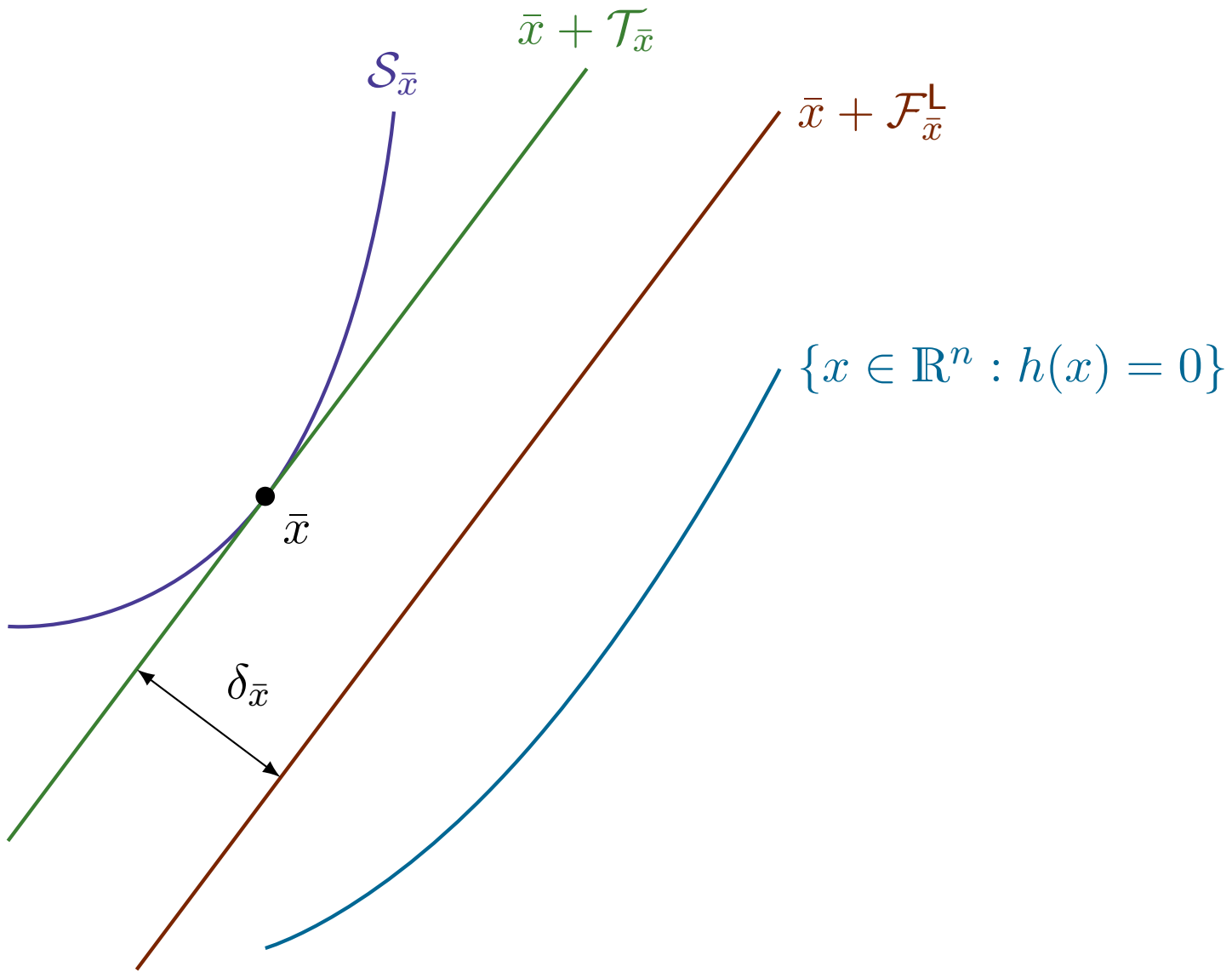

Figure 4.1: An illustration of $\mathcal{T}_{\bar{x}}$, $\mathcal{S}_{\bar{x}}$, $\mathcal{F}^{\mathsf{L}}_{\bar{x}}$, and the feasible set

Therefore, the SQP subproblem models the objective function in the tangent space of the contour surface of the constraint at the current iterate, and then minimizes this model in an affine space that shifts the tangent space towards improvement of feasibility. We observe that a smaller magnitude of $\nabla f(\bar{x})+\nabla h(\bar{x})^{\mathsf{T}}\bar{\lambda}$ implies a better accuracy of the model, as indicated by the error estimation in Theorem 4.1.1. If $(\bar{x},\bar{\lambda})$ is close to a KKT pair, then $\|\nabla f(\bar{x})+\nabla h(\bar{x})^{\mathsf{T}}\bar{\lambda}\|$ is small. On the other hand, if we defined $\nabla^2 Q$ by $\nabla^2 f(\bar{x})$ instead of $\nabla^2_{x,x}\mathcal{L}(\bar{x},\bar{\lambda})$, then under the assumptions of Theorem 4.1.1, we would have

$$\left|f\big(x(t)\big)-Q(x'(0)t)\right| \le \left(\nu t+\frac{1}{2}\left|x''(0)^{\mathsf{T}}\nabla f(\bar{x})\right|\right)t^2 \quad \text{for all } t\in(-\epsilon,\epsilon), \tag{4.1.15}$$

which can be obtained by setting $\bar{\lambda}=0$ in Theorem 4.1.1. However, since we are considering constrained optimization, we cannot expect $\|\nabla f(\bar{x})\|$ to be small even if $\bar{x}$ is close to a solution unless no constraint is taking effect at this solution. This explains why the second-order term of the SQP subproblem should be defined by the Hessian matrix of the Lagrangian, rather than that of $f$. A special case occurs if the Lagrangian and $f$ have the same Hessian, for example, when $h$ is linear. In such a case, (4.1.13) implies

$$x''(0)^{\mathsf{T}}\nabla h(\bar{x})^{\mathsf{T}}\bar{\lambda} = -\sum_{i=1}^{m}\bar{\lambda}_i x'(0)^{\mathsf{T}}\nabla^2 h_i(\bar{x})x'(0) = x'(0)^{\mathsf{T}}[\nabla^2 f(\bar{x})-\nabla^2_{x,x}\mathcal{L}(\bar{x},\bar{\lambda})]x'(0)=0,$$

and hence $x''(0)^{\mathsf{T}}\nabla f(\bar{x})$ in (4.1.15) equals $x''(0)^{\mathsf{T}}[\nabla f(\bar{x})+\nabla h(\bar{x})^{\mathsf{T}}\bar{\lambda}]$, leading to the same error estimation as in Theorem 4.1.1.

Geometrically speaking, when we regard $f$ as a function on the surface $\mathcal{S}_{\bar{x}}$ and approximate it in the tangent space $\mathcal{T}_{\bar{x}}$, the second-order Taylor expansion of $f$ generally does not attain the same accuracy as it does in $\mathbb{R}^n$, because neither this Taylor expansion nor $\mathcal{T}_{\bar{x}}$ can reflect the curvature information of $\mathcal{S}_{\bar{x}}$, unless the curvature is zero. In contrast, the objective function of the SQP subproblem can achieve a better accuracy since it includes such curvature information via $\nabla^2_{x,x}\mathcal{L}(\bar{x},\bar{\lambda})$.

It is worth mentioning that Theorem 4.1.1 still holds if $x$ satisfies

$$\bar{\lambda}^{\mathsf{T}}h\big(x(t)\big)=\bar{\lambda}^{\mathsf{T}}h(\bar{x}) \quad \text{for all } t\in\mathbb{R}, \quad \text{and} \quad x(0)=\bar{x} \tag{4.1.16}$$

instead of (4.1.8). Given (4.1.16), we will have

$$0=\frac{\mathrm{d}^2}{\mathrm{d}t^2}\bar{\lambda}^{\mathsf{T}}h\big(x(t)\big)\bigg|_{t=0}=\bar{\lambda}^{\mathsf{T}}\nabla h(\bar{x})x''(0)+\sum_{i=1}^{m}\bar{\lambda}_i x'(0)^{\mathsf{T}}\nabla^2 h_i(\bar{x})x'(0).$$

Plugging this into (4.1.12b), we obtain again (4.1.14) and hence, establish the same estimation of $|f(x(t))-Q(x'(0)t)|$. Therefore, $Q$ is a quadratic approximation of $f$ in

$$\mathcal{T}_{\bar{x},\bar{\lambda}}=\{d\in\mathbb{R}^n : \bar{\lambda}^{\mathsf{T}}\nabla h(\bar{x})d=0\},$$

which is linear space that contains $\mathcal{T}_{\bar{x}}$ as a subset. Note that the dimension of $\mathcal{T}_{\bar{x},\bar{\lambda}}$ is at least $n-1$, while that of $\mathcal{T}_{\bar{x}}$ is $n - \operatorname{rank}\big(\nabla h(\bar{x})\big)$.

### 4.1.4 Lagrangian and augmented Lagrangian of the SQP subproblem

In the sequel, we study the Lagrangian and the augmented Lagrangian of the SQP subproblem and observe their relations with those of the original optimization problem. We also point out that the augmented Lagrangian of the SQP subproblem is exactly the approximate augmented Lagrangian used in [143, 205].

For simplicity, instead of the problem (4.1.1), we consider the problem

$$\begin{aligned} \min_{x\in\mathbb{R}^n} \quad & f(x) && (4.1.17a)\\ \text{s.t.} \quad & h(x) = 0, && (4.1.17b)\\ & x \geq 0, && (4.1.17c) \end{aligned}$$

with $h : \mathbb{R}^n \to \mathbb{R}^m$. The problem (4.1.1) can be reformulated in this form[3]. Recall that the Lagrangian of (4.1.17) is

$$\mathcal{L}(x,\lambda) \overset{\text{def}}{=} f(x) + \lambda^{\mathsf{T}} h(x), \quad \text{for } x \geq 0 \text{ and } \lambda \in \mathbb{R}^m.$$

Let $x^k \geq 0$ and $\lambda^k \in \mathbb{R}^m$ be given. Correspondingly, the SQP subproblem of the problem (4.1.17) is

$$\begin{aligned} \min_{d\in\mathbb{R}^n} \quad & f(x^k) + \nabla f(x^k)^{\mathsf{T}} d + \frac{1}{2} d^{\mathsf{T}} H^k d && (4.1.18a)\\ \text{s.t.} \quad & h(x^k) + \nabla h(x^k) d = 0, && (4.1.18b)\\ & x^k + d \geq 0, && (4.1.18c) \end{aligned}$$

with $H^k \approx \nabla^2_{x,x}\mathcal{L}(x^k,\lambda^k)$. Note that the constant term $f(x^k)$ in (4.1.18a) may be excluded, as in (4.1.2a), but including it facilitates the discussion in the sequel.

Recall that the augmented Lagrangian [104, 154, 182, 183] of the problem (4.1.17) is

$$\mathcal{L}_{\mathsf{A}}(x,\lambda) \overset{\text{def}}{=} \mathcal{L}(x,\lambda) + \frac{\gamma}{2}\|h(x)\|^2, \quad \text{for } x \geq 0 \text{ and } \lambda \in \mathbb{R}^m, \qquad (4.1.19)$$

where $\gamma \geq 0$ is a penalty parameter. Denote by $\widetilde{\mathcal{L}}$ the Lagrangian of the SQP subproblem (4.1.18), i.e.,

$$\begin{aligned} \widetilde{\mathcal{L}}(d,\lambda) \overset{\text{def}}{=} & f(x^k) + \nabla f(x^k)^{\mathsf{T}} d + \frac{1}{2} d^{\mathsf{T}} H^k d \\ & + \lambda^{\mathsf{T}}[h(x^k) + \nabla h(x^k) d], \quad \text{for } d \geq -x^k \text{ and } \lambda \in \mathbb{R}^m, \end{aligned}$$

[3]In the reformulation, the dimension and the meaning of $x$ may be altered, but we do not change the notations since it does not lead to confusion.

and by $\widetilde{\mathcal{L}}_{\mathsf{A}}$ the augmented Lagrangian of the SQP subproblem (4.1.18), i.e.,

$$\widetilde{\mathcal{L}}_{\mathsf{A}}(d,\lambda) \stackrel{\text{def}}{=} \widetilde{\mathcal{L}}(d,\lambda) + \frac{\gamma}{2}\|h(x^k) + \nabla h(x^k)d\|^2, \quad \text{for } d \geq -x^k \text{ and } \lambda \in \mathbb{R}^m.$$

We now present some relations between $\mathcal{L}$ and $\widetilde{\mathcal{L}}$, and between $\mathcal{L}_{\mathsf{A}}$ and $\widetilde{\mathcal{L}}_{\mathsf{A}}$.

**Theorem 4.1.2.** *Assume that $f$ and $h$ are twice differentiable. If $H^k = \nabla^2_{x,x}\mathcal{L}(x^k, \lambda^k)$, then $\widetilde{\mathcal{L}}(d, \lambda^k)$ is the second-order Taylor expansion of $\mathcal{L}(x^k + d, \lambda^k)$ with respect to $d$ at $0$.*

*Proof.* This theorem can be verified by a straightforward calculation. □

**Theorem 4.1.3.** *Assume that $f$ and $h$ are twice differentiable. If*

$$H^k = \nabla^2 f(x^k) + \sum_{i=1}^{m}[\lambda_i^k + \gamma h_i(x^k)]\nabla^2 h_i(x^k), \tag{4.1.20}$$

*then $\widetilde{\mathcal{L}}_{\mathsf{A}}(d, \lambda^k)$ is the second-order Taylor expansion of $\mathcal{L}_{\mathsf{A}}(x^k + d, \lambda^k)$ with respect to $d$ at $0$.*

*Proof.* By direct calculations, we have

$$\nabla_x \mathcal{L}_{\mathsf{A}}(x^k, \lambda^k) = \nabla_x \mathcal{L}(x^k, \lambda^k) + \gamma \nabla h(x^k)^{\mathsf{T}} h(x^k),$$

and

$$\nabla^2_{x,x}\mathcal{L}_{\mathsf{A}}(x^k, \lambda^k) = \nabla^2_{x,x}\mathcal{L}(x^k, \lambda^k) + \gamma\left[\nabla h(x^k)^{\mathsf{T}}\nabla h(x^k) + \sum_{i=1}^{m} h_i(x^k)\nabla^2 h_i(x^k)\right].$$

Therefore, the second-order Taylor expansion of $\mathcal{L}_{\mathsf{A}}(x^k + d, \lambda^k)$ with respect to $d$ at $0$ is

$$\begin{aligned}
&\mathcal{L}_{\mathsf{A}}(x^k, \lambda^k) + \nabla_x \mathcal{L}_{\mathsf{A}}(x^k, \lambda^k)^{\mathsf{T}} d + \frac{1}{2} d^{\mathsf{T}} \nabla^2_{x,x}\mathcal{L}_{\mathsf{A}}(x^k, \lambda^k) d \\
&\quad = \mathcal{L}(x^k, \lambda^k) + \nabla_x \mathcal{L}(x^k, \lambda^k)^{\mathsf{T}} d + \frac{1}{2} d^{\mathsf{T}} H^k d + \frac{\gamma}{2}\|h(x^k) + \nabla h(x^k) d\|^2,
\end{aligned}$$

where $H^k$ is defined by (4.1.20). □

Intriguingly, an augmented Lagrangian method for solving the problem (4.1.17) updates traditionally the dual variable $\lambda^k$ by

$$\lambda^{k+1} = \lambda^k + \gamma h(x^k).$$

Therefore, the second-order Taylor expansion of $\mathcal{L}_{\mathsf{A}}$ can be interpreted as the augmented Lagrangian of the SQP subproblem (4.1.18) with $H^k = \nabla^2_{x,x}\mathcal{L}(x^k, \lambda^{k+1})$.

There is an interesting connection between the augmented Lagrangian of the SQP

subproblem (4.1.18) and the trust-region augmented Lagrangian methods studied by Yuan et al. [143, 205]. These methods employ the approximation

$$\mathcal{L}_{\mathsf{A}}(x^k+d,\lambda^k) \approx \mathcal{L}(x^k,\lambda^k)+\nabla_x\mathcal{L}(x^k,\lambda^k)^{\mathsf{T}}d+\frac{1}{2}d^{\mathsf{T}}H^kd+\frac{\gamma}{2}\|h(x^k)+\nabla h(x^k)d\|^2, \quad (4.1.21)$$

which is a quadratic approximation obtained by replacing $\mathcal{L}(x^k+d,\lambda^k)$ with a quadratic approximation, and replacing $h(x^k+d)$ in the penalty term with its first-order Taylor expansion. This approximation turns out to be the augmented Lagrangian of the SQP subproblem (4.1.18), as specified by Theorem 4.1.4.

**Theorem 4.1.4.** *Assume that $f$ and $h$ are differentiable. For any matrix $H^k \in \mathbb{R}^{n\times n}$, the right-hand side of (4.1.21) equals $\widetilde{\mathcal{L}}_{\mathsf{A}}(d,\lambda^k)$.*

*Proof.* As Theorem 4.1.2, this can be verified by a straightforward calculation. □

Several possibilities of $H^k$ are proposed in [143, 205]. For example, $H^k$ can be set to $\nabla^2_{x,x}\mathcal{L}(x^k,\lambda^k)$ or an approximation. As pointed out by [143, § 2.1], if $H^k$ is defined as in (4.1.20), then the right-hand side of (4.1.21) is the second-order Taylor expansion of $\mathcal{L}_{\mathsf{A}}(x^k+d,\lambda^k)$ with respect to $d$ at $0$, which agrees with Theorem 4.1.3. However, in the numerical experiments, [143, 205] choose $H^k=\nabla^2_{x,x}\mathcal{L}(x^k,\lambda^k)$.

Suppose that an algorithm defines a step $d^k$ based on the minimization of the right-hand side of (4.1.21). Then, Theorem 4.1.4 tells us that the algorithm can be regarded as an SQP method that approximately solves the SQP subproblem (4.1.18) by minimizing $\widetilde{\mathcal{L}}_{\mathsf{A}}(d,\lambda^k)$, i.e., by applying one single iteration of an augmented Lagrangian method.

## 4.2 Merit functions for the SQP method and the Maratos effect

In unconstrained optimization, a point $x \in \mathbb{R}^n$ is normally considered to be better than another point $y \in \mathbb{R}^n$ if $f(x) < f(y)$. However, this is not true in constrained optimization because the constraints must also be taken into account. This is usually done using *merit functions*. A merit function assesses the quality of a point by considering both $f$ and $c_i$, with $i \in \mathcal{I} \cup \mathcal{E}$. We present in what follows some classical merit functions.

We also mention that there exists an alternative to merit functions, known as the filter method [70–73, 178]. Roughly speaking, it attempts to minimize the objective function and the constraint violation simultaneously, as in biobjective optimization. However, we do not study the filter method in this thesis. Interested readers may refer to the references mentioned above for more details.

### 4.2.1 The Courant merit function

Perhaps the most classical merit function is the Courant function[4] [52], defined by

$$\varphi_\gamma(x) \stackrel{\text{def}}{=} f(x) + \gamma\Big(\sum_{i\in\mathcal{I}}[c_i(x)]_+^2 + \sum_{i\in\mathcal{E}} c_i(x)^2\Big), \quad \text{for } x \in \mathbb{R}^n \text{ and } \gamma \geq 0,$$

where $[\cdot]_+$ takes the positive part of a given number. The advantage of such a merit function is that it is differentiable if $f$ and $c_i$, for $i \in \mathcal{I} \cup \mathcal{E}$, are differentiable. However, normally, a global minimizer of $\varphi_\gamma$ is not a solution to (4.1.1) when $\gamma$ is finite. This phenomenon can be seen on the simple example of minimizing $x$ subject to $x \geq 0$.

### 4.2.2 Nonsmooth merit functions

Other examples of merit functions are the $\ell_p$-merit functions, defined by

$$\varphi_\gamma(x) \stackrel{\text{def}}{=} f(x) + \gamma\Big(\sum_{i\in\mathcal{I}}[c_i(x)]_+^p + \sum_{i\in\mathcal{E}} |c_i(x)|^p\Big)^{1/p}, \quad \text{for } x \in \mathbb{R}^n \text{ and } \gamma \geq 0,$$

Such merit functions enjoy the property of being exact under some mild assumptions. Roughly speaking, this means that a solution[5] to the constrained problem (4.1.1) can be obtained by minimizing $\varphi_\gamma$ when $\gamma$ is big enough. We will not discuss this concept further but only provide the following theorem for later reference. Interested readers may refer to [101, 133, 59].

**Theorem 4.2.1** ([44, Thm. 14.5.1])**.** *Assume that the functions $f$ and $c_i$ are twice continuously differentiable for all $i \in \mathcal{I} \cup \mathcal{E}$. Let $(x^*, \lambda^*)$ be a KKT pair to the problem (4.1.1) that satisfies the second-order sufficient condition of Theorem 1.3.3. If $\gamma \geq \|\lambda^*\|_q$, where $q$ is the Hölder conjugate of $p$, then $x^*$ satisfies the second-order sufficient condition for the minimization of $\varphi_\gamma$. Moreover, if $\gamma > \|\lambda^*\|_q$, then the two second-order sufficient conditions are equivalent.*

An inconvenience of the $\ell_p$-merit functions is that, even if $f$ and $c_i$ are differentiable for all $i \in \mathcal{I} \cup \mathcal{E}$, $\varphi_\gamma$ is likely not differentiable at the points $x \in \mathbb{R}^n$ where $c_i(x) = 0$ for all $i \in \mathcal{I} \cup \mathcal{E}$. However, smooth merit functions that enjoy exactness do exist, as shown in the next section.

[4]Courant proposed this merit function when dealing with boundary conditions of equilibrium and vibration problems. See [52, Pt. II, § 3] for details.

[5]Here, the term "solution" may be interpreted in different ways. It may be a global solution, a local solution, or a stationary point.

### 4.2.3 The augmented Lagrangian merit function

For simplicity, we assume here that $\mathcal{I} = \emptyset$. The augmented Lagrangian [104, 154, 182, 183] of the problem (4.1.1) is then

$$\mathcal{L}_{\mathsf{A}}(x, \lambda) \overset{\mathsf{def}}{=} \mathcal{L}(x, \lambda) + \frac{\gamma}{2} \sum_{i \in \mathcal{E}} c_i(x)^2, \quad \text{for } x \in \mathbb{R}^n \text{ and } \lambda = [\lambda_i]_{i \in \mathcal{E}}^{\mathsf{T}},$$

where $\mathcal{L}$ denotes the Lagrangian function of the problem (4.1.1). The augmented Lagrangian merit function is then defined as

$$\varphi_\gamma(x) \overset{\mathsf{def}}{=} \mathcal{L}_{\mathsf{A}}(x, \lambda^{\mathsf{LS}}(x)),$$

where $\lambda^{\mathsf{LS}}(x)$ denotes the least-norm solution to

$$\min_{\lambda} \left\| \nabla f(x) + \sum_{i \in \mathcal{E}} \lambda_i \nabla c_i(x) \right\|, \tag{4.2.1}$$

where $\lambda = [\lambda_i]_{i \in \mathcal{E}}^{\mathsf{T}}$. If $x^*$ is a solution to the problem (4.1.1), then $(x^*, \lambda(x^*))$ is a KKT pair. When $\mathcal{I} \neq \emptyset$, to achieve the same property, the complementary slackness conditions must be taken into account in (4.2.1).

A clear drawback of such a merit function is that it is expensive to evaluate, as one evaluation necessitates solving a linear least-squares problem that involves the gradients of $f$ and $c_i$, with $i \in \mathcal{E}$. Nonetheless, this merit function offers several advantages. First of all, if $f$ and $c_i$, for $i \in \mathcal{E}$, are differentiable, and if $\{\nabla c_i(x)\}_{i \in \mathcal{E}}$ are linearly independent for $x \in \mathbb{R}^n$, then $\lambda^{\mathsf{LS}}$ is differentiable at $x$ [44, Lem. 14.2.1], and hence, so is $\varphi_\gamma$. Moreover, under some mild assumptions, when $\gamma$ is large enough, the second-order sufficient conditions of (4.1.1) and of minimizing $\varphi_\gamma$ are equivalent [44, Thm. 14.6.1].

### 4.2.4 Maratos effect and second-order correction

In practice, to have global convergence, the SQP method needs to be globalized, normally using a merit function to decide whether to accept a step or not. However, the merit function may jeopardize the fast local convergence of the SQP method. This is known as the *Maratos effect* [129].

The reason behind the effect is that, for certain problems, the SQP method can generate a step that increases both the objective function and the constraint violation, no matter how close the current iterate is to a solution. Such a step would be rejected by any merit function that combines the objective function and the constraint violation so that it increases with respect to both. See examples of this phenomenon in [129, § 3.5] and [160].

One way to cope with the Maratos effect is to employ second-order correction steps. The idea behind the second-order correction is to modify the current step $d^k \in \mathbb{R}^n$ by a term $r^k \in \mathbb{R}^n$, so that

$$\begin{cases} [c_i(x^k + d^k + r^k)]_+ = \mathrm{o}(\|d^k\|^2), & i \in \mathcal{I}, \\ |c_i(x^k + d^k + r^k)| = \mathrm{o}(\|d^k\|^2), & i \in \mathcal{E}. \end{cases}$$

However, the step $d^k$ should not be substantially altered, so that $\|r^k\| = \mathrm{o}(\|d^k\|)$ is also imposed. A step that satisfies these conditions is commonly known as a *second-order correction step*. In practice, it is reasonable to make such a modification only when necessary (see [44, Alg. 15.3.1] and the discussion around). Perhaps the simplest of these steps, similar to [134, Eqs. (21) and (22)], is the least-squares solution of

$$\min_{r \in \mathbb{R}^n} \sum_{i \in \mathcal{I}} [c_i(x^k + d^k) + \nabla c_i(x^k + d^k)^\mathsf{T} r]_+^2 + \sum_{i \in \mathcal{E}} [c_i(x^k + d^k) + \nabla c_i(x^k + d^k)^\mathsf{T} r]^2.$$

Many other second-order correction steps can also be defined (see, e.g., [42, 43, 69, 77]).

Another possibility to handle the Maratos effect is properly choosing the merit function. For example, Powell and Yuan [172, § 4] showed that the Maratos effect could not occur when using the augmented Lagrangian merit function in a line-search SQP method. We refer to [190, § 7] for other techniques to cope with the Maratos effect.

## 4.3 The trust-region SQP method

As we mentioned previously, we need in practice a globalization strategy to make the method globally convergent. We present in what follows a method that uses the trust-region strategy in a derivative-based context. The new method that we will present in Chapter 5 uses a derivative-free variation of this trust-region SQP method.

### 4.3.1 Overview of the method

We now present the basic trust-region SQP method. The merit function we consider from now on is the $\ell_2$-merit function, defined by

$$\varphi_\gamma(x) \overset{\text{def}}{=} f(x) + \gamma \sqrt{\sum_{i \in \mathcal{I}} [c_i(x)]_+^2 + \sum_{i \in \mathcal{E}} |c_i(x)|^2}, \quad \text{for } x \in \mathbb{R}^n \text{ and } \gamma \geq 0,$$

where $\gamma$ is the penalty parameter. The algorithm we present below maintains a penalty parameter $\gamma^k \geq 0$ at the $k$th iteration, and we use $\varphi^k$ to denote the function $\varphi_{\gamma^k}$. We denote by $\widehat{\varphi}^k$ the $\ell_2$-merit function computed on the SQP subproblem (4.1.2), i.e.,

$$\widehat{\varphi}^k(d) \overset{\text{def}}{=} \nabla f(x^k)^\mathsf{T} d + \frac{1}{2} d^\mathsf{T} \nabla_{x,x}^2 \mathcal{L}(x^k, \lambda^k) d + \gamma^k \Phi(d), \quad \text{for } d \in \mathbb{R}^n,$$

where $\Phi$ is defined by

$$\Phi(d) \overset{\text{def}}{=} \sqrt{\sum_{i \in \mathcal{I}} [c_i(x^k) + \nabla c_i(x^k)^\mathsf{T} d]_+^2 + \sum_{i \in \mathcal{E}} [c_i(x^k) + \nabla c_i(x^k)^\mathsf{T} d]^2}, \quad \text{for } d \in \mathbb{R}^n.$$

Note that although $\varphi^k$ is a function of $x$, we define $\widehat{\varphi}^k$ as a function of $d$ for convenience. As mentioned in [44, § 15.3.2.1], such a function is a quadratic approximation of $\varphi^k$ at $x^k$, which is suitable for use in trust-region SQP method. The basic framework is given in Algorithm 4.2, where we do not include the second-order correction mechanism for simplicity. Note that this framework is conceptual, as many details are hidden, including the definition of "convergence" in the loop and the update of $\lambda^k$, which should be specified in a practical implementation. A similar framework for equality-constrained problems can be found in [173, § 2] and [144, Alg. 18.4].

### 4.3.2 Comments on the penalty parameter and merit function

To ensure that the requirement on Line 4 of Algorithm 4.2 can be satisfied, when choosing the trial step $d^k$, we make sure that there exists $\bar{\gamma} \geq 0$ such that for all $\gamma \geq \bar{\gamma}$, the property $\widehat{\varphi}^k(d^k) < \widehat{\varphi}^k(0)$ is guaranteed. According to the aforementioned definition of $\widehat{\varphi}^k$, this can be achieved if the step $d^k$ satisfies either

$$\sum_{i \in \mathcal{I}} [c_i(x^k) + \nabla c_i(x^k)^\mathsf{T} d^k]_+^2 + \sum_{i \in \mathcal{E}} [c_i(x^k) + \nabla c_i(x^k)^\mathsf{T} d^k]^2 < \sum_{i \in \mathcal{I}} [c_i(x^k)]_+^2 + \sum_{i \in \mathcal{E}} c_i(x^k)^2,$$

or if these two term equals, then

$$\nabla f(x^k)^\mathsf{T} d^k + \frac{1}{2} (d^k)^\mathsf{T} \nabla_{x,x}^2 \mathcal{L}(x^k, \lambda^k) d^k < 0.$$

In other words, the step $d^k$ either improves the feasibility of $x^k$ for the SQP subproblem or decreases the objective function of the SQP subproblem without worsening the aforementioned feasibility. In practice, it is possible that neither the constraint violation nor the objective function of the SQP subproblem can be reduced. A practical algorithm must handle such a case. However, from a theoretical point of view, it is possible to avoid such cases by making some mild assumptions (see, e.g., [173, Asm. 3.1]).

Note moreover that we imposed $\gamma^k \geq \|\lambda^k\|$. This is because Theorem 4.2.1 tells us that if $\gamma^k > \|\lambda^*\|$, where $\lambda^*$ is the Lagrange multiplier at the solution, then the merit function $\varphi^k$ enjoys an exactness property. Therefore, we attempt to fulfill this condition as much as possible using the information available so far. A scheme for updating the penalty parameter is given in Algorithm 4.3. Powell used in his solver COBYLA [161] a very similar scheme, with the constants $\upsilon_1 = 3/2$ and $\upsilon_2 = 2$. See [44, pp. 661–662] and [144, p. 542] for more discussions on the update of the penalty parameter.

**Algorithm 4.2:** Basic trust-region SQP method[†]

**Data:** Objective function $f$, constraint functions $\{c_i\}_{i\in\mathcal{I}\cup\mathcal{E}}$, initial guess $x^0 \in \mathbb{R}^n$, estimated Lagrange multiplier $\lambda^0 = [\lambda_i^0]_{i\in\mathcal{I}\cup\mathcal{E}}^\mathsf{T}$, initial trust-region radius $\Delta^0 > 0$, and parameters $0 < \eta_1 \leq \eta_2 < 1$ and $0 < \theta_1 < 1 < \theta_2$.

Set the penalty parameter $\gamma^{-1} \leftarrow 0$

**for** $k = 0, 1, \ldots$ *until convergence* **do**

Set the trial step $d^k$ to an approximate solution to

$$\min_{d\in\mathbb{R}^n} \quad \nabla f(x^k)^\mathsf{T} d + \frac{1}{2} d^\mathsf{T} \nabla^2_{x,x} \mathcal{L}(x^k, \lambda^k) d \tag{4.3.1a}$$

$$\text{s.t.} \quad c_i(x^k) + \nabla c_i(x^k)^\mathsf{T} d \leq 0, \; i \in \mathcal{I}, \tag{4.3.1b}$$

$$c_i(x^k) + \nabla c_i(x^k)^\mathsf{T} d = 0, \; i \in \mathcal{E}, \tag{4.3.1c}$$

$$\|d\| \leq \Delta^k, \tag{4.3.1d}$$

Pick a penalty parameter $\gamma^k \geq \max\{\gamma^{k-1}, \|\lambda^k\|\}$ providing $\widehat{\varphi}^k(d^k) < \widehat{\varphi}^k(0)$

Evaluate the trust-region ratio

$$\rho^k \leftarrow \frac{\varphi^k(x^k) - \varphi^k(x^k + d^k)}{\widehat{\varphi}^k(0) - \widehat{\varphi}^k(d^k)}$$

**if** $\rho^k \geq 0$ **then**

Update the trial point $x^{k+1} \leftarrow x^k + d^k$

**else**

Retain the trial point $x^{k+1} \leftarrow x^k$

**end if**

Estimate the Lagrange multiplier $\lambda^{k+1} = [\lambda_i^{k+1}]_{i\in\mathcal{I}\cup\mathcal{E}}^\mathsf{T}$

Update the trust-region radius

$$\Delta^{k+1} \leftarrow \begin{cases} \theta_1 \Delta^k, & \text{if } \rho^k \leq \eta_1, \\ \Delta^k, & \text{if } \eta_1 < \rho^k \leq \eta_2, \\ \theta_2 \Delta^k, & \text{otherwise} \end{cases}$$

**end for**

[†]Recall that superscripts are iteration counters rather than exponents.

**Algorithm 4.3:** Increasing the penalty parameter

**Data:** Merit function $\widehat{\varphi}_\gamma$ for $\gamma \geq 0$, trial step $d^k \in \mathbb{R}^n$, Lagrange multiplier estimate $\lambda^k$, previous penalty parameter $\gamma^{k-1} \geq 0$, and parameters $1 \leq \upsilon_1 \leq \upsilon_2$.

**Result:** Updated penalty parameter $\gamma^k \geq \gamma^{k-1}$.

```
Set γ̄ ← arg min{γ ≥ 0 : φ̂_γ(d^k) ≤ φ̂_γ(0)}
if γ^{k-1} ≤ υ_1 max{γ̄, ‖λ^k‖} then
    Set γ^k ← υ_2 max{γ̄, ‖λ^k‖}
else
    Set γ^k ← γ^{k-1}
end if
```

### 4.3.3 Composite steps for equality constrained optimization

The Line 3 of Algorithm 4.2 hides a crucial point: what should $d^k$ approximate if the trust-region subproblem (4.3.1) is infeasible? Since the constraints (4.3.1b) and (4.3.1c) are linear approximations of the constraints (4.1.1b) and (4.1.1c), we cannot guarantee the feasibility of the subproblem (4.3.1), even if the original problem is feasible.

We review the composite-step approach to coping with this difficulty in the equality-constrained case (i.e., $\mathcal{I} = \emptyset$). Most discussions here are covered in [44, § 15.4], except Proposition 4.3.1, which gives a simple observation on the Vardi approach introduced below.

The composite-step approach is to define $d^k$ as the sum of two steps, a *normal step* $n^k$ that aims at reducing the constraint violation, and a *tangential step* $t^k$ that aims at reducing the objective function (4.3.1a) without increasing the violation of the linearized constraints. More specifically, the normal step $n^k$ is an exact or approximate solution to the normal subproblem

$$\min_{d\in\mathbb{R}^n} \quad \sum_{i\in\mathcal{E}}[c_i(x^k) + \nabla c_i(x^k)^\mathsf{T} d]^2 \tag{4.3.2a}$$

$$\text{s.t.} \quad \|d\| \leq \zeta\Delta^k, \tag{4.3.2b}$$

for some $\zeta \in [0,1]$ (note that $\zeta$ should be positive unless $x^k$ is feasible). Further, the tangential step $t^k$ solves the tangential subproblem

$$\min_{d\in\mathbb{R}^n} \quad [\nabla f(x^k) + \nabla^2_{x,x}\mathcal{L}(x^k,\lambda^k)n^k]^\mathsf{T} d + \frac{1}{2}d^\mathsf{T}\nabla^2_{x,x}\mathcal{L}(x^k,\lambda^k)d \tag{4.3.3a}$$

$$\text{s.t.} \quad \|n^k + d\| \leq \Delta^k, \tag{4.3.3b}$$

$$d \in \mathcal{V}, \tag{4.3.3c}$$

either exactly or approximately, where $\mathcal{V}$ denotes some set such that for any $d \in \mathcal{V}$, $n^k + d$

does not violate the linearized constraints more than $n^k$ does. The exact formulation of $\mathcal{V}$ depends on how we measure the constraint violation.

In what follows, we describe three methods that are covered by this composite-step framework, namely the Byrd-Omojokun approach [30, 147], the Vardi approach [203], and the CDT approach [36]. It is noteworthy that Vardi did not explicitly formulate his method as a composite-step approach, but it fits well in the composite-step framework as explained in the sequel.

### Byrd-Omojokun approach

Even though the Byrd-Omojokun approach [30, 147] is not the earliest composite-step method, we first present it because it is the most straightforward implementation of the idea. It defines the normal step $n^k$ by solving (4.3.2) for some constant $\zeta \in (0, 1]$. Further, it defines the tangential step $t^k$ as a solution to (4.3.3) with

$$\mathcal{V} \stackrel{\text{def}}{=} \{d \in \mathbb{R}^n : \nabla c_i(x^k)^\mathsf{T} d = 0,\ i \in \mathcal{E}\}. \tag{4.3.4}$$

Note that $\zeta$ should not be set to one; otherwise, the feasible region of (4.3.3) might be a singleton. In practice, we can take for example $\zeta = 0.8$ (see, e.g., [144, Eq. (18.45)]).

### Vardi approach

The Vardi approach consists in replacing the constraints (4.3.1c) of the trust-region SQP subproblem by

$$\alpha c_i(x^k) + \nabla c_i(x^k)^\mathsf{T} d = 0,\ i \in \mathcal{E}, \tag{4.3.5}$$

for some $\alpha \in [0, 1]$, so that the trust-region subproblem after the replacement becomes feasible, this subproblem being

$$\min_{d \in \mathbb{R}^n} \quad \nabla f(x^k)^\mathsf{T} d + \frac{1}{2} d^\mathsf{T} \nabla^2_{x,x} \mathcal{L}(x^k, \lambda^k) d \tag{4.3.6a}$$

$$\text{s.t.} \quad \alpha c_i(x^k) + \nabla c_i(x^k)^\mathsf{T} d = 0,\ i \in \mathcal{E}, \tag{4.3.6b}$$

$$\|d\| \le \Delta^k. \tag{4.3.6c}$$

Note that if $\alpha = 0$, then (4.3.6) is clearly feasible. One can select $\alpha > 0$ under some assumptions, e.g., LICQ at $x^k$. The choice of $\alpha$ is important, and we refer to [44, § 15.4.1] for more discussions.

As explained in [44, § 15.4.1], the Vardi method can be understood as a composite-step approach as follows. Define the normal step $n^k$ as the least-norm solution to (4.3.5), which can be interpreted as an approximate solution to the normal subproblem (4.3.2) with a certain $\zeta \in [0, 1]$ (see Proposition 4.3.1). Further, define $t^k$ as a solution to (4.3.3) with $\mathcal{V}$ expressed in (4.3.4). It is then elementary to check that the composite step $n^k + t^k$ solves (4.3.6).

We now point out that the normal step $n^k$ is an approximate solution to (4.3.2) with $\zeta = \|n^k\|/\Delta^k$, which is at most one according to the definition of $n^k$.

**Proposition 4.3.1.** *Let $n^k$ be the least-norm solution to (4.3.5), and let $d^*$ be the least-norm global minimizer of (4.3.2a) without any constraint. Then $n^k$ solves*

$$\min_{d\in\mathbb{R}^n} \left\{ \sum_{i\in\mathcal{E}} [c_i(x^k) + \nabla c_i(x^k)^\mathsf{T} d]^2 : d = td^*,\ t \geq 0,\ \|d\| \leq \zeta\Delta^k \right\},$$

*with $\zeta = \|n^k\|/\Delta^k$.*

*Proof.* We let $A$ be the matrix whose rows are $\nabla c_i(x^k)^\mathsf{T}$ and $b$ the vector with components $c_i(x^k)$, with $i \in \mathcal{E}$, arranged in the same order. We have $d^* = -A^\dagger b$ and

$$n^k = \operatorname*{arg\,min}_{d\in\mathbb{R}^n}\{\|d\| : \alpha b + Ad = 0\} = -\alpha A^\dagger b = \alpha d^*,$$

where $A^\dagger$ denotes the Moore-Penrose pseudoinverse of $A$. Consequently,

$$\zeta\Delta^k = \|n^k\| = \alpha\|d^*\|.$$

Hence, we only need to show that $\alpha d^*$ is a least-squares solution to $b + Ad$ on the line segment $\{td^* : t \in [0, \alpha]\}$. This holds because $\alpha \leq 1$ and $\|b + Ad\|^2$ is nonincreasing when $d$ varies from $0$ to $d^*$ along a straight line, the monotonicity of $\|b + Ad\|^2$ coming from its convexity and the fact that $d^*$ is its global minimizer. □

Note that $d^*$ is a Newton step for (4.3.2a) at the trust region center. Therefore, the normal step $n^k$ of the Vardi approach is a truncated Newton step for the normal subproblem (4.3.2) if $\zeta = \|n^k\|/\Delta^k = \alpha\|d^*\|/\Delta^k$. This also reveals a relation between the Vardi and the Byrd-Omojokun approaches. The Vardi approach with $\alpha = \zeta\Delta^k/\|d^*\|$ can be regarded as an approximation to the Byrd-Omojokun approach because the former solves the normal subproblem (4.3.2) by a truncated Newton step, while the latter intends to solve this subproblem exactly, and the two approaches have no difference in the way they define the tangential step.

### Celis-Dennis-Tapia (CDT) approach

The CDT approach [35, 36] consists in replacing the constraints (4.3.1c) of the trust-region SQP subproblem with

$$\sum_{i\in\mathcal{E}} [c_i(x^k) + \nabla c_i(x^k)^\mathsf{T} d]^2 \leq \tau$$

for some $\tau \geq 0$, so that the trust-region subproblem after replacement becomes feasible, this subproblem being

$$\begin{aligned} \min_{d \in \mathbb{R}^n} \quad & \nabla f(x^k)^\mathsf{T} d + \frac{1}{2} d^\mathsf{T} \nabla^2_{x,x} \mathcal{L}(x^k, \lambda^k) d && \text{(4.3.7a)} \\ \text{s.t.} \quad & \sum_{i \in \mathcal{E}} [c_i(x^k) + \nabla c_i(x^k)^\mathsf{T} d]^2 \leq \tau, && \text{(4.3.7b)} \\ & \|d\| \leq \Delta^k. && \text{(4.3.7c)} \end{aligned}$$

The method selects the scalar $\tau$ in a composite-step fashion, as follows. It first takes a normal step $n^k$ by solving approximately the problem (4.3.2) with $\zeta = 1$, and then sets

$$\tau = \sum_{i \in \mathcal{E}} [c_i(x^k) + \nabla c_i(x^k)^\mathsf{T} n^k]^2.$$

The CDT approach does not contain a tangential step explicitly, but we can define it as $t^k = d^k - n^k$, which solves (4.3.3) with

$$\mathcal{V} = \left\{ d \in \mathbb{R}^n : \sum_{i \in \mathcal{E}} [c_i(x^k) + \nabla c_i(x^k)^\mathsf{T} (n^k + d)]^2 \leq \tau \right\}.$$

The set $\mathcal{V}$ consists of all the $d$ such that $n^k + d$ does not violate the linearized constraints more than $n^k$ in terms of the $\ell_2$ constraint violation.

Unlike the two former approaches, the tangential step of the CDT approach does not attempt to maintain or decrease the violation of each constraint but rather the global violation of all the constraints simultaneously. The advantage of the CDT approach over the Vardi and the Byrd-Omojokun approaches is that the feasible region of the tangential subproblem is usually wider given a certain normal step. Hence, a larger reduction in the objective function is usually obtained. However, solving such a subproblem in practice is challenging, because of the constraints (4.3.7b) and (4.3.7c). The subproblem, commonly known as the CDT subproblem, has enjoyed continued research interest since it was proposed (see, e.g., [39, 3, 21, 16, 214]).

### 4.3.4 Composite steps for general constrained optimization

In what follows, we discuss how to apply the composite-step approach to problems with inequality constraints, i.e., when $\mathcal{I} \neq \emptyset$.

For the normal step, as formulated in [44, § 15.4.4], we let $n^k$ solve

$$\begin{aligned} \min_{d \in \mathbb{R}^n} \quad & \sum_{i \in \mathcal{I}} [c_i(x^k) + \nabla c_i(x^k)^\mathsf{T} d]_+^2 + \sum_{i \in \mathcal{E}} [c_i(x^k) + \nabla c_i(x^k)^\mathsf{T} d]^2 && \text{(4.3.8a)} \\ \text{s.t.} \quad & \|d\| \leq \zeta \Delta^k, && \text{(4.3.8b)} \end{aligned}$$

for some $\zeta \in [0, 1]$, either exactly or approximately. Further, similar to the equality-constrained case, we let the tangential step $t^k$ solve (4.3.3) for some $\mathcal{V} \subseteq \mathbb{R}^n$ such that for any $d \in \mathcal{V}$, $n^k + d$ does not violate the linearized constraints more than $n^k$ does.

### Byrd-Omojokun approach

The Byrd-Omojokun approach has been originally introduced in a method that tackles both equality and inequality constraints. However, this method employs an active-set strategy, where the inequality constraints considered active are handled as equalities (see [147, § 3.2.6]). In other words, the Byrd-Omojokun approach has been defined from an equality-constrained perspective.

To extend the Byrd-Omojokun approach to the general case without employing an active-set method, we define the normal step $n^k$ by solving (4.3.8) for some $\zeta \in (0, 1]$. Further, we define the tangential step $t^k$ as a solution to (4.3.3) with the set $\mathcal{V}$ being all steps $d \in \mathbb{R}^n$ that satisfy

$$\begin{cases} \nabla c_i(x^k)^\mathsf{T} d \le \max\{-c_i(x^k) - \nabla c_i(x^k)^\mathsf{T} n^k, 0\}, \ i \in \mathcal{I}, & (4.3.9a) \\ \nabla c_i(x^k)^\mathsf{T} d = 0, \ i \in \mathcal{E}. & (4.3.9b) \end{cases}$$

Then we take the composite step $d^k = n^k + t^k$, which solves

$$\begin{aligned} \min_{d \in \mathbb{R}^n} \quad & \nabla f(x^k)^\mathsf{T} d + \frac{1}{2} d^\mathsf{T} \nabla^2_{x,x} \mathcal{L}(x^k, \lambda^k) d & (4.3.10a) \\ \text{s.t.} \quad & \nabla c_i(x^k)^\mathsf{T} d \le \max\{-c_i(x^k), \nabla c_i(x^k)^\mathsf{T} n^k\}, \ i \in \mathcal{I}, & (4.3.10b) \\ & \nabla c_i(x^k)^\mathsf{T} d = \nabla c_i(x^k)^\mathsf{T} n^k, \ i \in \mathcal{E}, & (4.3.10c) \\ & \|d\| \le \Delta^k, & (4.3.10d) \end{aligned}$$

This extension of the Byrd-Omojokun approach is different from the one presented in [44, § 15.4.4], which replaces the constraint (4.3.10b) with

$$\nabla c_i(x^k)^\mathsf{T} d \le \nabla c_i(x^k)^\mathsf{T} n^k, \ i \in \mathcal{I}. \tag{4.3.11}$$

Note that if $c_i(x^k) + \nabla c_i(x^k)^\mathsf{T} n^k < 0$, then the constraint in (4.3.11) is more restrictive than the original linearized constraint, which is $\nabla c_i(x^k)^\mathsf{T} d \le -c_i(x^k)$. Recall that the motivation of the composite-step approach is to relax the linearized constraints so that the trust-region subproblem becomes feasible. It seems, therefore, unnecessary to make a constraint become more restrictive. The feasible region of the subproblem (4.3.10) is usually wider than the one with (4.3.11) and hence, a larger reduction in the objective function (4.3.10a) can be expected. Note that the subproblem (4.3.10) is not more difficult to solve than its counterpart with (4.3.11). In the numerical experiments in Section 7.3.2, our Byrd-Omojokun approach performs much better than the one from [44, § 15.4.4].

### Vardi approach

The extension of the Vardi approach we propose is adapted from an idea proposed by Powell [157, Eqs. (2.7) and (2.8)] for a line-search SQP framework, which was published seven years earlier than Vardi's method.

We replace the constraints (4.3.1b) and (4.3.1c) of the trust-region SQP subproblem by

$$\begin{cases} c_i(x^k) + \nabla c_i(x^k)^\mathsf{T} d \le 0, \ i \in \mathcal{I} \setminus \mathcal{A}(x^k), & (4.3.12\text{a}) \\ \alpha c_i(x^k) + \nabla c_i(x^k)^\mathsf{T} d \le 0, \ i \in \mathcal{I} \cap \mathcal{A}(x^k), & (4.3.12\text{b}) \\ \alpha c_i(x^k) + \nabla c_i(x^k)^\mathsf{T} d = 0, \ i \in \mathcal{E}, & (4.3.12\text{c}) \end{cases}$$

for some $\alpha \in [0, 1]$, where $\mathcal{A}(x^k)$ designates the active set for the problem (4.1.1) at $x^k$ (see Definition 1.3.3), so that the trust-region subproblem after the replacement becomes feasible, this subproblem being

$$\min_{d \in \mathbb{R}^n} \quad \nabla f(x^k)^\mathsf{T} d + \frac{1}{2} d^\mathsf{T} \nabla^2_{x,x} \mathcal{L}(x^k, \lambda^k) d \qquad (4.3.13\text{a})$$

$$\text{s.t.} \quad c_i(x^k) + \nabla c_i(x^k)^\mathsf{T} d \le 0, \ i \in \mathcal{I} \setminus \mathcal{A}(x^k), \qquad (4.3.13\text{b})$$

$$\alpha c_i(x^k) + \nabla c_i(x^k)^\mathsf{T} d \le 0, \ i \in \mathcal{I} \cap \mathcal{A}(x^k), \qquad (4.3.13\text{c})$$

$$\alpha c_i(x^k) + \nabla c_i(x^k)^\mathsf{T} d = 0, \ i \in \mathcal{E}, \qquad (4.3.13\text{d})$$

$$\|d\| \le \Delta^k. \qquad (4.3.13\text{e})$$

Note that if $\alpha = 0$, then (4.3.13) is clearly feasible. The relaxation of the equality constraints (4.3.12c) is due to Vardi [203], but the relaxed inequality constraints (4.3.12a) and (4.3.12b) are adopted from Powell [157].

Remark that (4.3.12) does not modify a linearized inequality constraint if $x^k$ strictly satisfies the corresponding original constraint. Even more, if $x^k$ satisfies a constraint (whether inequality or equality), then its linearized counterpart is not modified. To be more specific, let $\mathcal{J}$ contain the indices of the constraints that are satisfied at $x^k$, i.e.,

$$\mathcal{J} \overset{\text{def}}{=} \{i \in \mathcal{I} \cup \mathcal{E} : c_i(x^k) \le 0, \ i \in \mathcal{I} \text{ and } c_i(x^k) = 0, \ i \in \mathcal{E}\}.$$

Then (4.3.12) does not modify the constraints in $\mathcal{J}$, and replaces the others by

$$\begin{cases} \alpha c_i(x^k) + \nabla c_i(x^k)^\mathsf{T} d \le 0, \ i \in \mathcal{I} \setminus \mathcal{J}, & (4.3.14\text{a}) \\ \alpha c_i(x^k) + \nabla c_i(x^k)^\mathsf{T} d = 0, \ i \in \mathcal{E} \setminus \mathcal{J}. & (4.3.14\text{b}) \end{cases}$$

This is because the difference between $(\mathcal{I} \cap \mathcal{A}(x^k)) \cup \mathcal{E}$ and $(\mathcal{I} \setminus \mathcal{J}) \cup (\mathcal{E} \setminus \mathcal{J})$ are the indices such that $c_i(x^k) = 0$, for which $\alpha c_i(x^k) = c_i(x^k)$.

We now explain why this method can be understood as a composite-step approach,

assuming that the system

$$
\begin{cases}
c_i(x^k) + \nabla c_i(x^k)^\mathsf{T} d \le 0, \ i \in \mathcal{I} \setminus \mathcal{A}(x^k), & \text{(4.3.15a)} \\
\alpha c_i(x^k) + \nabla c_i(x^k)^\mathsf{T} d = 0, \ i \in \mathcal{I} \cap \mathcal{A}(x^k), & \text{(4.3.15b)} \\
\alpha c_i(x^k) + \nabla c_i(x^k)^\mathsf{T} d = 0, \ i \in \mathcal{E}, & \text{(4.3.15c)}
\end{cases}
$$

admits a solution whose norm is at most $\Delta^k$, which ensures that (4.3.13) is feasible. Define the normal step $n^k$ as the least-norm solution to (4.3.15). Then $n^k$ can be regarded as an approximate solution to the normal subproblem (4.3.8) with $\zeta = 1$, because it satisfies $\|n^k\| \le \Delta^k$ and reduces (4.3.8a) to

$$
(1-\alpha)^2 \left[ \sum_{i \in \mathcal{I}} [c_i(x^k)]_+^2 + \sum_{i \in \mathcal{E}} c_i(x^k)^2 \right],
$$

being $(1-\alpha)^2$ times the corresponding value at the center of the trust region. Further, define $t^k$ as a solution to (4.3.3) with $\mathcal{V}$ being all steps $d \in \mathbb{R}^n$ that satisfy (4.3.9). In this setting, $n^k + t^k$ is a solution to (4.3.10), which turns out to be exactly (4.3.13), as specified by the following proposition.

**Proposition 4.3.2.** *Let $n^k$ be the least-norm solution to (4.3.15). Then the problems (4.3.10) and (4.3.13) are equivalent.*

*Proof.* We only need to show that the constraints (4.3.13b) to (4.3.13d) are equivalent to (4.3.10b) and (4.3.10c), the latter of which can be reformulated as

$$
\begin{cases}
\nabla c_i(x^k)^\mathsf{T} d \le \max\{-c_i(x^k), \nabla c_i(x^k)^\mathsf{T} n^k\}, \ i \in \mathcal{I} \setminus \mathcal{A}(x^k), & \text{(4.3.16a)} \\
\nabla c_i(x^k)^\mathsf{T} d \le \max\{-c_i(x^k), \nabla c_i(x^k)^\mathsf{T} n^k\}, \ i \in \mathcal{I} \cap \mathcal{A}(x^k), & \text{(4.3.16b)} \\
\nabla c_i(x^k)^\mathsf{T} d = \nabla c_i(x^k)^\mathsf{T} n^k, \ i \in \mathcal{E}. & \text{(4.3.16c)}
\end{cases}
$$

Since $\nabla c_i(x^k)^\mathsf{T} n^k = -\alpha c_i(x^k)$ for $i \in \mathcal{E}$, it is clear that (4.3.13d) is equivalent to (4.3.16c). Moreover, for $i \in \mathcal{I} \cap \mathcal{A}(x^k)$, we have $c_i(x^k) \ge 0$ and $\nabla c_i(x^k)^\mathsf{T} n^k = -\alpha c_i(x^k)$, so that (4.3.13c) is equivalent to (4.3.16b) because $\alpha \le 1$. Finally, for $i \in \mathcal{I} \setminus \mathcal{A}(x^k)$, since $n^k$ satisfies (4.3.15a), we have

$$
\max\{-c_i(x^k), \nabla c_i(x^k)^\mathsf{T} n^k\} = -c_i(x^k),
$$

so that (4.3.13b) is equivalent to (4.3.16a). The proof is complete. □

**Celis-Dennis-Tapia (CDT) approach**

The extension of the CDT approach is very direct and presented already in [44, § 15.4.4]. It consists in replacing the constraints (4.3.1b) and (4.3.1c) of the trust-region SQP subproblem by

$$\sum_{i\in\mathcal{I}}[c_i(x^k) + \nabla c_i(x^k)^\top d]_+^2 + \sum_{i\in\mathcal{E}}[c_i(x^k) + \nabla c_i(x^k)^\top d]^2 \leq \tau \tag{4.3.17}$$

for some $\tau \geq 0$, so that the trust-region subproblem after replacement becomes feasible. The method selects the scalar $\tau$ in a composite-step fashion as follows. It first takes a normal step $n^k$ by solving (approximately) the problem (4.3.8) with $\zeta = 1$, and then sets $\tau$ to the left-hand side in (4.3.17) with $d = n^k$.

### 4.3.5 Other methods for handling the trust-region SQP subproblem

There exist other methods to deal with the possible infeasibility of the trust-region SQP subproblem (4.3.1). In what follows, we briefly introduce the *sequential $\ell_p$-quadratic programming* (S$\ell_p$QP) and *sequential linear-quadratic programming* (SLQP) methods.

The S$\ell_p$QP method was first proposed by Fletcher [68] and further studied for example by Yuan [216]. Instead of the trust-region SQP subproblem (4.3.1), it solves

$$\min_{d\in\mathbb{R}^n} \quad \nabla f(x^k)^\top d + \frac{1}{2} d^\top \nabla^2_{x,x}\mathcal{L}(x^k, \lambda^k) d + \mu \Phi_p(d) \tag{4.3.18a}$$

$$\text{s.t.} \quad \|d\| \leq \Delta^k \tag{4.3.18b}$$

for some $\mu \geq 0$, where $\Phi_p$ is defined by the $\ell_p$-constraint violation of the linearized constraints, i.e.,

$$\Phi_p(d) \stackrel{\text{def}}{=} \left[\sum_{i\in\mathcal{I}}[c_i(x^k) + \nabla c_i(x^k)^\top d]_+^p + \sum_{i\in\mathcal{E}}|c_i(x^k) + \nabla c_i(x^k)^\top d|^p\right]^{1/p}, \quad \text{for } d \in \mathbb{R}^n.$$

When the trust region is defined by a polyhedral norm (e.g., the $\ell_1$- and $\ell_\infty$-norms), S$\ell_p$QP with $p \in \{1, \infty\}$ has the advantage that the subproblem (4.3.18) can be reformulated as a quadratic programming problem. Note that the original idea of Fletcher [68] is not to approximate the solution to the trust-region SQP subproblem (4.3.1) by solving (4.3.18). Instead, he proposed a trust-region method for minimizing composite nondifferentiable functions. When this method is applied to the $\ell_p$-exact penalty function of (4.1.1), the trust-region subproblem is (4.3.18). The idea of Yuan [216] is similar. For more information on S$\ell_p$QP, see [89, § 8.4.1].

SLQP [75, 31] computes an SQP step by a two-stage approach. The first stage minimizes approximately the linear part of the objective function of the SQP subproblem subject to the linearized constraints and a trust region. Similar to S$\ell_p$QP, it handles the linear constraints using an $\ell_1$-penalty function to avoid a possible infeasibility. With the approximate minimizer of the first stage, it estimates the active set at the solution, imposing the linear independence of the constraints in this set. The second stage minimizes the objective function of the SQP subproblem subject to the linearized constraints that are in the estimated active set, taking the inequalities as equalities. A major difference between [75] and [31] is that the former does not include a trust-region constraint in the second stage, while the latter does. In [75], the first-stage step is used as a safeguard that guarantees the global convergence, while [31] does not take such a safeguard, and the convergence is guaranteed by the trust-region framework.

## 4.4 Software based on the SQP method

Many well-known software packages are based on the SQP method. SNOPT [83, 84] for example is an implementation of the SQP method for inequality-constrained optimization. It approximates the Hessians of the Lagrangian function with a limited-memory quasi-Newton technique[6]. The SQP method is embedded in a line-search framework, which uses the augmented Lagrangian function as the merit function. It employs a reduced-Hessian algorithm (see [32] for more information about this class of methods) named SQOPT [85] for solving the subproblems. Moreover, the possible infeasibility of the SQP subproblem is handled using an $\ell_1$-penalty approach detailed in [84, § 1.2].

Another classical SQP solver is NLPQLP [189]. It is a line-search SQP method based on the augmented Lagrangian merit function. It approximates the Hessian of the Lagrangian function by the BFGS formula. The possible infeasibility of the SQP subproblem is handled by increasing the dimension of the decision variable by one to include a new variable that relaxes the infeasible constraints. Details on this mechanism are given in [188, Eq. (9)]. A remarkable feature of NLPQLP is its robustness with respect to inaccurate function and gradient evaluations. See [189] for more information.

Another instance is the ETR method [121] for equality-constrained optimization. It is a trust-region SQP method based on the $\ell_2$-merit function, which approximates the Hessians of the Lagrangian function by a limited-memory quasi-Newton method unless the user provides the Hessian matrices of the objective and constraint functions. The trust-region SQP subproblem is solved using a Byrd-Omojokun approach, described in Section 4.3.3.

We also mention CONDOR [201, 202], a trust-region method that extends Powell's UOBYQA [163] to constrained optimization. It is not an SQP method on its own but employs an SQP algorithm to solve its subproblem. This method assumes that derivatives of the objective function are unavailable, but both the values and the derivatives of the constraint functions can be evaluated inexpensively. It is a trust-region method based on the $\ell_1$-merit function, where at each iteration, a quadratic model of the objective function is minimized subject to the original constraints (rather than a model of them) and a trust-region constraint. The subproblem is solved with a "null-space, active-set approach" detailed in [201, § 9.3], and a line-search SQP method is applied to this subproblem in the presence of active nonlinear constraints.

[6] Although not studied in this thesis, the limited-memory quasi-Newton method [149, 193] is an important technique in numerical optimization. See the historical comments in [127, § 1, ¶ 3] for more details.

## 4.5 Summary and remarks

We presented in this chapter the classical SQP framework, and we got some insights by presenting several interpretations of its subproblem. To the best of our knowledge, the result in Theorem 4.1.1 is new. It shows that the objective function of the SQP subproblem is a natural approximation of the original objective function in the tangent space of a surface, this surface being the contour of the constraints at the current iterate. Moreover, the quality of this approximation depends on the stationary residual of the current approximate KKT pair. We also examined the augmented Lagrangian of the SQP subproblem. This augmented Lagrangian turns out to be indeed the objective function used in the augmented Lagrangian methods studied in [143, 205].

We then introduced the trust-region SQP framework and presented several methods to solve its subproblem. Interestingly, we derived in Proposition 4.3.1 a relation between two of these methods, namely the Byrd-Omojokun and the Vardi approaches, in the equality-constrained case. We also discussed how to extend these methods to handle inequality constraints, the extensions of the Byrd-Omojokun and the Vardi approaches being different from those presented in [44, § 15.4.4]. We will later in this thesis use the Byrd-Omojokun approach for our new optimization method presented in Chapter 5. As we will see in Section 7.3.2, our Byrd-Omojokun approach performs clearly better than the one in [44, § 15.4.4] for our new DFO software. The Vardi and the CDT approaches also need to be tested in the future.

*This page intentionally left blank*

# 5 COBYQA — a new derivative-free optimization method

This chapter introduces a new DFO method, named COBYQA. It is a derivative-free trust-region SQP method designed to tackle nonlinearly constrained optimization problems that include equality and inequality constraints. A particular feature of COBYQA is that it visits only points that respect the bound constraints, if any. This is useful because the objective functions of applications that admit bound constraints are often undefined when the bounds are violated.

Section 5.1 first formulates the problem to be tackled by COBYQA. Section 5.2 then presents the general framework of COBYQA and some important components. Section 5.3 provides further details on the management of the bound and linear constraints. Finally, Section 5.4 exhibits the merit function employed by COBYQA and the updates of its penalty parameter. Details on the methods used by COBYQA to solve its various subproblems will be provided in Chapter 6, and its Python implementation will be presented in Chapter 7.

## 5.1 Statement of the problem

The problems we consider in this chapter are of the form

$$\min_{x \in \mathbb{R}^n} \quad f(x) \tag{5.1.1a}$$

$$\text{s.t.} \quad c_i(x) \le 0, \; i \in \mathcal{I}, \tag{5.1.1b}$$

$$c_i(x) = 0, \; i \in \mathcal{E}, \tag{5.1.1c}$$

$$l \le x \le u, \tag{5.1.1d}$$

where $f$ and $c_i$ represent the objective and constraint functions, with $i \in \mathcal{I} \cup \mathcal{E}$ and the sets of indices $\mathcal{I}$ and $\mathcal{E}$ being finite and disjoint, but possibly empty. The lower bounds $l \in (\mathbb{R} \cup \{-\infty\})^n$ and the upper bounds $u \in (\mathbb{R} \cup \{\infty\})^n$ satisfy $l < u$. Note that the bound constraints (5.1.1d) are not included in the inequality constraints (5.1.1b), because they will be handled separately, as detailed in Section 5.3.

We will develop a derivative-free trust-region SQP method for the problem (5.1.1). The method, named COBYQA after *Constrained Optimization BY Quadratic Approximations*, does not use derivatives of the objective function or the nonlinear constraint functions, but models them using underdetermined interpolation based on the derivative-free symmetric Broyden update (see Section 2.4.2). This chapter presents the framework of the method, while the subproblems and the Python implementation will be discussed in Chapters 6 and 7, respectively.

## 5.2 The derivative-free trust-region SQP method

We presented in Chapter 4 the basic trust-region SQP method. We now adapt this framework to derivative-free settings.

### 5.2.1 Interpolation-based quadratic models

Recall that the basic trust-region SQP method presented in Algorithm 4.2 models the objective and constraint functions based on their gradients and Hessian matrices. Since we do not have access to such information, we use interpolation-based quadratic models of these functions. Specifically, we use quadratic models obtained by underdetermined interpolation based on the derivative-free symmetric Broyden update, detailed in Section 2.4.2.

At the $k$th iteration, we denote by $\hat{f}^k$ the quadratic model of $f$, and by $\hat{c}_i^k$ the quadratic model of $c_i$, for all $i \in \mathcal{I} \cup \mathcal{E}$. These models are built on an interpolation set $\mathcal{Y}^k \subseteq \mathbb{R}^n$, which is maintained by the method and updated along the iterations. It changes at most one point of $\mathcal{Y}^k$ at each iteration, and ensures that $x^k \in \mathcal{Y}^k$, where $x^k \in \mathbb{R}^n$ denotes the $k$th iterate.

The initial models $\hat{f}^0$ and $\hat{c}_i^0$, for $i \in \mathcal{I} \cup \mathcal{E}$, are built on the initial interpolation set $\mathcal{Y}^0 \subseteq \mathbb{R}^n$, defined as follows according to Powell [168]. The user provides a number $m$, satisfying

$$n + 2 \le m \le \frac{1}{2}(n + 1)(n + 2),$$

which will be the number interpolation points at each iteration. We are also provided with an initial guess $x^0 \in \mathbb{R}^n$ and an initial trust-region radius $\Delta^0 > 0$ that satisfies

$$\Delta^0 \le \frac{1}{2} \min_{1 \le i \le n} (u_i - l_i),$$

so that $u_i \ge l_i + 2\Delta^0$ for all $i \in \{1, 2, \dots, n\}$. Further, to avoid conflicts between the bounds and the interpolation points, $x^0$ is modified so that each component of $x^0$ is either on a bound or keeps a distance of at least $\Delta^0$ from the corresponding bounds, while being feasible.

The initial interpolation set $\mathcal{Y}^0 = \{y^1, y^2, \dots, y^m\}$ is constructed in the following way. The first $n + 1$ points of $\mathcal{Y}^0$ are defined by

$$y^i \stackrel{\text{def}}{=} \begin{cases} x^0 & \text{if } i = 1, \\ x^0 + \Delta^0 e_{i-1} & \text{if } 2 \le i \le n + 1 \text{ and } l_{i-1} \le x_{i-1}^0 < u_{i-1}, \\ x^0 - \Delta^0 e_{i-1} & \text{if } 2 \le i \le n + 1 \text{ and } x_{i-1}^0 = u_{i-1}, \end{cases}$$

and the following $\min\{n, m-n-1\}$ are defined by

$$y^i \overset{\text{def}}{=} \begin{cases} x^0 - \Delta^0 e_{i-n-1} & \text{if } n+2 \le i \le \min\{2n+1, m\} \text{ and } l_{i-n-1} < x^0_{i-n-1} \le u_{i-n-1}, \\ x^0 + 2\Delta^0 e_{i-n-1} & \text{if } n+2 \le i \le \min\{2n+1, m\} \text{ and } x^0_{i-n-1} = l_{i-n-1}, \\ x^0 - 2\Delta^0 e_{i-n-1} & \text{if } n+2 \le i \le \min\{2n+1, m\} \text{ and } x^0_{i-n-1} = u_{i-n-1}, \end{cases}$$

where $e_i$ denotes the $i$th standard coordinate vector of $\mathbb{R}^n$, i.e., the $i$th column of $I_n$. If $m > 2n+1$, for $i \in \{2n+2, 2n+3, \dots, m\}$, we set

$$y^i \overset{\text{def}}{=} y^{p(i)+1} + y^{q(i)+1} - x^0,$$

where $p$ and $q$ are defined by

$$p(i) \overset{\text{def}}{=} i - n - 1 - n\kappa(i) \quad \text{with} \quad \kappa(i) \overset{\text{def}}{=} \left\lfloor \frac{i-n-2}{n} \right\rfloor,$$

and

$$q(i) \overset{\text{def}}{=} \begin{cases} p(i) + \kappa(i) & \text{if } p(i) + \kappa(i) \le n, \\ p(i) + \kappa(i) - n & \text{otherwise.} \end{cases}$$

We point out that if $m \le 2n+1$ and $l < x^0 < u$, then the interpolation set $\mathcal{Y}^0$ is essentially the interpolation set studied in Section 2.5, which focus on $\Delta^0 = 1$. We have shown that this set is optimal for $m = 2n+1$, in a sense specified in Section 2.5. Therefore, as Powell did, we take the value $m = 2n+1$ by default.

An advantage of choosing such an initial interpolation set is that the coefficients of the minimum-norm quadratic interpolation model can be directly evaluated. See [168, § 9] for the detailed calculations.

Once the initial models are constructed, they will be maintained by the derivative-free symmetric Broyden update, in the same way as NEWUOA [166], BOBYQA [168], and LINCOA. See [164, 165], [166, § 4], and [168, § 4] for details.

### 5.2.2 A derivative-free trust-region SQP framework

We now present the derivative-free trust-region SQP framework employed by COBYQA. We denote for convenience by $\widehat{\mathcal{L}}^k$ the Lagrangian function evaluated on the models, i.e.,

$$\widehat{\mathcal{L}}^k(x, \lambda) \overset{\text{def}}{=} \hat{f}^k(x) + \sum_{i \in \mathcal{I} \cup \mathcal{E}} \lambda_i \hat{c}^k_i(x), \quad \text{for } x \in \mathbb{R}^n \text{ and } \lambda_i \in \mathbb{R}, \text{ with } i \in \mathcal{I} \cup \mathcal{E}.$$

As in Section 4.3, the merit function we consider is the $\ell_2$-merit function, defined for a given penalty parameter $\gamma^k \ge 0$ by

$$\varphi^k(x) \overset{\text{def}}{=} f(x) + \gamma^k \sqrt{\sum_{i \in \mathcal{I}} [c_i(x)]_+^2 + \sum_{i \in \mathcal{E}} |c_i(x)|^2}, \quad \text{for } x \in \mathbb{R}^n.$$

We denote by $\widehat{\varphi}^k$ the $\ell_2$-merit function computed on the SQP subproblem, i.e.,

$$\widehat{\varphi}^k(d) \overset{\text{def}}{=} \nabla \hat{f}^k(x^k)^\mathsf{T} d + \frac{1}{2} d^\mathsf{T} \nabla^2_{x,x} \widehat{\mathcal{L}}^k(x^k, \lambda^k) d + \gamma^k \Phi(d), \quad \text{for } d \in \mathbb{R}^n, \tag{5.2.1}$$

where $\Phi$ is defined by

$$\Phi(d) \overset{\text{def}}{=} \sqrt{\sum_{i \in \mathcal{I}} [\hat{c}_i^k(x^k) + \nabla \hat{c}_i^k(x^k)^\mathsf{T} d]_+^2 + \sum_{i \in \mathcal{E}} [\hat{c}_i^k(x^k) + \nabla \hat{c}_i^k(x^k)^\mathsf{T} d]^2}. \tag{5.2.2}$$

The framework is given in Algorithm 5.1, which hides many details, including the definition of "convergence" in the loop, the calculation of $d^k$, the update of $\lambda^k$, the update of $\Delta^k$, and the geometry improvement of the interpolation set, which will be specified in Sections 5.2.3 to 5.2.5 and 5.2.7, respectively.

In the implementation of Line 16 in Algorithm 5.1, $x^{k+1}$ is always set to a best point in $\mathcal{Y}^{k+1}$ according to $\varphi^k$, in the sense that it minimizes $\varphi^k$ in $\mathcal{Y}^{k+1}$. If $x^k$ turns out to be optimal in $\mathcal{Y}^{k+1}$ according to $\varphi^k$, then we choose $x^{k+1} = x^k$, even if there are other optimal points. Note however that we cannot naively set $x^{k+1}$ to either $x^k + d^k$ if $\rho^k > 0$, or $x^k$ otherwise. This is because $x^k$ is a best point in $\mathcal{Y}^k$ according to $\varphi^{k-1}$, but it may not be the best one according to $\varphi^k$.

### 5.2.3 Solving the trust-region SQP subproblem

To solve the trust-region SQP subproblem (5.2.3), we employ the Byrd-Omojokun approach that we presented in Section 4.3.4. It first generates a normal step $n^k$ as an approximate solution to

$$\min_{d \in \mathbb{R}^n} \quad \sum_{i \in \mathcal{I}} [\hat{c}_i^k(x^k) + \nabla \hat{c}_i^k(x^k)^\mathsf{T} d]_+^2 + \sum_{i \in \mathcal{E}} [\hat{c}_i^k(x^k) + \nabla \hat{c}_i^k(x^k)^\mathsf{T} d]^2 \tag{5.2.4a}$$

$$\text{s.t.} \quad l \le x^k + d \le u, \tag{5.2.4b}$$

$$\|d\| \le \zeta \Delta^k, \tag{5.2.4c}$$

for some $\zeta \in (0, 1)$. In the implementation, we choose $\zeta = 0.8$. If $x^k$ is feasible for the original problem (5.1.1), we choose $n^k = 0$. Further, it generates a tangential step $t^k$ by solving approximately

$$\min_{d \in \mathbb{R}^n} \quad [\nabla \hat{f}^k(x^k) + \nabla^2_{x,x} \widehat{\mathcal{L}}^k(x^k, \lambda^k) n^k]^\mathsf{T} d + \frac{1}{2} d^\mathsf{T} \nabla^2_{x,x} \widehat{\mathcal{L}}^k(x^k, \lambda^k) d \tag{5.2.5a}$$

$$\text{s.t.} \quad \nabla \hat{c}_i^k(x^k)^\mathsf{T} d \le \max\{-\hat{c}_i^k(x^k) - \nabla \hat{c}_i^k(x^k)^\mathsf{T} n^k, 0\}, \ i \in \mathcal{I}, \tag{5.2.5b}$$

$$\nabla \hat{c}_i^k(x^k)^\mathsf{T} d = 0, \ i \in \mathcal{E}, \tag{5.2.5c}$$

$$l \le x^k + n^k + d \le u, \tag{5.2.5d}$$

$$\|n^k + d\| \le \Delta^k, \tag{5.2.5e}$$

**Algorithm 5.1:** Derivative-free trust-region SQP method[†]

**Data:** Objective function $f$, constraint functions $\{c_i\}_{i \in \mathcal{I} \cup \mathcal{E}}$, initial guess $x^0 \in \mathbb{R}^n$, and initial trust-region radius $\Delta^0 > 0$.

Set the penalty parameter $\gamma^{-1} \leftarrow 0$
Build the initial interpolation set $\mathcal{Y}^0 \subseteq \mathbb{R}^n$ described in Section 5.2.1
Redefine $x^0$ to a solution to $\min_{y \in \mathcal{Y}^0} f(y)$
Estimate the initial Lagrange multiplier $\lambda^0 = [\lambda_i^0]_{i \in \mathcal{I} \cup \mathcal{E}}^\mathsf{T}$
**for** $k = 0, 1, \dots$ *until convergence* **do**
    Update $\hat{f}^k$ and $\hat{c}_i^k$ for $i \in \mathcal{I} \cup \mathcal{E}$ as mentioned in Section 5.2.1
    Set the trial step $d^k \in \mathbb{R}^n$ to an approximate solution to

$$
\begin{aligned}
\min_{d \in \mathbb{R}^n} \quad & \nabla \hat{f}^k(x^k)^\mathsf{T} d + \frac{1}{2} d^\mathsf{T} \nabla_{x,x}^2 \widehat{\mathcal{L}}^k(x^k, \lambda^k) d && (5.2.3a)\\
\text{s.t.} \quad & \hat{c}_i^k(x^k) + \nabla \hat{c}_i^k(x^k)^\mathsf{T} d \le 0, \; i \in \mathcal{I}, && (5.2.3b)\\
& \hat{c}_i^k(x^k) + \nabla \hat{c}_i^k(x^k)^\mathsf{T} d = 0, \; i \in \mathcal{E}, && (5.2.3c)\\
& l \le x^k + d \le u, && (5.2.3d)\\
& \|d\| \le \Delta^k && (5.2.3e)
\end{aligned}
$$

    Pick a penalty parameter $\gamma^k \ge \max\{\gamma^{k-1}, \|\lambda^k\|\}$ providing $\widehat{\varphi}^k(d^k) < \widehat{\varphi}^k(0)$
    Evaluate the trust-region ratio

$$
\rho^k \leftarrow \frac{\varphi^k(x^k) - \varphi^k(x^k + d^k)}{\widehat{\varphi}^k(0) - \widehat{\varphi}^k(d^k)}
$$

    **if** $\rho^k > 0$ **then**
        Choose a point $\bar{y} \in \mathcal{Y}^k$ to remove from $\mathcal{Y}^k$
    **else**
        Choose a point $\bar{y} \in \mathcal{Y}^k \setminus \{x^k\}$ to remove from $\mathcal{Y}^k$
    **end if**
    Update the interpolation set $\mathcal{Y}^{k+1} \leftarrow (\mathcal{Y}^k \setminus \{\bar{y}\}) \cup \{x^k + d^k\}$
    Update the current iterate $x^{k+1}$ to a solution to $\min_{y \in \mathcal{Y}^{k+1}} \varphi^k(y)$
    Estimate the Lagrange multiplier $\lambda^{k+1} = [\lambda_i^{k+1}]_{i \in \mathcal{I} \cup \mathcal{E}}^\mathsf{T}$
    Update the trust-region radius from $\Delta^k$ to $\Delta^{k+1}$
    Improve the geometry of $\mathcal{Y}^{k+1}$ if necessary
**end for**

[†]Recall that superscripts are iteration counters rather than exponents.

and then sets the composite step $d^k = n^k + t^k$. Details on the calculations of the tangential and normal steps are provided in Sections 6.2 and 6.3, respectively.

As mentioned in Section 4.3.4, our Byrd-Omojokun approach is different from the one proposed in [44, § 15.4.4], which replaces the constraint (5.2.5b) with

$$\nabla \hat{c}_i^k(x^k)^\mathsf{T} d \le 0, \ i \in \mathcal{I}.$$

Our approach provides evidently better performance in practice for COBYQA, as is demonstrated by the numerical experiments in Section 7.3.2.

If we have $\mathcal{I} \cup \mathcal{E} = \emptyset$, we then necessarily have $n^k = 0$ because our algorithm guarantees that $l \le x^k \le u$. The tangential subproblem is then identical to the original subproblem (5.2.3). In COBYQA, such a subproblem is approximately solved using the TCG method that Powell employed in BOBYQA [168].

When $\mathcal{I} \cup \mathcal{E} \ne \emptyset$, the tangential subproblem (5.2.5) is approximately solved using the TCG method that Powell designed for LINCOA [171]. This method requires that

1. the origin $d = 0$ is feasible, and
2. the trust region is centered on the origin.

The first requirement is satisfied by (5.2.5). To meet the second one, following the work of [121, Eq. (2.10)] and [44, Eq. (15.4.3)], we can solve instead

$$\begin{aligned}
\min_{d \in \mathbb{R}^n} \quad & [\nabla \hat{f}^k(x^k) + \nabla_{x,x}^2 \widehat{\mathcal{L}}^k(x^k, \lambda^k) n^k]^\mathsf{T} d + \frac{1}{2} d^\mathsf{T} \nabla_{x,x}^2 \widehat{\mathcal{L}}^k(x^k, \lambda^k) d && (5.2.6a)\\
\text{s.t.} \quad & \nabla \hat{c}_i^k(x^k)^\mathsf{T} d \le \max\{-\hat{c}_i^k(x^k) - \nabla \hat{c}_i^k(x^k)^\mathsf{T} n^k, 0\}, \ i \in \mathcal{I}, && (5.2.6b)\\
& \nabla \hat{c}_i^k(x^k)^\mathsf{T} d = 0, \ i \in \mathcal{E}, && (5.2.6c)\\
& l \le x^k + n^k + d \le u, && (5.2.6d)\\
& \|d\| \le \sqrt{(\Delta^k)^2 - \|n^k\|^2}, && (5.2.6e)
\end{aligned}$$

the only difference with (5.2.5) being the trust-region constraint (5.2.6e). In doing so, we would have

$$\|d^k\| = \|n^k + t^k\| \le \max_{d \in \mathbb{R}^n} \left\{ \|d\| + \sqrt{(\Delta^k)^2 - \|d\|^2} : \|d\| \le \Delta^k \right\} = \sqrt{2}\Delta^k.$$

Therefore, to maintain the property $\|d^k\| \le \Delta^k$, we replace $\Delta^k$ in the constraints (5.2.4c) and (5.2.6e) with $\Delta^k/\sqrt{2}$.

### 5.2.4 Estimating the Lagrange multiplier

We now introduce the strategy employed by COBYQA to estimate the Lagrange multiplier on Line 17 of Algorithm 5.1.

Similar to [44, § 15.2, p. 626] and [144, Eq. (18.22)], we let $\lambda^{k+1}$ be a *least-squares*

*Lagrange multiplier* (see also [61, § 3.3]), i.e., a solution to

$$\min_{\lambda} \quad \left\| \nabla \hat{f}^k(x^k) + \sum_{i \in \mathcal{I} \cup \mathcal{E}} \lambda_i \nabla \hat{c}_i^k(x^k) \right\| \tag{5.2.7a}$$

$$\text{s.t.} \quad \lambda_i = 0, \ i \in \{j \in \mathcal{I} : \hat{c}_j^k(x^k) < 0\} \tag{5.2.7b}$$

$$\lambda_i \geq 0, \ i \in \mathcal{I}, \tag{5.2.7c}$$

where $\lambda = [\lambda_i]_{i \in \mathcal{I} \cup \mathcal{E}}^{\mathsf{T}}$. Note that for all $i \in \mathcal{I}$, $\hat{c}_i^k(x^k) = c_i(x^k)$, because $\hat{c}_i^k$ interpolates $c_i$ on $\mathcal{Y}^k$ and $x^k \in \mathcal{Y}^k$. The conditions (5.2.7b) correspond to the complementary slackness conditions of the inactive constraints at $x^k$. Even though the complementary slackness conditions at the solution are $\lambda_i^* c_i(x^*) = 0$ for $i \in \mathcal{I}$, it is unreasonable to force $\lambda_i^{k+1} = 0$ if $c_i(x^k) > 0$ and $i \in \mathcal{I}$, because the corresponding constraint is currently active. Details on the calculations of the least-squares Lagrange multiplier are provided in Section 6.4.

The initial Lagrange multiplier $\lambda^0$ on Line 4 of Algorithm 5.1 is evaluated in a very similar fashion. It is set to a solution to (5.2.7) with $k = 0$, so that it equals $\lambda^1$.

We also mention that the least-squares multiplier can be regarded as an approximation of the multiplier of the SQP subproblem, known as the QP multiplier. See [144, pp. 538–539] for more details.

### 5.2.5 Managing the trust-region radius

We now detail the trust-region radius management strategy employed by COBYQA, corresponding to Line 18 of Algorithm 5.1. The strategy does not only maintain the trust-region radius $\Delta^k$ but also a lower bound $\delta^k$ of it. This strategy is taken from UOBYQA [163], NEWUOA [166], BOBYQA [168], and LINCOA, which were briefly described in Sections 3.2.3 to 3.2.6.

The update of $\Delta^k$ is given in Algorithm 5.2, which is typical for trust-region methods, except that the lower bound $\delta^k$ is imposed. The idea behind $\delta^k$ will be explained later. The parameters chosen in COBYQA are $\eta_1 = 0.1$, $\eta_2 = 0.7$, $\theta_1 = 0.5$, $\theta_2 = 1.4$, $\theta_3 = \sqrt{2}$, and $\theta_4 = 2$. Note that it is crucial to ensure that $\theta_2 < \theta_3$. Otherwise, if $\Delta^k = \delta^k$ and the trial step $d^k$ performs very well (i.e., $\rho^k > \eta_2$), then the method would set $\Delta^{k+1} = \Delta^k$, but $\Delta^{k+1} > \Delta^k$ is expected if $\|d^k\|$ is close to $\Delta^k$.

The value of $\delta^k$ can be regarded as an indicator of the *resolution* of the algorithm at the current iteration. If we do not impose $\Delta^k \geq \delta^k$, the trust-region radius $\Delta^k$ may be reduced to a value that is too small for the current precision of the algorithm, making the interpolation points concentrate too much. The value of $\delta^k$ is never increased and is decreased when the algorithm decides that the work for the current value of $\delta^k$ is finished. This is reflected in practice by

1. the fact that $\Delta^k$ reaches its lower bound $\delta^k$,
2. a poor performance of the trust-region trial step $d^k$, and

**Algorithm 5.2:** Updating the trust-region radius

**Data:** Current lower bound on the trust-region radius $\delta^k > 0$, current trust-region radius $\Delta^k \geq \delta^k$, current trust-region ratio $\rho^k \in \mathbb{R}$, current trial step $d^k \in \mathbb{R}^n$, and parameters $0 < \eta_1 \leq \eta_2 < 1$, and $0 < \theta_1 < 1 \leq \theta_2 < \theta_3 < \theta_4$.

**Result:** Updated trust-region radius $\Delta^{k+1}$.

Update the trust-region radius

$$\Delta^{k+1} \leftarrow \begin{cases} \theta_1 \Delta^k & \text{if } \rho^k \leq \eta_1, \\ \max\{\theta_1 \Delta^k, \|d^k\|\} & \text{if } \eta_1 < \rho^k \leq \eta_2, \\ \min\{\theta_3 \Delta^k, \max\{\theta_1 \Delta^k, \theta_4 \|d^k\|\}\} & \text{otherwise} \end{cases}$$

**if** $\Delta^{k+1} \leq \theta_2 \delta^k$ **then**

    $\Delta^{k+1} \leftarrow \delta^k$

**end if**

3. the fact that the interpolation points are all close enough to $x^{k+1}$.

More specifically, similar to Powell's solvers, COBYQA reduces the value of $\delta^k$ if

$$\delta^k = \Delta^k, \quad \rho^k \leq \eta_1, \quad \text{and} \quad \max_{y \in \mathcal{Y}^{k+1}} \|y - x^{k+1}\| \leq 2\delta^k.$$

There is also another exceptional case where $\delta^k$ is reduced, which will be detailed at the end of Section 5.2.7. When COBYQA decides to reduce $\delta^k$, it invokes Algorithm 5.3 to do so, with the parameters $\eta_3 = 16$, $\eta_4 = 250$, and $\theta_5 = 0.1$. Algorithm 5.3 and these parameters are adopted from UOBYQA [163], NEWUOA [166], BOBYQA [168], and LINCOA.

If Line 2 of Algorithm 5.3 is reached, then the algorithm has finished the work for the smallest desired resolution, and hence, we stop the computations. This is, in fact, the definition of "convergence" we mentioned in Algorithm 5.1, which is also the stopping criterion used in Powell's solvers.

### 5.2.6 Updating the interpolation set

We now detail how COBYQA updates the interpolation set $\mathcal{Y}^k$. To this end, we only need to detail the choice of $\bar{y}$ in Lines 11 and 13 of Algorithm 5.1. Then, $\mathcal{Y}^{k+1}$ will be set to $(\mathcal{Y}^k \setminus \{\bar{y}\}) \cup \{x^k + d^k\}$.

Let us first consider the case $\rho^k \geq 0$, so that $\bar{y}$ is chosen from $\mathcal{Y}^k$. We denote by $\bar{x}$ the closest point from $x^k$ that solves

$$\min_{y \in \mathcal{Y}^k} \varphi^k(y).$$

**Algorithm 5.3:** Reducing the lower bound on the trust-region radius

**Data:** Final trust-region radius $\delta^\infty > 0$, current lower bound on the trust-region radius $\delta^k \geq \delta^\infty$, updated trust-region radius $\Delta^{k+1} \geq \delta^k$, and parameters $1 \leq \eta_3 < \eta_4$ and $0 < \theta_5 < 1$.

**Result:** Reduced lower bound on trust-region radius $\delta^{k+1}$ and modified trust-region radius $\Delta^{k+1}$.

```
if δ^k = δ^∞ then
    Terminate the optimization method
else
    Update the lower bound on the trust-region radius
```

$$\delta^{k+1} \leftarrow \begin{cases} \theta_5 \delta^k & \text{if } \eta_4 < \delta^k/\delta^\infty, \\ \sqrt{\delta^k \delta^\infty} & \text{if } \eta_3 < \delta^k/\delta^\infty \leq \eta_4, \\ \delta^\infty & \text{otherwise} \end{cases}$$

```
    Update the trust-region radius Δ^{k+1} ← max{Δ^{k+1}, δ^{k+1}}
end if
```

This point is not necessarily $x^k$, because $\varphi^k$ may differ from $\varphi^{k-1}$ due to the update of the penalty parameter. If several such points exist, we select any of them. We then set $\bar{y}$ to a solution to

$$\max_{y \in \mathcal{Y}^k} |\omega(y)| \|y - \bar{x}\|^4, \tag{5.2.8}$$

where $\omega(y)$ is a scaling factor, detailed later. The term $\|y - \bar{x}\|$ in (5.2.8) intends to set $\bar{y}$ to a point far away from $\bar{x}$, because the interpolation points must be reasonably close in order to achieve a good precision of the models. To explain the scalar factor $\omega(y)$, recall from Section 2.4.3 that the inverse of the coefficient matrix of the interpolation problem is maintained. Whenever an interpolation point is replaced with a new one, this inverse will be updated by the formula detailed in [165, Eq. (2.12)], which has a scalar denominator. We set $\omega(y)$ to this denominator corresponding to the replacement of $y \in \mathcal{Y}^k$ with $x^k + d^k$. This is because we want this denominator to be far away from zero to avoid numerical difficulties. The choice of the formula (5.2.8) is taken from LINCOA, the most recent DFO solver of Powell. He also made very similar choices for NEWUOA [166, Eq. (7.4)] and BOBYQA [168, Eq. (6.1)].

The strategy for the case $\rho^k < 0$ is quite similar. We set $\bar{y}$ to a solution to

$$\max_{y \in \mathcal{Y}^k \setminus \{x^k\}} |\omega(y)| \|y - \bar{x}\|^4.$$

This case enforces $\bar{y} \neq x^k$ to avoid removing $x^k$ from $\mathcal{Y}^k$, as this point is better than $x^k + d^k$ according to $\varphi^{k-1}$.

Theoretically, it is tempting to set $\bar{y}$ to the point that optimizes the $\Lambda$-poisedness (see Section 2.4.2) of $\mathcal{Y}^{k+1}$ in some set, say

$$\{x \in \mathbb{R}^n : \|x - \bar{x}\| \le \delta^k\}. \tag{5.2.9}$$

This idea has, however, a major drawback: it is expensive to implement due to the definition of the $\Lambda$-poisedness. Hence, from a computational perspective, we decide not to employ this technique in COBYQA.

### 5.2.7 Improving the geometry of the interpolation set

We now detail the geometry-improving procedure, corresponding to Line 19 of Algorithm 5.1. As we mentioned in Chapter 2, the interpolation set has to be reasonably well-poised to ensure the accuracy of the interpolation models. However, it is known that model-based DFO methods tend to lose the poisedness of their interpolation set as the iterations progress. To cope with this difficulty, such methods usually include a geometry-improving mechanism. This topic is well studied. See, e.g., [45, 46, 65, 185].

We adapt the geometry-improving strategies designed by Powell for BOBYQA [168] and LINCOA. Briefly speaking, a point $\bar{y} \in \mathcal{Y}^{k+1}$ is chosen to be replaced with another one $x^{k+1} + r^{k+1}$, where $r^{k+1} \in \mathbb{R}^n$ is a *model step* to improve the geometry of $\mathcal{Y}^{k+1}$.

First of all, the point $\bar{y}$ is chosen to be a solution to

$$\max_{y \in \mathcal{Y}^{k+1}} \|y - x^{k+1}\|. \tag{5.2.10}$$

Further, similar to [166, Eq. (6.6)], the step $r^{k+1} \in \mathbb{R}^n$ should provide a large value of $|L_{\bar{y}}(x^{k+1} + r^{k+1})|$, where $L_{\bar{y}}$ denotes the minimum Frobenius norm Lagrange polynomials for $\mathcal{Y}^{k+1}$ associated with $\bar{y}$ (see Definition 2.4.2). The reason for taking such a step is as follows. Recall that we need to update the inverse of the coefficient matrix of the interpolation problem to replace $\bar{y}$ with $x^{k+1} + r^{k+1}$ in $\mathcal{Y}^{k+1}$. As already mentioned, the updating formula has a scalar denominator. This denominator is always lower bounded by $L_{\bar{y}}(x^{k+1} + r^{k+1})^2$ as shown in [166, Eq. (6.5)]. Therefore, by selecting a model step $r^{k+1}$ that provides a large value of $|L_{\bar{y}}(x^{k+1} + r^{k+1})|$, we keep the denominator away from zero, avoiding numerical difficulties. In COBYQA, the model step $r^{k+1}$ is selected from two alternatives $r_B^{k+1}$ and $r_L^{k+1}$, detailed below.

The first step $r_B^{k+1}$, adopted from BOBYQA, is an approximate solution to

$$\max_{d \in \mathbb{R}^n} \quad |L_{\bar{y}}(x^{k+1} + d)| \tag{5.2.11a}$$

$$\text{s.t.} \quad l \le x^{k+1} + d \le u, \tag{5.2.11b}$$

$$\|d\| \le \bar{\Delta}, \tag{5.2.11c}$$

where, as in LINCOA, we set $\bar{\Delta} = \max\{\Delta^{k+1}/10, \delta^{k+1}\}$. Note that the linearized con-straints are not taken into account in this subproblem. We will detail how to solve the subproblem (5.2.11) in Section 6.5.

The second step $r_L^{k+1}$, adopted from LINCOA, is defined as follows. Let $\mathcal{W}$ be some approximation of the active set of the linearized constraints at $x^{k+1}$. In COBYQA, we set $\mathcal{W}$ to be the working set returned by the active-set method that evaluates the tangential step (5.2.6) (see Section 6.2 for details). Then $r_L^{k+1}$ is a Cauchy step for $|L_{\bar{y}}(x^{k+1}+d)|$ in the null space of the linearized constraints in $\mathcal{W}$ subject to $\|d\| \le \bar{\Delta}$. In our experiments, considering this step turns out to be crucial for the performance of COBYQA.

The selection between $r_B^{k+1}$ and $r_L^{k+1}$ is done as follows. We set $r^{k+1}$ to $r_L^{k+1}$ if

1. the point $x^{k+1} + r_L^{k+1}$ is feasible with respect to the linearized constraints and the bounds, and
2. the denominator of the updating formula corresponding to $r_L^{k+1}$ is at least a tenth of the one corresponding to $r_B^{k+1}$ in absolute value.

Otherwise, we set $r^{k+1}$ to $r_B^{k+1}$. When checking the feasibility of $x^{k+1} + r_L^{k+1}$, we allow small infeasibility of the linearized constraints, to include a tolerance for contributions from computer rounding errors.

Now that the general framework of the geometry-improving phase is set up, we need to decide when to entertain such a mechanism. It is triggered when

1. the trust-region trial step $d^k$ performed poorly, and
2. some interpolation points are far away from $x^{k+1}$.

More specifically, COBYQA modifies the interpolation set $\mathcal{Y}^{k+1}$ if

$$\rho^k \le \eta_1 \quad \text{and} \quad \max_{y \in \mathcal{Y}^{k+1}} \|y - x^{k+1}\| > \max\{\Delta^{k+1}, 2\delta^{k+1}\}.$$

When this condition holds, COBYQA entertains the geometry-improving phase given in Algorithm 5.4.

**Algorithm 5.4:** Geometry-improving phase

**Data:** Bounds $l$ and $u$, current interpolation set $\mathcal{Y}^{k+1} \subseteq \mathbb{R}^n$, current lower bound on the trust-region radius $\delta^{k+1}$, current trust-region radius $\Delta^{k+1}$, and current best point so far $x^{k+1} \in \mathcal{Y}^{k+1}$.

**Result:** Updated interpolation set $\mathcal{Y}^{k+1}$ and updated iterate $x^{k+1}$.

Set $\bar{y}$ to be a solution to (5.2.10)
Compute the model step $r^{k+1} \in \mathbb{R}^n$
Update the interpolation set $\mathcal{Y}^{k+1} \leftarrow (\mathcal{Y}^{k+1} \setminus \{\bar{y}\}) \cup \{x^{k+1} + r^{k+1}\}$
Update the current iterate $x^{k+1}$ to a solution to $\min_{y \in \mathcal{Y}^{k+1}} \varphi^k(y)$

There is an exceptional case where we jump from Line 7 to Line 19 in Algorithm 5.1, after setting the values $\gamma^k = \gamma^{k-1}$, $\mathcal{Y}^{k+1} = \mathcal{Y}^k$, $x^{k+1} = x^k$, and $\lambda^{k+1} = \lambda^k$. Algo-rithm 5.1 omits this case for simplicity, and we detail it now. This jump is done if the trust-region trial step $d^k$ generated at Line 7 of Algorithm 5.1 is small compared to the current trust-region radius, more precisely, if $\|d^k\| < \Delta^k/2$. In this case, we first reduce $\Delta^{k+1}$ to $\Delta^{k+1}/2$, setting it further to $\delta^k$ if this reduced value is below $\theta_2\delta^k$. We then invoke Algorithm 5.4 if

$$\max_{y \in \mathcal{Y}^{k+1}} \|y - x^{k+1}\| \geq \Delta^{k+1}.$$

However, if $\|d^k\| < \Delta^k/2$ has happened at five consecutive iterations, or $\|d^k\| < \Delta^k/10$ at three consecutive iterations, then the geometry-improving phase is not invoked. Instead, we decrease $\delta^k$ by invoking Algorithm 5.3, and continue with the next trust-region iteration. This mechanism has been developed by Powell in his solver LINCOA, and similar strategies exist in UOBYQA [163], NEWUOA [166], and BOBYQA [168].

We also employ the following mechanism from NEWUOA, which aims at improving the performance of the quadratic models. We denote for convenience $\tilde{f}^k$ the least Frobenius norm model of $f$ on $\mathcal{Y}^k$. When evaluating the models at Line 6 of Algorithm 5.1, if at three consecutive iterations we have

$$\Delta^{k-1} = \delta^{k-1}, \quad \rho^{k-1} \leq \eta_5, \quad \text{and} \quad \|\nabla \hat{f}^k(x^k)\| \geq \nu \|\nabla \tilde{f}^k(x^k)\|,$$

with $\eta_5 \in (0, 1)$ and $\nu \geq 0$, then $\hat{f}^k$ is replaced with $\tilde{f}^k$, and the models of the nonlinear constraints are also replaced with their least Frobenius norm counterparts. We choose in COBYQA $\eta_5 = 0.01$ and $\nu = 10$. Since COBYQA maintains the inverse of the coefficient matrix of the interpolation problem, the least Frobenius norm models are easy to evaluate. See Section 2.4.3 for details. As will be demonstrated in Section 7.3.3, this procedure improves the performance of COBYQA.

## 5.3 Management of bound and linear constraints

COBYQA accepts three types of constraints, namely bound, linear, and nonlinear constraints. They are handled separately, for the following reasons.

First of all, bound constraints often represent inalienable physical or theoretical restrictions. In many applications for which COBYQA is designed, the objective function is not defined if the bound constraints are violated. For example, see the parameter tuning problem given in Section 1.2.2. For this reason, we implement COBYQA so that it visits only points that respect the bounds, as is also the case for the BOBYQA method, presented in Section 3.2.5.

The linear constraints are usually less restrictive. In applications, the objective function is often well-defined even at points that violate linear constraints. Therefore, we do not enforce COBYQA to always respect the linear constraints. However, when setting up a model $\hat{c}_i^k$ of a linear constraint $c_i$, we take directly $\hat{c}_i^k = c_i$ instead of constructing $\hat{c}_i^k$ by interpolation. In exact arithmetic, the derivative-free symmetric Broyden update described in Section 2.4.2 automatically ensures that $\hat{c}_i^k = c_i$ if $c_i$ is linear. However, setting $\hat{c}_i^k = c_i$ is obviously better, because it reduces the computational expense and avoid rounding errors of constructing these models.

## 5.4 Merit function and update of the penalty parameter

We now discuss the merit function we use in COBYQA. Recall that we decided to use the $\ell_2$-merit function. This is because the $\ell_2$-merit function $\widehat{\varphi}^k$ computed on the SQP subproblem, given in (5.2.1), satisfies Proposition 5.4.1. In this proposition, we denote by $\widehat{\varphi}_\gamma$ the $\ell_2$-merit function evaluated on the SQP subproblem, i.e.,

$$\widehat{\varphi}_\gamma(d) \overset{\text{def}}{=} \nabla \hat{f}^k(x^k)^\mathsf{T} d + \frac{1}{2} d^\mathsf{T} \nabla_{x,x}^2 \widehat{\mathcal{L}}^k(x^k, \lambda^k) d + \gamma \Phi(d), \quad \text{for } d \in \mathbb{R}^n,$$

where $\Phi$ is defined by (5.2.2) and where $\gamma$ is the corresponding penalty parameter. The same as $\widehat{\varphi}^k$, we do not include the bound constraints in the definition of $\widehat{\varphi}_\gamma$, because the points visited by COBYQA always satisfy the bound constraints (5.1.1d).

**Proposition 5.4.1.** *Let $d^k = n^k + t^k$ be a composite step generated by the Byrd-Omojokun approach. Assume that $n^k$ satisfies either $\Phi(n^k) < \Phi(0)$ or $n^k = 0$, and that $t^k$ satisfies*

$$\nabla \hat{f}^k(x^k)^\mathsf{T} t^k + \frac{1}{2} (t^k)^\mathsf{T} \nabla_{x,x}^2 \widehat{\mathcal{L}}^k(x^k, \lambda^k) t^k \le 0. \tag{5.4.1}$$

*Then there exists $\bar{\gamma} \ge 0$ such that for all $\gamma \ge \bar{\gamma}$, we have $\widehat{\varphi}_\gamma(d^k) \le \widehat{\varphi}_\gamma(0)$. Moreover, the least possible such value is*

$$\bar{\gamma} = \begin{cases} 0 & \text{if } \Phi(d^k) = \Phi(0), \\ \left[ \dfrac{\nabla \hat{f}^k(x^k)^\mathsf{T} d^k + (d^k)^\mathsf{T} \nabla_{x,x}^2 \widehat{\mathcal{L}}^k(x^k, \lambda^k) d^k / 2}{\Phi(0) - \Phi(d^k)} \right]_+ & \text{otherwise.} \end{cases}$$

*Proof.* Since $t^k$ satisfies (5.2.6b) and (5.2.6c), we have $\Phi(n^k + t^k) \le \Phi(n^k)$ so that it is necessary that $\Phi(d^k) \le \Phi(0)$. If $\Phi(d^k) = \Phi(0)$, then $\Phi(n^k) = \Phi(0)$ and hence, $n^k = 0$. In such a case, the result is clear, because of (5.4.1). Otherwise, if $\Phi(d^k) < \Phi(0)$, the result is clear because $\widehat{\varphi}_\gamma$ is linear with respect to $\gamma$. □

The strategy of COBYQA to set the penalty parameter at Line 8 of Algorithm 5.1 is given in Algorithm 4.3. Note that this strategy imposes $\gamma^k \ge \gamma^{k-1}$. However, in the implementation of COBYQA, we do include a procedure that attempts to reduce $\gamma^k$, preventing it from becoming large, and avoiding computational difficulties. This procedure is not included in Algorithm 5.1 for simplicity, and is triggered whenever $\delta^k$ is decreased, i.e., immediately after Algorithm 5.3 is entertained. Algorithm 5.5 presents this procedure. For simplicity, we assume in this algorithm that $\mathcal{E} = \emptyset$. If $\mathcal{E} \neq \emptyset$, COBYQA will view each equality constraint as two inequality constraints when invoking Algorithm 5.5. This procedure is the same as the one used by Powell in his solver COBYLA [161] for a similar purpose.

---

**Algorithm 5.5:** Reducing the penalty parameter‡

---

**Data:** Objective function $f$, constraint functions $\{c_i\}_{i \in \mathcal{I} \cup \mathcal{E}}$, current interpolation set $\mathcal{Y}^{k+1} \subseteq \mathbb{R}^n$, and current penalty parameter $\gamma^k \geq 0$.

**Result:** Modified penalty parameter $\gamma^k \geq 0$.

1 Set the indices of the "important" constraints

$$\mathcal{I}^* \overset{\text{def}}{=} \Big\{i \in \mathcal{I} : \min_{y \in \mathcal{Y}^{k+1}} c_i(y) < \max_{y \in \mathcal{Y}^{k+1}} 2c_i(y)\Big\}$$

2 **if** $\mathcal{I}^* = \emptyset$ **then**

3     Set $\gamma^k \leftarrow 0$

4 **else**

5     Replace $\gamma^k$ by

$$\min \left\{ \frac{\max_{y \in \mathcal{Y}^{k+1}} f(y) - \min_{y \in \mathcal{Y}^{k+1}} f(y)}{\min_{i \in \mathcal{I}^*} \{\max_{y \in \mathcal{Y}^{k+1}} c_i(y) - \min_{y \in \mathcal{Y}^{k+1}} [c_i(y)]_-\}}, \gamma^k \right\} \tag{5.4.2}$$

6 **end if**

---

‡We assume that $\mathcal{E} = \emptyset$; otherwise, each equality constraint is viewed as two inequality constraints.

The notation $[\cdot]_-$ in Algorithm 5.5 means to take the negative part of a given number. This update of $\gamma^k$ is detailed in [161, Eq. (12) and (13)] and the discussion around. The ratio in (5.4.2) can be understood as follows. Its numerator corresponds to a typical change in the objective function, while the denominator's term

$$\max_{y \in \mathcal{Y}^{k+1}} c_i(y) - \min_{y \in \mathcal{Y}^{k+1}} [c_i(y)]_-$$

represents a typical change in the $i$th constraint if $\min_{y \in \mathcal{Y}^{k+1}} c_i(y) \leq 0$, and some distance to feasibility otherwise. See [161, § 4] for more discussions about (5.4.2).

## 5.5 Summary and remarks

In this chapter, we presented the general framework of our new model-based DFO method, named COBYQA, and highlighted several important ingredients.

In a nutshell, COBYQA is a DFO trust-region SQP method. It builds quadratic models of the objective and constraint functions using the derivative-free symmetric Broy-den update, estimates its Lagrange multiplier by the least-squares Lagrange multiplier, and solves its trust-region SQP subproblem with a Byrd-Omojokun approach. It uses an $\ell_2$-merit function to assess the quality of each trial step.

COBYQA updates its trust-region radius in a way that is typical for trust-region methods, except that an adaptive lower bound is imposed on the trust-region radius. Being updated along the iterations, this lower bound is used as an indicator of the current resolution of the algorithm, and prevents the trust-region radius from becoming too small compared to the resolution. COBYQA never increases this lower bound, and decreases it if the algorithm decides that the work for the current resolution is finished. This strategy is adopted from Powell's DFO solvers.

Recall that model-based DFO methods usually entertain a geometry-improving procedure, because they tend to lose the poisedness of their interpolation set as the iterations progress. The geometry improvement of COBYQA excludes from the interpolation set the furthest point from the current iterate, and includes a point that maximizes the Lagrange polynomial associated with the excluded point. This is adopted from [166, 168]. The objective of this procedure is to improve the conditioning of the interpolation system. As a consequence, the geometry improvement also makes the interpolation problem numerically more stable.

An important feature of COBYQA is that it never attempts to visit points lying outside of the bounds. All interpolation points and trial points always respect the bounds. The motivation for us to design COBYQA in this way is that, in many applications, the objective function is undefined when the bound constraints are violated.

Finally, we note that there exists another derivate-free trust-region SQP method, named SQPDFO [199, 96, 200]. This method reformulates the inequality constraints using slack variables, so that the reformulated problem has only equality constraints and bound constraints. This increases the dimension of the problem, which may lead to inefficiency, especially if the amount of inequality constraints is large. COBYQA does not have this problem because it handles the inequality constraints directly, without introducing slack variables. Moreover, the current (August 2022) publicly available implementation of SQPDFO[1] does not always respect the bound constraints, particularly during the geometry-improving phase. In contrast, as mentioned above, COBYQA always respects the bounds, which is crucial for many engineering and industrial applications (see Section 5.3).

[1] See `https://github.com/DLR-SC/sqpdfo`.

*This page intentionally left blank*

# 6 COBYQA — solving the subproblems

This chapter introduces the methods employed by COBYQA to solve its various subproblems. We first present the classical Steihaug-Toint *truncated conjugate gradient* (TCG) method in Section 6.1. We then introduce in Sections 6.2 and 6.3 active-set variations of the TCG method to solve the tangential subproblem (5.2.6) and the normal subproblem (5.2.4), respectively. Section 6.4 adapts the NNLS algorithm [125, Alg. 23.10] to solve (5.2.7) for estimating the Lagrange multiplier in COBYQA. Finally, Section 6.5 details the procedure that COBYQA employs to solve the geometry-improving subproblem (5.2.11).

## 6.1 The truncated conjugate gradient (TCG) method

The classical trust-region subproblem for unconstrained optimization is of the form

$$\min_{s \in \mathbb{R}^n} \quad Q(s) \tag{6.1.1a}$$

$$\text{s.t.} \quad \|s\| \le \Delta, \tag{6.1.1b}$$

with $Q \in \mathbb{Q}(\mathbb{R}^n)$ being the trust-region model and $\Delta > 0$ being the trust-region radius. The usual convergence results for trust-region methods do not require us to solve (6.1.1) exactly. Rather, we only need to find an inexact solution $s^*$ that satisfies

$$Q(0) - Q(s^*) \ge c\|\nabla Q(0)\| \min\left\{\frac{\|\nabla Q(0)\|}{\|\nabla^2 Q\|}, \Delta\right\},$$

for some $c > 0$, where we assume that $\|\nabla Q(0)\|/\|\nabla^2 Q\| = \infty$ if $\nabla^2 Q = 0$. This is achieved for example by the Cauchy step, i.e., the step that minimizes $Q$ along $-\nabla Q(0)$, subject to the trust-region constraint [156, Thm. 4].

A well-known method for finding such an inexact solution to (6.1.1) is the Steihaug-Toint TCG method [196, 198], presented in Algorithm 6.1. Briefly speaking, this algorithm entertains conjugate gradients iterations, and stops the computations if a step reaches the boundary of the trust region.

In practice, $\alpha_\Delta^k$ in Line 5 of Algorithm 6.1 is obtained by solving the quadratic equation

$$\|s^k + \alpha p^k\|^2 = \Delta^2,$$

**Algorithm 6.1:** Steihaug-Toint TCG method

**Data:** Quadratic function $Q$ and trust-region radius $\Delta > 0$.

**Result:** An approximate solution to (6.1.1).

1 Set $s^0 \leftarrow 0$<br>
2 Set the search direction $p^0 \leftarrow -\nabla Q(s^0)$<br>
3 **for** $k = 0, 1, \dots$ *until* $\nabla Q(s^k)^\mathsf{T} p^k \geq 0$ **do**<br>
4     Set $\alpha_Q^k \leftarrow \operatorname{arg\,min} \{Q(s^k + \alpha p^k) : \alpha \geq 0\}$<br>
5     Set $\alpha_\Delta^k \leftarrow \operatorname{arg\,max} \{\alpha \geq 0 : \|s^k + \alpha p^k\| \leq \Delta\}$<br>
6     Set the step size $\alpha^k \leftarrow \min \{\alpha_Q^k, \alpha_\Delta^k\}$<br>
7     Update the iterate $s^{k+1} \leftarrow s^k + \alpha^k p^k$<br>
8     **if** $\alpha^k < \alpha_\Delta^k$ **then**<br>
9         Set $\beta^k \leftarrow \|\nabla Q(s^{k+1})\|^2 / \|\nabla Q(s^k)\|^2$<br>
10         Update $p^{k+1} \leftarrow -\nabla Q(s^{k+1}) + \beta^k p^k$<br>
11     **else**<br>
12         Break<br>
13     **end if**<br>
14 **end for**

which provides

$$\alpha_\Delta^k = \frac{\sqrt{[(s^k)^\mathsf{T} p^k]^2 + \|p^k\|^2 (\Delta^2 - \|s^k\|^2)} - (s^k)^\mathsf{T} p^k}{\|p^k\|^2}. \tag{6.1.2}$$

For $\alpha_Q^k$ in Line 4, we have the explicit formula

$$\alpha_Q^k = \begin{cases} \dfrac{-\nabla Q(s^k)^\mathsf{T} p^k}{(p^k)^\mathsf{T} (\nabla^2 Q) p^k}, & \text{if } (p^k)^\mathsf{T} (\nabla^2 Q) p^k > 0, & (6.1.3a) \\ \infty, & \text{otherwise}, & (6.1.3b) \end{cases}$$

which is valid due to the stopping criterion $\nabla Q(s^k)^\mathsf{T} p^k \geq 0$ in Line 3 of the algorithm. As detailed in Proposition 6.1.1, this stopping criterion is indeed equivalent to $\|p^k\| = 0$, which is usually used when presenting the TCG algorithm.

**Proposition 6.1.1.** *In Line 3 of Algorithm 6.1, the stopping criterion $\nabla Q(s^k)^\mathsf{T} p^k \geq 0$ is satisfied if and only if $\|p^k\| = 0$.*

*Proof.* It is obvious that $\|p^k\| = 0$ implies $\nabla Q(s^k)^\mathsf{T} p^k \geq 0$. For the converse, we assume without loss of generality that $k \geq 1$. When the stopping criterion in Line 3 is tested, we must have $\alpha^{k-1} = \alpha_Q^{k-1}$, and $\alpha_Q^{k-1}$ is evaluated by (6.1.3a); otherwise, $\alpha^{k-1}$ would equal $\alpha_\Delta^{k-1}$, and the algorithm would have stopped in Line 12. By (6.1.3a) and the fact that $s^k = s^{k-1} + \alpha_Q^{k-1} p^{k-1}$, we then have $\nabla Q(s^k)^{\mathsf{T}} p^{k-1} = 0$. Therefore, due to the definition of $p^k$, we see that $\nabla Q(s^k)^{\mathsf{T}} p^k \geq 0$ implies $\nabla Q(s^k) = 0$, which further leads to $\beta^{k-1} = 0$, and hence $\|p^k\| = 0$. □

The TCG algorithm enjoys many good properties. In particular, Yuan [217] showed that the TCG step provides at least half the reduction provided by the exact solution to (6.1.1) if $\nabla^2 Q$ is positive definite. The tangential and normal subproblem solvers employed by COBYQA are constrained variations of the TCG method, which will be presented in Sections 6.2 and 6.3.

## 6.2 Solving the tangential subproblem

We now present the constrained variations of the TCG method we use in COBYQA to solve the tangential subproblem (5.2.6). We introduce two variations, one for the bound-constrained case (when $\mathcal{I} \cup \mathcal{E} = \emptyset$), and the other for the linearly constrained case (when $\mathcal{I} \cup \mathcal{E} \neq \emptyset$).

### 6.2.1 Bound-constrained case

We introduce here the bound-constrained variation of the TCG method designed by Powell for his solver BOBYQA [168]. The method is described in [168, § 3], and we present it in Algorithm 6.2 for ease of reference. It is an active-set variation of the TCG method in Algorithm 6.1.

This method is employed by COBYQA to solve the tangential subproblem (5.2.6) when only bound constraints exist. In such a case, (5.2.6) is of the form

$$\min_{s \in \mathbb{R}^n} \quad Q(s) \tag{6.2.1a}$$

$$\text{s.t.} \quad l \leq s \leq u, \tag{6.2.1b}$$

$$\|s\| \leq \Delta, \tag{6.2.1c}$$

where the objective function $Q$ is quadratic, the lower bound $l \in (\mathbb{R} \cup \{-\infty\})^n$ and the upper bound $u \in (\mathbb{R} \cup \{\infty\})^n$ satisfy $l \leq 0$, $u \geq 0$, and $l < u$. Note that the lower bound $l$ in (6.2.1) indeed corresponds to the bound $l - x^k$ in (5.2.6), but we abuse the notation $l$ for simplicity. The situation for $u$ is similar.

**Description of the working set**

The method maintains a working set of active bounds. At the $k$th iteration, if a bound constraint is in the working set, then the corresponding coordinate of $s^k$ is fixed at this bound, and it will not be changed in the subsequent iterations. In the subspace of the coordinates that are not fixed, a TCG step is taken. If a new bound is hit by the TCG step, the bound is added to the working set, and the procedure is restarted. The working set is only enlarged through the iterations, and hence, the restarting happens only finitely many times.

The initial working set is a subset of the active bounds at the starting point $s^0 = 0$, which is the trust-region center. More specifically, the initial working set is

$$\mathcal{W}^0 \overset{\text{def}}{=} \left\{ i \in \{1, \dots, n\} : l_i = s_i^0 \text{ and } \frac{\partial Q}{\partial s_i}(s^0) \geq 0, \text{ or } u_i = s_i^0 \text{ and } \frac{\partial Q}{\partial s_i}(s^0) \leq 0 \right\}, \quad (6.2.2)$$

which excludes the active bounds such that a step along $-\nabla Q(s^0)$ would depart from the bounds.

### Description of the TCG method

Algorithm 6.2 details the TCG method. Given a vector $z \in \mathbb{R}^n$, we use $\Pi^k(z) \in \mathbb{R}^n$ to denote the vector whose $i$th component is $z_i$ if $i \notin \mathcal{W}^k$, and zero otherwise.

---

**Algorithm 6.2:** Bound-constrained TCG method [168, § 3]

---

**Data:** Quadratic function $Q$, bounds $l < u$, and trust-region radius $\Delta > 0$.

**Result:** An approximate solution to (6.2.1).

1 Set $s^0 \leftarrow 0$ and $\mathcal{W}^0$ according to (6.2.2)
2 Set the search direction $p^0 \leftarrow -\Pi^0(\nabla Q(s^0))$
3 **for** $k = 0, 1, \dots$ *until* $\nabla Q(s^k)^\mathsf{T} p^k \geq 0$ **do**
4     Set $\alpha_Q^k \leftarrow \arg\min \{Q(s^k + \alpha p^k) : \alpha \geq 0\}$ using (6.1.3)
5     Set $\alpha_\Delta^k \leftarrow \arg\max \{\alpha \geq 0 : \|s^k + \alpha p^k\| \leq \Delta\}$ using (6.1.2)
6     Set $\alpha_B^k \leftarrow \arg\max \{\alpha \geq 0 : l \leq s^k + \alpha p^k \leq u\}$
7     Set the step size $\alpha^k \leftarrow \min \{\alpha_Q^k, \alpha_\Delta^k, \alpha_B^k\}$
8     Update the iterate $s^{k+1} \leftarrow s^k + \alpha^k p^k$
9     **if** $\alpha^k < \min\{\alpha_\Delta^k, \alpha_B^k\}$ **then**
10         Set $\mathcal{W}^{k+1} \leftarrow \mathcal{W}^k$
11         Set $\beta^k \leftarrow \|\Pi^{k+1}(\nabla Q(s^{k+1}))\|^2 / \|\Pi^{k+1}(\nabla Q(s^k))\|^2$
12         Update $p^{k+1} \leftarrow -\Pi^{k+1}(\nabla Q(s^{k+1})) + \beta^k p^k$
13     **else if** $\alpha^k < \alpha_\Delta^k$ **then**
14         Add the first bound hit by $s^{k+1}$ to $\mathcal{W}^k$, obtaining $\mathcal{W}^{k+1}$
15         Set $p^{k+1} \leftarrow -\Pi^{k+1}(\nabla Q(s^{k+1}))$
16     **else**
17         Break
18     **end if**
19 **end for**

---

Note that $\alpha_B^k$ in Line 6 of Algorithm 6.2 is $\infty$ if all points in $\{s^k + \alpha p^k : \alpha \geq 0\}$ lie

between the bounds $l$ and $u$.

### Stopping criteria

Similar to Proposition 6.1.1, we observe the following for the stopping criterion in Line 3 of Algorithm 6.2.

> **Proposition 6.2.1.** *In Line 3 of Algorithm 6.2, the stopping criterion $\nabla Q(s^k)^\mathsf{T} p^k \geq 0$ is satisfied if and only if $\|p^k\| = 0$.*

*Proof.* It is obvious that $\|p^k\| = 0$ implies $\nabla Q(s^k)^\mathsf{T} p^k \geq 0$. Suppose that $\nabla Q(s^k)^\mathsf{T} p^k \geq 0$ in Line 3 of Algorithm 6.2, and assume without loss of generality that $k \geq 1$. If Line 3 is reached from Line 15, then

$$\|p^k\|^2 = \|\Pi^k(\nabla Q(s^k))\|^2 = Q(s^k)^\mathsf{T} \Pi^k(\nabla Q(s^k)) = -Q(s^k)^\mathsf{T} p^k \leq 0,$$

and hence $\|p^k\| = 0$. If Line 3 is reached from Line 12, then we can establish $\|p^k\| = 0$ by an argument similar to the proof of Proposition 6.1.1. □

Powell [168, § 3] proposed two additional stopping criteria for the TCG algorithm, in order to avoid unworthy computations. Algorithm 6.2 omits them for simplicity, and we specify them now. First, the algorithm is stopped if the move along $p^k$ is expected to give a relatively small reduction in $Q$, namely if

$$\big\|\Pi^k(\nabla Q(s^k))\big\|\Delta \leq \nu[Q(s^0) - Q(s^k)] \tag{6.2.3}$$

for some constant $\nu > 0$. The second additional stopping criterion is not mentioned in [168] but implemented in BOBYQA, where the algorithm is also terminated if the reduction in $Q$ provided by the recent iteration is tiny compared to the reduction so far, that is if

$$Q(s^k) - Q(s^{k+1}) \leq \nu[Q(s^0) - Q(s^{k+1})]. \tag{6.2.4}$$

COBYQA includes both of the stopping criteria with $\nu = 0.01$ when implementing Algorithm 6.2.

### Moving round the trust-region boundary

COBYQA employs an adapted version of Algorithm 6.2 for solving its tangential subproblem (5.2.6) if $\mathcal{I} \cup \mathcal{E} = \emptyset$. The adaptation, also proposed by Powell [168, § 3] and implemented in BOBYQA, is based on the following observation. Let $s^*$ be the point returned by Algorithm 6.2. If $s^*$ satisfies $\|s^*\| = \Delta$, it is likely that the objective function in (6.2.1a) can be further decreased by moving this point round the trust-region boundary. The method employed by COBYQA then attempts to further reduce the objective function by refining $s^*$ with a step that approximately solves

$$\begin{aligned}
\min_{s \in \mathbb{R}^n} \quad & Q(s^* + s) && \text{(6.2.5a)}\\
\text{s.t.} \quad & l \le s^* + s \le u, && \text{(6.2.5b)}\\
& \|s^* + s\| = \Delta, && \text{(6.2.5c)}\\
& s \in \operatorname{span}\{\Pi^*(s^*), \Pi^*(\nabla Q(s^*))\}, && \text{(6.2.5d)}
\end{aligned}$$

where $\Pi^*$ corresponds initially to the operator $\Pi^k$ at the last iteration of Algorithm 6.2. Note that a substantial reduction in $Q$ is unlikely if either $\Pi^*(\nabla Q(s^*))$ or the angle between $\Pi^*(s^*)$ and $-\Pi^*(\nabla Q(s^*))$ is tiny. Therefore, the method designed by Powell attempts to refine $s^*$ round the trust-region boundary only if

$$\|\Pi^*(s^*)\|^2 \|\Pi^*(\nabla Q(s^*))\|^2 - [\Pi^*(s^*)^{\mathsf{T}} \Pi^*(\nabla Q(s^*))]^2 > \xi [Q(s^0) - Q(s^*)]^2,$$

for some $\xi \in (0, 1)$. In COBYQA, we choose $\xi = 10^{-4}$. When this refinement is entertained, the method builds an orthogonal basis $\{\Pi^*(s^*), w\}$ of $\operatorname{span}\{\Pi^*(s^*), \Pi^*(\nabla Q(s^*))\}$ by computing the vector $w \in \mathbb{R}^n$ that satisfies

$$w^{\mathsf{T}} \Pi^*(s^*) = 0, \quad w^{\mathsf{T}} \Pi^*(\nabla Q(s^*)) < 0, \quad \text{and} \quad \|w\| = \|\Pi^*(s^*)\|.$$

Further, the method considers the function

$$s(\theta) \stackrel{\text{def}}{=} [\cos(\theta) - 1] \Pi^*(s^*) + \sin(\theta) w, \quad \text{for } 0 \le \theta \le \pi/2,$$

finds an approximate solution $\theta^*$ to the univariate problem

$$\begin{aligned}
\min_{\theta \in \mathbb{R}} \quad & Q(s^* + s(\theta)) && \text{(6.2.6a)}\\
\text{s.t.} \quad & l \le s^* + s(\theta) \le u, && \text{(6.2.6b)}\\
& 0 \le \theta \le \pi/2, && \text{(6.2.6c)}
\end{aligned}$$

and uses $s(\theta^*)$ as an approximate solution to (6.2.5). Note that the choice of $w$ ensures that $\|s^* + s(\theta^*)\| = \|s^*\| = \Delta$.

If the value of $\theta^*$ is restricted by (6.2.6b), then the first active bound is added to the working set $\mathcal{W}^*$, and the refining procedure specified above is invoked again after $s^*$ is updated to $s^* + s(\theta^*)$, and $\Pi^*$ is updated according to $\mathcal{W}^*$. Since the working set is only enlarged, this procedure is invoked at most $n - \operatorname{card}(\mathcal{W}^*)$ times, where $\mathcal{W}^*$ denotes the working set when Algorithm 6.2 terminates.

### 6.2.2 Linearly constrained case

We now introduce the linearly constrained variation of the TCG method designed by Powell [171] for his solver LINCOA. This algorithm is described in [171, § 3, § 5]. For ease of reference, we present it in Algorithm 6.3 and briefly explain the idea behind the algorithm in the sequel. The original method in [171] is devised for linear inequality constraints only, and Algorithm 6.3 adapts it to handle linear equality constraints as well.

This method is employed by COBYQA when $\mathcal{I} \cup \mathcal{E} \neq \emptyset$. In such a case, the tangential subproblem (5.2.6) is of the form

$$\min_{s \in \mathbb{R}^n} \quad Q(s) \tag{6.2.7a}$$

$$\text{s.t.} \quad As \leq b, \tag{6.2.7b}$$

$$Cs = 0, \tag{6.2.7c}$$

$$\|s\| \leq \Delta, \tag{6.2.7d}$$

where the objective function $Q$ is quadratic, $A \in \mathbb{R}^{m_1 \times n}$, $C \in \mathbb{R}^{m_2 \times n}$, and $b \in \mathbb{R}^{m_1}$ with $b \geq 0$. In this form, we include the bound constraints in (6.2.7b). Since the method is a feasible method, the bounds will be respected. In what follows, we denote the $j$th column of $A^{\mathsf{T}}$ by $a_j$, and that of $C^{\mathsf{T}}$ by $c_j$.

The method uses an active-set strategy. It maintains a working set of the linear *inequality* constraints. As in Algorithm 6.2, if the working set is modified, then the TCG procedure is restarted. However, this method allows constraints to be removed from the working set, while Algorithm 6.2 does not.

### Description of the working set

The working set is not directly the set of the active inequality constraints at the current iterate, for the following reason. Assume that $b_j$ is positive and tiny for a certain $j$, so that the $j$th constraint is almost active at $s^0 = 0$. If $j$ does not belong to the initial working set and if we have $a_j^{\mathsf{T}} \nabla Q(s^0) < 0$, then it is likely that the first step $s^1$ generated by the TCG method has a small norm, because a tiny step along the search direction $-\nabla Q(s^0)$ may reach the boundary of the feasible set. Moreover, this tiny step would lead to a change of the working set, and as observed by Powell [171, § 3], small steps that cause changes in the working set are numerically expensive, and hence, should be avoided. Therefore, we must consider constraints with small residuals when defining the working set.

The method simultaneously defines the initial search direction and the initial working set. More precisely, for any $s \in \mathbb{R}^n$ that is feasible for (6.2.7), we set

$$\mathcal{J}(s) \stackrel{\text{def}}{=} \{j \in \{1, \dots, m_1\} : b_j - a_j^{\mathsf{T}} s \leq \sigma \Delta \|a_j\|\},$$

for some $\sigma \in (0, 1)$. We take $\sigma = 0.2$ in COBYQA. The method defines the initial search direction $p^0$ as the solution to

$$\min_{s \in \mathbb{R}^n} \quad \|s + \nabla Q(s^0)\| \tag{6.2.8a}$$
$$\text{s.t.} \quad a_j^\mathsf{T} s \le 0, \; j \in \mathcal{J}(s^0), \tag{6.2.8b}$$
$$Cs = 0. \tag{6.2.8c}$$

In other words, $p^0$ is the closest direction from $-\nabla Q(s^0)$ such that any step from $s^0$ along $p^0$ moves no closer to the boundary of the constraints in $\mathcal{J}(s^0)$. The initial working set is chosen to be

$$\mathcal{W}^0 = \{j \in \mathcal{J}(s^0) : a_j^\mathsf{T} p^0 = 0\}. \tag{6.2.9}$$

It excludes the inequality constraints with $j \in \mathcal{J}(s^0)$ and $a_j^\mathsf{T} p^0 < 0$, because any step from $s^0$ along $p^0$ moves further away from these constraints.

After defining $p^0$ and $\mathcal{W}^0$, the method then considers the restriction of (6.2.7) to the orthogonal complement[1] of

$$\{a_j : j \in \mathcal{W}^0\} \cup \{c_j : 1 \le j \le m_2\}, \tag{6.2.10}$$

and entertains CG iterations in this space, using $p^0$ as the first search direction. Note that $p^0$ is the projection of $-\nabla Q(s^0)$ onto this space (see Lemma 6.2.1). When a linear inequality constraint or the trust-region constraint is hit by a CG iterate, the CG step is truncated and the CG iterations are terminated. However, upon this termination, if the current iterate $s^{k+1}$ is not on the trust-region boundary, then the procedure is restarted, changing the next search direction $p^{k+1}$ to the unique solution to

$$\min_{s \in \mathbb{R}^n} \quad \|s + \nabla Q(s^{k+1})\| \tag{6.2.11a}$$
$$\text{s.t.} \quad a_j^\mathsf{T} s \le 0, \; j \in \mathcal{J}(s^{k+1}), \tag{6.2.11b}$$
$$Cs = 0. \tag{6.2.11c}$$

Further, the working set is modified to be

$$\mathcal{W}^{k+1} = \{j \in \mathcal{J}(s^{k+1}) : a_j^\mathsf{T} p^{k+1} = 0\}. \tag{6.2.12}$$

#### Description of the TCG method

Algorithm 6.3 details the linearly constrained TCG method designed by Powell [171]. In the algorithm, we denote by $\Pi^k$ the projection onto the orthogonal complement of

$$\{a_j : j \in \mathcal{W}^k\} \cup \{c_j : 1 \le j \le m_2\}.$$

Algorithm 6.3 takes almost the same form as Algorithm 6.2. The role of $\alpha_L^k$ in Al-

[1]The orthogonal complement of (6.2.10) is the null space of the corresponding constraints.

**Algorithm 6.3:** Linearly constrained TCG method [171, § 3, § 5]

**Data:** Quadratic function $Q$, matrices $A$ and $C$, right-hand side $b \geq 0$, trust-region radius $\Delta > 0$, and constant $\sigma \in (0, 1)$.

**Result:** An approximate solution to (6.2.7).

Set $s^0 \leftarrow 0$
Set $p^0$ to the solution to (6.2.8)
Set $\mathcal{W}^0$ as described in (6.2.9)
**for** $k = 0, 1, \ldots$ *until* $\nabla Q(s^k)^\mathsf{T} p^k \geq 0$ **do**
  Set $\alpha_Q^k \leftarrow \arg\min\{Q(s^k + \alpha p^k) : \alpha \geq 0\}$ using (6.1.3)
  Set $\alpha_\Delta^k \leftarrow \arg\max\{\alpha \geq 0 : \|s^k + \alpha p^k\| \leq \Delta\}$ using (6.1.2)
  Set $\alpha_L^k \leftarrow \arg\max\{\alpha \geq 0 : A(s^k + \alpha p^k) \leq b\}$
  Set the step size $\alpha^k \leftarrow \min\{\alpha_Q^k, \alpha_\Delta^k, \alpha_L^k\}$
  Update the iterate $s^{k+1} \leftarrow s^k + \alpha^k p^k$
  **if** $\alpha^k < \min\{\alpha_\Delta^k, \alpha_L^k\}$ **then**
    Set $\mathcal{W}^{k+1} \leftarrow \mathcal{W}^k$
    Set $\beta^k \leftarrow \|\Pi^{k+1}(\nabla Q(s^{k+1}))\|^2 / \|\Pi^{k+1}(\nabla Q(s^k))\|^2$
    Update $p^{k+1} \leftarrow -\Pi^{k+1}\big(\nabla Q(s^{k+1})\big) + \beta^k p^k$
  **else if** $\alpha^k < \alpha_\Delta^k$ **then**
    Set $p^{k+1}$ to the solution to (6.2.11)
    Set $\mathcal{W}^{k+1}$ as described in (6.2.12)
  **else**
    Break
  **end if**
**end for**

gorithm 6.3 is identical to that of $\alpha_B^k$ in Algorithm 6.2. The only essential difference between the two algorithms is how they define the working set and how they define the searching direction when the working set is updated.

Line 14 of Algorithm 6.3 differs from Powell's algorithm described in [171, § 3, § 5]. Powell's algorithm replaces $\alpha^k < \alpha_\Delta^k$ in Line 14 with $\|s^k\| \le (1-\sigma)\Delta$ (see the end of [171, § 3]), which is a stronger condition because $\|s^k\| < \Delta$ implies $\alpha^k < \alpha_\Delta^k$. In our numerical experiments with COBYQA, the condition $\alpha^k < \alpha_\Delta^k$ works slightly better than $\|s^k\| \le (1-\sigma)\Delta$.

Problems (6.2.8) and (6.2.11) are solved using the Goldfarb-Idnani method for quadratic programming [87]. It is an active-set method designed for minimizing positive definite quadratic function subject to linear constraints.

### Stopping criteria

Similar to Propositions 6.1.1 and 6.2.1, we will show that the stopping criterion in Line 4 of Algorithm 6.3 is also equivalent to $\|p^k\| = 0$. To this end, we need the following simple observation taken from [171, § 3]. Its proof is *trivial*, but we include it for completeness.

**Lemma 6.2.1.** *For $g \in \mathbb{R}^n$ and $\mathcal{J} \subseteq \{1, \dots, m_1\}$, let $\tilde{p}$ be the unique solution to*

$$\min_{s \in \mathbb{R}^n} \quad \|s - g\| \tag{6.2.13a}$$

$$\text{s.t.} \quad a_j^\mathsf{T} s \le 0,\ j \in \mathcal{J}, \tag{6.2.13b}$$

$$Cs = 0. \tag{6.2.13c}$$

*Then $\tilde{p}$ is the projection of $g$ onto the orthogonal complement of*

$$\big\{a_j : j \in \mathcal{J},\ a_j^\mathsf{T} \tilde{p} = 0\big\} \cup \big\{c_j : 1 \le j \le m_2\big\}.$$

*Proof.* Define $\tilde{\mathcal{J}} = \{j \in \mathcal{J} : a_j^\mathsf{T} \tilde{p} = 0\}$. Then it suffices to prove that $\tilde{p}$ solves

$$\min_{s \in \mathbb{R}^n} \quad \|s - g\| \tag{6.2.14a}$$

$$\text{s.t.} \quad a_j^\mathsf{T} s = 0,\ j \in \tilde{\mathcal{J}}, \tag{6.2.14b}$$

$$Cs = 0, \tag{6.2.14c}$$

whose solution is unique. By the definitions of $\mathcal{J}$ and $\tilde{\mathcal{J}}$, there exists a neighborhood $\mathcal{N}$ of $\tilde{p}$ such that

$$\mathcal{N} \cap \{s \in \mathbb{R}^n : a_j^\mathsf{T} s = 0,\ j \in \tilde{\mathcal{J}},\ Cs = 0\} \subseteq \{s \in \mathbb{R}^n : a_j^\mathsf{T} s \le 0,\ j \in \mathcal{J},\ Cs = 0\}. \tag{6.2.15}$$

Since $\tilde{p}$ is a global solution to (6.2.13), we obtain from (6.2.15) that $\tilde{p}$ is a local solution to (6.2.14). This finishes the proof because (6.2.14) is convex. □

**Proposition 6.2.2.** *In Line 4 of Algorithm 6.3, the stopping criterion $\nabla Q(s^k)^\mathsf{T} p^k \ge 0$ is satisfied if and only if $\|p^k\| = 0$.*

*Proof.* It is obvious that $\|p^k\| = 0$ implies $\nabla Q(s^k)^\mathsf{T} p^k \ge 0$. Suppose that $\nabla Q(s^k)^\mathsf{T} p^k \ge 0$ in Line 4 of Algorithm 6.3, and assume without loss of generality that $k \ge 1$. If Line 4 is reached from Line 16, then we have $p^k = -\Pi^k(\nabla Q(s^k))$ according to Lemma 6.2.1. Thus,

$$\|p^k\|^2 = \|\Pi^k(\nabla Q(s^k))\|^2 = Q(s^k)^\mathsf{T} \Pi^k(\nabla Q(s^k)) = -Q(s^k)^\mathsf{T} p^k \le 0,$$

and hence $\|p^k\| = 0$. If Line 4 is reached from Line 13, then we can establish $\|p^k\| = 0$ by an argument similar to the proof of Proposition 6.1.1. □

Similar to the bound-constrained case, Powell [171, § 2] proposed two additional stopping criteria for the TCG algorithm, in order to avoid unworthy computations. Algorithm 6.3 omits them for simplicity, and we specify them now. First of all, the algorithm is stopped if the move along $p^k$ is expected to give a relatively small reduction in $Q$, namely if

$$\alpha_\Delta^k |\nabla Q(s^k)^\mathsf{T} p^k| \le \nu[Q(s^0) - Q(s^k)] \tag{6.2.16}$$

for some constant $\nu > 0$. The condition (6.2.16) can be regarded as a refined and weakened version of (6.2.3). Moreover, the algorithm also tests (6.2.4) and stop when this condition holds. As demonstrated by Powell [171, § 5], these two stopping criteria guarantee that the algorithm terminates after a finite number of iterations. COBYQA includes both of the stopping criteria with $\nu = 0.01$ when implementing Algorithm 6.3.

## 6.3 Solving the normal subproblem

In this section, we present the method employed by COBYQA to approximately solve its normal subproblem (5.2.4). Unlike for the tangential subproblem, the objective function of the normal subproblem is not quadratic, but piecewise quadratic. More specifically, problem (5.2.4) is of the form

$$\begin{aligned} \min_{z \in \mathbb{R}^n} \quad & \|[\hat{A}z - \hat{b}]_+\|^2 + \|\check{A}z - \check{b}\|^2 && \text{(6.3.1a)} \\ \text{s.t.} \quad & l \le z \le u, && \text{(6.3.1b)} \\ & \|z\| \le \Delta, && \text{(6.3.1c)} \end{aligned}$$

where $\hat{A} \in \mathbb{R}^{m_1 \times n}$, $\check{A} \in \mathbb{R}^{m_2 \times n}$, $\hat{b} \in \mathbb{R}^{m_1}$, $\check{b} \in \mathbb{R}^{m_2}$, $l \in (\mathbb{R} \cup \{-\infty\})^n$, and $u \in (\mathbb{R} \cup \{\infty\})^n$. We assume that the bounds $l$ and $u$ satisfy $l \le 0$, $u \ge 0$, and $l < u$. Note that these bounds are not the same as in (5.1.1). In (6.3.1), $[\cdot]_+$ takes the positive part of a number.

### Reformulation of the problem

The main difficulty in solving (6.3.1), even approximately, is the piecewise quadratic term $\|[\hat{A}z - \hat{b}]_+\|^2$ in its objective function. To handle this term, we introduce a new variable $v \in \mathbb{R}^{m_1}$ and reformulate (6.3.1) as

$$
\begin{aligned}
\min_{(z,v) \in \mathbb{R}^n \times \mathbb{R}^{m_1}} \quad & \|\check{A}z - \check{b}\|^2 + \|v\|^2 \\
\text{s.t.} \quad & \hat{A}z - v \le \hat{b}, \\
& v \ge 0, \\
& l \le z \le u, \\
& \|z\| \le \Delta.
\end{aligned}
$$

This problem takes the form

$$
\begin{aligned}
\min_{s \in \mathbb{R}^{n+m_1}} \quad & Q(s) && \text{(6.3.2a)} \\
\text{s.t.} \quad & As \le b, && \text{(6.3.2b)} \\
& \|s\|_D \le \Delta, && \text{(6.3.2c)}
\end{aligned}
$$

with $s$ being the vector that concatenates $z$ and $v$, and

$$
Q(s) \stackrel{\mathsf{def}}{=} -2 \begin{bmatrix} \check{b}^\mathsf{T} \check{A} & 0 \end{bmatrix} s + s^\mathsf{T} \begin{bmatrix} \check{A}^\mathsf{T} \check{A} & 0 \\ 0 & I_{m_1} \end{bmatrix} s, \quad A \stackrel{\mathsf{def}}{=} \begin{bmatrix} \hat{A} & -I_{m_1} \\ 0 & -I_{m_1} \\ I_n & 0 \\ -I_n & 0 \end{bmatrix}, \quad \text{and} \quad b \stackrel{\mathsf{def}}{=} \begin{bmatrix} \hat{b} \\ 0 \\ u \\ -l \end{bmatrix}.
$$

In (6.3.2), $\|\cdot\|_D$ denotes the seminorm defined by

$$
\|s\|_D \stackrel{\mathsf{def}}{=} \sqrt{s^\mathsf{T} D s},
$$

with $D$ being the diagonal matrix whose first $n$ diagonal elements are one, and the remaining $m_1$ components are zero.

Compared with the original problem (6.3.1), problem (6.3.2) increases the dimension of the problem, and introduce additional linear inequality constraints. However, it has the advantage that its objective function is quadratic. In addition, this reformulation enables us to apply a TCG method, as specified below.

### Description of the TCG method

We observe that problem (6.3.2) takes almost the same form as (6.2.7) with $C = 0$, the only difference being that the trust-region constraint (6.3.2c) is defined using the seminorm $\|\cdot\|_D$ instead of the $\ell_2$-norm. Based on this observation, we propose to approximately solve (6.3.2) using Algorithm 6.4, which is a straightforward adaptation of Algorithm 6.3 that replaces $\|\cdot\|$ with $\|\cdot\|_D$ in the trust region.

**Algorithm 6.4:** TCG method for approximately solving (6.3.2)

**Data:** Quadratic function $Q$, matrices $A$ and $D$, right-hand side $b$, trust-region radius $\Delta > 0$, and constant $\sigma \in (0, 1)$.

**Result:** An approximate solution to (6.3.2).

Identical to Algorithm 6.3 with $C = 0$, except that the initial guess $s^0$ is not $0$ but a feasible point, and the formula in Line 6 is replaced with

$$\alpha_\Delta^k \leftarrow \arg\max\{\alpha \ge 0 : \|s^k + \alpha p^k\|_D \le \Delta\}.$$

Note that Algorithm 6.4 requires the initial guess to be feasible; otherwise, $\alpha_\Delta^0$ or $\alpha_L^0$ may not exist (see Lines 6 and 7 of Algorithm 6.3). In practice, we take

$$s^0 = \begin{bmatrix} 0 \\ [-\hat{b}]_+ \end{bmatrix}.$$

## 6.4 Evaluating the least-squares Lagrange multiplier

We present in this section the method employed by COBYQA to solve the least-squares problem (5.2.7). It is adapted from the NNLS algorithm [125, Alg. 23.10] as follows.

Problem (5.2.7) takes the form

$$\min_{s \in \mathbb{R}^n} \quad \mu(s) \stackrel{\text{def}}{=} \frac{1}{2}\|As - b\|^2 \tag{6.4.1a}$$

$$\text{s.t.} \quad s_i \ge 0, \; i = 1, 2, \dots, n_0, \tag{6.4.1b}$$

where $A \in \mathbb{R}^{m \times n}$, $b \in \mathbb{R}^m$, and $n_0$ is a nonnegative integer with $n_0 \le n$. We observe that if $n_0 = 0$, then (6.4.1) is a simple unconstrained least-squares problem, which can be solved using traditional methods. Our major concern in this section is the case $n_0 \ge 1$.

In order to solve problem (6.4.1) when $n_0 \ge 1$, we construct an active-set method based on [125, Alg. 23.10], that maintains a working set of the nonnegativity constraints. It is a particular case of the BVLS algorithm [195]. The method is described in Algorithm 6.5. Given a working set $\mathcal{W} \subseteq \{1, 2, \dots, n_0\}$, the algorithm uses $\varrho(\mathcal{W})$ to denote the least-norm solution to

$$\min_{s \in \mathbb{R}^n} \frac{1}{2}\|A_{\mathcal{W}} s - b\|^2$$

where $A_{\mathcal{W}}$ is the matrix whose $i$th column is that of $A$ if $i \notin \mathcal{W}$, and zero otherwise.

The KKT conditions at Line 2 of Algorithm 6.5 are easy to verify in practice. A point $s^k \in \mathbb{R}^n$ satisfies the KKT conditions for (6.4.1) if and only if $\nabla\mu(s^k) \ge 0$ and

$$\frac{\partial \mu}{\partial s_i}(s^k) = 0 \quad \text{if } i > n_0 \text{ or } s_i^k > 0.$$

**Algorithm 6.5:** NonNegative Least-Squares method

**Data:** Matrix $A$, vector $b$, and integer $n_0 \geq 1$.

**Result:** A solution to (6.4.1).

Set $s^0 \leftarrow 0$, $\mathcal{W}^{-1} \leftarrow \{1, 2, \dots, n_0\}$, and $k \leftarrow 0$

**while** *the KKT conditions do not hold at* $s^k$ **do**

  Build $\mathcal{W}^k$ by removing from $\mathcal{W}^{k-1}$ the smallest $i$ that solves

$$\min_{i \in \mathcal{W}^k} \frac{\partial \mu}{\partial s_i}(s^k)$$

  Set the trial point $p^k \leftarrow \varrho(\mathcal{W}^k)$

  **while** $p_i^k \leq 0$ *for some* $i \notin \mathcal{W}^k$ *with* $i \leq n_0$ **do**

    Set the step size

$$\alpha^k \leftarrow \min \left\{ \frac{s_i^k}{s_i^k - p_i^k} : p_i^k \leq 0,\ i \notin \mathcal{W}^k,\ i \leq n_0 \right\}$$

    Update $s^{k+1} \leftarrow s^k + \alpha^k (p^k - s^k)$ and increment $k$

    Expand the working set $\mathcal{W}^k \leftarrow \mathcal{W}^{k-1} \cup \{i \leq n_0 : s_i^k = 0\}$

    Set the trial point $p^k \leftarrow \varrho(\mathcal{W}^k)$

  **end while**

  Update $s^{k+1} \leftarrow p^k$ and increment $k$

**end while**

For more discussion on this algorithm, see [125, § 23.3] and [195].

## 6.5 Solving the geometry-improving subproblem

We now detail how COBYQA approximately solves the geometry-improving subproblem (5.2.11). Such a subproblem is of the form

$$\max_{s \in \mathbb{R}^n} \quad |Q(s)| \tag{6.5.1a}$$

$$\text{s.t.} \quad l \le s \le u, \tag{6.5.1b}$$

$$\|s\| \le \Delta, \tag{6.5.1c}$$

where $Q$ a is quadratic function, the bounds $l \in (\mathbb{R} \cup \{-\infty\})^n$ and $u \in (\mathbb{R} \cup \{\infty\})^n$ satisfy $l \le 0$, $u \ge 0$, and $l < u$. Note that these bounds are not the same as in (5.1.1).

The method employed by COBYQA is adopted from BOBYQA [168]. It computes two alternative approximations, and select the one that provides the larger absolute value of the denominator of the updating formula (see Section 5.2.7).

The first alternative is a constrained Cauchy step. More specifically, the method evaluates two steps, one for the minimization of $Q(s)$, one for the minimization of $-Q(s)$, and selects the one that provides the largest value of $|Q(s)|$. For the minimization of $Q$, we first define a direction $s^c$ by solving

$$\min_{s \in \mathbb{R}^n} \quad Q(0) + \nabla Q(0)^{\mathsf{T}} s \tag{6.5.2a}$$

$$\text{s.t.} \quad l \le s \le u, \tag{6.5.2b}$$

$$\|s\| \le \Delta. \tag{6.5.2c}$$

The considered step is then obtained by minimizing $Q$ along $s^c$ subject to the bounds and the trust-region constraint. For the minimization of $-Q$, we define the direction by solving an analog of (6.5.2), with $Q$ changed to $-Q$, and then obtain the step by minimizing $-Q$ along this direction subject to the constraints.

The second alternative is as follows. Let $\mathcal{Y}$ be the current interpolation set. Without any loss of generality, assume that $0 \in \mathcal{Y}$ is the best point so far. The method minimizes $|Q(s)|$ along the straight lines through 0 and the other interpolation points. In other words, it computes

$$\begin{aligned} \max_{s \in \mathbb{R}^n} \quad & |Q(s)| \\ \text{s.t.} \quad & l \le s \le u, \\ & \|s\| \le \Delta, \\ & s \in \{\alpha y : \alpha \in \mathbb{R},\ y \in \mathcal{Y}\}. \end{aligned}$$

This subproblem can be reformulated as $\operatorname{card}(\mathcal{Y}) - 1$ maximization problems in the variable $\alpha$, which can be directly solved.

Another method that the author tried to approximately solve (6.5.1) is to employ Algorithm 6.3 on both $Q$ and $-Q$, taking the solution that provides the larger absolute value of $Q$. However, in our numerical experiments on COBYQA, this strategy performs much worse than the one explained above. This concurs with the observations of Powell [167] made on BOBYQA. Therefore, we decide not to use this strategy in the implementation of COBYQA.

## 6.6 Summary and remarks

We presented in this chapter the detailed methods employed by COBYQA to solve its various subproblems.

To solve the tangential subproblem (5.2.6) and the normal subproblem (5.2.4), we apply active-set variations of the TCG algorithm. In particular, the tangential subproblem (5.2.6) is either solved using the subproblem solver of BOBYQA [168] or a that of LINCOA [171] with a slight modification, depending on whether it admits linear constraints or not. Moreover, we adapted the subproblem solver of LINCOA to approximately solve a reformulation of the normal subproblem.

We then presented a method for solving the least-squares problem (5.2.7) that defines the least-squares Lagrange multiplier. This solver is a straightforward adaptation of the NNLS algorithm [125, Alg. 23.10] for nonnegative least-squares problems.

Finally, we introduced the method employed by COBYQA for solving the geometry-improving subproblem (5.2.11). This method is adopted from BOBYQA [168]. Two approximate solutions are calculated, the better one (according to some criterion) being returned.

# 7 COBYQA — implementation and experiments

In this chapter, we provide details on the implementation of COBYQA. In particular, Section 7.1 presents in detail the management of the user's inputs, describes every stopping criterion, and exhibits a mechanism for returning the best iterate, which is not necessarily the last one. Section 7.2 introduces the Python implementation, which is open-source and publicly available, and it also provides some examples of use. Finally, we present in Section 7.3 some numerical experiments, comparing in different scenarios COBYQA with NEWUOA, BOBYQA, LINCOA, and COBYLA. These numerical experiments confirm the good performance of COBYQA. In particular, it evidently outperforms COBYLA on all experiments, often by a significant margin.

## 7.1 Additional algorithmic details

### 7.1.1 Preprocessing of the inputs

When presenting the framework of COBYQA in Chapter 5, we assumed that the bound constraints satisfy $l < u$, and that the initial trust-region radius satisfies

$$\Delta^0 \le \frac{1}{2} \min_{1 \le i \le n} (u_i - l_i). \tag{7.1.1}$$

However, the inputs received from the user may not fulfill these conditions. In the implementation of COBYQA, we preprocess the inputs to satisfy them as long as possible.

First of all, if $l_i > u_i$ for some $i \in \{1, 2, \dots, n\}$, then no computation is attempted since the problem is infeasible. If there exist bound constraints such that $l_i$ and $u_i$ takes the same value, then we fix the corresponding component of $x$ at this value and perform the computations on the problem restricted to the remaining components of $x$. After this restriction, if (7.1.1) does not hold, then the value of $\Delta^0$ is reduced to the right-hand side of (7.1.1).

### 7.1.2 Additional stopping criteria

The only stopping criterion presented so far is that we terminate the algorithm when the lower-bound $\delta^k$ needs to be reduced to a value below a threshold $\delta^\infty$. However, in the implementation, we also stop the computations if

1. the number of function evaluations reaches a maximal value,
2. the number of iterations reaches a maximal value,

3. a target value on the objective function is reached by a feasible iterate,

The only remaining stopping criterion is when the computations must stop due to computer rounding errors. In practice, this happens if the denominator of the updating formula described in Chapter 5 becomes zero, which should not happen in exact arithmetic.

### 7.1.3 Returning the best iterate

It is important to remark that the current iterate may not be the best point encountered so far by COBYQA. This depends on the penalty parameters in the merit function that we use to assess the quality of the iterates. Indeed, a point that is removed from a previous interpolation set may be better than the current iterate according to the current merit function. Therefore, it is not always reasonable to return the last iterate as the approximate solution when COBYQA terminates.

Instead, COBYQA returns a point that is selected in the following way from all the points visited upon the termination. First of all, we consider only the points whose constraint violations are at most twice as large as the least one. Hence, if a feasible point has been produced, only the feasible points will be under consideration. Among all these points, COBYQA selects the one that minimizes the merit function value, with the penalty parameter taking the value of the last iteration. If several minimizers exist, we choose the one with the least constraint violation; if multiple choices are still possible, we select one with the least objective value.

In order to perform the selection specified above, the current version of COBYQA stores all the iterates, together with the corresponding objective function values and constraint violations. This is expensive in terms of memory. Indeed, it is not necessary to save a point if both its objective function value and its constraint violation are larger than those of another point visited by COBYQA. We will improve the implementation of COBYQA in this way in future versions.

## 7.2 Description of the Python implementation

In this section, we provide details on the Python implementation of COBYQA and the way we distribute it.

### 7.2.1 Choice of programming languages

We chose to develop COBYQA in Python because it is a simple open-source language, widely used by both researchers and industrial practitioners. However, the subproblem solvers described in Chapter 6 are not written in pure Python, but in Cython [14], which is a compiled language that aims at improving the performance of Python code.

The motivation is that these solvers are the bottleneck for the execution time of the linear algebra calculations made by COBYQA. Indeed, the initial version of COBYQA was entirely written in Python, but the execution time was prohibitively long on problems of dimensions $n = 20$ and above.

Another motivation for us to develop the COBYQA in Python is the extensive availability of libraries. COBYQA does not only rely on Cython, but also on NumPy [103] and SciPy [204]. These packages provide array structures, together with basic mathematical operations and basic scientific computing tools. In particular, SciPy provides interfaces for BLAS [19] and LAPACK [4] subroutines in Cython, which are intensively used by the subproblem solvers.

### 7.2.2 Distribution of the package

The Python implementation of COBYQA is open-source and distributed under the BSD-3-Clause license. The source code of COBYQA is available at

`https://github.com/ragonneau/cobyqa`.

A complete documentation for COBYQA is available at

`https://www.cobyqa.com/`.

A C/C++ compiler is needed to build the subproblem solvers of COBYQA written in Cython. To ease the efforts of users, we made wheel distributions available for

1. Windows, both 32-bit and 64-bit operating systems,
2. macOS, both x64 and ARM64 architectures, and
3. Linux, both x64 and ARM64 architectures.

The current version of COBYQA is available for Python 3.7, 3.8, 3.9, and 3.10, i.e., all the Python version supported as of August 2022.

### 7.2.3 Illustrations of usage

The package COBYQA provides a function `minimize`, which takes the following arguments.

1. The objective function to be minimized,
2. the initial guess,
3. bound constraints (optional),
4. linear constraints (optional),
5. nonlinear constraints (optional), and
6. a dictionary of options (listed in Table 7.1, each being optional).

Table 7.1: Options of COBYQA

| | |
|---|---|
| `rhobeg` | Initial trust-region radius |
| `rhoend` | Final trust-region radius |
| `npt` | Number of interpolation points |
| `maxfev` | Maximum number of objective and constraint function evaluations |
| `maxiter` | Maximum number of iterations |
| `target` | Target value on the objective function |
| `disp` | Whether to print pieces of information on the execution of the solver |
| `debug` | Whether to make debugging tests during the execution |

In the following, we will provide example to illustrate how to specify the aforementioned arguments.

#### Examples of usage

For example, consider the unconstrained minimization of the chained Rosenbrock function

$$f(x) = \sum_{i=1}^{n-1} 100(x_{i+1} - x_i^2)^2 + (1 - x_i)^2, \quad \text{for } x \in \mathbb{R}^n.$$

It only serves as an illustration. In practice, it is clearly unreasonable to employ a DFO method to solve such a problem. Listing 7.1 shows how to solve such a problem using COBYQA with $n = 5$, starting from the initial guess $x^0 = [1.3, 0.7, 0.8, 1.9, 1.2]^\mathsf{T}$.

Listing 7.1: Solving the Rosenbrock problem using COBYQA

```
>>> from scipy.optimize import rosen
>>> from cobyqa import minimize
>>>
>>> x0 = [1.3, 0.7, 0.8, 1.9, 1.2]  # starting point
>>> res = minimize(rosen, x0)
>>> res.x
array([1., 1., 1., 1., 1.])
```

To see how to supply to COBYQA bound and linear constraints, consider the Exam-ple 16.4 of [144], defined as

$$\begin{aligned}
\min_{x \in \mathbb{R}^2} \quad & (x_1 - 1)^2 + (x_2 - 2.5)^2 \\
\text{s.t.} \quad & -x_1 + 2x_2 \le 2, \\
& x_1 + 2x_2 \le 6, \\
& x_1 - 2x_2 \le 2, \\
& x_1 \ge 0, \\
& x_2 \ge 0.
\end{aligned}$$

Listing 7.2 shows how to solve this problem using COBYQA starting from the initial guess $x^0 = [2, 0]^\mathsf{T}$.

Listing 7.2: An example of COBYQA with linear constraints

```
>>> from cobyqa import minimize
>>>
>>> def obj(x):
...     return (x[0] - 1.0) ** 2.0 + (x[1] - 2.5) ** 2.0
>>>
>>> x0 = [2.0, 0.0]
>>> xl = [0.0, 0.0]
>>> Aub = [[-1.0, 2.0], [1.0, 2.0], [1.0, -2.0]]
>>> bub = [2.0, 6.0, 2.0]
>>> res = minimize(quadratic, x0, xl=xl, Aub=Aub, bub=bub)
>>> res.x
array([1.4, 1.7])
```

Finally, to see how to supply nonlinear constraints to COBYQA, consider the Problem G of [161], defined by

$$\begin{aligned}
\min_{x \in \mathbb{R}^3} \quad & x_3 \\
\text{s.t.} \quad & -5x_1 + x_2 - x_3 \le 0, \\
& 5x_1 + x_2 - x_3 \le 0, \\
& x_1^2 + x_2^2 + 4x_2 - x_3 \le 0.
\end{aligned}$$

Listing 7.3 shows how to solve such a problem using COBYQA, starting from the initial guess $x^0 = [1, 1, 1]^\mathsf{T}$. Note that we can supply the linear constraints are nonlinear constraints. However, in doing so, the code cannot detect that the constraint is linear, and hence, will build quadratic models of these constraints.

Listing 7.3: An example of COBYQA with nonlinear constraints

```
>>> from cobyqa import minimize
>>>
>>> def obj(x):
...     return x[2]
>>>
>>> def cub(x):
...     return x[0] ** 2.0 + x[1] ** 2.0 + 4.0 * x[1] - x[2]
>>>
>>> x0 = [1.0, 1.0, 1.0]
>>> Aub = [[-5.0, 1.0, -1.0], [5.0, 1.0, -1.0]]
>>> bub = [0.0, 0.0]
>>> res = minimize(obj, x0, Aub=Aub, bub=bub, cub=cub)
>>> res.x
array([ 0., -3., -3.])
```

**Details of the returned structure**

The function `minimize` returns a structure (that we named `res` in the examples above) whose attributes are detailed in Table 7.2.

Table 7.2: Structure for the result of the optimization algorithm

| | |
|---|---|
| `x` | Solution point |
| `status` | Termination status of the optimization solver reflecting which stopping criterion has been reached |
| `message` | Description of the termination status |
| `fun` | Value of the objective function at the solution point |
| `jac` | Approximate gradient of the objective function at the solution point |
| `nfev` | Number of objective and constraint function evaluations |
| `nit` | Number of iterations |
| `maxcv` | $\ell_\infty$-constraint violation at the solution point, if applicable |

## 7.3 Numerical experiments

To conclude this chapter, we present extensive numerical experiments for COBYQA. We will make comparisons using performance and data profiles obtained on CUTEst problems [90], with at most $50$ variables and at most $1{,}000$ constraints, excluding bound constraints. These problems are listed in Appendices A.1 to A.4, according to the types of their constraints.

As introduced in Section 1.5, the performance and data profiles require us to define a merit function. Let $v_\infty(x) \in \mathbb{R}$ be the $\ell_\infty$-constraint violation at $x \in \mathbb{R}^n$. The merit function $\varphi$ is defined by

$$\varphi(x) = \begin{cases} f(x), & \text{if } v_\infty(x) \le \nu_1, \\ \infty, & \text{if } v_\infty(x) \ge \nu_2, \\ f(x) + \gamma v_\infty(x), & \text{otherwise}, \end{cases}$$

for $x \in \mathbb{R}^n$. We choose in the experiments $\nu_1 = 10^{-10}$, $\nu_2 = 10^{-5}$, and $\gamma = 10^5$. Note that $v_\infty$ is always zero for unconstrained problems.

For all the experiments, we present the performance and data profiles corresponding to $\tau \in \{10^{-1}, 10^{-3}, 10^{-5}, 10^{-7}\}$ in the convergence test (1.5.1). The starting points of the problems are set to the default ones provided in CUTEst. The initial trust-region radius $\Delta^0$ is set to one, and the final value $\delta^\infty$ for the lower bound on the trust-region radius is set to $10^{-6}$. The maximal number of function evaluations is $500n$, where $n$ denotes the dimension of the problem being solved. For the methods that employ underdetermined quadratic interpolation models, the number of interpolation points is set to $2n + 1$.

### 7.3.1 Testing the Hessian of the SQP subproblem

Recall that the objective function of the SQP subproblem takes a Hessian matrix that approximates $\nabla^2_{x,x} \mathcal{L}(x^k, \lambda^k)$. In COBYQA, this Hessian is set to $\nabla^2_{x,x} \widehat{\mathcal{L}}^k(x^k, \lambda^k)$, defined in Section 5.2.2. It is erroneous to set this Hessian to $\nabla^2 \hat{f}^k(x^k)$, unless there is no nonlinear constraints (see Section 4.1). To demonstrate this fact, we compare COBYQA with a modification that sets this Hessian to $\nabla^2 \hat{f}^k(x^k)$.

Figures 7.1 and 7.2 provide the performance and data profiles for these two methods, based on the nonlinearly constrained problems listed in Appendix A.4. According to these profiles, COBYQA performs much better than the modified method, which is expected. It solves more problems, and uses less function evaluations than the modified method on most of the problems. These results demonstrate the effectiveness of the SQP methodology of COBYQA. It also shows that the estimated Lagrange multipliers we chose in Section 5.2.4 are reasonable.

### 7.3.2 Testing our Byrd-Omojokun approach for inequality constraints

To handle inequality constraints, COBYQA takes an extension of the Byrd-Omojokun that we devised in Section 4.3.4. This extension is different from the one presented in [44, § 15.4.4]. To demonstrate the advantage of our extension, we compare COBYQA with a modified version that employs the Byrd-Omojokun approach in [44, § 15.4.4].

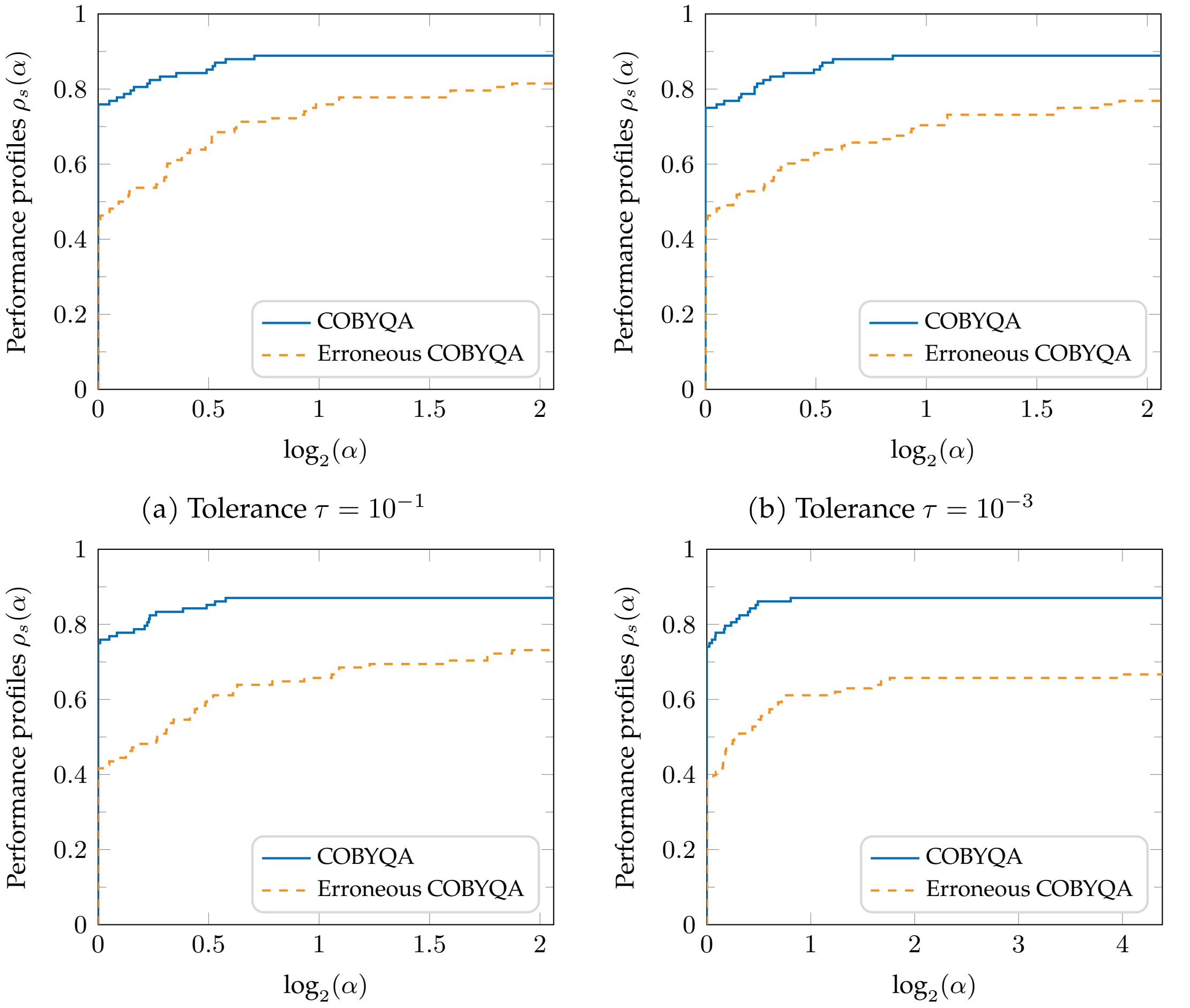

(a) Tolerance $\tau = 10^{-1}$

(b) Tolerance $\tau = 10^{-3}$

(c) Tolerance $\tau = 10^{-5}$

(d) Tolerance $\tau = 10^{-7}$

Figure 7.1: Performance profiles of COBYQA and the modified version with erroneous Hessian on nonlinearly constrained problems

In this comparison, we consider only linearly and nonlinearly constrained problems with at least one inequality constraint. The problems are taken from Appendices A.3 and A.4 accordingly.

Figures 7.3 and 7.4 presents the corresponding performance and data profiles of COBYQA and the modified version specified above. According to these profiles, our extension to the Byrd-Omojokun approach evidently outperforms the one presented in [44, § 15.4.4]. As mentioned in Section 4.3.4, this can be explained by the fact that our Byrd-Omojokun approach provides a larger feasible region to the tangential steps. Hence, a larger reduction in the objective function of the SQP subproblem can be expected.

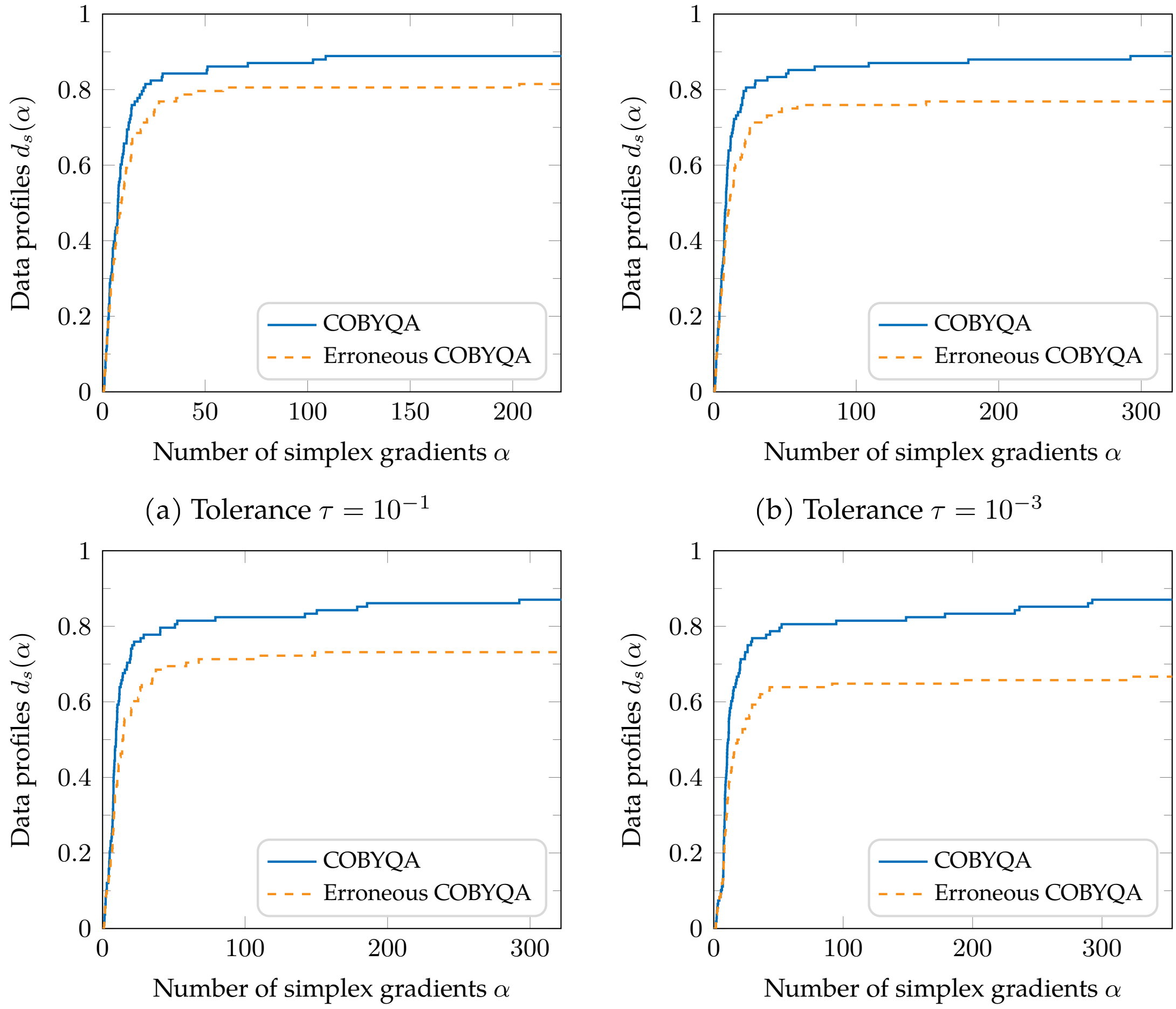

(a) Tolerance $\tau = 10^{-1}$

(b) Tolerance $\tau = 10^{-3}$

(c) Tolerance $\tau = 10^{-5}$

(d) Tolerance $\tau = 10^{-7}$

Figure 7.2: Data profiles of COBYQA and the modified version with erroneous Hessian on nonlinearly constrained problems

### 7.3.3 Performance of the derivative-free symmetric Broyden update

Recall that COBYQA employs quadratic models based on the derivative-free symmetric Broyden update (see Section 2.4.2). Moreover, it includes a strategy to replace these models with the minimum Frobenius norm models when the method considers that the current models are inefficient, as detailed at the end of Section 5.2.7.

To justify our choice for the models, we now compare the performance of COBYQA with two modifications. The first one, referred to as COBYQA PSB, always uses the models based on the derivative-free symmetric Broyden update. The second one, referred to as COBYQA MNH, always uses the minimum Frobenius norm models.

Figures 7.5 and 7.6 provide the performance and data profiles for these three methods, based on all the problems listed in Appendices A.1 to A.4. These profiles show that COBYQA outperforms COBYQA PSB, which in turn outperforms COBYQA MNH. These results confirm that our choice for the quadratic models is reasonable.

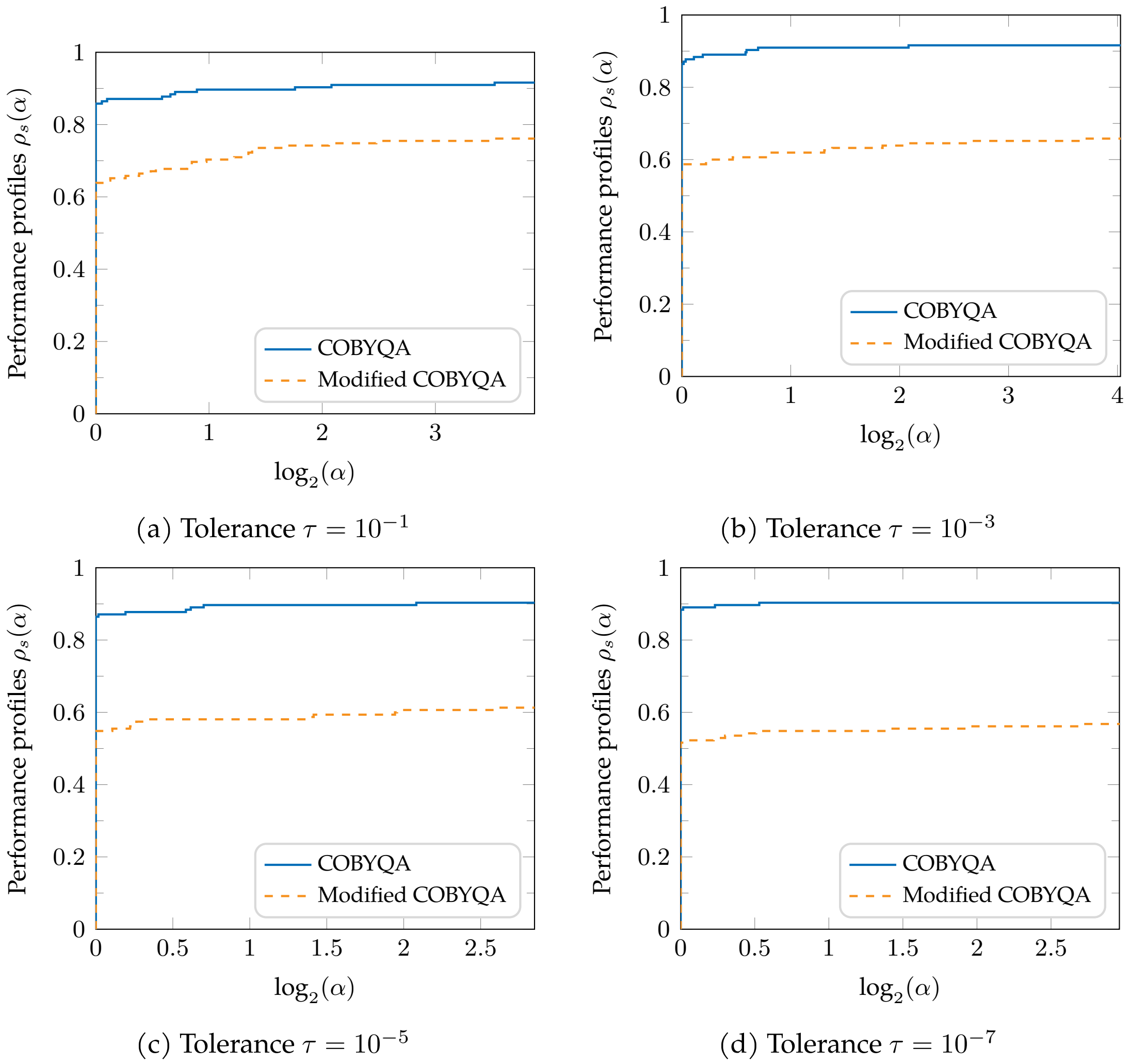

(a) Tolerance $\tau = 10^{-1}$

(b) Tolerance $\tau = 10^{-3}$

(c) Tolerance $\tau = 10^{-5}$

(d) Tolerance $\tau = 10^{-7}$

Figure 7.3: Performance profiles of COBYQA and a modified version on and nonlinearly constrained problems with at least one inequality constraint

### 7.3.4 Performance of COBYQA on unconstrained problems

In the remaining of this chapter, we compare the performance of COBYQA with the Powell's DFO methods that are available through PDFO. In all the comparisons, we include both COBYQA and COBYLA, because both are capable of handling all types of problems, and because COBYQA is expected to be a successor to COBYLA.

We start with the comparison of COBYQA, NEWUOA, and COBYLA on unconstrained problems. Figures 7.7 and 7.8 present the performance and data profiles of these solvers, based on all the problems listed in Appendix A.1. These are all the unconstrained problems of the CUTEst library with at most 50 variables.

We first observe that COBYQA outperforms COBYLA in general, the margin being large for the tolerance $\tau$ is low. According to the performance profiles in Figure 7.7, with $\tau = 10^{-1}$, $10^{-3}$, $10^{-5}$, and $10^{-7}$, COBYQA solves about 98 %, 90 %, 85 %, and 79 % of the problems, respectively, whereas the percentages for COBYLA are about 89 %, 67 %, 50 %, and 39 %. Moreover, for $\tau \in \{10^{-3}, 10^{-5}, 10^{-7}\}$, COBYQA takes the least number of function evaluations to converge on about 35 % of the problems, while the percentages for COBYLA are between 15 % and 30 %.

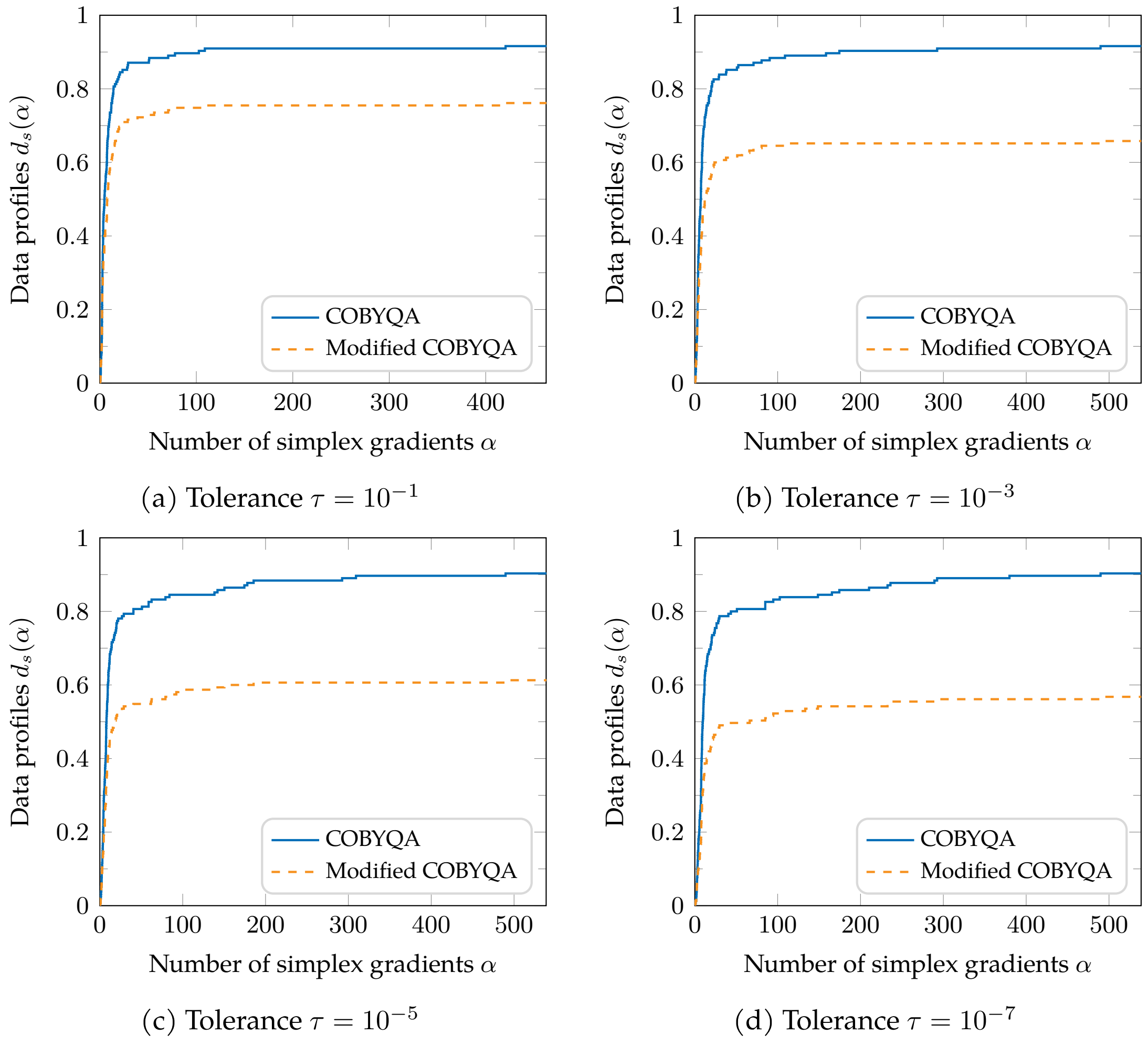

Figure 7.4: Data profiles of COBYQA and a modified version on linearly and nonlinearly constrained problems with at least one inequality constraint

However, for $\tau = 10^{-1}$, COBYLA uses the smallest number of function evaluations to converge on about 45 % of problems, while the percentage for COBYQA is lower, being about 25 %. This is expected, because COBYLA necessitates $n+1$ function evaluations to build the initial model for an $n$-dimensional problem, while COBYQA requires $2n + 1$. Consequently, COBYLA takes less function evaluations to make progress at the early stage and hence, performs better for a low-precision comparison.

We also observe that COBYQA and NEWUOA perform comparably in this experiment, which is expected because they share very similar algorithms. NEWUOA slightly outperforms COBYQA, but the difference is marginal. This is likely due to minor dif-ferences in the algorithms, in their implementations, and computer rounding errors. In future releases of COBYQA, we will try to close the tiny gap between NEWUOA and COBYQA by refining the implementation.

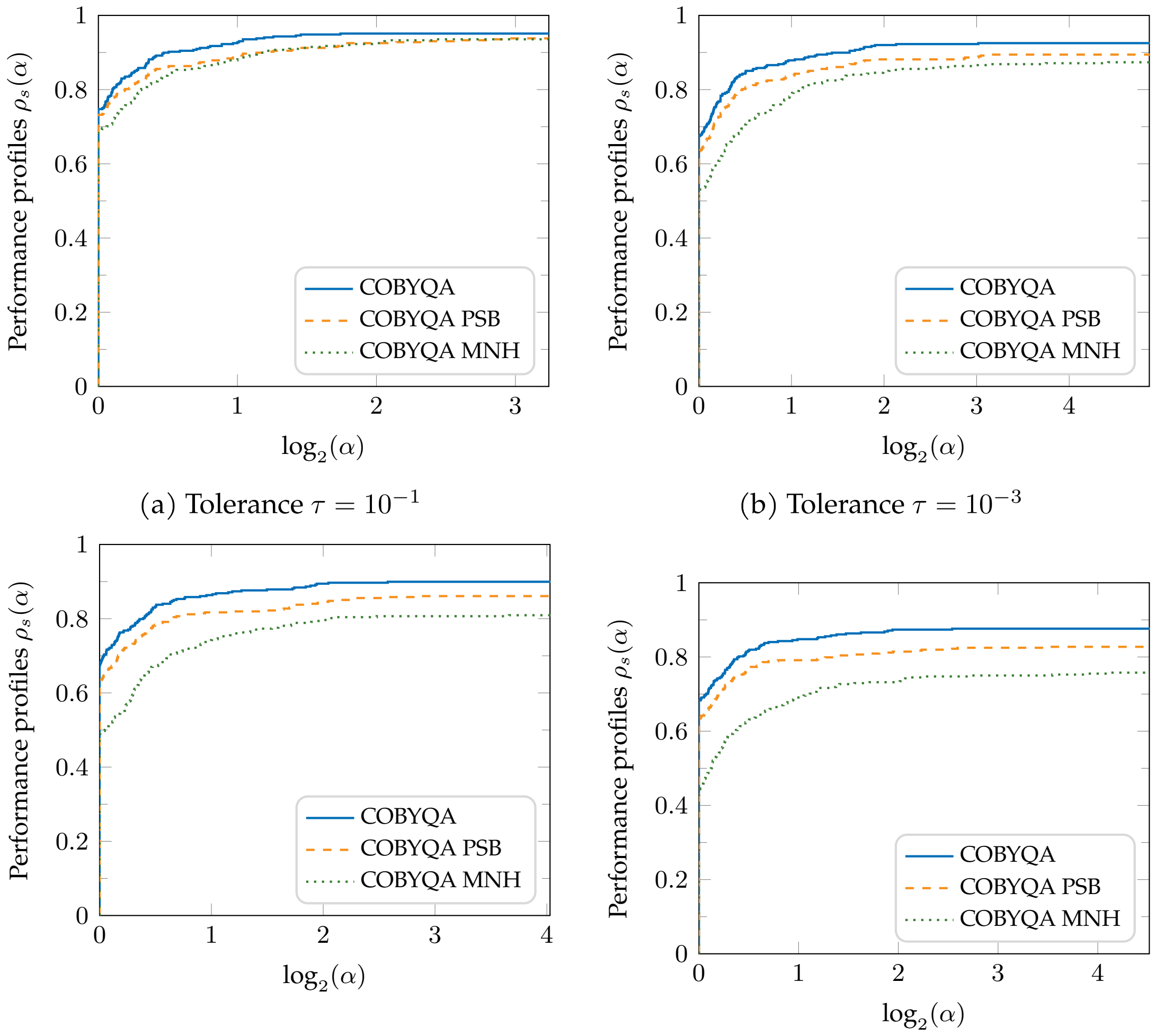

(a) Tolerance $\tau = 10^{-1}$

(b) Tolerance $\tau = 10^{-3}$

(c) Tolerance $\tau = 10^{-5}$

(d) Tolerance $\tau = 10^{-7}$

Figure 7.5: Performance profiles of COBYQA, COBYQA PSB, and COBYQA MNH on all problems

### 7.3.5 Performance of COBYQA on bound-constrained problems

We now compare the performance of COBYQA, BOBYQA, and COBYLA on bound-constrained problems. Figures 7.9 and 7.10 present the performance and data profiles of these solvers, based on all the problems listed in Appendix A.2. These are all the bound-constrained problems of the CUTEst library with at most 50 variables.

The comparison of COBYQA with COBYLA is quite similar to the one in the unconstrained case. More specifically, for $\tau \in \{10^{-3}, 10^{-5}, 10^{-7}\}$, COBYQA outperforms COBYLA by a large margin. For the tolerance $\tau = 10^{-1}$, COBYQA solves more prob-lems than COBYLA, the percentages being about 98 % and 91 %. However, for this low-precision comparison, according to the performance profiles in Figure 7.9, COBYLA solves more problems than COBYQA using the least number of function evaluations, the percentages being about 53 % and 35 %. This is expected for the same reason as the unconstrained case.

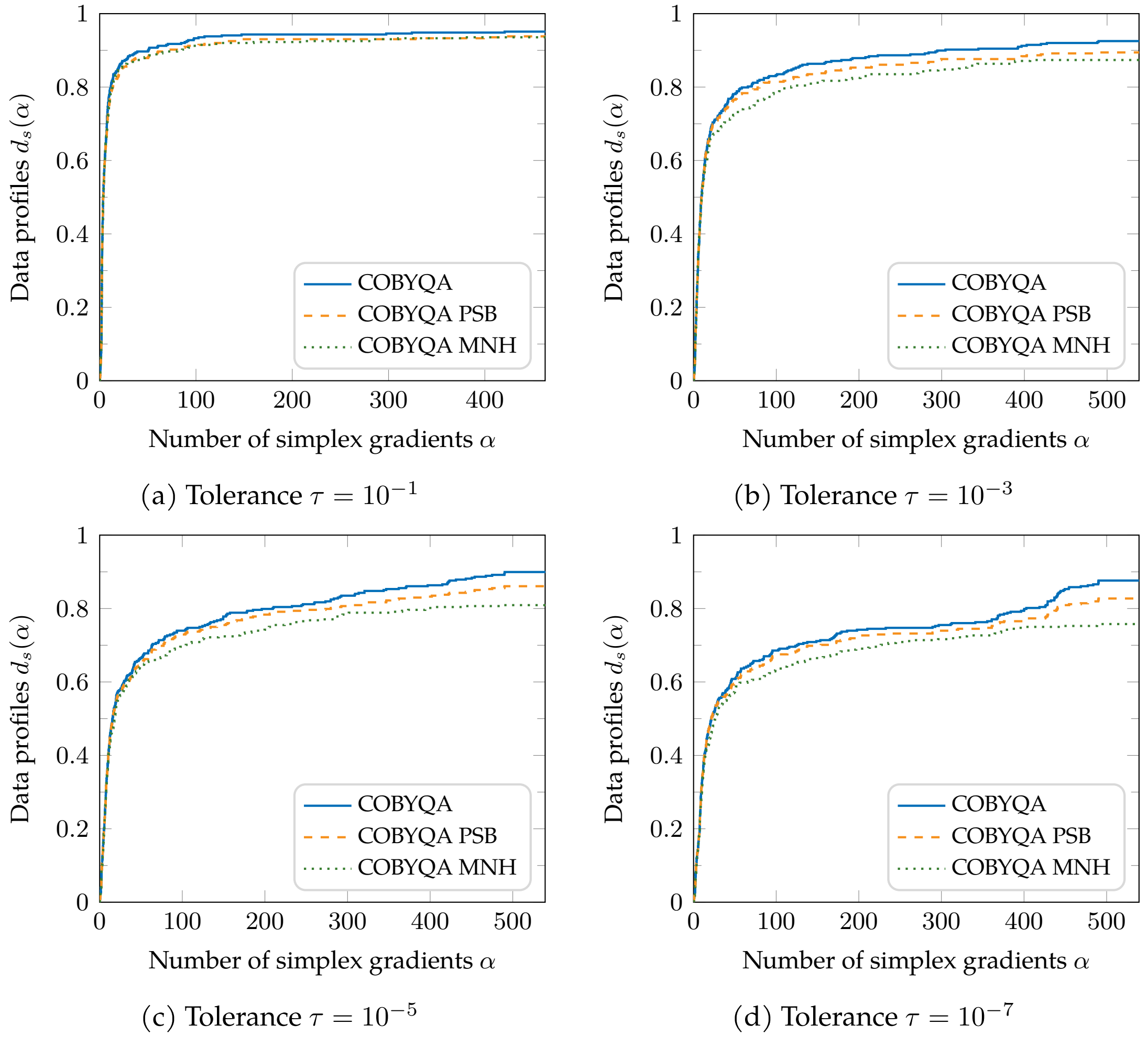

(a) Tolerance $\tau = 10^{-1}$

(b) Tolerance $\tau = 10^{-3}$

(c) Tolerance $\tau = 10^{-5}$

(d) Tolerance $\tau = 10^{-7}$

Figure 7.6: Data profiles of COBYQA, COBYQA PSB, and COBYQA MNH on all problems

Moreover, we observe that COBYQA outperforms BOBYQA in general. For the high-precision test with $\tau = 10^{-7}$, the advantage of COBYQA over BOBYQA is considerable. When $\tau = 10^{-3}$ and $10^{-5}$, COBYQA slightly outperforms BOBYQA. Finally, for the tolerance $\tau = 10^{-1}$, the two methods perform similarly. It is interesting that COBYQA outperforms BOBYQA on bound-constrained problems, because BOBYQA is particularly designed for this type of problems, while COBYQA is a general-purpose solver.

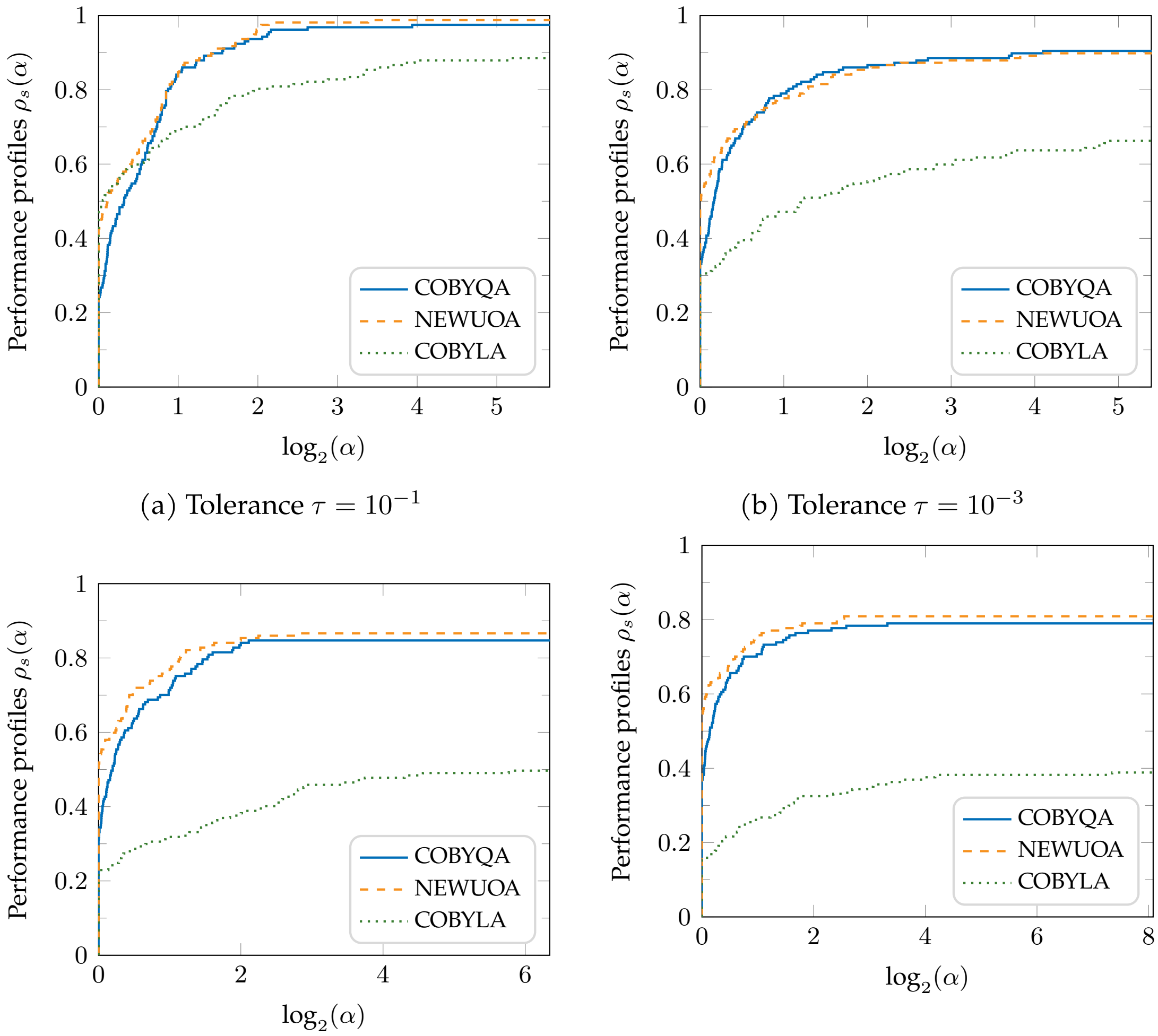

(a) Tolerance $\tau = 10^{-1}$

(b) Tolerance $\tau = 10^{-3}$

(c) Tolerance $\tau = 10^{-5}$

(d) Tolerance $\tau = 10^{-7}$

Figure 7.7: Performance profiles of COBYQA, NEWUOA, and COBYLA on unconstrained problems

### 7.3.6 Performance of COBYQA on linearly constrained problems

We now compare the performance of COBYQA, LINCOA, and COBYLA on linearly constrained problems. Recall that one of the critical features of COBYQA is that is always respect bound constraints, if any. This is important in many applications, because the objective function may be undefined when bounds are violated (see Section 5.3). However, LINCOA and COBYLA do not always respect bound constraints. To observe different behaviors of these methods with respect to bounds we consider two scenarios in the experiments. In the first one, we allow the solvers to violate bound constraints, if any; the objective functions of the problems are well-defined even if bound constraints are violated. In the second one, we assume that the bounds are inviolable, and we set the objective function values to $\infty$ if any bound constraint is violated. As will be shown, COBYQA outperforms both LINCOA and COBYLA in the second scenario with large

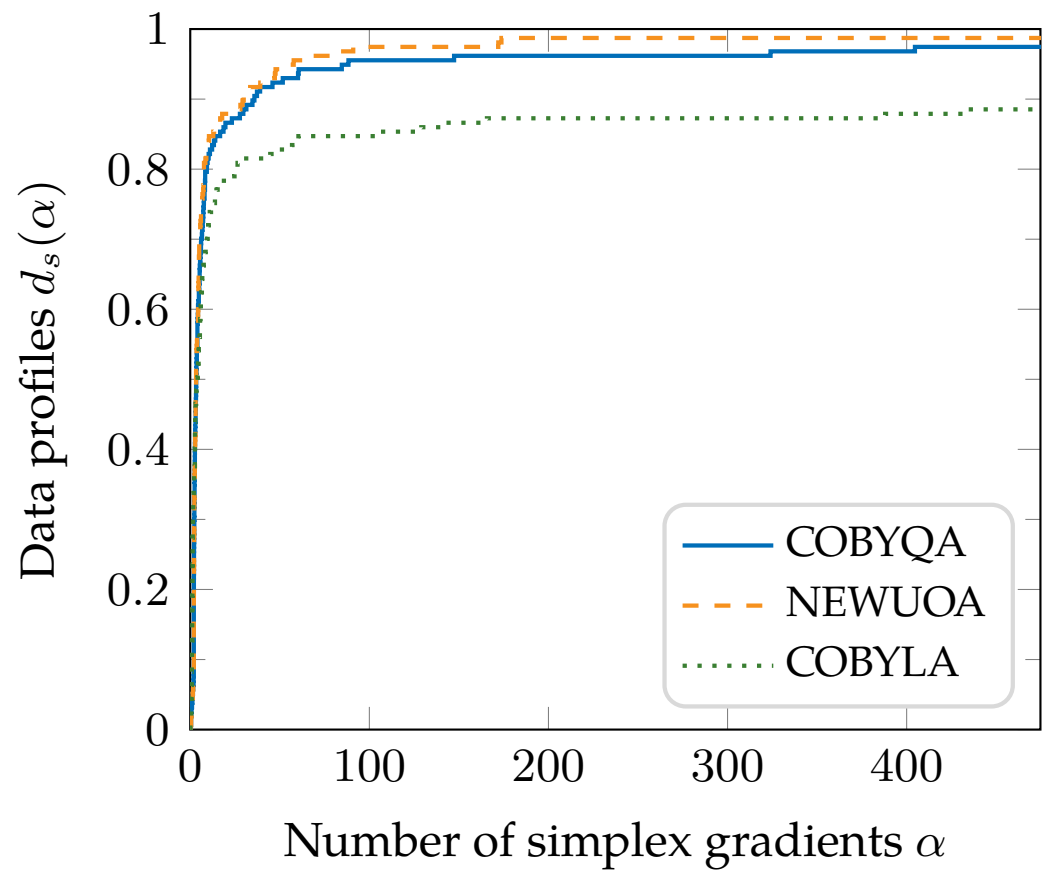

(a) Tolerance $\tau = 10^{-1}$

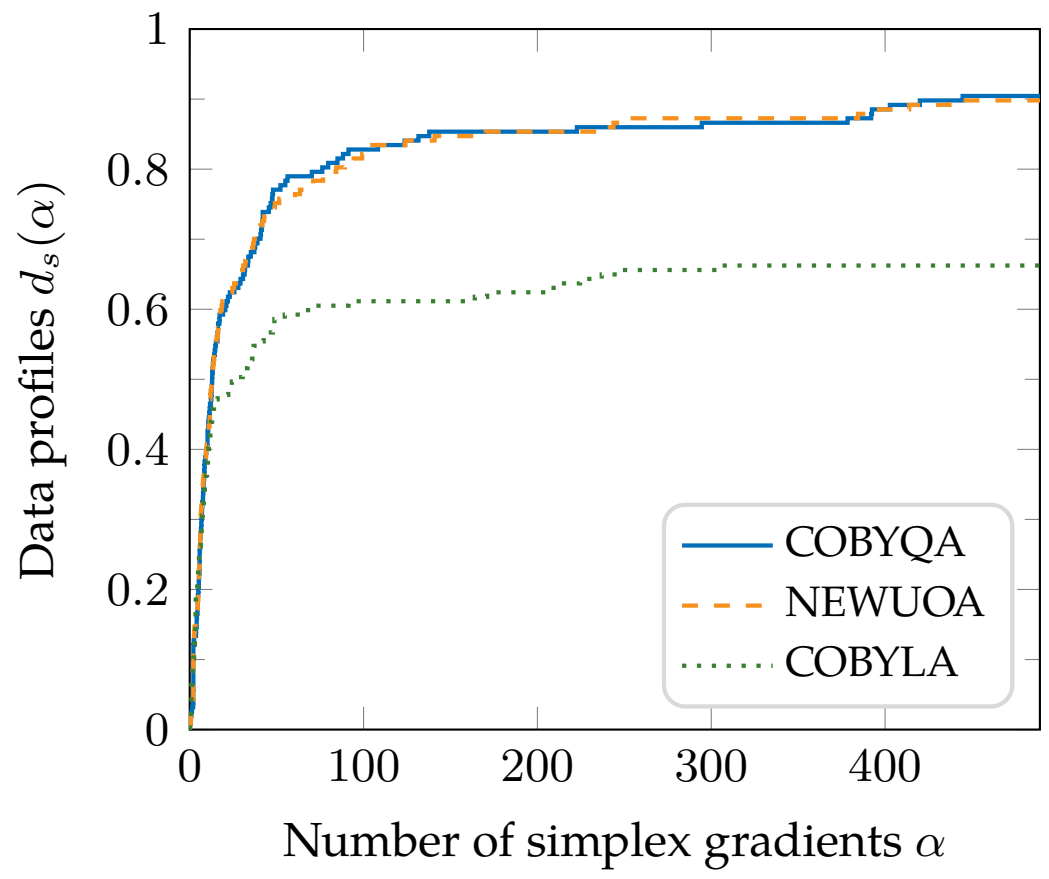

(b) Tolerance $\tau = 10^{-3}$

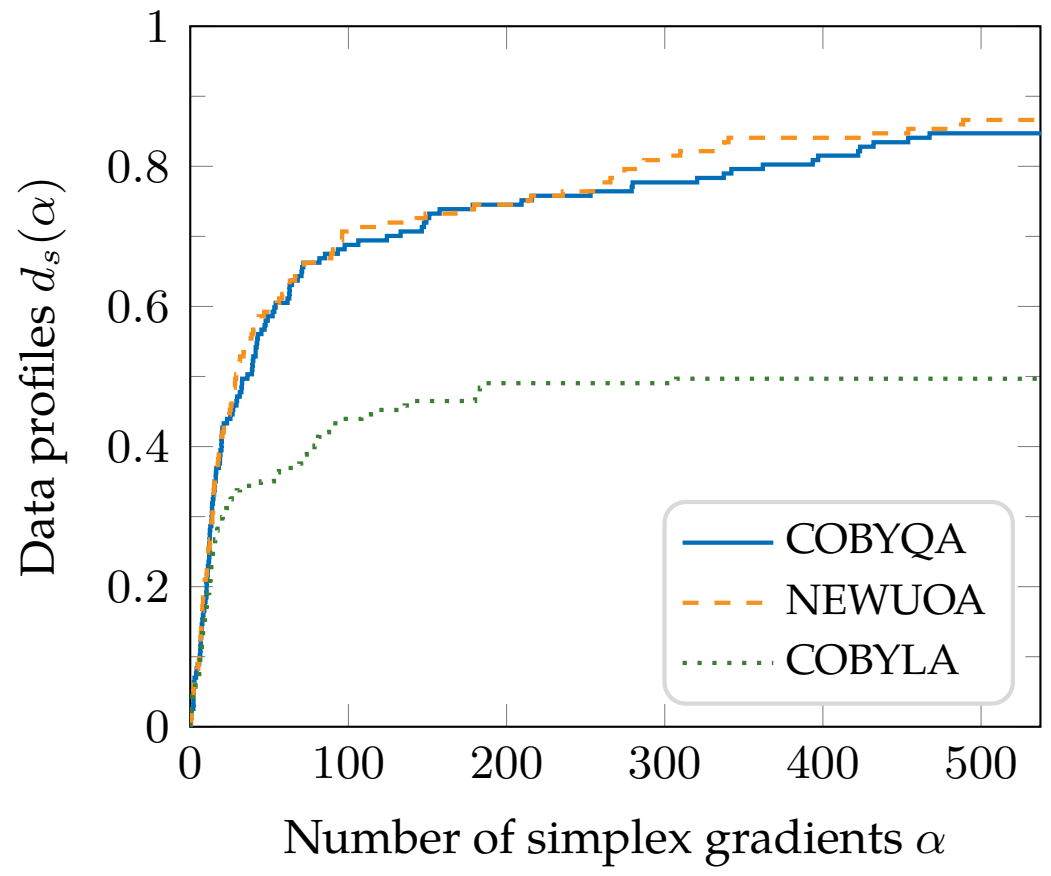

(c) Tolerance $\tau = 10^{-5}$

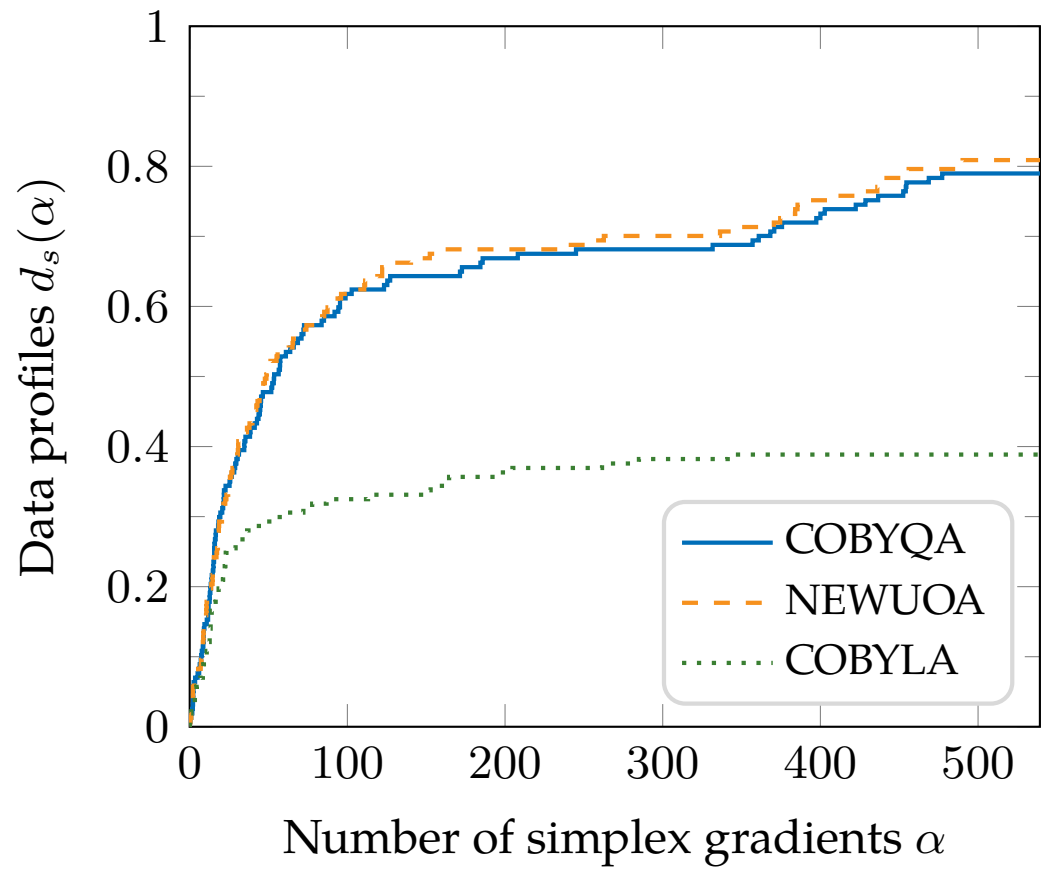

(d) Tolerance $\tau = 10^{-7}$

Figure 7.8: Data profiles of COBYQA, NEWUOA, and COBYLA on unconstrained problems

margins.

#### Comparison when the bounds are violable

Figures 7.11 and 7.12 present the performance and data profiles of these solvers, based on all the problems listed in Appendix A.3. These are all the linearly constrained problems of the CUTEst library with at most 50 variables and at most 1,000 linear constraints.

COBYQA outperforms COBYLA when $\tau = 10^{-5}$ and $10^{-7}$ with a fair margin. They perform similarly when $\tau = 10^{-3}$. However, when $\tau = 10^{-1}$, COBYLA solves more problems than both COBYQA and LINCOA using the least number of function evaluations. The advantage of COBYLA on this low-precision test is expected, as explained in Section 7.3.4.

However, COBYQA is always outperformed by LINCOA for all values of $\tau$. Note that LINCOA is specifically designed for linearly constrained problems, while COBYQA is a general-purpose solver. On of the possible reasons for the advantage of LINCOA over COBYQA is that it uses a more careful strategy to decide when to replace the default models with the minimum Frobenius norm ones. In further releases of COBYQA, we will investigate the differences between COBYQA and LINCOA, and try to improve the performance of COBYQA on linearly constrained problems.

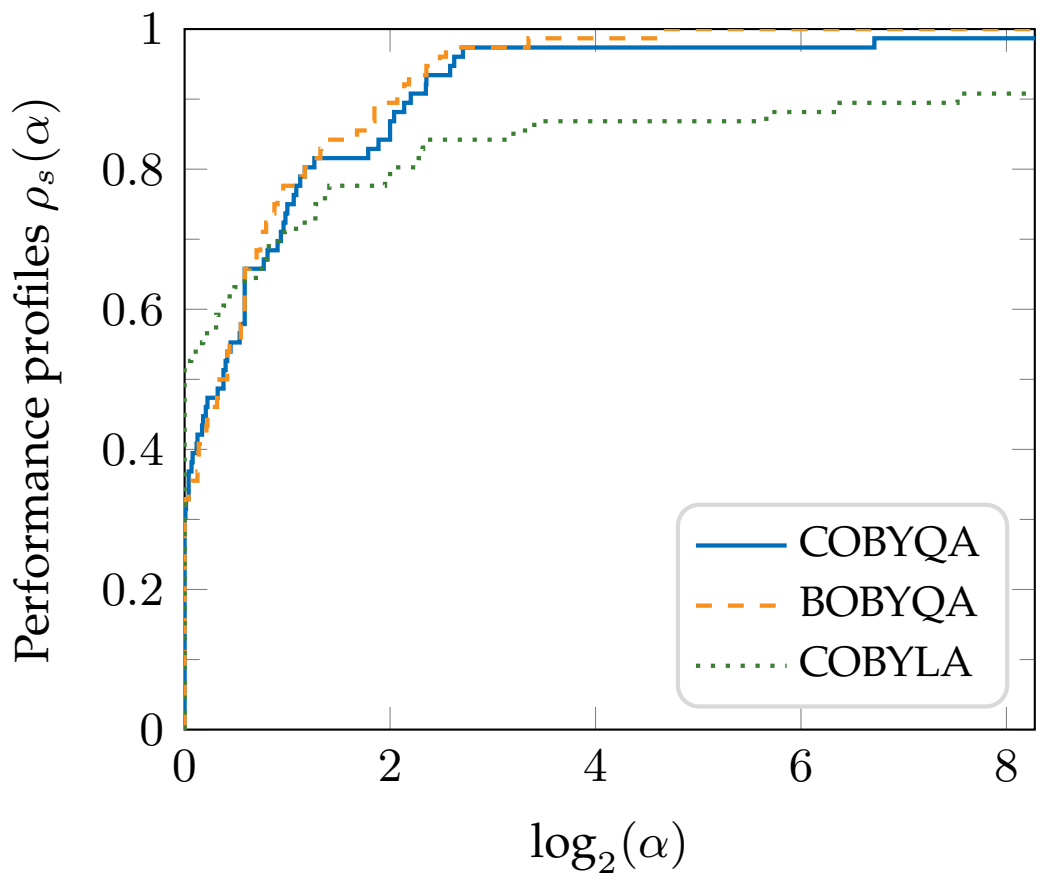

(a) Tolerance $\tau = 10^{-1}$

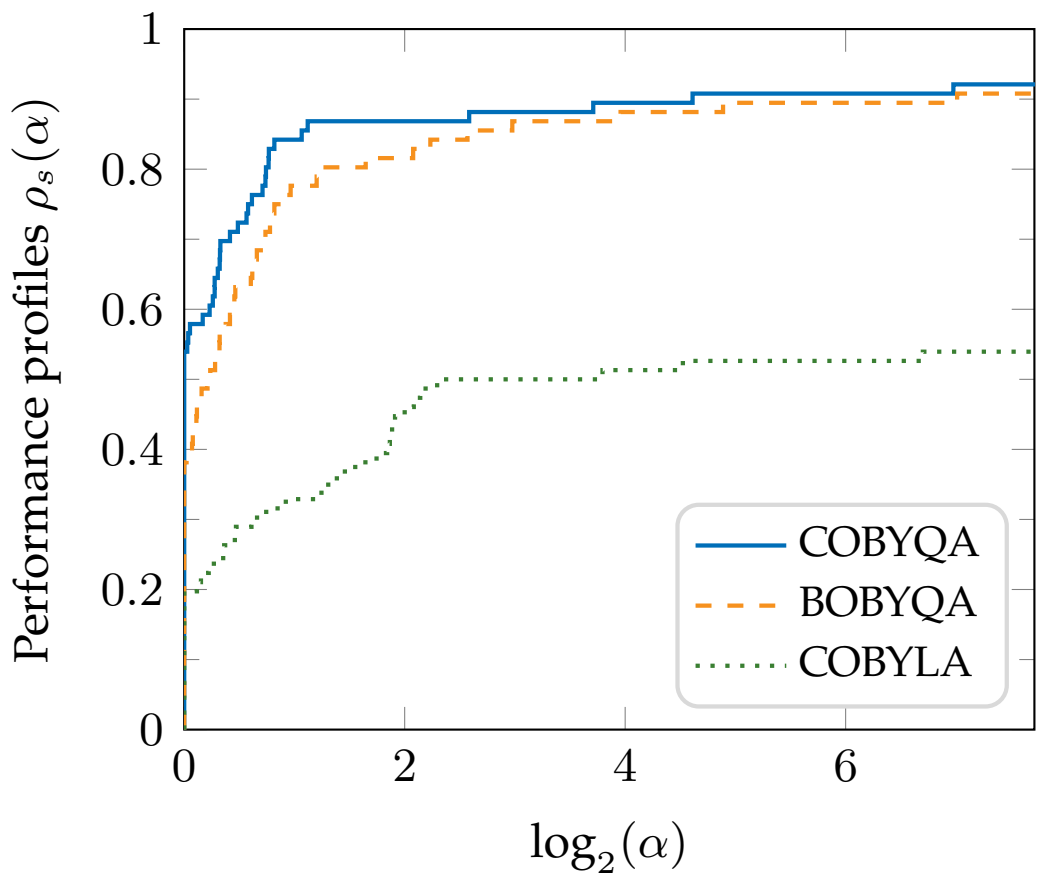

(b) Tolerance $\tau = 10^{-3}$

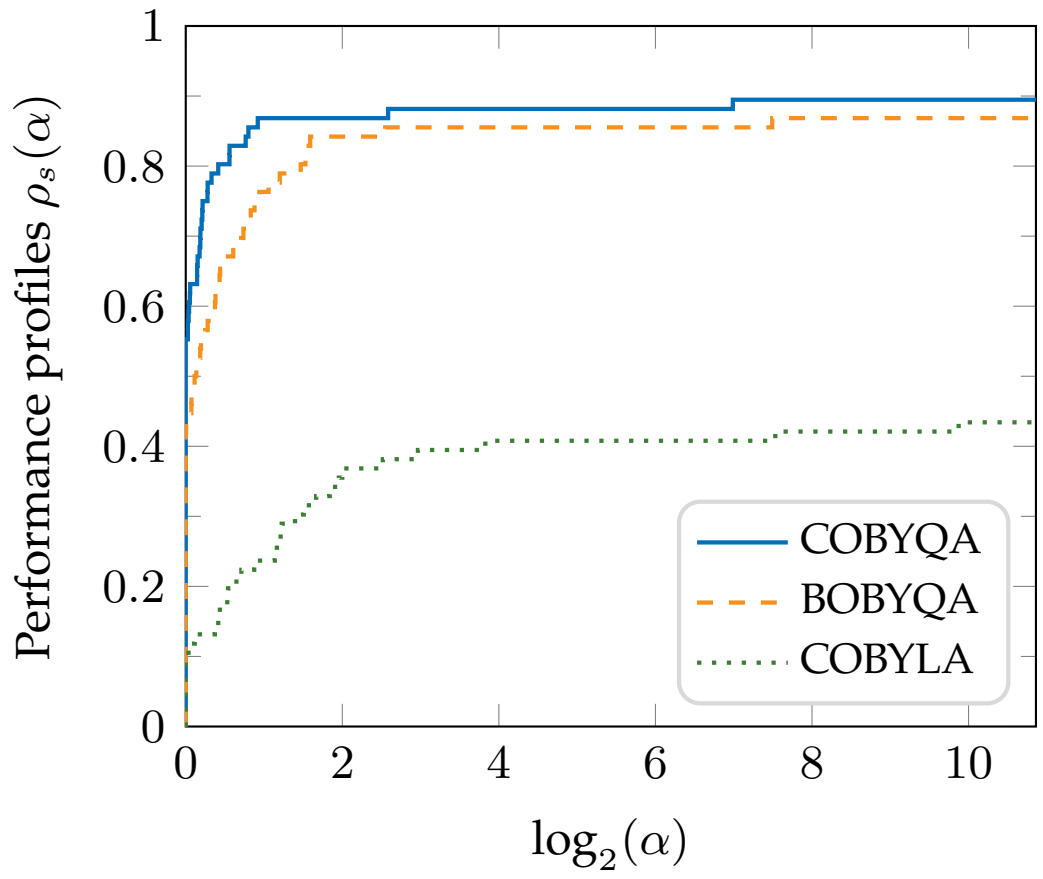

(c) Tolerance $\tau = 10^{-5}$

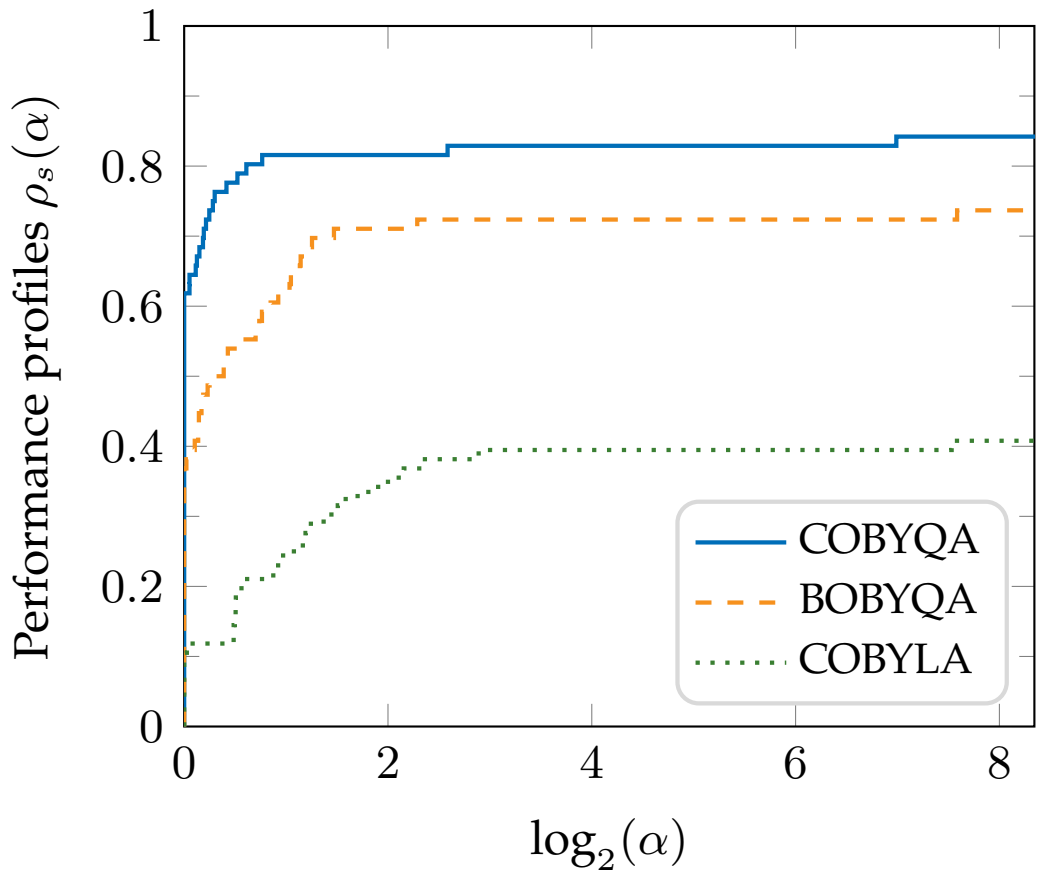

(d) Tolerance $\tau = 10^{-7}$

Figure 7.9: Performance profiles of COBYQA, BOBYQA, and COBYLA on bound-constrained problems

#### Comparison when the bounds are inviolable

We now consider the experiment when the bound constraints are inviolable. Among all the linearly constrained problems of the CUTEst library with at most 50 variables and at most 1,000 linear constraints listed in Appendix A.3, we test only those with at least one bound constraint. Figures 7.13 and 7.14 present the performance and data profiles of COBYQA, LINCOA, and COBYLA based on these problems.

On this experiment, COBYQA outperforms LINCOA and COBYQA, the margins be-ing significant in the high-precision tests with $\tau = 10^{-5}$ and $10^{-7}$. Generally speaking, COBYQA solves more problems, and uses the least amount of function evaluations on most of them.

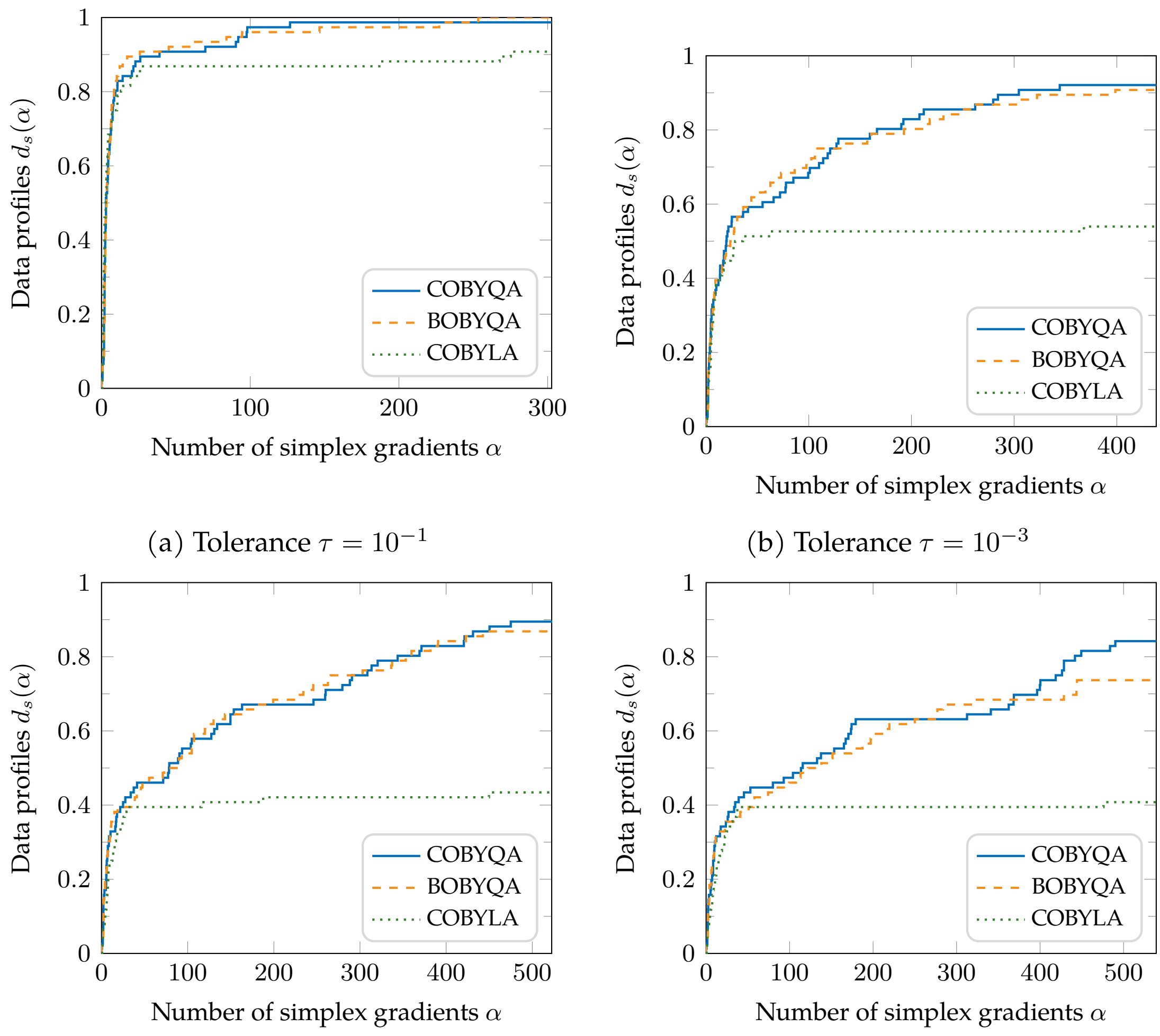

(a) Tolerance $\tau = 10^{-1}$

(b) Tolerance $\tau = 10^{-3}$

(c) Tolerance $\tau = 10^{-5}$

(d) Tolerance $\tau = 10^{-7}$

Figure 7.10: Data profiles of COBYQA, BOBYQA, and COBYLA on bound-constrained problems

It is noteworthy that COBYQA performs much better than LINCOA on these linearly constrained problems with inviolable bounds, whereas the comparison is the opposite when the bounds are violable. Consequently, we expect that COBYQA is a better solver than LINCOA for applications that contain inviolable bound constraints.

### 7.3.7 Performance of COBYQA compared with COBYLA

We now compare the performances of COBYQA and COBYLA. We conduct three experiments. The first two experiments compare the methods on nonlinearly constrained problems, with violable and inviolable bounds, respectively. In the second experiment, if a bound is violated, we set the objective function value to $\infty$, leaving the constraint values untouched. The last experiment compares COBYQA and COBYLA on all types of problems, constrained and unconstrained. As will be shown, COBYQA outperforms COBYLA in all scenarios.

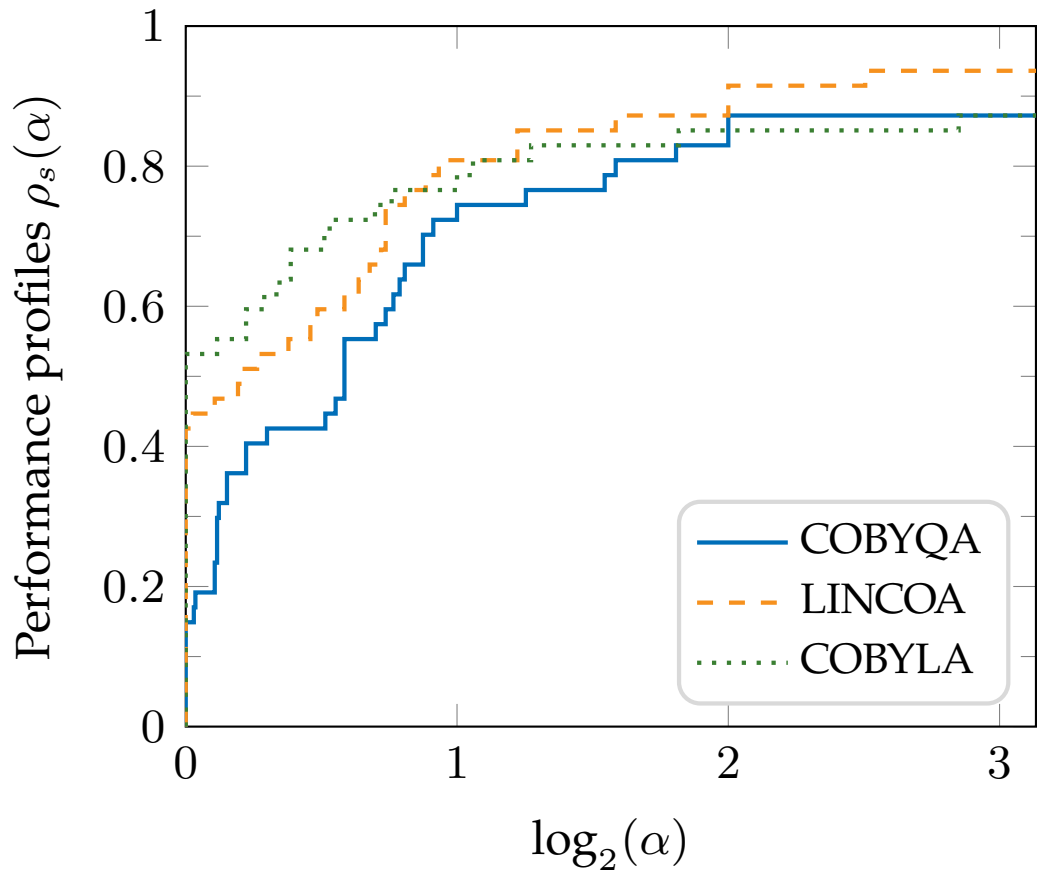

(a) Tolerance $\tau = 10^{-1}$

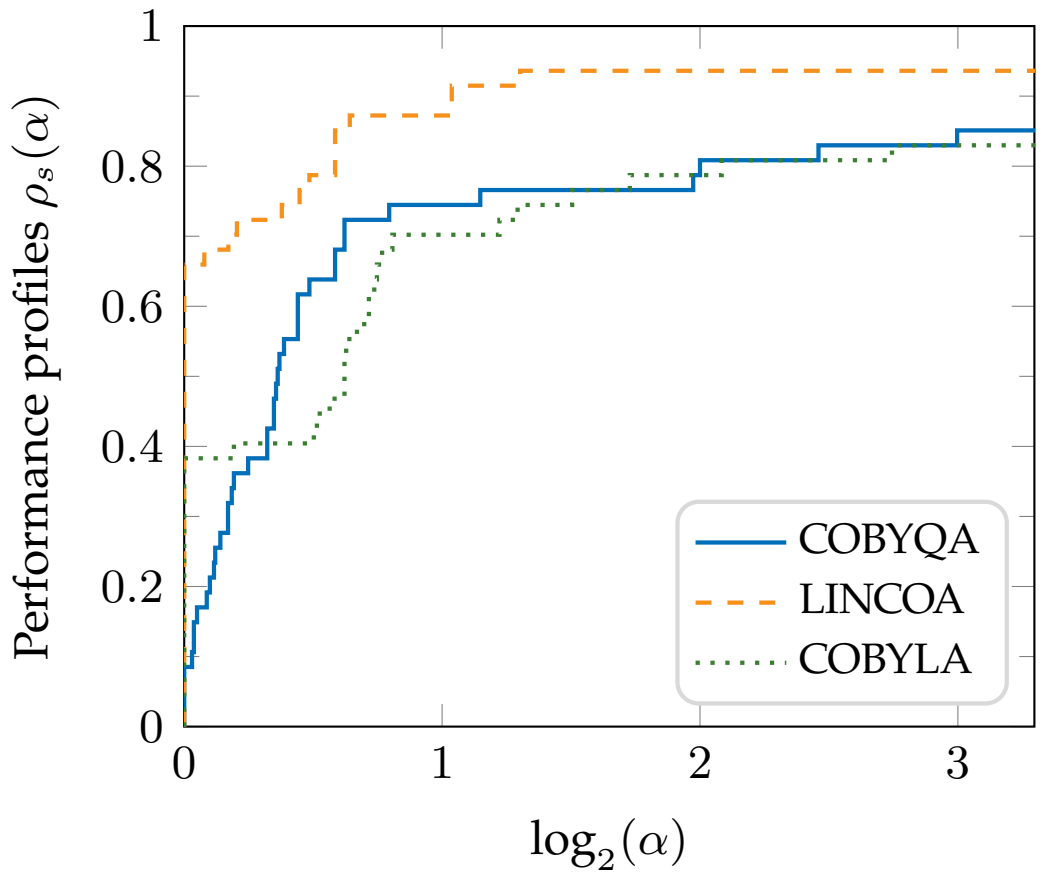

(b) Tolerance $\tau = 10^{-3}$

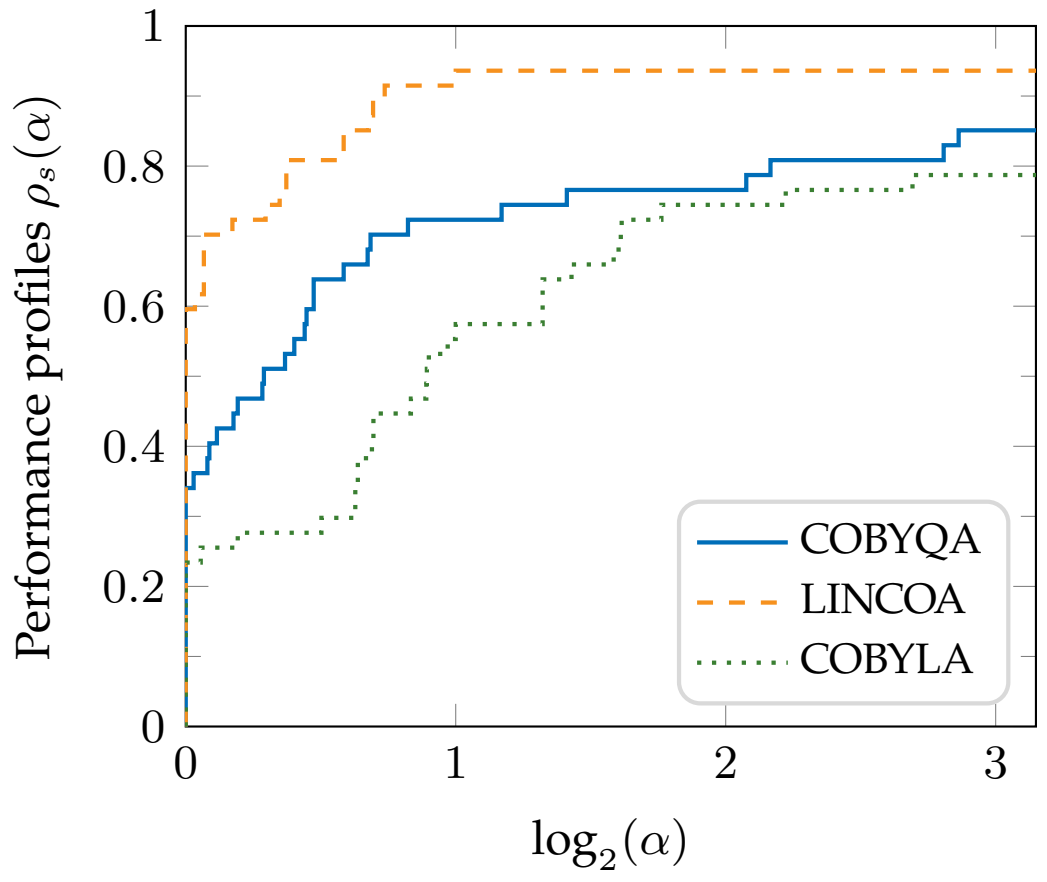

(c) Tolerance $\tau = 10^{-5}$

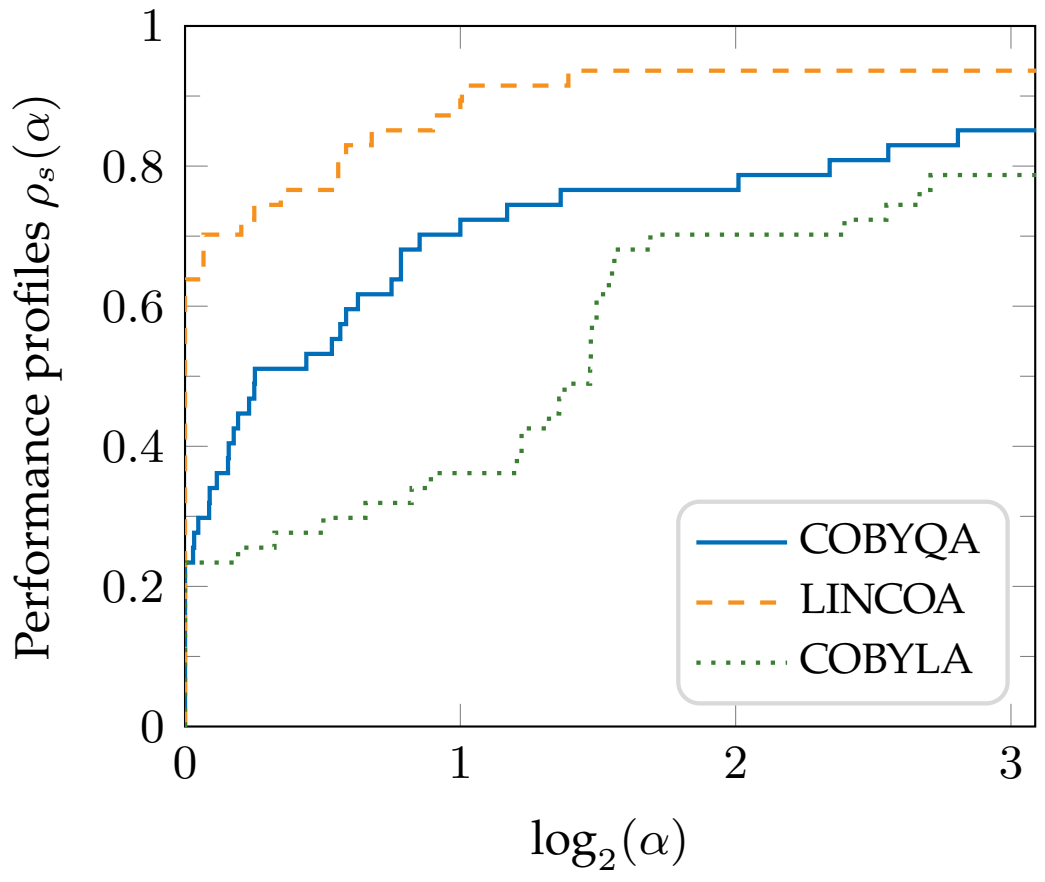

(d) Tolerance $\tau = 10^{-7}$

Figure 7.11: Performance profiles of COBYQA, LINCOA, and COBYLA on linearly constrained problems

**Comparison on nonlinearly constrained problem when the bounds are violable**

Figures 7.15 and 7.16 present the performance and data profiles of these solvers, based on all the problems listed in Appendix A.4. These are all the nonlinearly constrained problems of the CUTEst library with at most 50 variables and at most 1,000 constraints, excluding bound constraints.

We observe that COBYQA provides better performance than COBYLA for all considered tolerances. More specifically, COBYQA uses no more function evaluations than COBYLA to converge on about 60 % of the problems. They eventually solve a similar

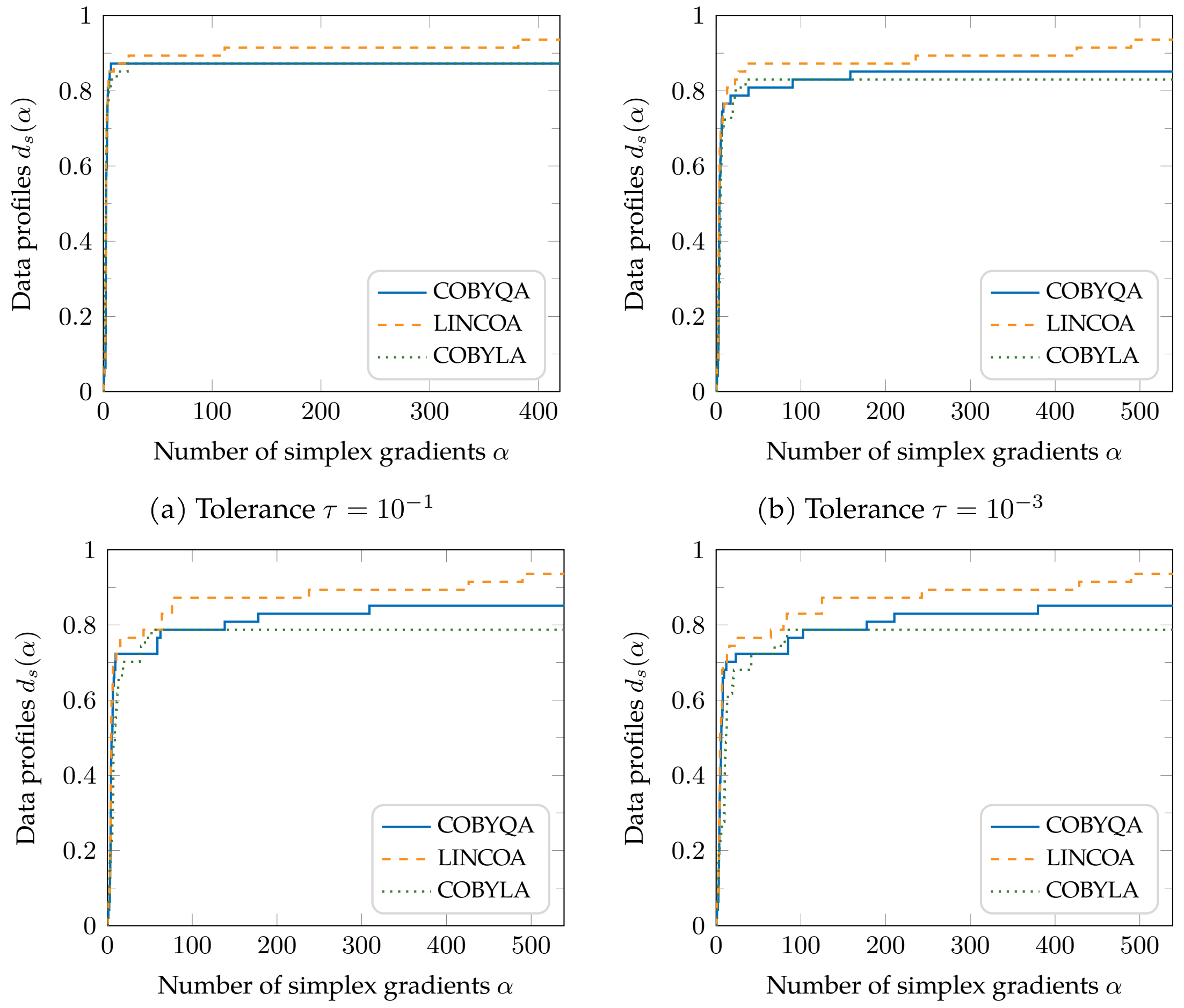

(a) Tolerance $\tau = 10^{-1}$

(b) Tolerance $\tau = 10^{-3}$

(c) Tolerance $\tau = 10^{-5}$

(d) Tolerance $\tau = 10^{-7}$

Figure 7.12: Data profiles of COBYQA, LINCOA, and COBYLA on linearly constrained problems

amount of problems, namely around 80 %.

### Comparison on nonlinearly constrained problem when the bounds are inviolable

We now consider the experiment when the bound constraints are inviolable. Among all the nonlinearly constrained problems of the CUTEst library with at most 50 variables and at most 1,000 constraints, excluding bound constraints, listed in Appendix A.4, we test only those with at least one bound constraint. Figures 7.17 and 7.18 present the performance and data profiles of COBYQA and COBYLA based on these problems.

The results are similar to the ones in Figures 7.15 and 7.16, although the advantage of COBYQA over COBYLA becomes even more visible. This is also expected, because COBYQA is designed to always respect the bound constraints, while COBYLA is not.

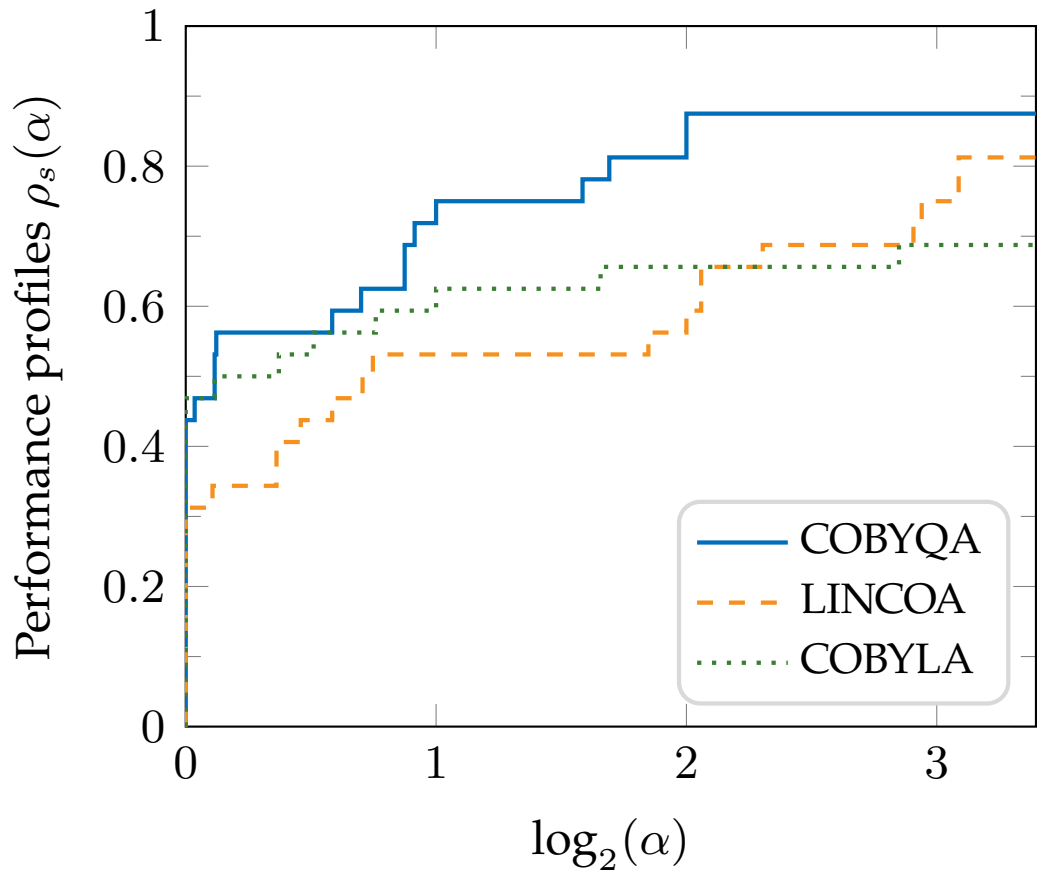

(a) Tolerance $\tau = 10^{-1}$

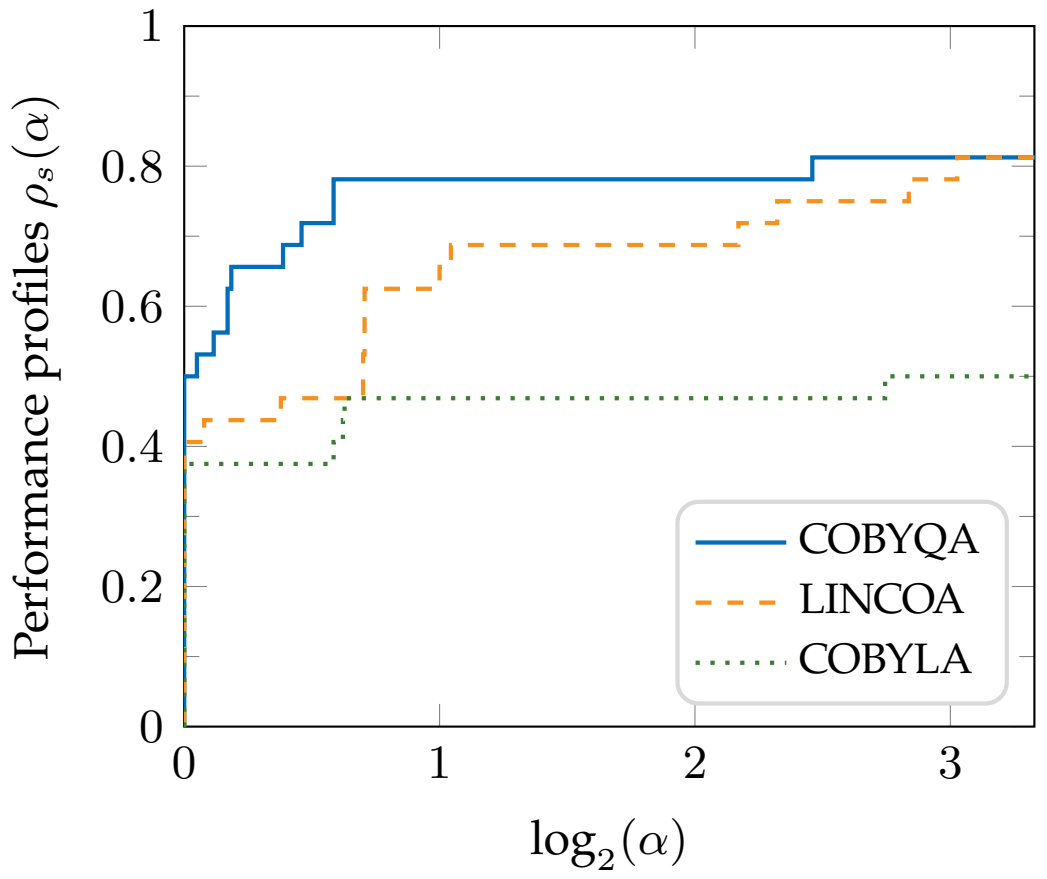

(b) Tolerance $\tau = 10^{-3}$

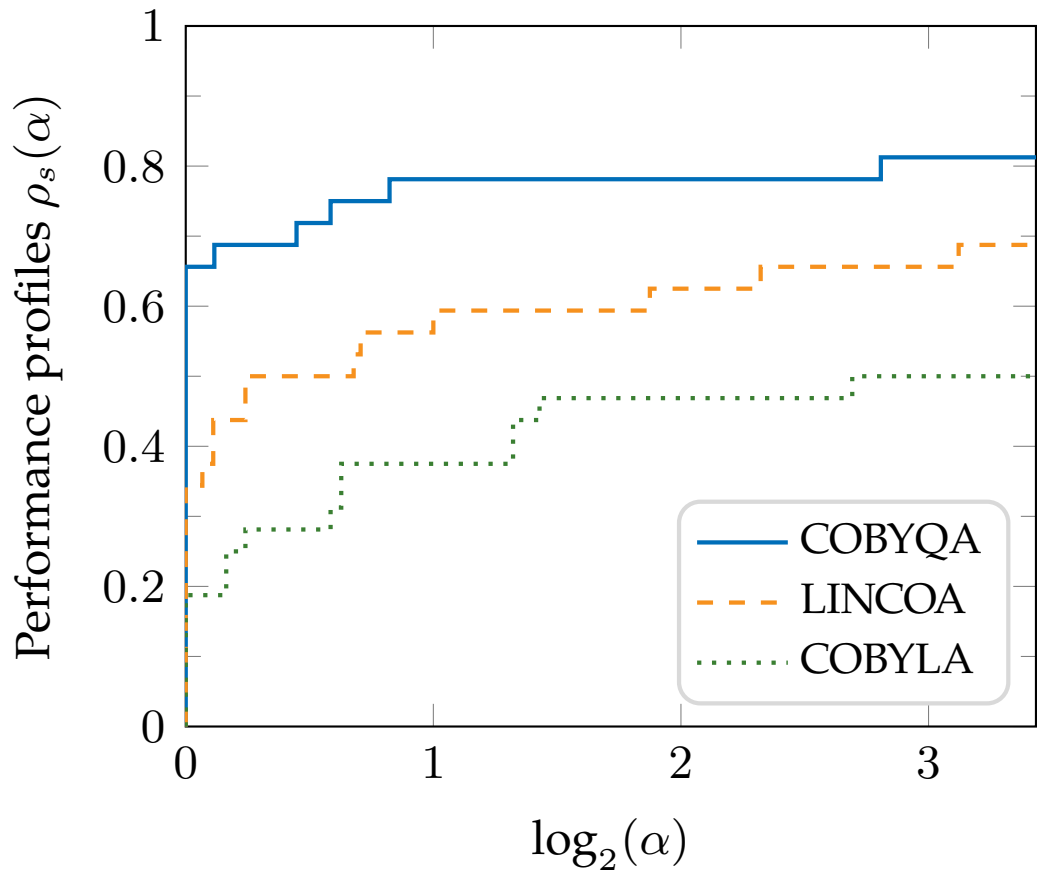

(c) Tolerance $\tau = 10^{-5}$

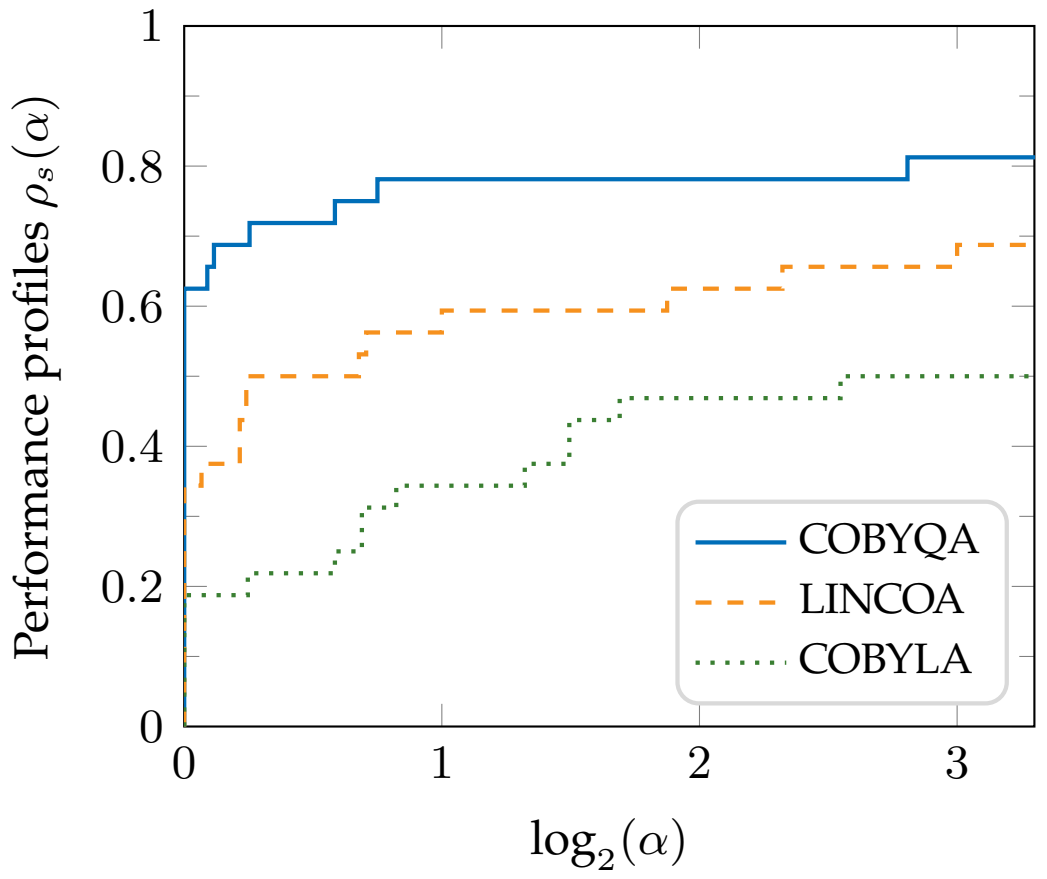

(d) Tolerance $\tau = 10^{-7}$

Figure 7.13: Performance profiles of COBYQA, LINCOA, and COBYLA on linearly constrained problems with inviolable bounds

### Comparison on all types of problems

One of the motivations for designing COBYQA is to provide a successor for COBYLA. To demonstrate that it is reasonable to replace COBYLA with COBYQA in general, we compare their performances on all types of problems, constrained and unconstrained. Figures 7.19 and 7.20 provide the performance and data profiles on all the CUTEst problems of dimension at most 50 listed in Appendices A.1 to A.4.

The superiority of COBYQA over COBYLA is demonstrated evidently by these profiles. COBYQA outperforms COBYLA with a substantial margin, and the margin is larger when the precision of the comparison is higher. When $\tau = 10^{-3}$, $10^{-5}$, and $10^{-7}$, COBYQA solves much more problems than COBYLA, while taking less function evaluations on most of the problems solved. These results are not surprising, because COBYLA uses linear models of the objective and constraint functions, whereas COBYQA uses quadratic models, which approximate these functions better in general.

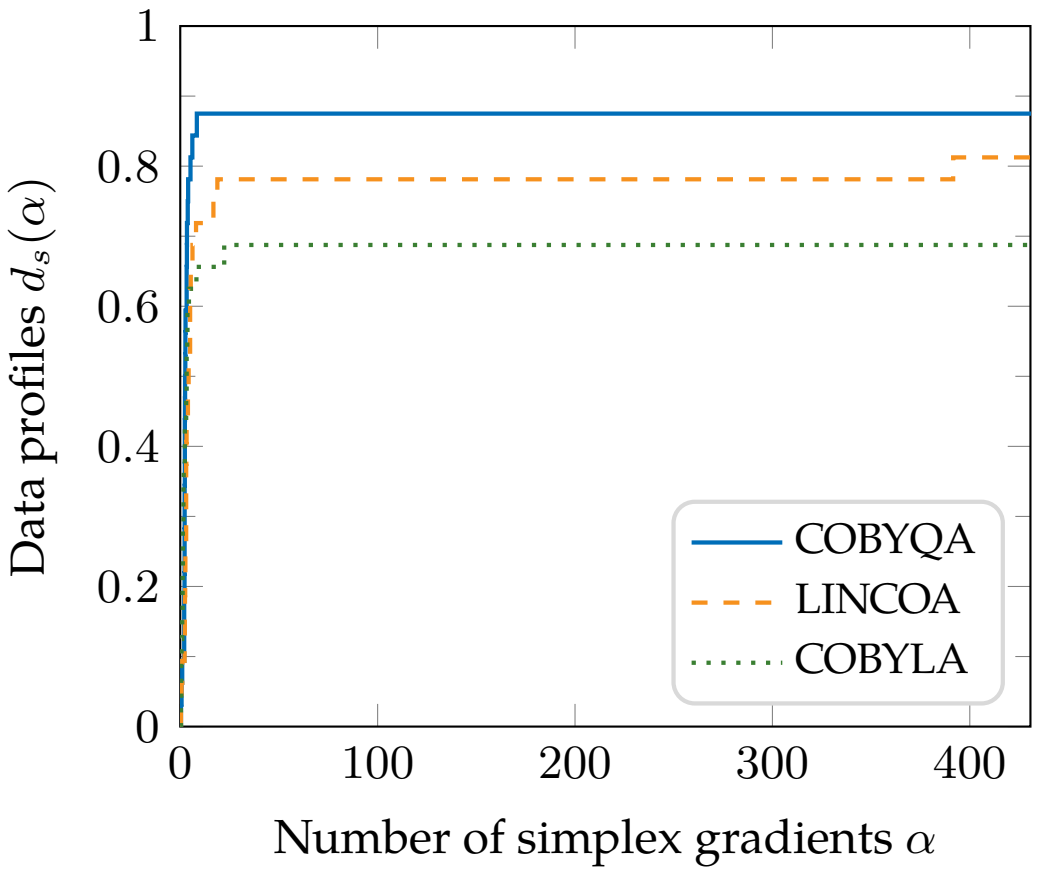

(a) Tolerance $\tau = 10^{-1}$

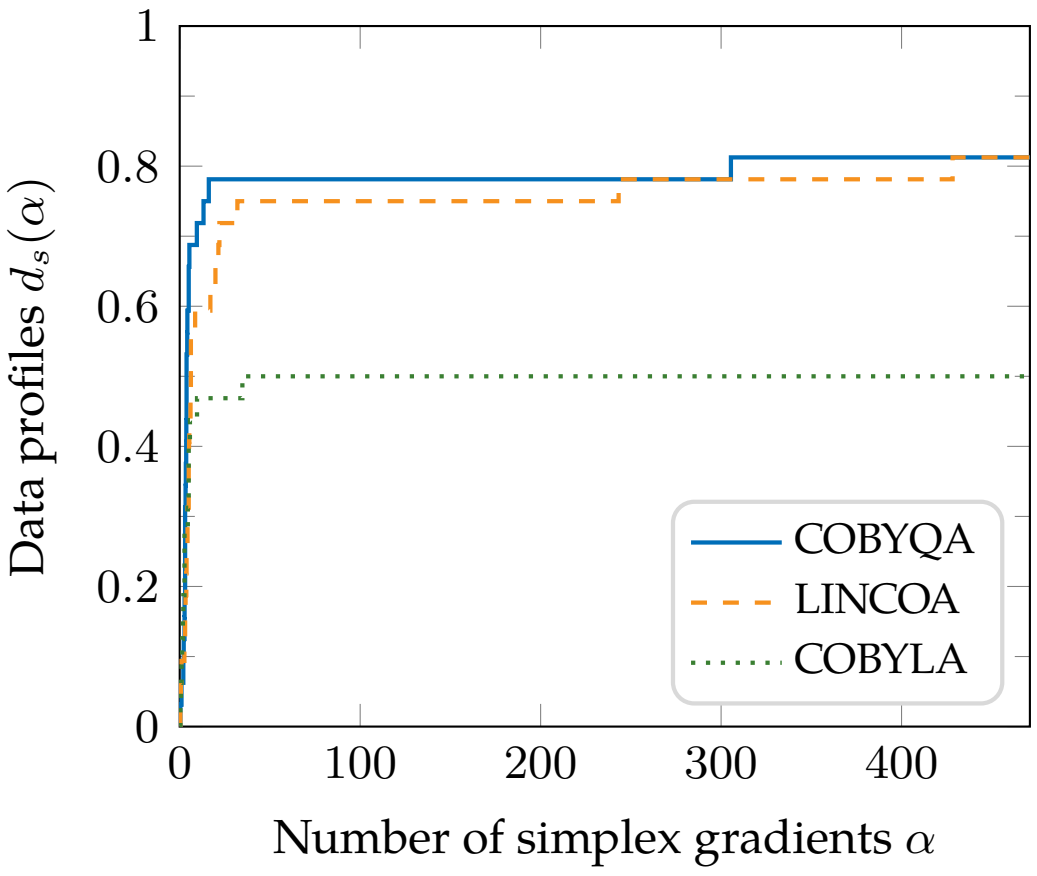

(b) Tolerance $\tau = 10^{-3}$

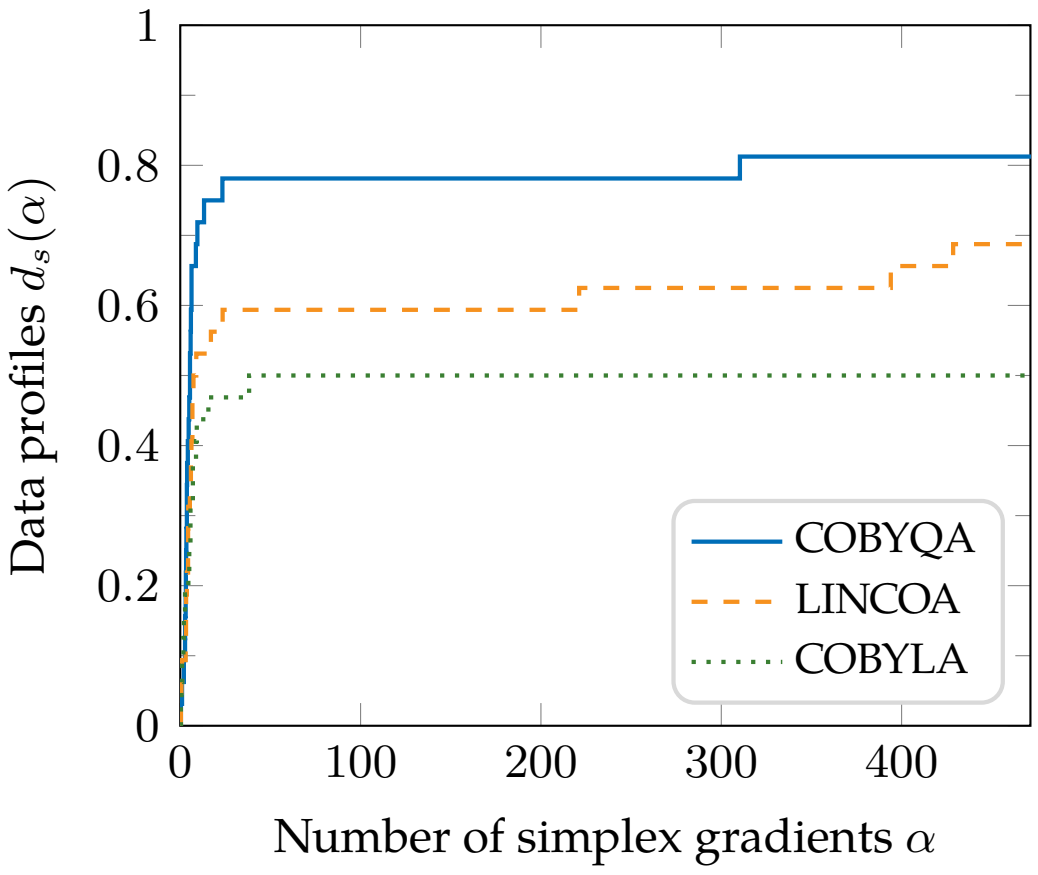

(c) Tolerance $\tau = 10^{-5}$

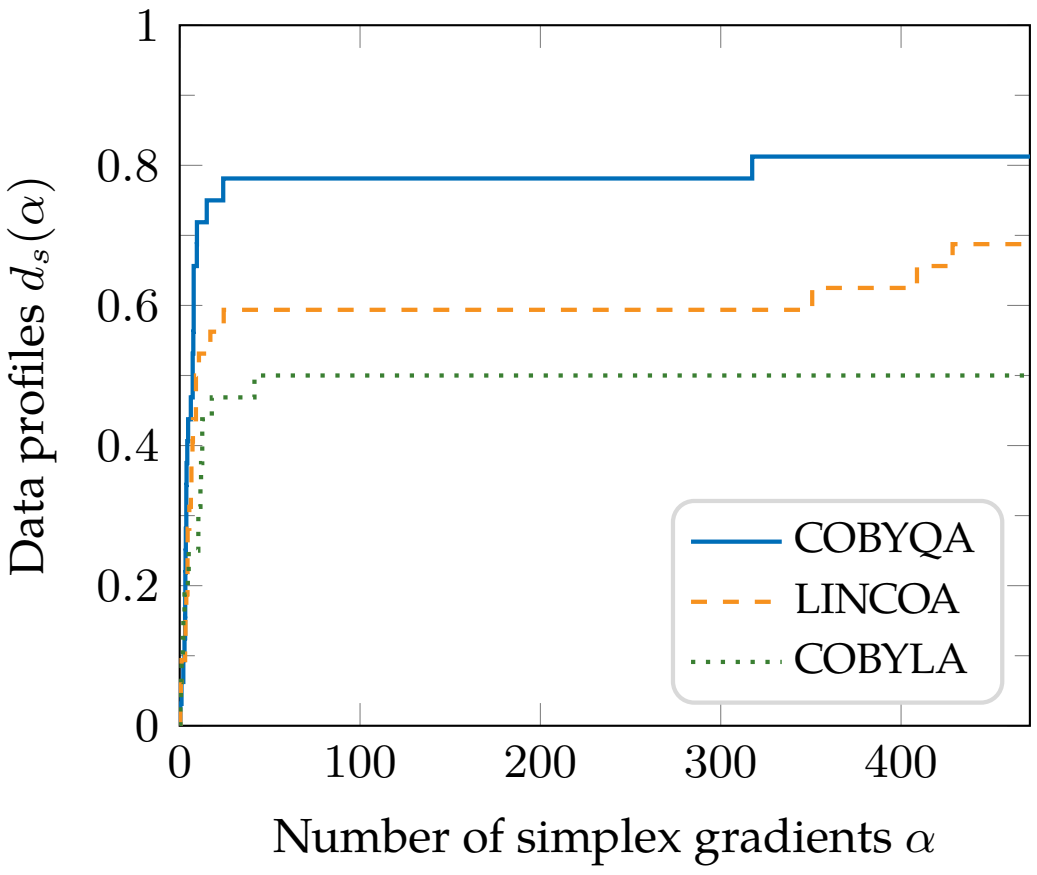

(d) Tolerance $\tau = 10^{-7}$

Figure 7.14: Data profiles of COBYQA, LINCOA, and COBYLA on linearly constrained problems with inviolable bounds

We also conducted a similar experiment when the bound constraints are inviolable. The results are similar to those presented in Figures 7.19 and 7.20 and hence omitted, although the margins are even larger, as expected.

To conclude, COBYQA performs convincingly better than COBYLA, and hence, is a great successor for COBYLA as a general-purpose solver.

## 7.4 Summary and remarks

In this chapter, we presented the Python implementation of COBYQA, and extensive numerical experiments.

We first provided additional algorithmic details, which were omitted in the general presentation of COBYQA in Chapters 5 and 6.

We then detailed the Python implementation of COBYQA, which is available at

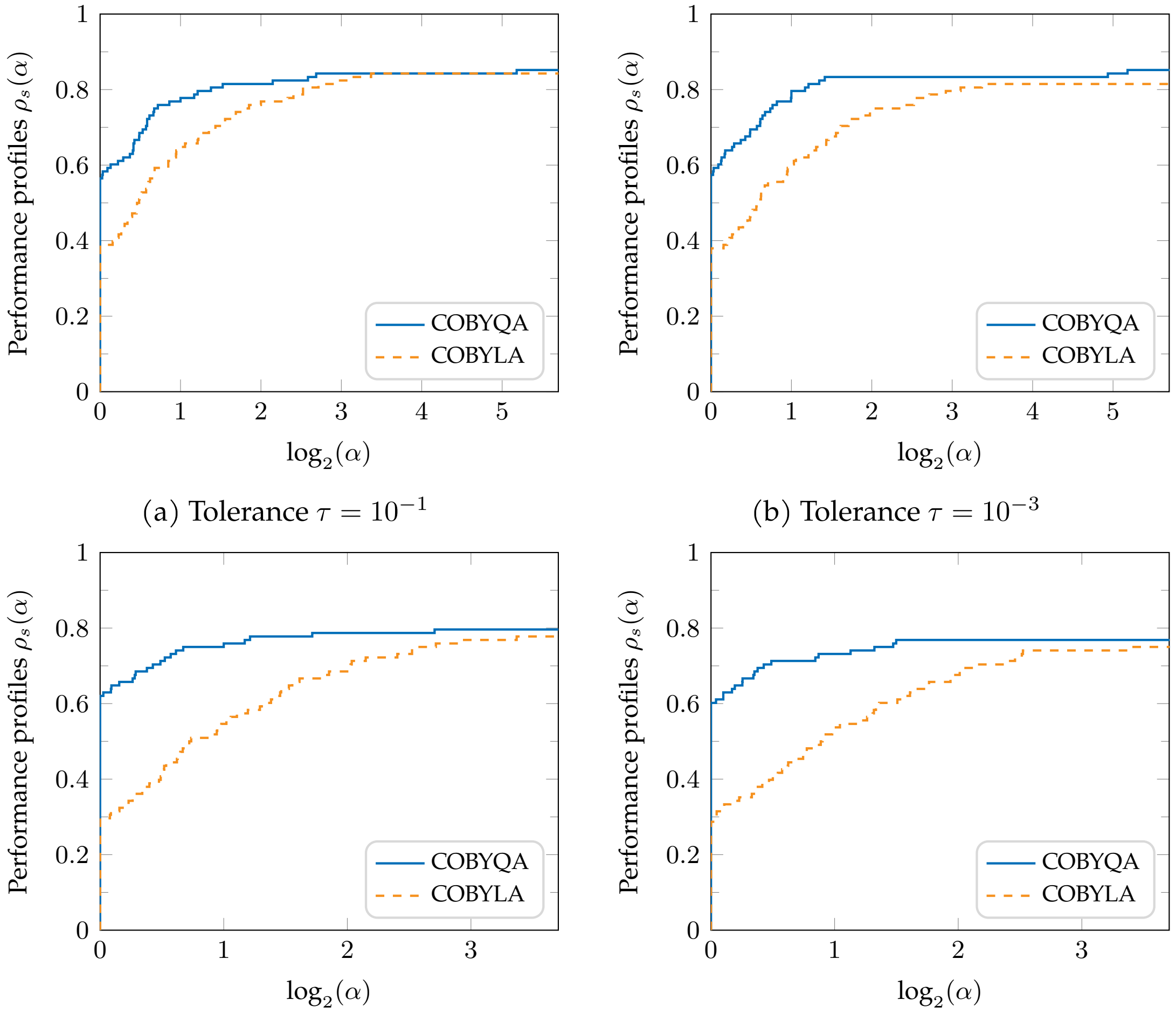

(a) Tolerance $\tau = 10^{-1}$ (b) Tolerance $\tau = 10^{-3}$

(c) Tolerance $\tau = 10^{-5}$ (d) Tolerance $\tau = 10^{-7}$

Figure 7.15: Performance profiles of COBYQA and COBYLA on nonlinearly constrained problems

`https://www.cobyqa.com/`. We briefly introduced some programming aspects of the implementation, provided an overview of the signature of the main Python function, and illustrated the usage by several simple examples.

Finally, we presented extensive numerical experiments on COBYQA. First, we demonstrated the effectiveness of several critical strategies of COBYQA, including the SQP methodology, our extension to the Byrd-Omojokun approach in the inequality case, and the derivative-free symmetric Broyden update of Powell. Second, we tested the performance of COBYQA compared with Powell's DFO solvers, which are standard benchmarks for the assessment of DFO methods. The performance of COBYQA is comparable to NEWUOA on unconstrained problems, and it outperforms BOBYQA on bound-constrained ones, while also being able to tackle more general problems. In contrast to LINCOA and COBYLA, a strength of COBYQA is that it always respects bound constraints. On linearly constrained problems, COBYQA clearly outperforms LINCOA if the problems contain inviolable bound constraints. Most importantly, COBYQA has much better performance than COBYLA on all types of problems, constrained and unconstrained, no matter whether bound constraints (if any) can be violated or not.

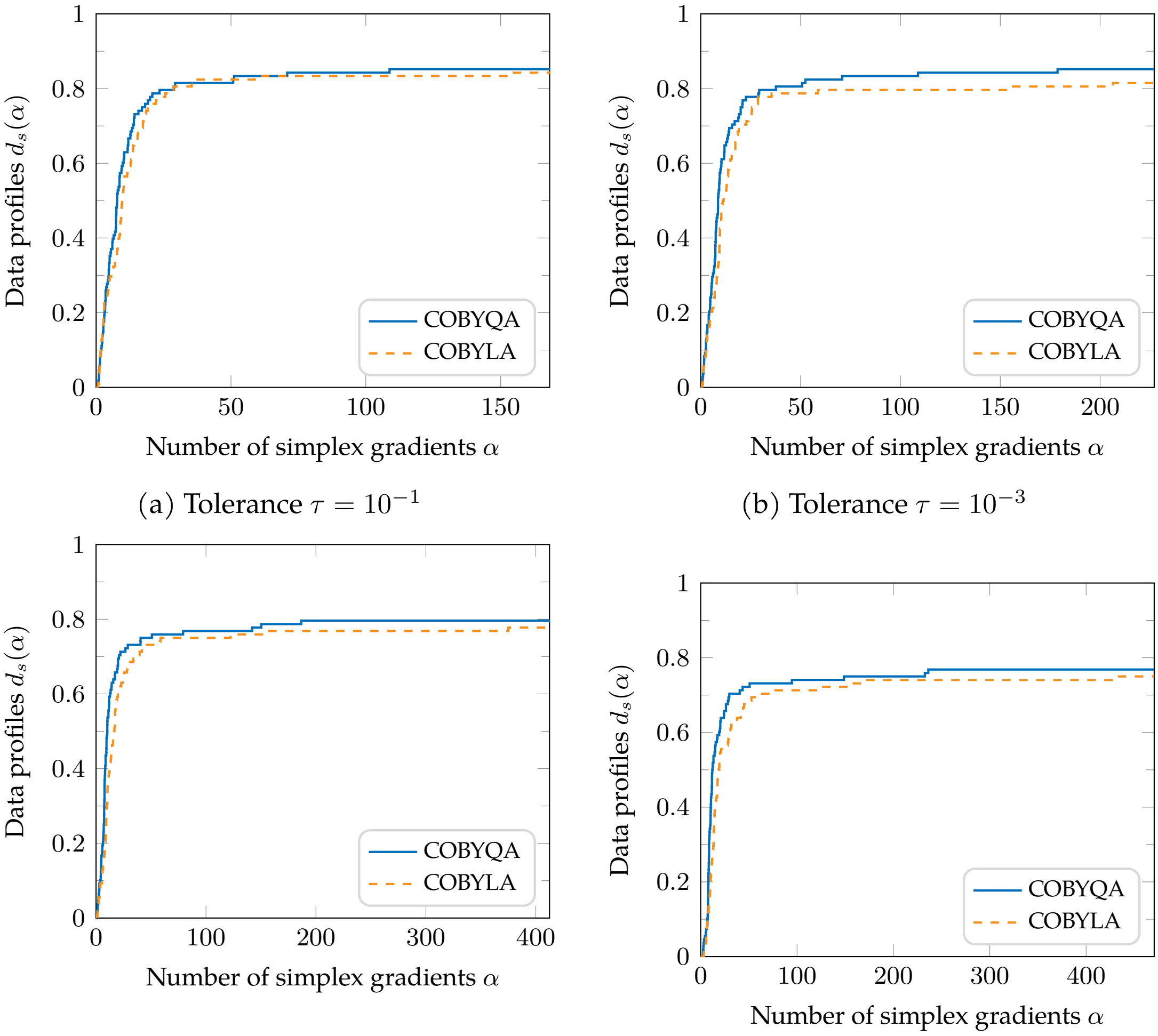

(a) Tolerance $\tau = 10^{-1}$

(b) Tolerance $\tau = 10^{-3}$

(c) Tolerance $\tau = 10^{-5}$

(d) Tolerance $\tau = 10^{-7}$

Figure 7.16: Data profiles of COBYQA and COBYLA on nonlinearly constrained problems

The implementation of COBYQA has been highly challenging, as is the case for all model-based DFO solvers[1]. We endeavored to develop COBYQA in order to provide a successor to COBYLA as a general-purpose DFO solver. Finally, the numerical experiments have demonstrated that our goal is achieved.

[1]Recall that Powell [166] mentioned that the development of NEWUOA was time-consuming, and "was very frustrating," so was the implementation of COBYQA.

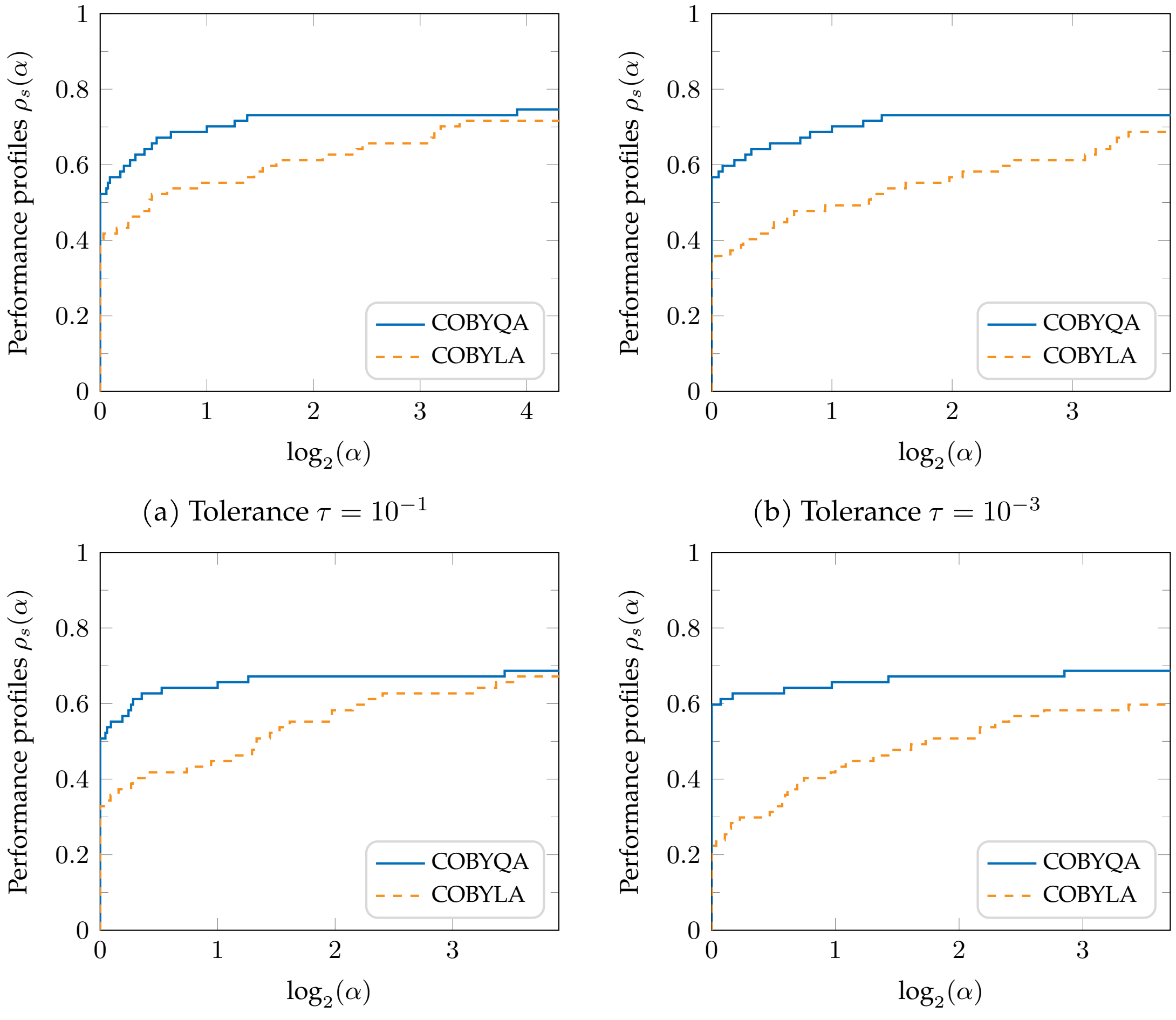

(a) Tolerance $\tau = 10^{-1}$

(b) Tolerance $\tau = 10^{-3}$

(c) Tolerance $\tau = 10^{-5}$

(d) Tolerance $\tau = 10^{-7}$

Figure 7.17: Performance profiles of COBYQA and COBYLA on nonlinearly constrained problems with inviolable bounds

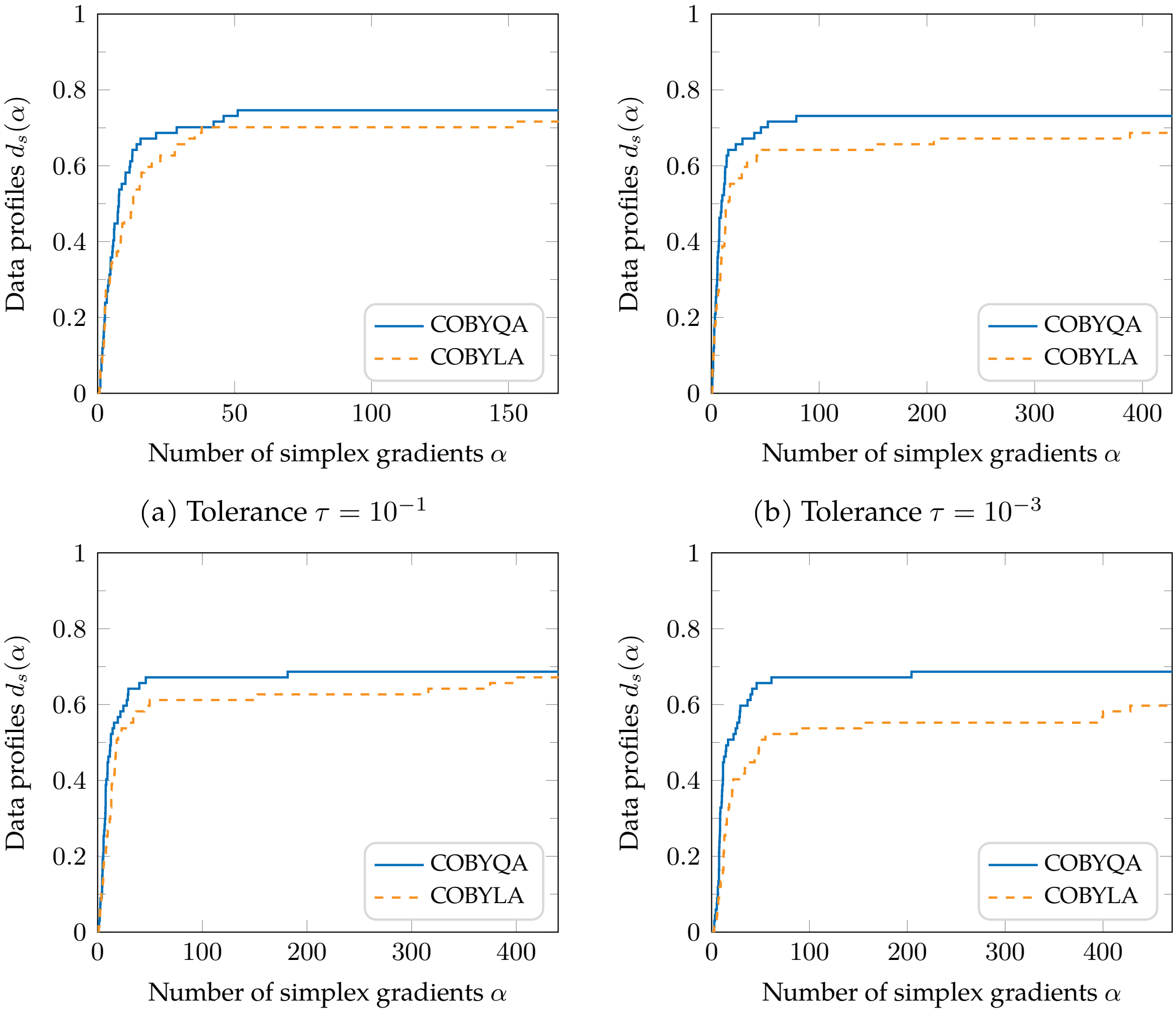

(a) Tolerance $\tau = 10^{-1}$

(b) Tolerance $\tau = 10^{-3}$

(c) Tolerance $\tau = 10^{-5}$

(d) Tolerance $\tau = 10^{-7}$

Figure 7.18: Data profiles of COBYQA and COBYLA on nonlinearly constrained problems with inviolable bounds

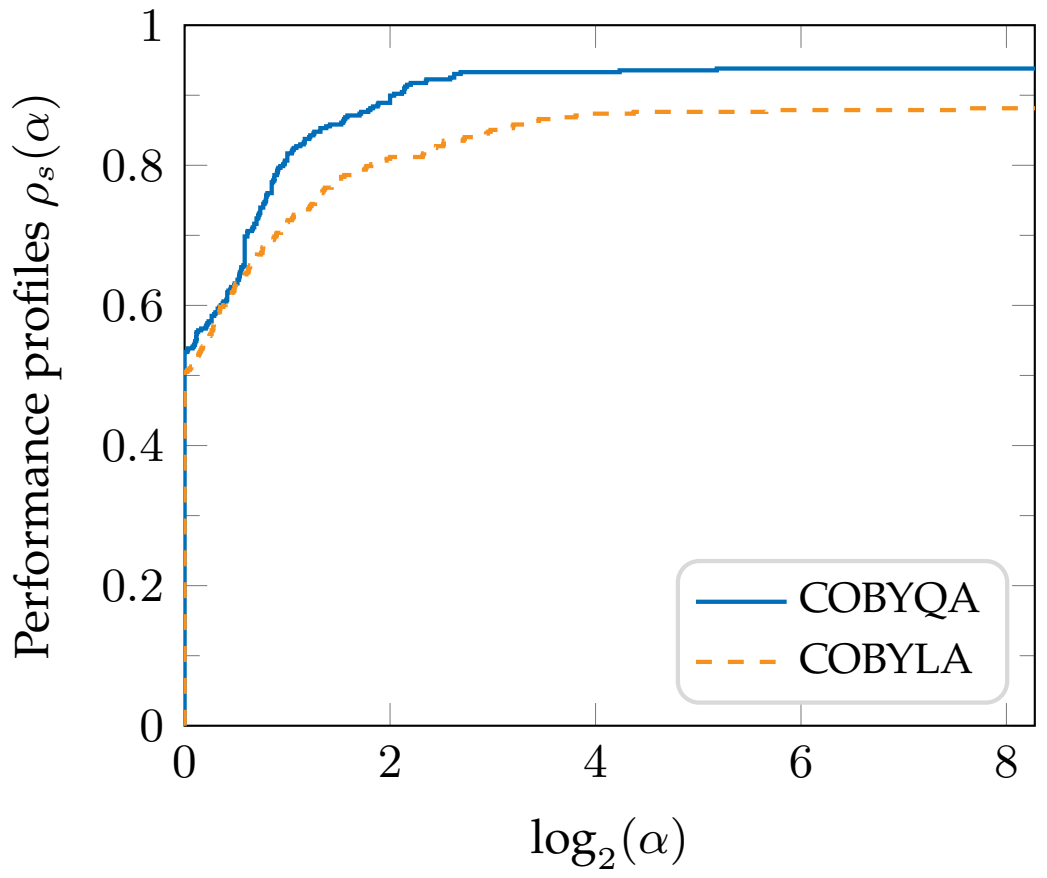

(a) Tolerance $\tau = 10^{-1}$

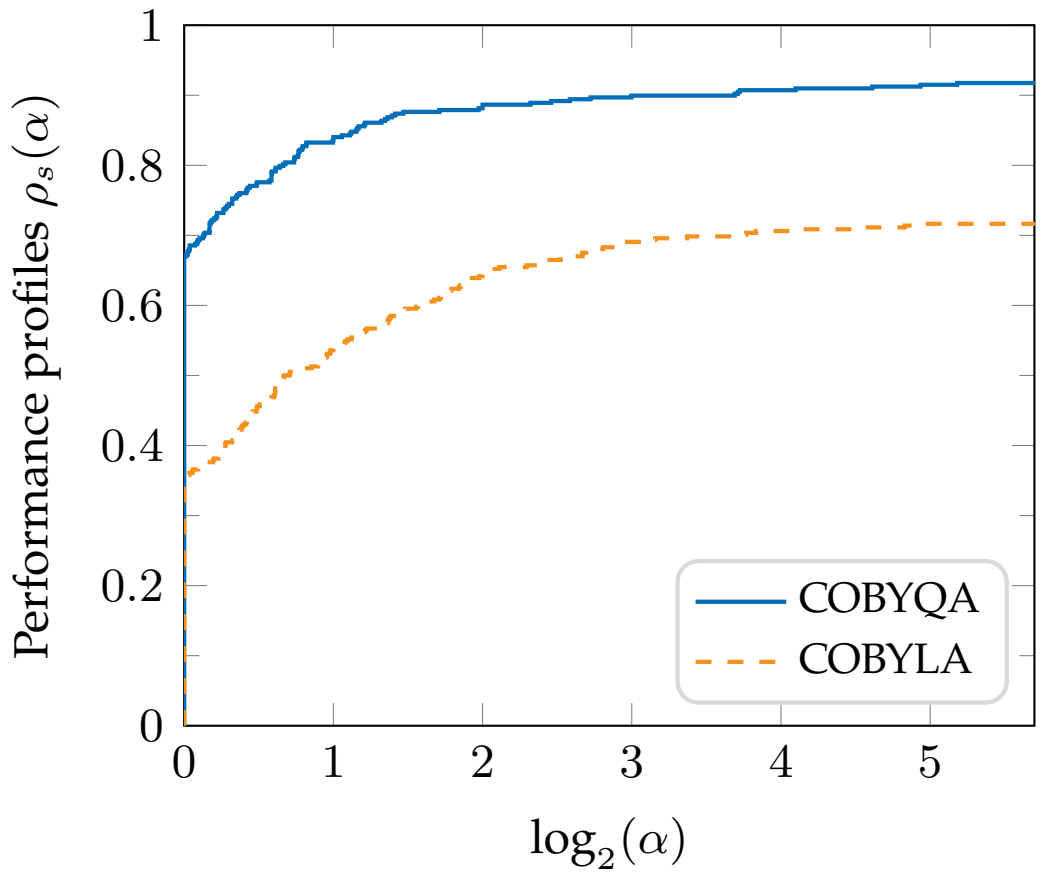

(b) Tolerance $\tau = 10^{-3}$

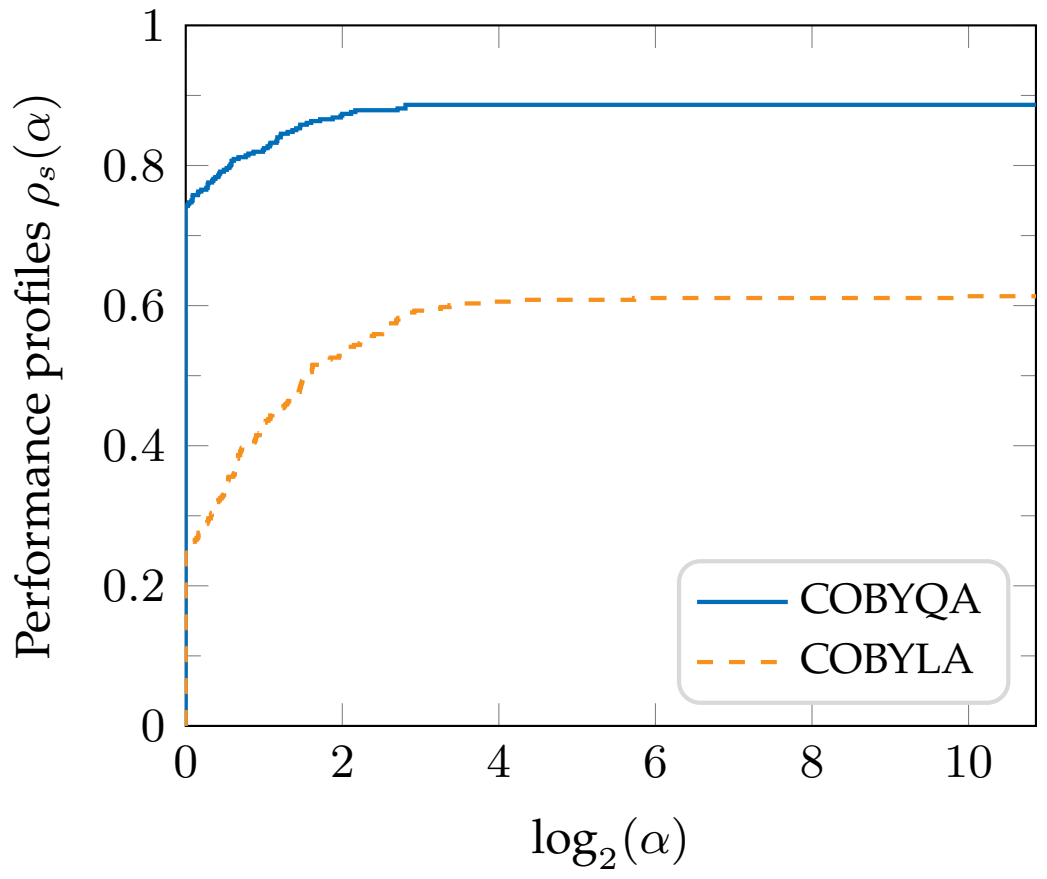

(c) Tolerance $\tau = 10^{-5}$

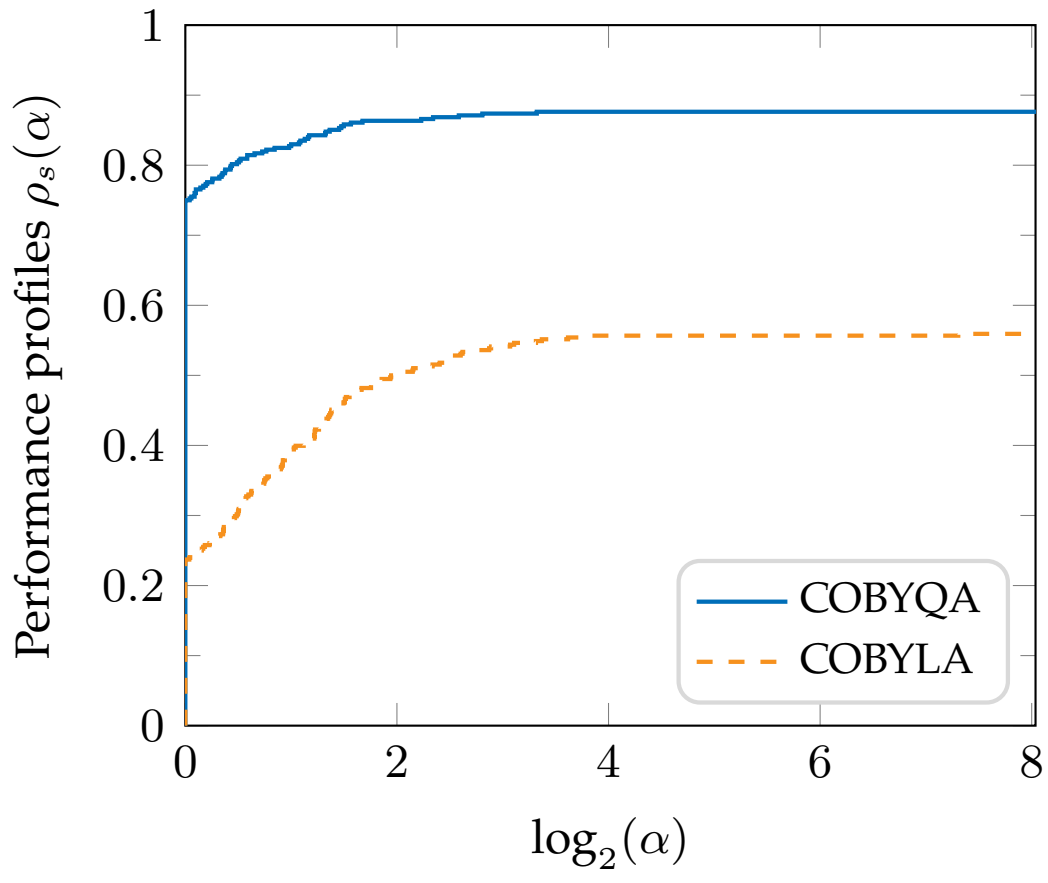

(d) Tolerance $\tau = 10^{-7}$

Figure 7.19: Performance profiles of COBYQA and COBYLA on all types of problems

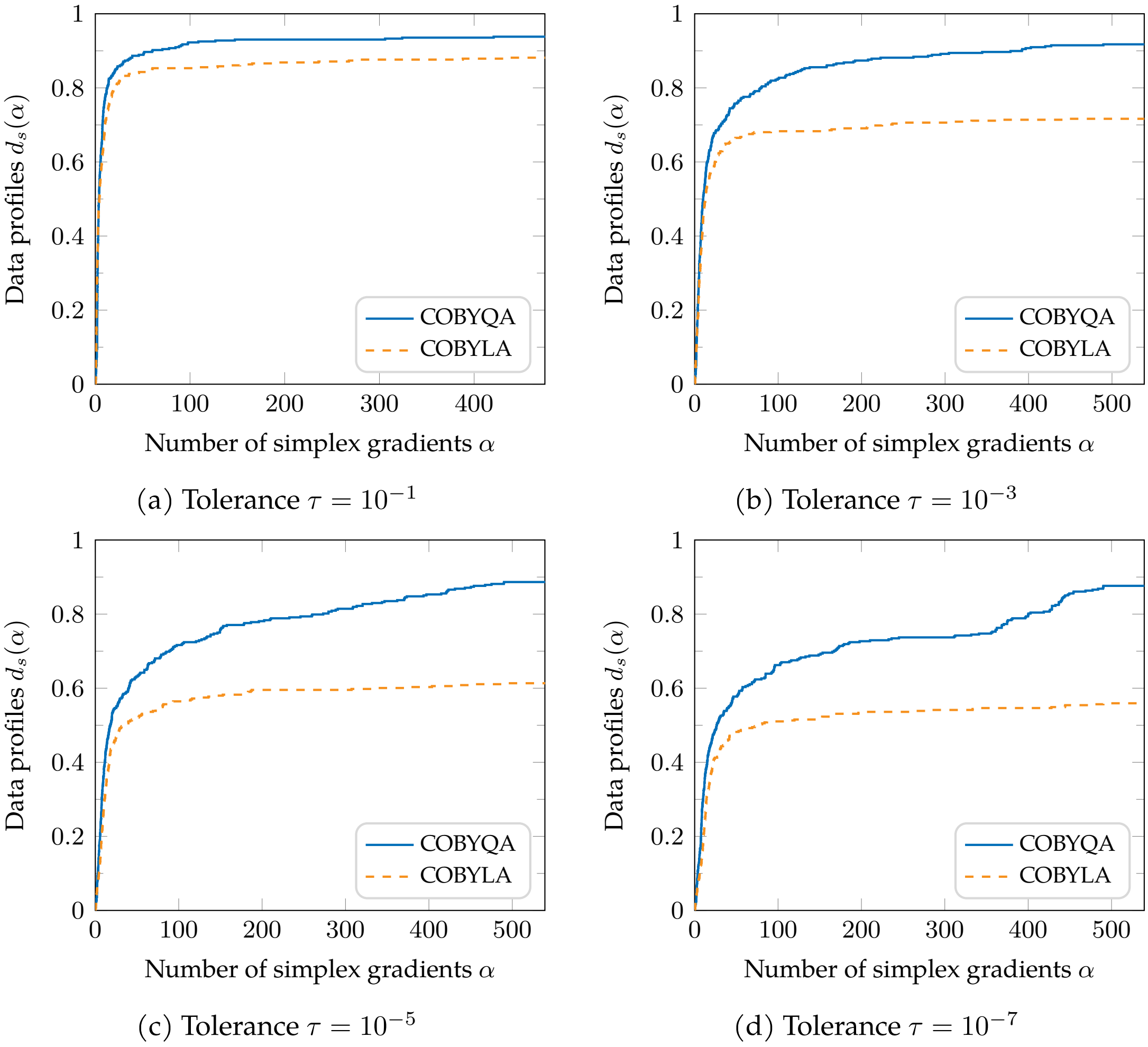

(a) Tolerance $\tau = 10^{-1}$

(b) Tolerance $\tau = 10^{-3}$

(c) Tolerance $\tau = 10^{-5}$

(d) Tolerance $\tau = 10^{-7}$

Figure 7.20: Data profiles of COBYQA and COBYLA on all types of problems

*This page intentionally left blank*

# 8 Conclusion and future research directions

## 8.1 Conclusion

In this thesis, we investigated model-based methods for DFO, covering both theoretical and practical aspects.

On the theoretical side, we first studied properties of an interpolation set for underdetermined interpolation employed by Powell in NEWUOA. We showed that this set is optimal when composed of $2n + 1$ points, $n$ being the dimension of the problem. This supports the recommendation made by Powell, and it also guides us in devising the initial interpolation set of our new DFO method.

We also examined several theoretical aspects of the SQP method. In particular, we highlighted a new interpretation of the SQP subproblem, regarding its objective function as a quadratic approximation of the original objective function in the tangent space of a surface. Moreover, we established a relation between the Vardi and the Byrd-Omojokun approaches in the equality-constrained case, interpreting the former as an approximation of the latter. We also briefly studied the Lagrangian and augmented Lagrangian functions of the SQP subproblem, establishing connections between the SQP method and certain augmented Lagrangian methods.

Writing a thesis on model-based DFO methods, we allocated a major part of our efforts to practices and implementations of these methods, which are essential in DFO. In particular, we developed in Chapter 3 a cross-platform package named PDFO, providing user-friendly MATLAB and Python interfaces for employing late Prof. M. J. D. Powell's DFO solvers. The package also patches several bugs contained in the Fortran 77 code. As of August 2022, this package has been downloaded more than $30{,}000$ times, and has been included in GEMSEO, an engine for MDO.

The aforementioned theoretical and practical study lays a solid groundwork for our development of a new derivative-free method for nonlinearly constrained optimization. Named COBYQA, it is a trust-region SQP method that models the objective and constraint functions using quadratic models based on the derivative-free symmetric Broyden update. The development of this method is fully documented in Chapters 5 to 7, covering the mathematical basis, the detailed descriptions of the subproblem solvers, and the Python implementation. In this development, we summarized, adapted, and extended various tools for the practices of SQP and DFO. For example, we proposed a new extension of the Byrd-Omojokun approach for the inequality-constrained case, which worked remarkably well in our tests. Finally, extensive numerical experiments demonstrated the excellent performance of COBYQA, confirming that we have achieved our initial goal of developing a general-purpose DFO solver, as a successor to COBYLA.

## 8.2 Future research directions

The study of this thesis opens many interesting directions for further investigations. We briefly mention four of them.

1. When developing COBYQA, we followed Powell's philosophy reflected in his work on DFO, which is to first produce an efficient and robust solver that performs well in practice, instead of setting up a conceptual framework without a practical implementation. Now that COBYQA is computationally successful, we will proceed to analyze its convergence properties, including global convergence and worst-case complexity. We expect that such an analysis can be done based on the existing techniques for analyzing DFO methods in the unconstrained case, and the SQP method in the derivative-based case.
2. We will further improve the implementation of COBYQA to enhance its performance, particularly in the linearly constrained case, recalling that the current version does not perform as well as LINCOA on linearly constrained problems with violable bounds. In addition, we plan to implement COBYQA in modern Fortran (Fortran 2018 or above), because the current Python version can be slow for high-dimensional problems (recall, however, the dominating expense of DFO in practice is the function evaluations rather than the execution of the solver's own code). We will also make COBYQA available in other languages, such as MATLAB, Julia, and R. The implementation of PDFO in Julia and R are also under consideration.
3. Our new interpretation of the SQP subproblem sheds lights on the approximation of nonlinear functions on manifolds. It will be interesting to investigate the implication of our interpretation on manifold optimization.
4. The connections between the SQP and the augmented Lagrangian methods that we mentioned in Section 4.1.4 also deserve further investigations. We will study the implications of these connections on the theory and practice of the SQP and the augmented Lagrangian methods, trying to develop theories that bridge these two methodologies and algorithms that combine their strengths.

*This page intentionally left blank*

# A List of testing problems

In this appendix, we provide the lists of CUTEst problems that are tested in our experiments. We employed PyCUTEst[1] to load the CUTEst problems in Python, and we select all the available problems with a dimension at most 50 and a number of constraints at most 1,000, excluding bound constraints. Some experiments require the problems to admit at least one bound constraint, so that we extract the problems from these lists accordingly. For detailed information on these problems, we refer to the official GitHub repository of CUTEst at `https://github.com/ralna/CUTEst`.

## A.1 Unconstrained problems

Table A.1: List of unconstrained problems

| | | | |
|---|---|---|---|
| AKIVA | ALLINITU | ARGLINA | ARGLINB |
| ARGLINC | ARGTRIGLS | BARD | BEALE |
| BENNETT5LS | BIGGS6 | BOX | BOX3 |
| BOXBODLS | BOXPOWER | BRKMCC | BROWNAL |
| BROWNBS | BROWNDEN | BROYDN3DLS | BROYDNBDLS |
| BRYBND | CHNROSNB | CHNRSNBM | CHWIRUT1LS |
| CHWIRUT2LS | CLIFF | COSINE | CUBE |
| DENSCHNA | DENSCHNB | DENSCHNC | DENSCHND |
| DENSCHNE | DENSCHNF | DIXON3DQ | DJTL |
| DQDRTIC | DQRTIC | ECKERLE4LS | EDENSCH |
| ENGVAL1 | ENGVAL2 | ENSOLS | ERRINROS |
| ERRINRSM | EXPFIT | EXTROSNB | FBRAIN3LS |
| FLETBV3M | FLETCBV2 | FLETCBV3 | FLETCHBV |
| FLETCHCR | FREUROTH | GAUSSIAN | GBRAINLS |
| GENHUMPS | GENROSE | GROWTHLS | GULF |
| HAHN1LS | HAIRY | HATFLDD | HATFLDE |

[1] See `https://jfowkes.github.io/pycutest/`.

| | | | |
|---|---|---|---|
| HATFLDFL | HEART6LS | HEART8LS | HELIX |
| HIELOW | HILBERTA | HILBERTB | HIMMELBB |
| HIMMELBF | HIMMELBG | HIMMELBH | HUMPS |
| INDEFM | JENSMP | KIRBY2LS | KOWOSB |
| LANCZOS1LS | LANCZOS2LS | LANCZOS3LS | LIARWHD |
| LOGHAIRY | LSC1LS | LSC2LS | MANCINO |
| MARATOSB | MEXHAT | MEYER3 | MGH09LS |
| MGH10LS | MISRA1BLS | MISRA1DLS | MOREBV |
| NCB20B | NONCVXU2 | NONCVXUN | NONDIA |
| OSBORNEB | OSCIGRAD | OSCIPATH | PALMER1C |
| PALMER1D | PALMER2C | PALMER3C | PALMER4C |
| PALMER5C | PALMER5D | PALMER6C | PALMER7C |
| PALMER8C | PARKCH | PENALTY1 | PENALTY2 |
| POWELLBSLS | POWELLSG | POWER | QUARTC |
| RAT42LS | ROSENBR | ROSENBRTU | ROSZMAN1LS |
| S308 | SBRYBND | SCHMVETT | SCOSINE |
| SCURLY10 | SENSORS | SINEVAL | SINQUAD |
| SISSER | SNAIL | SPARSINE | SPARSQUR |
| SSBRYBND | SSCOSINE | SSI | STRATEC |
| STREG | THURBERLS | TOINTGOR | TOINTGSS |
| TOINTPSP | TOINTQOR | TQUARTIC | TRIDIA |
| VARDIM | VAREIGVL | VESUVIALS | VESUVIOLS |
| VESUVIOULS | VIBRBEAM | WATSON | YFITU |
| ZANGWIL2 | | | |

## A.2 Bound-constrained problems

Table A.2: List of bound-constrained problems

| | | | |
|---|---|---|---|
| 3PK | ALLINIT | BLEACHNG | BQP1VAR |
| BQPGABIM | BQPGASIM | CAMEL6 | CHEBYQAD |
| CVXBQP1 | EG1 | FBRAIN2LS | FBRAINLS |
| HART6 | HATFLDA | HATFLDB | HATFLDC |
| HIMMELP1 | HS1 | HS110 | HS2 |
| HS25 | HS3 | HS38 | HS3MOD |
| HS4 | HS45 | HS5 | LINVERSE |
| LOGROS | MAXLIKA | MCCORMCK | MDHOLE |
| NCVXBQP1 | NCVXBQP2 | NCVXBQP3 | NONSCOMP |
| OSLBQP | PALMER1 | PALMER1A | PALMER1B |
| PALMER1E | PALMER2 | PALMER2A | PALMER2B |
| PALMER2E | PALMER3 | PALMER3A | PALMER3B |
| PALMER3E | PALMER4 | PALMER4A | PALMER4B |
| PALMER4E | PALMER5A | PALMER5B | PALMER5E |
| PALMER6A | PALMER6E | PALMER7A | PALMER7E |
| PALMER8A | PALMER8E | PENTDI | PFIT1LS |
| PFIT2LS | PFIT3LS | PFIT4LS | PROBPENL |
| PSPDOC | S368 | SANTALS | SIM2BQP |
| SIMBQP | SINEALI | WEEDS | YFIT |

## A.3 Linearly constrained problems

Table A.3: List of linearly constrained problems

| | | | |
|---|---|---|---|
| AVGASA | AVGASB | AVION2 | BIGGSC4 |
| BT3 | DALLASS | EQC | EXPFITA |
| EXPFITB | EXPFITC | FCCU | GENHS28 |
| HATFLDH | HS105 | HS118 | HS21 |
| HS21MOD | HS24 | HS268 | HS28 |
| HS35 | HS35I | HS35MOD | HS36 |
| HS37 | HS41 | HS44 | HS44NEW |
| HS48 | HS49 | HS50 | HS51 |
| HS52 | HS53 | HS54 | HS76 |
| HS76I | HS9 | HUBFIT | LSQFIT |
| ODFITS | PENTAGON | QC | QCNEW |
| S268 | WATER | ZECEVIC2 | |

## A.4 Nonlinearly constrained problems

Table A.4: List of nonlinearly constrained problems

| | | | |
|---|---|---|---|
| ACOPP14 | ACOPR14 | ALLINITA | ALLINITC |
| ALSOTAME | BATCH | BT1 | BT11 |
| BT12 | BT2 | BT4 | BT5 |
| BT6 | BT7 | BT8 | BURKEHAN |
| DEMBO7 | DIPIGRI | FLETCHER | FLT |
| GILBERT | HIMMELP2 | HIMMELP3 | HIMMELP4 |
| HIMMELP5 | HIMMELP6 | HS100 | HS100LNP |
| HS100MOD | HS104 | HS108 | HS11 |
| HS111 | HS111LNP | HS113 | HS114 |
| HS117 | HS12 | HS14 | HS15 |
| HS16 | HS17 | HS18 | HS20 |
| HS22 | HS23 | HS26 | HS27 |
| HS29 | HS30 | HS31 | HS32 |
| HS33 | HS40 | HS42 | HS43 |
| HS46 | HS47 | HS56 | HS57 |
| HS59 | HS6 | HS60 | HS61 |
| HS65 | HS68 | HS69 | HS7 |
| HS70 | HS71 | HS74 | HS75 |
| HS77 | HS78 | HS79 | HS80 |
| HS81 | HS83 | HS90 | HS91 |
| HS93 | HS99 | HS99EXP | JANNSON3 |
| JANNSON4 | LAUNCH | LEWISPOL | LOOTSMA |
| MARATOS | MATRIX2 | MESH | MISTAKE |
| MWRIGHT | OPTCNTRL | OPTPRLOC | ORTHREGB |
| ROBOT | S316-322 | S365 | S365MOD |
| SYNTHES1 | SYNTHES2 | SYNTHES3 | TRY-B |
| TWOBARS | ZECEVIC3 | ZECEVIC4 | ZY2 |

# Index